\documentclass[12pt, psamsfonts, reqno]{amsart}

\usepackage[T1]{fontenc}
\usepackage[latin1]{inputenc}
\usepackage[frenchb, english]{babel}
\usepackage{amsthm}
\usepackage{amscd}
\usepackage{tensor}

\usepackage{hyperref}
\usepackage[noabbrev]{cleveref}

\usepackage{amssymb,amsmath,bbm}
\usepackage{tikz-cd}
\usepackage[all]{xy}
\usepackage{mathrsfs}
\usepackage[left=0.0cm,right=0.0cm,top=1.7cm,bottom=1.7cm]{geometry}
\usepackage{hyperref}
\usepackage{comment}
\usepackage{enumitem}

\numberwithin{equation}{subsection}

\theoremstyle{plain}
\newtheorem{theorem}{Theorem}[subsection]

\newtheorem{lemma}[theorem]{Lemma}
\newtheorem{proposition}[theorem]{Proposition}
\newtheorem{prop}[theorem]{Proposition}
\newtheorem{corollary}[theorem]{Corollary}
\newtheorem{cor}[theorem]{Corollary}
\newtheorem*{theorem*}{Theorem}

\theoremstyle{definition}
\newtheorem{definition}[theorem]{Definition}
\newtheorem{example}[theorem]{Example}

\newtheorem{hypothesis}[theorem]{Hypothesis}

\newtheorem{notation}[theorem]{Notation}
\newtheorem*{claim}{Claim}
\newtheorem{construction}[theorem]{Construction}

\newtheorem{warning}[theorem]{Warning}

\theoremstyle{remark}
\newtheorem{remark}[theorem]{Remark}

\usepackage{tikz}
\usetikzlibrary{calc}

\def\Q{{\bf Q}}

\def\Z{{\bb Z}}

\def\N{{\bb N}}

\def\F{{\bf F}}

\def\ev{ev}

\def\Fil{\mathrm{Fil}}
\def\Hom{\mathrm{Hom}}
\def\Ext{\mathrm{Ext}}

\def\Sym{\mathrm{Sym}}

\def\epsilon{\varepsilon}
\def\det{\mathrm{det}}

\def\GL{\mathbf{GL}}

\def\M{\mathbf{M}}

\def\Cond{\mathsf{Cond}}
\def\cond{\mathrm{cond}}
\def\Solid{\mathsf{Solid}}
\def\solid{\mathrm{solid}}
\def\iHom{\underline{\mathrm{Hom}}}

\def\sol{{\square}}

\DeclareMathOperator{\Mod}{Mod}

\DeclareMathOperator{\Spa}{Spa}
\DeclareMathOperator{\Spec}{Spec}

\DeclareMathOperator{\Spf}{Spf}
\DeclareMathOperator{\gr}{gr}

\DeclareMathOperator{\Top}{Top}
\DeclareMathOperator{\Sch}{Sch}

\DeclareMathOperator{\Sets}{Sets}

\DeclareMathOperator{\id}{id}

\DeclareMathOperator{\Ab}{Ab}
\DeclareMathOperator{\LMod}{LMod}

\DeclareMathOperator{\Map}{Map}

\DeclareMathOperator{\Ind}{Ind}

\DeclareMathOperator{\End}{End}

\DeclareMathOperator{\CAlg}{CAlg}

\DeclareMathOperator{\AnRing}{AnRing}
\DeclareMathOperator{\Ani}{Ani}

\DeclareMathOperator{\AniRing}{AniRing}

\DeclareMathOperator{\Shv}{Shv}
\DeclareMathOperator{\PShv}{PShv}
\DeclareMathOperator{\pr}{pr}

\DeclareMathOperator{\Cat}{Cat}

\DeclareMathOperator{\Fun}{Fun}
\DeclareMathOperator{\Hub}{Hub}
\DeclareMathOperator{\CondRing}{CondRing}

\DeclareMathOperator{\AnSpec}{AnSpec}

\DeclareMathOperator{\Corr}{Corr}

\DeclareMathOperator{\Nil}{Nil}
\DeclareMathOperator{\QCoh}{QCoh}
\DeclareMathOperator{\Coh}{Coh}
\DeclareMathOperator{\Pro}{Pro}
\DeclareMathOperator{\Fin}{Fin}
\DeclareMathOperator{\Sm}{Sm}

\DeclareMathOperator{\Prof}{Prof}
\DeclareMathOperator{\Set}{Set}
\DeclareMathOperator{\CondSet}{CondSet}
\DeclareMathOperator{\CondAb}{CondAb}
\DeclareMathOperator{\CondAni}{CondAni}
\DeclareMathOperator{\Ring}{Ring}
\DeclareMathOperator{\AnStack}{AnStk}
\DeclareMathOperator{\AnStk}{AnStk}
\DeclareMathOperator{\CHaus}{CHaus}
\DeclareMathOperator{\Null}{Null}
\DeclareMathOperator{\iNull}{\underline{Null}}
\DeclareMathOperator{\Shift}{Shift}
\DeclareMathOperator{\sInd}{sInd}
\DeclareMathOperator{\fib}{fib}

\newcommand{\n}[1]{\mathcal{#1}}
\newcommand{\ob}[1]{\mathrm{#1}}
\newcommand{\bb}[1]{\mathbb{#1}}
\newcommand{\bbf}[1]{\mathbf{#1}}
\newcommand{\f}[1]{\mathfrak{#1}}
\newcommand{\cat}[1]{\mathsf{#1}}
\newcommand{\s}[1]{\mathscr{#1}}

\DeclareMathOperator{\op}{\begin{scriptsize}
op
\end{scriptsize}}

\DeclareMathOperator{\red}{\begin{scriptsize}
red
\end{scriptsize}}
\DeclareMathOperator{\aff}{\begin{scriptsize}
aff
\end{scriptsize}}

\DeclareMathOperator{\an}{\begin{scriptsize}
an
\end{scriptsize}}
\DeclareMathOperator{\st}{\begin{scriptsize}
st
\end{scriptsize}}

\DeclareMathOperator{\light}{\begin{scriptsize}
light
\end{scriptsize}}
\DeclareMathOperator{\Zar}{\begin{scriptsize}
Zar
\end{scriptsize}}
\DeclareMathOperator{\gas}{\begin{scriptsize}
gas
\end{scriptsize}}
\DeclareMathOperator{\bnd}{\begin{scriptsize}
bnd
\end{scriptsize}}
\DeclareMathOperator{\finpres}{\begin{scriptsize}
finpres
\end{scriptsize}}

\title{Notes on solid geometry}

\author{Juan Esteban Rodr\'iguez Camargo}

\makeatletter
\def\@tocline#1#2#3#4#5#6#7{\relax
  \ifnum #1>\c@tocdepth 
  \else
    \par \addpenalty\@secpenalty\addvspace{#2}%
    \begingroup \hyphenpenalty\@M
    \@ifempty{#4}{%
      \@tempdima\csname r@tocindent\number#1\endcsname\relax
    }{%
      \@tempdima#4\relax
    }%
    \parindent\z@ \leftskip#3\relax \advance\leftskip\@tempdima\relax
    \rightskip\@pnumwidth plus4em \parfillskip-\@pnumwidth
    #5\leavevmode\hskip-\@tempdima
      \ifcase #1
       \or\or \hskip 1em \or \hskip 2em \else \hskip 3em \fi%
      #6\nobreak\relax
    \dotfill\hbox to\@pnumwidth{\@tocpagenum{#7}}\par
    \nobreak
    \endgroup
  \fi}
\makeatother

\begin{document}

\begin{abstract}
These are expanded notes of a seminar held in Columbia university during the Spring and Fall of 2024 about the  theory of analytic  stacks of Clausen and Scholze, with a focus in the theory of solid mathematics. The seminar is inspired from the Lecture Series of Analytic Stacks of  Clausen and Scholze during the winter semester of 2023. All the theory of light condensed mathematics, analytic stacks and the proof of Serre duality must be attributed  to Clausen and Scholze, any mistake or misconception is totally due to the author. The only original work in these notes is the discussion of smooth, \'etale and finite presentation morphisms of solid Huber rings in Section 7.
\end{abstract}

\maketitle

\tableofcontents

\section{Introduction}
\label{SectionIIntroTalk}

 Different geometric theories appear all across mathematics:  differentiable manifolds, complex and real analytic varieties, rigid analytic spaces, adic spaces, Berkovich spaces, algebraic varieties and schemes, formal schemes,  etc.  The aim of ``Analytic Stacks'' is to provide a general ecosystem where the previous (and many more!) ``theories of analytic and algebraic geometries'' cohabit and interact with each other. To motivate the content of these notes, let us make explicit the obstructions that mathematicians have face all over the years when dealing with analytic geometry, and how condensed mathematics and analytic stacks have solved them. 
 
\subsection{Light condensed sets}
\label{SubsecI:LightCondensed}

The building blocks in theories such as algebraic varieties or schemes consist simply of commutative rings  satisfying some additional algebraic properties. This leads to a pleasant treatment of geometry that is studied in purely algebraic terms. However, in other theories such as differentiable manifolds and complex  or rigid analytic varieties, the building blocks turn out to be some sort of topological rings, more often Banach or Fr\'echet rings. Then, any  form of ``analytic geometry'' that has a similar formalism as algebraic geometry must be built up over an algebraic theory of ``topological rings''. However, history has shown  that the datum of a topology does not mixes very well with that of an algebraic structure. A very simple and clever solution to this issue is provided by condensed mathematics \cite{ClausenScholzeCondensed2019}, where ``topology'' is changed by the topos of (light) condensed sets. Therefore, our first replacement for topological ``\textit{preferred algebraic structure}'' (eg. ring/module/abelian group/monoid) will be a condensed ``\textit{preferred algebraic structure}''.

The idea behind condensed mathematics follows the philosophy of Grothendieck saying that a space $X$ must be studied by looking at maps $Y\to X$ from some ``test objects'' $Y$. For this approach to be useful, one needs to choose the ``test objects'' wisely. In our situation, we want to study (reasonable) topological spaces, which certainly should include compact Hausdorff spaces. It turns out that compact Hausdorff spaces can be reconstructed from a certain class of ``very acyclic''  spaces.   Concretely, let $\Prof$ be the category of profinite sets/totally disconnected compact Hausdorff spaces. We endow $\Prof$ with the Grothendieck topology whose covers are given by finitely many jointly surjective maps. As justification for this choice, recall that any surjective map of compact Hausdorff spaces is a quotient map, and that any compact Hausdorff space $X$ admits a surjection from a profinite set. For instance, the closed interval $[0,1]$ admits a surjective map from the Cantor set $\prod_{\bb{N}}\{0,1\}\to [0,1]$ by sending  a sequence $(a_n)$ to the real number written in binary decimals 
\[
(a_n)\mapsto 0.a_1a_2a_2\cdots.
\]

\begin{definition} A \textit{condensed set} is a sheaf $T:\Prof^{\op}\to \Set$ (modulo some set-theoretical technicalities, namely,  it must be accessible), we let $\CondSet$ denote the category of condensed sets. For  $X$ a (reasonable) topological space (eg. Hausdorff) compactly generated, we  define its condensification $\underline{X}\in \CondSet$ by taking as $S$-points for $S\in \Prof$
\[
\underline{X}(S)=\Map(S,X)
\]
the space of continuous maps from $S$ to $X$. 
\end{definition}

Most of the spaces we care of in topology (such as countably generated CW complexes), geometry (eg. manifolds), and analysis (eg. Banach, Fr\'echet spaces) are endowed with a topology for which understanding converging sequences is often enough. More precisely, the most interesting\footnote{Take this claim with a grain of salt, here \textit{interesting} is from the point of view of analysis. Spectral spaces are quite interesting by themselves, examples arising from spectrum of rings, but they are not well suited to do analysis.} topological spaces are (locally) metrizable.  Hence, a good balance  between capturing all the relevant topological information and  avoiding  unnecessary technicalities is given by light condensed sets:

\begin{definition}
A \textit{light profinite set} is a metrizable profinite set, we $\Prof^{\light}$ be the category of light profinite sets. A \textit{light condensed set} is a sheaf $T: \Prof^{\light,\op} \to \Set$, we let $\CondSet^{\light}$ denote the category of light condensed sets. 
\end{definition}

Finally, for any algebraic structure $\s{C}$ (aka. a category with small limits and colimits), its condensification $\Cond(\s{C})$ is the category of sheaves $T:\Prof^{\light,\op}\to \s{C}$ from light profinite sets to $\s{C}$. For example, we can talk about (light) condensed abelian groups, rings, monoids, etc. The notion of light condensed ``\textit{preferred algebraic structure}'' is the replacement we   use for its topological analogue.

\subsection{Analytic rings}
\label{SubsecI:AnalyticRings}

 As it was mentioned before,  part of the datum of the building blocks in a general theory of analytic geometry must involve some kind of topological (aka condensed) ring. On the other hand, the most fundamental invariant of a space $X$ in both analytic or algebraic geometry is its category of (quasi-)coherent sheaves $\QCoh(X)$. In classical algebraic geometry, this category is obtained by gluing, using ``Zariski descent'', the category of modules $\Mod(A)$ of commutative rings $A$.  However, in the case of complex and rigid geometries, the best that one can (classically) do in a systematic and algebraic manner is to built up the category of coherent modules $\Coh(X)$, imposing in this way finiteness conditions to the sheaves living over $X$. In particular, for a general morphism $f:X\to Y$ of rigid or complex analytic spaces, the sheaf $f_* \s{O}_X$ does not belong to the category attached to $X$. On the other hand, even though condensed rings are some kind of topological rings, in analytic geometries we often want to have some kind of ``completed tensor product'' and a category of ``complete modules''. It turns out that if $A$ and $B$ are two condensed rings, then the underlying ring of $A\otimes_{\bb{Z}} B$ is just the algebraic tensor $A(*)\otimes_{\bb{Z}} B(*)$, showing that we still need to do something else.

The notion of analytic ring appears as a solution to the previous problematic. The datum of an analytic ring $A$ consists of a condensed ring $A^{\triangleright}$ and a stable $\infty$-category $\ob{D}(A)$ of ``complete $A$-modules''. Before we list the features of $\ob{D}(A)$, let us do a brief detour explaining this jump from an abelian category of modules to a stable $\infty$-category:   in classical algebraic geometry, the category $\QCoh(X)$ of quasi-coherent sheaves is endowed with a symmetric tensor product $\otimes_{\s{O}_X}$. Thanks to  this tensor product one can construct fiber products $X\times_Y Z$ of    (affine) schemes by simply taking the (affine) scheme represented by the tensor product of rings. However, when dealing with cohomological invariants of algebraic varieties, it is natural to enter the world of derived categories.  In this realm the ``correct fiber product'' $X\times_Y Z$ is not longer constructed using the ``abelian'' tensor product of rings, but instead the ``derived tensor product''. Thanks to the modern foundations of higher category theory and higher algebra, eg. \cite{HigherTopos,HigherAlgebra,LurieSpectralAlg}, we have nowadays the categorical tools to develop theories of ``derived algebraic geometries'' as in \cite{LurieDerivedAlgebraic,LurieSpectralAlg}.  Nevertheless,  \Cref{TheoAnRingStructure} says that the category of complete modules $\ob{D}(A)$ only depends on the  abelian category of its heart, so strictly speaking   the additional ``analytic structure'' in the definition of analytic ring does not require the use of higher algebra. However, when gluing analytic rings into analytic stacks, the jump to higher category theory is unavoidable.

More precisely, in the former theories of analytic geometry, classical abelian or triangulated categories of quasi-coherent sheaves are not enough to obtain descend and glue to more general spaces (a reason for this  is the lack of ``complete'' flatness even for some simple maps such as open immersions of rigid or complex analytic spaces). Instead,  stable $\infty$-categories are perfectly suited for these purposes. As  consequence of the previous explanation, the general theory of analytic rings depends in higher categorical foundations (eg. the underlying condensed ring $A^{\triangleright}$ should be an animated or  a condensed $\bb{E}_{\infty}$-ring), even though the most fundamental examples still can be explained in the world of abelian categories. For the reader that is not comfortable   with the language of higher category theory, we recommend to consider $\ob{D}(A)$ as a classical derived category in a first approach, and to accept some  features of $\infty$-derived categories for granted such as the existence of arbitrary (small) limits of $\infty$-categories \cite[\S 3.3.3]{HigherTopos}, or  the adjoint functor theorem \cite[Corollary 5.5.2.9]{HigherTopos}.

Going back to the category $\ob{D}(A)$,  it ought to satisfy the following properties:

\begin{enumerate}

\item It should be a full subcategory $\ob{D}(A)\subset \ob{D}(A^{\triangleright})$ of the derived $\infty$-category of condensed $A^{\triangleright}$-modules stable under all limits and colimits, and it should be ``tensored over condensed abelian groups''. These are the basic requirements for doing homological algebra over $A$, and being ``$A$-complete'' should be a property and not extra structure. 

\item There is a ``completion functor'' $A\otimes_{A^{\triangleright}}-:\ob{D}(A^{\triangleright}) \to \ob{D}(A)$, left adjoint to the natural inclusion (note that we have dropped derived decorations in the tensor). Moreover, $\ob{D}(A)$ can be uniquely promoted to a symmetric monoidal category such that $A\otimes_{A^{\triangleright}}-$ is a symmetric monoidal functor. Similarly as for schemes, we require our category of modules to be endowed with a ``completed tensor product'' that will generalize those  of classical theories of analytic geometry. 

\item The completion functor $A\otimes_{A^{\triangleright}}-$ should preserve connective objects: $A\otimes_{A^{\triangleright}}-:\ob{D}(A^{\triangleright})_{\geq 0}\to \ob{D}(A^{\triangleright})_{\geq 0}$. This will endow $\ob{D}(A)$ with a $t$-structure arising from condensed $A^{\triangleright}$-modules. 

\item We have $A^{\triangleright}\in \ob{D}(A)$ (we want our  topological ring to be complete!).

\end{enumerate}

\begin{definition}
An \textit{analytic ring $A$} is a pair $(A^{\triangleright}, \ob{D}(A))$ consisting of a light condensed animated ring $A^{\triangleright}$, and a full subcategory $\ob{D}(A)\subset \ob{D}(A^{\triangleright})$ of ``\textit{complete modules}'' satisfying properties (1)-(4) above. A morphism of analytic rings $f:A\to B$ is a morphism of condensed animated rings $A^{\triangleright} \to B^{\triangleright}$ such that the forgetful functor $f_*:\ob{D}(B^{\triangleright})\to \ob{D}(A^{\triangleright})$ sends $\ob{D}(B)$ to $\ob{D}(A)$.  We let $\AnRing$ denote the $\infty$-category of analytic rings. 
\end{definition}

It turns out that $\AnRing$ a is a presentable $\infty$-category\footnote{This is claimed in the course of Analytic Stacks  \cite{AnalyticStacks}, unfortunately I do not know a proof of this fact. In practice, however, there are ways to go around this as we will see in the constructions of categories of analytic stacks.} (cf. \cite[\S 5.5]{HigherTopos} for the notion of presentability), in particular it admits all (small) colimits (cf.  \cite[Proposition 12.12]{ClauseScholzeAnalyticGeometry} and \cite[Proposition 2.3.15]{MannSix}). Analytic rings are  the bulding blocks in the theory of analytic stacks.

\subsection{Analytic stacks}
\label{SubsecI:AnalyticStacks}

Let $\Ring$ be the category of rings. Schemes are constructed out from $\Ring$ by gluing using the Zariski topology. In particular, a scheme can be seen as an object in $\Shv_{\Zar}(\Ring^{\op},\Set)$, i.e. a sheaf for the Zariski topology in the opposite category of rings (aka affine schemes). Similarly, algebraic spaces (resp. Artin stacks) are obtained by ``gluing affine schemes'' along \'etale or smooth maps, they then define sheaves in more refined Grothendieck topologies such as the \'etale or flat topologies.  Moreover, when defining stacks in derived algebraic geometry \cite{LurieDerivedAlgebraic}, it is mandatory to not just consider functors with values in sets but in anima $\Ani$ (aka $\infty$-groupoids or spaces).

For the theory of analytic stacks, we want to define a suitable Grothendieck topology $G$ on $\AnRing$ such that ``analytic stacks'' are given by (``good'' hyper)sheaves 
\[
\AnStack=\Shv_{G}(\AnRing^{\op},\Ani).
\]
The key question is what Grothendieck topology  we should  consider. Well, by definition analytic rings \textbf{are not just} its underlying condensed ring but its category of modules. Indeed, an analytic ring is (essentially) completely determined by its category of modules! Thus, whatever Grotendieck topology we choose, the functor $A\mapsto \ob{D}(A)$ should certainly satisfy descent. On the other hand, we want a  refine enough Grothendieck topology that explains already existing ``identifications'' in classical analytic geometries that show up in nature: 

Let $\bb{Q}_p$ be the field of $p$-adic numbers, and consider the projective space $\bb{P}^1_{\bb{Q}_p}$. There are different ways to construct $\bb{P}^1_{\bb{Q}_p}$. First, we have the algebraic-geometric manner that glues the (spectrum of the) rings $\bb{Q}_p[T]$ and $\bb{Q}_p[T^{-1}]$ along the intersection $\bb{Q}_p[T^{\pm 1}]$. On the other hand, we have rigid geometry and we can construct $\bb{P}^1_{\bb{Q}_p}$ by gluing the (adic spectrum of the) Tate algebras $\bb{Q}_p\langle T \rangle$ and $\bb{Q}_p\langle T^{-1} \rangle$ along the intersection $\bb{Q}_p\langle T^{\pm 1} \rangle$. Thus, we want the theory of analytic stacks to be able to identify these both constructions of $\bb{P}^1_{\bb{Q}_p}$ as the same space, obtaining  as a result a geometric version of GAGA theorems. 

In later sections, we shall introduce the  definition of the Grothendieck topology used for defining analytic stacks. A key tool in its definition will be the abstract theory of six functor formalisms \cite{HeyerMannSix}.

\subsection{Examples}
\label{SubsecI:Examples}

Throughout the introduction of light condensed sets, analytic rings and analytic stacks, we shall study in more detail some examples arising from algebraic geometry and the theory of adic spaces (solid theory).  In this introduction we will  shortly mention the existence of archimedean and global examples of analytic rings (liquid and gaseous theory). 

\textit{Solid abelian groups}. Let $\CondAb^{\light}$ denote the category of light condensed abelian groups. We shall define the subcategory of (light) solid abelian groups $\Solid \subset \CondAb^{\light}$ by imposing a condition extracted from the idea that  ``converging sequences in non-archimedean analysis are precisely the null sequences''. The category of solid abelian groups is endowed with a tensor product that we denote by $\otimes_{\sol}$, it has $\bb{Z}$ as unit, and so it defines an analytic ring $\bb{Z}_{\sol}$ that we call the \textit{solid integers}.  The category $\Solid$ has a compact projective generator $\prod_{\bb{N}} \bb{Z}$ that is flat for $\otimes_{\sol}$, and satisfies 
\[
\prod_{I} \bb{Z} \otimes_{\sol} \prod_{J} \bb{Z} =  \prod_{I\times J} \bb{Z}
\] 
for countable sets $I,J$. This category is completely disjoint from archimedean analysis, namely, the solidification of $\bb{R}$ is just $0$.  Examples of solid abelian groups are discrete groups, $p$-adically complete modules, $\bb{Q}_p$-Banach and Fr\'echet spaces, etc. Moreover, most of the completed tensor products appearing in non-archimedean analysis coincide with $\otimes_{\sol}$ (eg. $p$-complete tensor products of Banach spaces, projective tensor product of nuclear Fr\'echet spaces). The category of solid abelian groups is also the smallest full subcategory of condensed abelian groups containing $\mathbb{Z}$, and stable under limits and colimits.

\textit{Liquid vector spaces}.  Let $q\in (0,1]$.  The analytic ring of liquid real vector spaces was constructed in \cite{ClauseScholzeAnalyticGeometry}. The construction of this analytic ring requires a lot of effort due to the non-locally convex functional analysis involved. For instance, if $\bb{R}_{<q}$ denotes the analytic ring of $(<q)$-liquid real numbers, and $S$ is a profinite set, then the free liquid real vector space $\bb{R}_{<q}[S]$ is not the naive guess of signed Radon measures on $S$, but a certain space of $(<q)$-convex Radon measures. The liquid tensor product agrees with the projective tensor product for nuclear Fr\'echet spaces, as well as for their duals, see \cite[IV]{CondensedComplex}.

\textit{Gaseous rings}. As we shall see later, one of the main advantages of the new foundations for the theory of analytic rings, based on light condensed sets, is that it is much easier to construct analytic rings out from inverting some concrete maps of modules. The difficulty is then translated in computing the functors of ``measures'' $A[S]$ for $S\in \Prof^{\light}$. The initial gaseous ring with fixed pseudo-uniformizer is defined via some universal property in the category of analytic rings. It specializes in both solid and liquid rings, and its underlying ring $\bb{Z}[\widehat{q}]^{\gas}\subset \bb{Z}[[q]]$ consists on power series of at most polynomial growth\footnote{This is a highly non trivial computation that is also claimed in \cite{AnalyticStacks}. Unfortunately I haven't been able to fully understand the proof of this.}:
\[
\bb{Z}(\widehat{q})^{\gas}(*)=\{\sum_{n>>{-\infty}}  a_n q^n : \; \exists \; m,k>0  \mbox{ such that } \lim_{n\to \infty} |a_n| (n+m)^{-k} =0  \}.
\]
The gaseous ring was partially motivated from the construction of Tate's elliptic curve $\bb{G}^{\an}_{m,A}/q^{\bb{Z}}$ in an universal way, and it gives a conceptual explanation  why the coefficients of the algebraic equation defining the Tate curve are power series with polynomial growth.

\subsection{Overview of the document}

 The first three sections are taken from the course of analytic stacks of Clausen and Scholze \cite{AnalyticStacks}.   In \Cref{Section:LightCondensed} we introduce the set up for light condensed mathematics, we define light profinite sets, light condensed sets and light condensed abelian groups. 
 
  In \Cref{SectionLightSolid} we define the category of solid abelian groups as those condensed abelian groups for which all null sequence is summable in a suitable strong sense. We study in detail the category of solid abelian groups, and give some examples of solid tensor products. 
  
   \Cref{SectionAnalyticRings} reviews the general theory of  analytic rings; we explain how to compute colimits of analytic rings, the dependence of analytic ring structures only on the underlying abelian category, the invariance of analytic ring structures with respect to higher homotopy groups or nilpotent thickenings, and a general way to construct examples of analytic rings by imposing summability conditions on null sequences.  
 
 In \Cref{SubsecMoreSolidRings} we begin the study of solid analytic rings, we follow the lecture notes \cite{ClausenScholzeCondensed2019}, Mann's thesis \cite{MannSix}, and the course \cite{AnalyticStacks}.  We start with a discussion of the smashing spectrum of a symmetric monoidal stable category; this general formalism is very useful for the intuition behind the geometry of the six functor formalism of analytic stacks. We study solid analytic ring structures of finite type algebras over $\Z$ and stablish some basic properties of their categories of complete modules. Then, we discuss different categories of sheaves for schemes (classical quasi-coherent sheaves and solid sheaves), and categories of solid sheaves for discrete adic spaces. 

Next, in \Cref{s:AnStacks}, we introduce the category of analytic stacks following \cite{AnalyticStacks} as well as some corrections that has been addressed in \cite{dRFF} and \cite{ScholzeGestalted}. We briefly discuss the abstract theory of six functor formalisms and the category of kernels following  \cite{HeyerMannSix} and \cite{SixFunctorsScholze}. We construct the six functor formalism for analytic rings and from this, following a rather general procedure, construct categories of analytic stacks. This is the most technical section of the notes, and it freely uses  the theory of linear presentable  categories and higher algebra. We finish by producying two main examples of analytic stacks, those comming from algebraic geometry (called algebraic stacks in these notes), and those comming from condensed anima (called Betti stacks and studied in detail in \cite{HeyerMannSix}).

We end with \Cref{SectionSerreDualitySolid} with the proof of Serre duality following the strategy of \cite{CondensedComplex}. We introduce the category of solid Huber rings; a generalization of the classical category of complete Huber pairs. We construct the solid spectrum of a solid Huber ring, a spectral space that generalizes the classical adic spectrum of Tate Huber rings, define solid Huber stacks and solid adic spaces.  We discuss some basic facts about the theory of cotangent complexes of analytic rings. We introduce the notions of morphism (locally) of solid finite presentation, solid \'etale and solid smooth morphisms  between solid Huber rings (and solid adic spaces), and show some basic stability properties. Finally, we state and prove Serre duality for (locally) solid smooth morphisms of solid adic spaces; we deduce the classical statement of Serre duality for schemes, and deduce Serre duality for rigid varieties.

\subsection*{Acknowledgements} I heartily thank  all the participants of the seminar on solid geometry during 2024 at Columbia university for their motivation and interest in learning the theory of condensed mathematics and analytic geometry. I thank Johannes Ansch\"utz,  Arthur-C\'esar Le Bras, Guido Bosco, Yutaro Mikami and Semen Slobodianuk, for many comments and corrections in previous versions of this document. I am specially greatful to Dustin Clausen and Peter Scholze for many discussions and explanations of their beautiful theory, and for letting me present part of their work in these notes.


\section{Light condensed mathematics}\label{Section:LightCondensed}

In this talk we will study the basics in light condensed mathematics, this  involves light profinite sets, light condensed sets and light condensed abelian groups.

\subsection{Light profinite sets}

Condensed mathematics proposes a  better algebraic  framework that replaces topological spaces, namely condensed sets. The building blocks of condensed sets are profinite sets which we recall down below:

\begin{prop}
\label{PropositionEquivalentProfinite}

The following categories are equivalent.

\begin{enumerate}

\item The  pro-category of finite sets $\Pro(\Fin)$ where maps are given by 
\[
\Map(\varprojlim_{i} S_i,\varprojlim_j T_j)= \varprojlim_{j} \varinjlim_i \Map(S_i,T_i).
\]

\item The category of  totally disconnected compact Hausdorff spaces with continuous maps. 

\item The opposite category of Boolean algebras.

\end{enumerate}

We let $\Prof$ denote the category of profinite sets, considered as in (1) or (2) above. 

\end{prop}
\begin{proof}

We just construct the equivalences. From (1) to (2) we take a projective system $\{S_i\}_{i}$ and pass to the topological space $S=\varprojlim_i S_i$ endowed with the limit topology.  From (2) to (1) we take a totally disconnected compact Hausdorff space and consider the projective system $\{S_i\}_{i\in I}$  of finite quotients of $S$, equivalently, the projective system of partitions of $S$ in clopen subspaces. From (2) to (3) we take a totally disconnected compact Hausdorff space $S$ and consider the Boolean algebra $A=C(S,\bb{F}_2)$ of continuous functions from $S$ to $\bb{F}_2$. From (3) to (2) we take a Boolean algebra $A$ and consider its spectrum $\Spec A$ as a topological space. 
\end{proof}

A delicate issue when working with the category of all profinite sets is that it is not essentially small, i.e. there is not a set of isomorphism classes of objects. On the other hand, all the spaces we actually care about that appear in geometry, topology or analysis (such as manifolds,  CW complexes, Banach or Fr\'echet spaces) admit a metric, and can be recovered within a \textbf{set} of smaller profinite sets. 

\begin{prop}
\label{PropLightProfinite}
Let $S$ be a profinite set, the following are equivalent: 

\begin{enumerate}

\item $C(S,\bb{Z})$ is countable,

\item $S$ is metrizable,

\item $S$ is $2$-countable,

\item $S$ can be written as a sequential limit of finite sets. 

\end{enumerate}

\end{prop}
\begin{proof}
Urysohn's metrization theorem implies that a compact Hausdorff space is metrizable if and only if it is 2-countable, this shows (2) $\Leftrightarrow$ (3).

 (3) $\Leftrightarrow$  (4). By  \Cref{PropositionEquivalentProfinite}, the passage from a totally disconnected compact Hausdorff space $S$ to a projective system of finite sets is made by taking the system of partitions of $S$ into clopen subspaces, since $S$ is 2-countable this projective system is countable. Conversely, if $S=\varprojlim_{\bb{N}} S_n$ is a sequencial limit of finite sets, taking the fibers of the maps $S\to S_n$ defines a  countable basis for the topology of   $S$. 
 
(4) $ \Rightarrow$ (1). If $S=\varprojlim_{\bb{N}} S_n$, then $C(S,\bb{Z})= \varinjlim_{n} C(S_n, \bb{Z})$ which is countable. 

(1) $\Rightarrow$ (3). Finally, if $C(S,\bb{Z})$ is countable, then $C(S,\bb{F}_2)$ is countable, and $S= \Spec C(S,\bb{F}_2)$ has at most countably many clopen subspaces,  proving that $S$ is $2$-countable. Indeed, clopen subspaces of $S$ are in bijection with the elements of $C(S,\bb{F}_2)$.
\end{proof}

\begin{definition}
\label{DefinitionLightprof}
A profinite set is \textit{light} if it satisfies the equivalent conditions of  \Cref{PropLightProfinite}. We let $\Prof^{\light}$ denote the category of light profinite sets. 
\end{definition}

One has the basic stability properties of light profinite sets under countable limits.

\begin{prop}
\label{PropSequentialLimitSurj}
The category of light profinite sets admits countable limits. Moreover, sequential limits of surjections is a surjection. 
\end{prop}
\begin{proof}
Stability under countable limits follows from  \Cref{PropLightProfinite} and the fact that a countable limit of 2-countable topological spaces is 2-countable. Let $S=\varprojlim_n S_n$ be a sequencial limits of surjections, then the  map $S\to S_n$ is  surjection, namely, given $x_n\in S_n$ take lifts $x_{n+m}\in S_{n+m}$, inductively such that $x_{n+m+1}$ maps to $x_{m}$.
\end{proof}

Next, we prove some nice features that are special to the category of light profinite sets.

\begin{prop}
\label{PropOpenLightDisjointUnion}
Let $S$ be a light profinite set and let $U\subset S$ be an open subspace. Then $U$ is a countable disjoint union of light profinite sets. 
\end{prop}
\begin{proof}
Let us write $S=\varprojlim_{n} S_n$ and let $Z=S\backslash U$.    Then $Z=\varprojlim_n Z_n$ with $Z_n\subset S_n$  the image of $Z$ in $S_n$. Let $\pi_n: S\to S_n$ and $\pi_{m,n}:S_m\to S_n$ denote the projection maps. We define  $Y_0=S_0\backslash Z_0$ and for $n\geq 1$ we let $Y_n= S_n\backslash (Z_n \cup \pi^{-1}_{n,n-1} Y_{n-1}\cup \cdots \cup \pi^{-1}_{n,0}(Y_0))$. Then
\[
U=\bigsqcup_{n\in \bb{N}}  \pi_n^{-1}(Y_n).
\]
\end{proof}

\begin{prop}
\label{PropInjectivelight}
Let $S$ be a light profinite set. Then $S$ is an injective object in $\Prof^{\light}$. 
\end{prop}
\begin{proof}
Let $f:X\to Y$ be an injection of light profinite sets and let $g:X\to S$ be a map.   The map $f$ is a closed immersion, then we can write it as a sequential limit $\varprojlim_{n} (f_n:X_n\to Y_n)$ of injective maps of  finite sets.  We can write the map $g$ as a sequential limit of finite sets $\varprojlim_{n} (g_n:X_{k_n}\to S_n)$ with $(k_n)_n$ being  some increasing sequence. After taking a subsequence we can assume that $k_n=n$. Then, we can always find a map $h_0: Y_0 \to  S_0$ extending $g_0$, and provided the extension $h_n:Y_n\to S_n$, we can always find a map $h_{n+1}:Y_{n+1}\to S_{n+1}$ extending $g_{n+1}$ that reduces to $h_n$ in the $n$-th step. Taking the limit $h=\varprojlim_{n}h_n$ we get the desired map $h:Y \to S$ extending $g$. 
\end{proof}

\begin{prop}[{\cite[Theorem 5.4]{ClausenScholzeCondensed2019}}]
\label{PropFreeContinuousFunctions}
Let $S$ be a light profinite set, then the space of continuous functions $C(S,\bb{Z})$ is a free $\bb{Z}$-module. 
\end{prop}
\begin{proof}
Let us write $S=\varprojlim_{n} S_n$ as a sequential limit with surjective maps. We can find compatible sections 
\[
S_0\to S_1\to S_2\to \cdots
\]
and then inductively find compatible sections $S_0\to S$, $S_1\to S$, $\cdots$. Then, we know that 
\[
C(S,\bb{Z})=\varinjlim_{n} C(S_n,\bb{Z}),
\] 
and we just found compatible sections of $C(S_n,\bb{Z})\to C(S,\bb{Z})$, since the modules $C(S_n,\bb{Z})$ are free, this shows that $C(S,\bb{Z})$ is also free. 
\end{proof}

\begin{example}
\label{ExamplesLightProfinite}
The two examples of light profinite sets that will be the most relevant for us are the following:

\begin{enumerate} 

\item The one point compactification of $\bb{N}$, namely, $\bb{N}\cup \{\infty\}$. It can be written as 
\[
\bb{N}\cup \{\infty\} = \varprojlim_{n} \{1,2,\ldots, n,\infty\}
\]
where for $m\geq n$ the map $\{1,2,\ldots, m,\infty\}\to \{1,2,\ldots, n,\infty\}$ sends all the elements $k\geq n+1$  to $\infty$. 

\item The Cantor set $S=\prod_{\bb{N}}\{0,1\}$, it admits a surjective map onto the interval $[0,1]$ by taking binary decimal expansions. 

\end{enumerate}

\end{example}

The relevance of the Cantor set is explained in the following proposition. 

\begin{prop}
\label{ProplightQuotientCantor}
A profinite set is light if and only if it admits a surjective map from the Cantor set. 
\end{prop}
\begin{proof}
Let $S=\varprojlim_{n} S_n$ be a light profinite set, and let us suppose that $S\to S_n$ is surjective for all $n$. Then, we can always find a sequence of non-negative integers $(k_n)_{n\in \bb{N}}$ and compatible surjection maps for varying $n$
\[
\prod_{m=0}^{k_n}\{0,1\} \to S_n.
\]
Taking the limit we get the desired surjection from the Cantor set. The converse follows from the fact that the quotient of a $2$-countable topological space is still $2$-countable. 
\end{proof}

\subsection{Light condensed sets}

After the previous preparations of light profinite sets, we can finally define  light condensed sets (cf. \cite[Definition 1.2]{ClausenScholzeCondensed2019}):

\begin{definition}
\label{DefinitionLightCondensed}
A \textit{light condensed set} is a sheaf in the category of light profinite sets for the Grothendieck topology generated by finite disjoint unions of jointly surjective maps. More concretely, a condensed set is a functor $T:\Prof^{\light,\op}\to \Set$ such that 
\begin{enumerate}

\item $T(\emptyset)=*$.

\item $T(S_1\sqcup S_2 )=T(S_1)\times T(S_2)$.

\item  For all surjective map $S_1\to S_2$ we have 
\[
T(S_2)=\mathrm{eq}(T(S_2)\rightrightarrows T(S_2\times_{S_1} S_2) ).
\]

\end{enumerate} 

 We let $\CondSet^{\light}$ denote the category of light condensed sets.   
\end{definition}

\begin{remark}
\label{RemarkRepleteTopos}
By  \Cref{PropSequentialLimitSurj}, sequential limits of covers in $\Prof^{\light}$ are covers. In particular, the topos of condensed sets is replete in the sense of \cite[\S 3]{bhatt2014proetale}, namely, sequencial limits $T=\varprojlim_{n} T_n$ of condensed sets with surjective maps are still surjective. Indeed, by definition of the Grothendieck topology, given $S_0\to T_0$ an $S_0$-point of $T_0$ there is a surjective map $S_1\to S_0$ and a lift $S_1\to T_1$. Repeating this process we find a compatible sequence of points $S_n\to T_n$ with $S_{n+1}\to S_n$ a surjective map. Then, taking limits $S=\varprojlim_n S_n\to T$, we get a lift of $S_0\to T_0$ to $S\to T$ and the map $S\to S_0$ is a cover in the Grothendieck topology  being surjective by  \Cref{PropSequentialLimitSurj}.
\end{remark}

\begin{example}
\label{ExampleToptoCond}
\begin{enumerate}

\item Let $T$ be a light condensed set, then the set $T(\bb{N}\sqcup \{\infty\})$ is heuristically the space of convergent sequences with fixed limit, namely, this is exactly the case when $T$ arises from the condensification of a topological space. If $T=\underline{X}$ arises from a Hausdorff space then the set of convergent sequences are determined by its restriction to $\bb{N}$, i.e. the map $T(\bb{N}\sqcup\{\infty\})\to T(\bb{N})$ is injective.  In general, a convergent sequence can have different limits, so the map $T(\bb{N}\sqcup \{\infty\})\to T(\bb{N})$ is not necessarily injective. 

\item Let $\Top$ denote the category of topological spaces. We define the condensification functor 
\[
\underline{(-)}: \Top\to \CondSet^{\light}
\]
mapping a topological space $X$ to the condensed set $\underline{X}:S\mapsto C(S,X)$ for $S\in \Prof^{\light}$.

\item The Yoneda embedding $\Prof^{\light}\hookrightarrow  \CondSet^{\light}$ maps a profinite set $S$ to its condensification $\underline{S}$.  Since $\Prof^{\light}$ is a small category,  any condensed set can be written as a colimit of light profinite sets. More precisely, we have that 
\[
T=\varinjlim_{S\to T} \underline{S}
\]
as a condensed set.  From now we will not make further distinction between $S$ and $\underline{S}$ for $S$  a light profinite set.

\end{enumerate}
\end{example}

As we saw in the previous example, there is a natural functor from topological spaces to light condensed sets by taking mapping spaces. The following proposition shows that this functor is fully faithful in a reasonable  subcategory of topological spaces (cf. \cite[Proposition 1.7]{ClausenScholzeCondensed2019})

\begin{prop}
\label{PropCondensification}
The condensification functor has a left adjoint called the ``underlying topological space'', mapping a condensed set $T$ to the topological space given by 
\[
T(*)_{\mathrm{top}}=\varinjlim_{S\to T} S
\]
where the colimit is taken in the category of topological spaces.  More precisely, $T(*)_{\mathrm{top}}$ has underlying set $T(*)$ and  topology determined by the set of maps 
\[
\bigsqcup_{S\to T} S\to T(*). 
\]
In particular, the functor $\underline{(-)}$ is fully faithful in metrizably  compactly generated  spaces (eg.  metrizable compact Hausdorff spaces).

\end{prop}
\begin{proof}
Since $T= \varinjlim_{S\to T} S$ as a condensed set, the statement reduces to the fact that for a profinite set $S$ and a topological space $X$ we have 
\[
\underline{X}(S)=C(S,X). 
\]
\end{proof}

\begin{remark}
\label{RemarkSurjectionCondensed}
Let us make more explicit what means to be an epimorphism for topological spaces when considered as condensed sets.  Let $X\to Y$ be a map of topological spaces such that their condensification $\underline{X}\to \underline{Y}$ is an epimorphism. This means that for any light profinite set $S$ and any map $f:S\to Y$, there is a surjection from a light profinite set $S'\to S$ and a map $S'\to X$ lifting $S$. For instance, if $X\to Y$ is a surjection of compact Hausdorff spaces then so is its condensification. However, this property does not hold true in general, for example, consider the map $\bigsqcup_n B_n\to \varinjlim_n B_n$, where  $( B_n)_{n}$  is an inductive system of Banach spaces with injective transitions maps (LB spaces);  for this map  to be surjective is sufficient that  the transition maps are of compact type, and one can produce counter examples otherwise. In other words, the condensification of the topological colimit  $\varinjlim_n B_n$ is not necessarily the colimit of the condensification of the Banach spaces $B_n$ unless one adds some hypothesis in the transitions maps.  

 As a concrete example, consider $B_n=\ell_{\infty}(\N)$ the Banach space over a non-discrete normed field $K$ (eg. $\bb{Q}_p$ or $\bb{C}$) given by sequences $(a_0,a_1,a_2,\ldots)$ with $a_n\in K$ with $\sup_{n}|a_n|<\infty$ endowed with the $\sup$ norm. Consider the transition map $f_{n}\colon B_n\to B_{n+1}$ given by shifting the sequence $(a_0,a_1,a_2,\ldots)\mapsto (a_1,a_2,a_3\ldots)$. The transition maps are not of compact type and the classical colimit of Banach spaces $B:=\varinjlim_{n} B_n$ is pathological as it gives rise to a non-Hausdorff space. On the other hand, the colimit seen as a condensed abelian group is a perfectly well behaved space from the point of view of homological algebra and it actually captures some topological information of $K$; te underlying set of the condensed abelian group $B$ is the space of tails of bounded sequences in $K$. 
\end{remark}

In any topos there is a notion of quasi-compact and quasi-separated objects, in the case of light condensed abelian groups these properties can be stated in more concrete terms.

\begin{definition}
\label{DefinitionQCQS}
A condensed set $T$ is \textit{quasi-compact} if there is a surjection $S\to T$ from a profinite set. A condenset set $T$ is \textit{quasi-separated} if for every two maps from profinite sets $S\to T \leftarrow S'$, the fiber product $S\times_T S'$ is quasi-compact. 
\end{definition}

\begin{remark}
By definition, the Grothendieck topology of $\Prof^{\light}$ is finitary, this makes the profinite sets quasi-compact objects in the topos of condensed sets. Moreover, since light profinite sets are stable under countable limits, they are stable under pullbacks and so they are  quasi-separated. This makes $\CondSet$ a coherent topos.   On the other hand, if $T$ is a condensed set and $S,S'\to T$ are maps from profinite sets to $T$, then $S\times_T S'$ is a subobject of $S\times S$, therefore $T$ is quasi-separated if and only if for all $S,S'$ as before $S\times_T S'$ is also profinite. 
\end{remark}

We can describe concretely the qcqs objects in $\CondSet$. 

\begin{prop}
\label{PropCompactHausQCQS}
Let $\CHaus^{\light}$ be the category of metrizable compact Hausdorff spaces. Then the condensification functor induces an equivalence from $\CHaus^{\light}$ to the category of qcqs  light condensed sets. Moreover, the category of quasi-separated light condensed sets is equivalent to the ind-category  with injective transition maps of metrizable compact Hausdorff spaces  $\Ind_{\mathrm{inj}}(\CHaus^{\light})$. 
\end{prop}
\begin{proof}
First, we claim that a quasi-compact subobject of a light profinite set is necessarily profinite.  For this, let $f:S\to S'$ be a map of light  profinite sets, we want to see that the image of $f$ is a closed subspace of $S'$. Let $\mathrm{Im}(f)\subset S'$ be the image as topological space,  it is profinite and we know that  $f$ factors through the condensification of  $\mathrm{Im}(f)$. Then, we are left to show that if $f$ is a surjection of light profinite sets, then it is an epimorphism as condensed sets, but this is clear by the definition of the Grothendieck topology of $\Prof^{\light}$.

  Let $T$ be a qcqs object in $\CondSet$, then there is a surjection $S\to T$ from a light profinite set such that $S\times_{T} S$ is also profinite. Then, $T$ arises as the quotient of a light profinite set by a light profinite equivalence relation, making $T(*)_{\mathrm{top}}$ a metrizable compact Hausdorff space, and the natural map $\underline{T(*)_{\mathrm{top}}} \to T$ from the adjunction is an equivalence by  \Cref{RemarkSurjectionCondensed}. 
  
  Conversely, let $X$ be a metrizable compact Hausdorff space and fix a countable basis $\f{U}$ of $X$ consisting on compact subspaces.  Let $I$  denote the category consisting on tuples $(U_s)_{s\in S_i}$ indexed by  finite sets $S_i$ where the $\{U_s: \; s\in S_i\}$ form a compact cover of $X$ by elements in $\f{U}$. A map $F: (U_s)_{s\in S_j}\to (V_{s})_{s\in S_i}$ is a map of finite sets $f: S_j\to S_i$ and inclusions $U_{s}\subset V_{f(s)}$ for all $s\in S_j$.  The category $I$ is cofiltered, namely, given two  finite covers of $X$ in $\f{U}$ one can always find a common refinement, and given two maps $F,G: (U_s)_{s\in S_j}\to (V_s)_{s\in S_i}$ one can find a refinement $H:(W_s)_{s\in S_k}\to (U_s)_{s\in S_j}$ such that $F\circ H=G\circ H$. Then,  $S=\varprojlim_{i} S_i$ is a light profinite set.  We can define a map  $f:S\to X$ by sending a compatible system of compact subsets $x=\{U_{s_j}\}_{j\in I}$ to its intersection $f(x)=\bigcap_{j} U_{s_j}$ which is necessarily a point. The map $f$ is  continuous and  a surjection from a light profinite set onto $X$.

   By  \Cref{RemarkSurjectionCondensed}, the map of condensed sets $S\to \underline{X}$ is surjective, and the fiber product $S\times_{\underline{X}} S$ is the condensification of the topological fiber product which is a light profinite set, this shows that $\underline{X}$ is qcqs as wanted. 

Finally, let $T$ be a quasi-separated light condensed set, and let $S\to T$ be a map from a profinite set $S$. Then the image $X$ of $S$ in $T$ is qcqs since $S\times_X S= S\times_T S$ is profinite. This shows that $T$ can be written as a union of qcqs condensed sets by injective maps, which produces an object in $\Ind_{\mathrm{inj}}(\CHaus^{\light})$, furthermore, since qcqs condensed sets are compact objects in $\CondSet$ this map is fully faithful. Conversely, given a cofiltered diagram $\{X_i\}_{i}$ of light compact Hausdorff spaces with injective transition maps, the colimit $T=\varinjlim_{i} \underline{X_i}$ of condensed sets is quasi-separated, namely, given any two maps from profinite sets $S,S'\to T$   there is some $i$ such that $S,S'$ factor through $X_i$, and $S\times_T S'=S\times_{\underline{X_i}} S'$ is profinite. 
\end{proof}

\subsection{Light condensed abelian groups}

Next, we define light condensed abelian groups and prove some of its most important features.

\begin{definition}
\label{DefinitionCondAb}
The category of \textit{light condensed abelian groups} $\CondAb^{\light}$ is the category of abelian group  objects in $\CondSet^{\light}$. Equivalently,  it is the category of abelian sheaves on light profinite sets. 
\end{definition}

\begin{example}
\label{ExampleFreeCondensed}
\begin{enumerate}
\item The forgetful functor
\[
\CondAb^{\light} \to \CondSet
\]
has a left adjoint $T\mapsto \bb{Z}[T]$ given  by the free abelian group generated by a condensed set. The condensed abelian group $\bb{Z}[T]$ is given by the sheafification of the functor mapping a light profinite set $S$ to the free abelian group $\bb{Z}[T(S)]$.  

\item  Let $A$ be a topological abelian group, then $\underline{A}$ has a natural structure of light condensed abelian group. Indeed, the condensification functor preserves finite limits and the structure of an abelian group for $A$ is encoded in some diagrams such as $+:A\times A\to A$. 

\item Let $\bb{R}$ be the real numbers endowed with the addition and its natural topology, then $\underline{\bb{R}}$ is a condensed abelian group. On the other hand, if $\bb{R}^{\delta}$ is endowed with the discrete topology then $\underline{\bb{R}^{\delta}}$ is another condensed abelian group with same underlying group as $\underline{\bb{R}}$. There is an inclusion $\underline{\bb{R}}^{\delta}\subset \underline{\bb{R}}$ which is not an isomorphism. Indeed, for a light profinite set $S$ we have 
\[
\underline{\bb{R}}/\underline{\bb{R}}^{\delta} (S) = C(S,\bb{R})/C^{lc}(S,\bb{R}),
\]
where $C^{lc}(S,\bb{R})$ is the space of locally constant functions from $S$ to $\bb{R}$. 

\end{enumerate}
\end{example}

\begin{theorem}
\label{TheoStructureCondAb}
The category $\CondAb^{\light}$ is a Grothendieck abelian category endowed with a natural symmetric monoidal structure and an internal $\Hom$.  Moreover, it has the following properties 
\begin{enumerate}

\item Countable products are exact (countable AB4*) and satisfy countable (AB6).

\item Sequential limits of surjective maps are surjective. 

\item The object $\bb{Z}[\bb{N}\sqcup \{\infty\}]$ is internally projective. 

\end{enumerate}
\end{theorem}
\begin{proof}
The fact that $\CondAb^{\light}$ is a Grothendieck abelian category is a general fact about sheaves on abelian groups on a (essentially small) site. It also has a natural tensor product given by the sheafification of  the tensor product of presheaves.  The internal $\Hom$ is just the right adjoint of the tensor product.  Points (1) and (2) hold from the fact that the topos of condensed sets is replete, and the analogue AB4* and AB6 properties for abelian groups (see \cite[Theorem 2.2]{ClausenScholzeCondensed2019} for the non-light version of (1)). 

We now prove (3).  It suffices  to show that the space of null sequences $P=\bb{Z}[\bb{N}\cup\{\infty\}]/(\infty )$  is internally projective. We want to show that for a surjection $A\to B$ of light condensed abelian groups, a light profinite set $S$, and a map $g:\bb{Z}[S]\otimes N\to B$,  there is a dashed arrow making the following diagram commutative
\[
\begin{tikzcd}
& A \ar[d] \\ 
\bb{Z}[S]\otimes P\ar[r,"g"] \ar[ru,dashed]& B
\end{tikzcd}
\]
after possibly  replacing $S$ by a cover. We have that $\bb{Z}[S]\otimes P= \bb{Z}[S\times (\bb{N}\times \{\infty\})]/(\bb{Z}[S\times \{\infty\}])$. Then the map $g$ is the same as a map $S\times(\bb{N}\times \{\infty\}) \to B$ sending $S\times \{\infty\}$ to $0$.  By hypothesis, there is a surjection $f:S'\to S\times (\bb{N}\cup\{\infty\})$ and a map $S'\to A$ lifting $g$.   For $n\in \bb{N}$ let $S_{n}'$ be the fiber over $S\times \{n\}$ (which is still a surjection).   By  \Cref{PropInjectivelight} we can find retractions $r_n:S'\to S'_{n}\subset S'$, and construct the following diagram of locally profinite sets
\[
\begin{tikzcd}
S'  \times \bb{N}  \ar[r, "\bigsqcup_{n} r_n"]  \ar[rd,"\bigsqcup_{n} f\circ r_n"'] & S'  \ar[d,"f"] \\
 &  S\times (\bb{N}\cup \infty).
\end{tikzcd}
\]
We can find a light profinite compactification $S''$ of $S'\times \bb{N}$ such that $S\times \bb{N} \to S'$ extends to  $S''\to S'$  (Exercise, construct one of such compactifications).   Let $D$ be the boundary of $S''$, by  \Cref{PropInjectivelight} we can find another retraction $r:S''\to D$. Let $h:S''\to S'\to A$ be the composite map, then $h-h\circ r$ induces a map  
\[
\bb{Z}[S'']/\bb{Z}[D]=\bb{Z}[S']\otimes P\to A
\]
that lifts $g$ proving what we wanted. 
\end{proof}

\begin{remark}
\label{RemarkOnZNull}
It is surprising that the object $\bb{Z}[\bb{N}\cup \{\infty\}]$ is internally projective in $\CondAb^{\light}$. This does not happens at the level of profinite sets, for example the map $(2\bb{N}\cup\{\infty\})\bigsqcup (2\bb{N}+1 \cup\{\infty\})\to \bb{N}\cup \{\infty\} $ does not admit a split.  This condensed abelian group will be key in the construction of examples on analytic rings. 
\end{remark}

The following lemma shows that the tensor product of condensed abelian groups is nothing but a condensed enhancement of the tensor product on underlying abelian groups. 

\begin{lemma}\label{Lemma:evaluationPointSymmetricMonoidal}
The evaluation map $\ev(*)\colon \CondAb\to \Ab$ is symmetric monoidal. In other words, if $A$ and $B$ are condensed abelian groups then the natural map $A(*)\otimes B(*) \xrightarrow{\sim} (A\otimes B)(*)$ is an isomorphism.
\end{lemma}
\begin{proof}
Since the point $*$ is a compact projective profinite set, the map $\ev(*)$ sends surjections of condensed abelian groups to surjections of abelian groups, as well as filtered colimits. Therefore, by writting a condensed abelian group as a quotient of filtered colimits of free condensed abelian groups on light profinite sets, it suffices to show that if $S_1$ and $S_2$ are light profinite sets then the natural map 
\begin{equation}\label{eqwo0jaownroawr}
\Z[S_1](*)\otimes \Z[S_2](*)\xrightarrow{\sim} \Z[S_1\times S_2](*)
\end{equation}
is an isomorphism. By \cite[Proposition 2.1]{ClauseScholzeAnalyticGeometry} we have, for $S=\varprojlim_i S_i$ a light profinite set written as a limit of finite sets, the presentation as condensed set
\[
\Z[S]=\bigcup_{n} \varprojlim_i \Z[S_i]_{\leq n}
\]
where $\Z[S_i]_{\leq n}\subset \Z[S_i]$ is the finite subset of sums $\sum_{s\in S_i} a_s s$ with $\sum_{s\in S_i} |a_{s}|\leq n$.   From this presentation one sees directly that $\Z[S](*) = \Z[S(*)]$ is the free abelian group generated by the underlying set of $S$. This implies the isomorphism of \eqref{eqwo0jaownroawr} proving what we wanted. 
\end{proof}

We can define the condensed cohomology as follows:

\begin{definition}
Let $T\in \CondSet^{\light}$ be a light condensed set and $M$ a discrete abelian group, we define the \textit{condensed cohomology of $T$ with values in $M$} to be 
\[
R\Gamma_{\cond}(T, M):= R\Hom(\bb{Z}[T],M). 
\]

\end{definition}

Condensed cohomology behaves as expected in good cases. 

\begin{prop}[{\cite[Theorem 3.2]{ClausenScholzeCondensed2019}}]
\label{PropCohomologyProfinites}
Let $S$ be a profinite set and $M$ a discrete abelian group, then 
\[
R\Gamma_{\cond}(S,M)= C(S,M)
\]
is the space of continuous (eq. locally constant) functions from $S$ to $M$. 
\end{prop}
\begin{proof}
It is clear that  $H^0_{\cond}(S,M)$ is just the space of continuous maps from $S$ to $M$. To show that the higher cohomology groups vanish, it suffices to show that for a cover $S'\to S$ with \v{C}ech nerve $(S^{',\times_S {n_1}})_{[n]\in \Delta^{\op}}$ the \v{C}ech cohomology complex
\begin{equation}
\label{eqCheck}
0\to C(S',M)\to C(S'\times_S S',M) \to \cdots
\end{equation}
is acyclic in cohomological degrees $\geq 1$. For this, we can write the surjection $S'\to S$ as a sequential limit of finite sets with surjective maps $\varprojlim_n (S_n'\to S_n)$. Then the \v{C}ech complex \eqref{eqCheck} is the colimit of the \v{C}ech complexes of the surjections $S_n'\to S_n$, which are acyclic in degrees $\geq 1$ since any surjection of finite sets splits. 
\end{proof}

\begin{prop}[{\cite[Theorem 3.2]{ClausenScholzeCondensed2019}}]
\label{PropCondensedAsSheafCoho}
Let $X$ be a light compact Hausdorff space and $M$ a discrete abelian group, then there is a natural isomorphism
\[
R\Gamma_{\cond}(X,M)= R\Gamma(X,M)
\]
between condensed  and \v{C}ech cohomology. 
\end{prop}
\begin{proof}
Since $X$ is compact Hausdorff we can formally reduce to the case $M=\bb{Z}$. Let $X_{\Prof}:= \Prof^{\light}_{/X}$ be the site of light profinite sets over $X$. Then condensed cohomology of $X$ is the same as the cohomology in $X_{\Prof}$. Let $X_{\mathrm{top}}$ be the site consisting on  closed subspaces of $X$ with coverings given by finite unions of  closed subspaces admitting a refinement by an open cover. Then \v{C}ech cohomology of $X$ is the same as the cohomology on $X_{\mathrm{top}}$. We have a natural morphism of sites
\[
\eta:X_{\Prof}\to X_{\mathrm{top}}.
\] 
It suffices to show that the natural map $\bb{Z}\to R\eta_* \bb{Z}$ is an isomorphism. This can be proved at  stalks, so let $x\in X$, then the stalk $R\eta_*\bb{Z}|_x$ is the same as the pushforward of the fiber over $x$, which is nothing but the condensed cohomology of a point which is $\bb{Z}$.
\end{proof}


\section{Light solid abelian groups}
\label{SectionLightSolid}

The theory of solid abelian groups was introduced in \cite{ClausenScholzeCondensed2019}, it plays a fundamental role in non-archimedean analytic geometries and non-archimedean analysis. The category $\Solid$ of solid abelian groups is a full subcategory of $\CondAb$, stable under limits, colimits and extensions, and containing $\bb{Z}$; it is  actually  the smallest category satisfying those properties. In its ``classical construction'' \footnote{If we are allowed to call classical a construction just made around 2018-2019; five to six years ago compared to the time these notes are being written.}, the theory of locally compact abelian groups and its extensions as condensed abelian groups play a key role. However, within the new framework of light condensed mathematics, the theory of solid abelian groups can be formally developed from the more intuitive idea that the ``summable sequences'' in non-archimedean analysis are precisely the ``null-sequences''. In the following we will explain how this very simple idea naturally guides us to the correct definition of Solid.

\subsection{Null-sequences and summability}
\label{SubsecNullseq}

Let $K$ be a local field and $V$ a Banach space over $K$. Recall that a \textit{null-sequence} in $V$ is a sequence $(v_n)_{n\in \bb{N}}$ converging to $0$. Similarly, a \textit{summable sequence} is a sequence $(v_n)_{n\in \bb{N}}$ such that the partial sums $\sum_{i=0}^n v_i$ converge to an element in  $v$ that we denote by $\sum_{n} v_n$. One of the first properties   that we learn in a course of analysis is that a summable sequence $(v_n)$ has tails $w_n=\sum_{i\geq n} v_n$ converging to $0$. In other words, we have a map 
\[
\{\mbox{summable sequences}\} \to \{\mbox{null sequences}\}: (v_n)\mapsto (w_n). 
\]
On the other hand, given a null sequence $\{w_n\}_{n\in \bb{N}}$ we can form the sequence $x_n:= w_{n}-w_{n+1}$ which turns out to be summable in $V$, namely, 
\[
v_n:=\sum_{i=0}^n x_n = w_0-w_{n+1}
\]
and $(v_n)_{n}$ converges to $w_0$ as $n\to \infty$.  Thus, we get a bijection 
\[
\{\mbox{null sequences}\} \to \{\mbox{summable sequences }\} : (w_n)_n\mapsto (x_n)_{n}=(w_{n}-w_{n+1}).
\]
In general, any summable sequence in $V$ is also a null-sequence. The converse does not hold for archimedean fields (eg. $(1/n)_{n}$), but it does hold for non-archimedean fields thanks to the ultrametric inequality. 

 Therefore, a way to isolate non-archimedean analysis from condensed abelian groups is by asking that any null-sequence is summable, namely, that the map 
\[
1-S:\{\mbox{null sequences}\} \to \{\mbox{null sequences}\},
\]
where $S$ is the shift map $(v_n)\mapsto (v_{n+1})$, is a bijection. 

In order to formalize this idea, first we need to be able to talk about null-sequences of  condensed abelian groups.

\begin{definition}
\label{DefinitionNullSequence}
We let $P:=\bb{Z}[\bb{N}\cup\{\infty \}]/ (\infty)$. Given a condensed abelian group $A$ its \textit{space of null sequences} is given by $\Null(A)=\Hom(P,A)$, we also let $\iNull(A):= \iHom(P,A)$. 
\end{definition}

\begin{example}
\label{ExampleNullExtraDatum}
We continue in the spirit of  \Cref{ExampleToptoCond} (1). For a quasi-separated condensed abelian group $A$ being a null-sequence is an actual property of the underlying sequence, namely, the map 
\[
\Null(A)\to \Map(\bb{N}, A) = \prod_{\bb{N}} A(*)
\] 
is injective. However, for general condensed abelian groups null-sequences are not properties on sequences but additional structure you put in the condensed abelian group. As example, let $\bb{R}$ be the real numbers with the usual topology, and let $\bb{R}^{\delta}$ be the real numbers with the discrete topology. Then $\bb{R}/\bb{R}^{\delta}$, if  scary as topological abelian group, it is a well defined condensed abelian group, and for any light profinite set $S$ we have that
\[
\bb{R}/ \bb{R}^{\delta}(S) =C(S,\bb{R})/C^{lc}(S,\bb{R})
\]
is the quotient of continuous maps from $S\to \bb{R}$ modulo locally constant maps from $S$ to $\bb{R}$. Applying this to $S=\bb{N}\cup\{\infty\}$ we get that $\bb{R}/\bb{R}^{\delta}(S)$ is a non-zero space of null-sequences  while $\bb{R}/\bb{R}^{\delta}(*)=0$, this shows that a null-sequence in that non quasi-separated quotient remembers the tails of the virtually zero sequences. 
\end{example}

An additional feature of $P$ is that it has a natural structure of algebra making $\bb{Z}[T]=\bb{Z}[\bb{N}]\to P$ an algebra morphism. 

\begin{prop}
\label{PropNullSequenceAlgebra}
The  addition map
\[
\bb{N}\times \bb{N}\to \bb{N}
\]
induces an algebra structure on $P$, we shall denote this algebra by $\bb{Z}[\widehat{q}]$. 
\end{prop}

To prove   \Cref{PropNullSequenceAlgebra}, it will suffices to show the following lemma

\begin{lemma}
\label{LemmaCompactificationProfinites}
Consider a surjective map of light profinite sets $S\to S'$ and let $U\subset S'$ be an open subspace such that $S\times_{S'} U \to U$ is an homeomorphism. Let $D$ and $D'$ be the complements of $U$ in $S$ and $S'$ respectively. Then we have a pushout square in $\CondSet$
\[
\begin{tikzcd}
D\ar[r] \ar[d] & S \ar[d]\\
D' \ar[r] & S'.
\end{tikzcd}
\] 
\end{lemma}
\begin{proof}
We have a surjection of condensed sets $S\to S'$  whose \v{C}ech fiber is given by  $S\times_{S'} S = \Delta S \cup (D\times_{D'} D)\subset S\times S$.  Since $S\to S'$ is surjective, we have that $S'=S/(S\times_{S'} S)= S/(\Delta S \cup (D\times_{D'} D))$, which is exactly the pushout $S\bigsqcup_{D} D'$. 
\end{proof}

\begin{definition}
\label{DefinitionNullmeasures}
Let $U$ be a light locally profinite set, i.e. a countable disjoint union of light profinite sets. We let $P_{U}:= \bb{Z}[U\cup \{\infty \}]/(\infty)$ be the space of \textit{measures on $U$ vanishing at $\infty$}.  
\end{definition}

\begin{prop}
\label{PropVanishingMeasuresInvariance}
Let $U$ be a light locally profinite set, let $S$ be any  compactification of $U$ in a light profinite set and let $D$ be the boundary, then there is a natural isomorphism of condensed abelian groups $P_U=\bb{Z}[S]/\bb{Z}[D]$.
\end{prop}
\begin{proof}
By \Cref{LemmaCompactificationProfinites} we have a pushout diagram
\[
\begin{tikzcd}
D\ar[r]\ar[d] & S\ar[d] \\
* \ar[r]& U\cup\{\infty \},
\end{tikzcd}
\]
applying the left adjoint $\bb{Z}[-]$ we get a push out diagram at the level of free abelian groups, which induces the isomorphism 
\[
\bb{Z}[S]/\bb{Z}[D]=\bb{Z}[U\cup\{\infty\}]/(\infty). 
\]
\end{proof}

\begin{proof}[Proof of  \Cref{PropNullSequenceAlgebra}]

We can endow $\bb{N}\cup\{\infty\}$ with a structure of additive monoid by declaring $\infty+a=\infty$. Then, $\bb{Z}[\bb{N}\cup \{\infty\}]$ has a natural algebra structure such that $\bb{Z}[\infty]$ is an ideal, this endows $P$ with an algebra structure. More explicitly, consider the addition map
\[
(\bb{N}\cup\{\infty\})\times (\bb{N}\cup\{\infty\})\to \bb{N}\cup\{\infty\},
\] 
it sends the boundary of $(\bb{N}\cup\{\infty\})\times( \bb{N}\cup\{\infty\})$ to the boundary of $\bb{N}\cup\{\infty\}$, and by    \Cref{PropVanishingMeasuresInvariance}  it defines a map 
\[
P\otimes P\to P,
\]
compatible with the multiplication map $\bb{Z}[T]\otimes \bb{Z}[T]\to \bb{Z}[T]$. It is easy to check that this defines an algebra structure on $P$. 
\end{proof}

\subsection{Solid abelian groups form an analytic ring}
\label{SubsecSolidAb}

Now we define the category of solid abelian groups, for this, note that the condensed abelian group $P$ parametrizing null sequences has an endomorphism $\Shift:P\to P$ which is induced from the map of profinite sets $\bb{N}\cup \{\infty\}\to \bb{N}\cup \{\infty\}$ mapping $\infty$ to $\infty$ and $n$ to $n+1$, we call $\Shift$ the \textit{shift map}. 

\begin{definition}
\label{DefinitionSolidAb}
Consider the map $1-\Shift:P\to P$.  A light condensed abelian group $A$  is said \textit{solid} if the natural map
\[
\iHom(P,A) \xrightarrow{1-\Shift^*} \iHom(P,A)
\]
is an isomorphism. We let $\Solid\subset \CondAb^{\light}$ denote the full subcategory of (light) solid  abelian groups. 

More generally, given $C\in \ob{D}(\CondAb^{\light})$ an object in the derived ($\infty$-)category of condensed abelian groups, we say that $C$ if solid if the natural map
\[
R\iHom(P,C)\xrightarrow{1-\Shift^*} R\iHom(P,C)
\]
is an equivalence. We let $\ob{D}(\CondAb)^{\sol}\subset \ob{D}(\CondAb^{\light})$ be the full subcategory of solid objects. 
\end{definition}

\begin{remark}
By  \Cref{TheoStructureCondAb}, the object $P$ is internally projective in the category of light condensed abelian groups, in particular there is no difference between the derived or non derived Hom space $\iHom(P,A)$. This shows that  $\Solid \subset \ob{D}(\CondAb)^{\sol}$. 
\end{remark}

\begin{remark}\label{RemDefSolidAb2}
Let $\Z[\widehat{q}]$ be the object $P$ seen as a condensed ring. The shift map corresponds to the multiplication by $q$,  and $\Z[\widehat{q}]$ is clearly an integral domain. Thus, a condensed abelian group $A$ is solid if and only if the derived internal Hom space $R\iHom(\Z[\widehat{q}]/(1-q),A)=0$ vanishes. 
\end{remark}

The main theorem regarding the category of solid abelian groups is the following:

\begin{theorem}
\label{TheoSolidAb}
The category $\Solid$ is a Grothendieck abelian category stable under limits, colimits and extensions in $\CondAb$. Furthermore, the following properties hold: 
\begin{enumerate}

\item $\bb{Z}\in \Solid$.

\item There is a left adjoint $(-)^{\sol}:\CondAb\to \Solid$ for the inclusion that we call the solidification functor.

\item  There is a unique symmetric monoidal structure $\otimes_{\sol}$ on $\Solid$ making $(-)^{\sol}$ symmetric monoidal. 

\item $\bb{R}^{\sol}=0$ (solid abelian groups kill the archimedean theory).

\end{enumerate}

Moreover, $\ob{D}(\CondAb)^{\sol}$ is a presentable full subcategory of $\ob{D}(\CondAb)$ stable under limits and colimits, and the following properties are satisfied:

\begin{enumerate}
\setcounter{enumi}{4}
\item  The inclusion $\ob{D}(\CondAb)^{\sol}\to \ob{D}(\CondAb)$ has a left adjoint $(-)^{L\sol}$.

\item An object $C\in \ob{D}(\CondAb)$ is solid if and only if $H^i(C)\in \Solid$ for all $i\in \bb{Z}$, i.e. the natural $t$-structure on $\ob{D}(\CondAb)$ induces a $t$-structure on $\ob{D}(\CondAb)^{\sol}$.

\item For $C\in \ob{D}(\CondAb)^{\sol}$ and $M\in  \ob{D}(\CondAb)$ we have $R\iHom(M,C)\in \ob{D}(\CondAb)^{\sol}$. 

\item The category $\ob{D}(\CondAb)^{\sol}$ has a unique symmetric monoidal structure $\otimes^L_{\sol}$ making $(-)^{L\sol}$ symmetric monoidal.

\item The natural map $\ob{D}(\Solid)\to \ob{D}(\CondAb)$ of derived categories is fully faithful, and has essential image  $\ob{D}(\CondAb)^{\sol}$.

\item The functor $(-)^{L\sol}$ is the left derived functor of $(-)^{\sol}$.

\item The functor $\otimes^L_{\sol}$ is the left derived functor of $\otimes_{\sol}$.

\item For $S=\varprojlim_{n}S_n$ a light profinite set, there is a natural equivalence
\[
\bb{Z}_{\sol}[S]:=(\bb{Z}[S])^{L\sol} \xrightarrow{\sim} \varprojlim_{n} \bb{Z}[S_n] \cong \prod_{J} \bb{Z}
\]
with $J$ a countable set. In particular,  $\bb{Z}_{\sol}[S]$ is a compact projective solid abelian group, and if $S$ is infinite $\bb{Z}_{\sol}[S]$ is a compact projective generator of $\Solid$. 

\item For $I$ and $J$ countable sets we have 
\[
\prod_{I} \bb{Z} \otimes^L_{\sol} \prod_{J} \bb{Z} = \prod_{I\times J} \bb{Z}. 
\]

\item The object $\prod_{\bb{N}}\ \bb{Z}$ is flat in $\Solid$.

\end{enumerate}

\end{theorem}

In \cite{ClausenScholzeCondensed2019} a lot of effort is made in order to prove  \Cref{TheoSolidAb} and the only  obvious property was point (12), this is because solid abelian groups were constructed by first defining the functor of measures $S\mapsto \bb{Z}_{\sol}[S]$.  Furthermore, property (14) is not true in arbitrary solid abelian groups (with a counter example due to Effimov).  It turns out that with  \Cref{DefinitionSolidAb} most of the theorem is immediate. 

\begin{prop}
\label{PropFirstPartTheoSolid}
The category $\Solid$ is a Grothendieck abelian  category. Furthermore, points (1)-(8) hold. Moreover,  property (12) implies (9) and (10), and property (13) implies (11). 

\end{prop}
\begin{proof}
Recall that the category $\Solid$ is defined as the full subcategory of condensed abelian groups $A$ such that the map $1-\Shift^*$ on $\iHom(P,A)$ is an isomorphism. Since $P$ is internally projective, this condition is clearly stable under limits, colimits and extensions in $\CondAb$, making $\Solid$  an abelian category.   The same argument shows that $\ob{D}(\CondAb)^{\sol}$ is stable under limits and colimits in $\ob{D}(\CondAb)$. It is left to show that $\Solid$ and $\ob{D}(\CondAb)^{\square}$ are presentable, for this, consider $Q=\mathrm{cone}(P\to \varinjlim_{1-\Shift} P)$, then an object $C$ is solid if and only if $R\iHom(Q,C)=0$. Presentability then follows from \cite[Theorem 5.5.3.18]{HigherTopos} (one could also just evoke the general reflection principle of \cite{ReflectionPrinciple} as these categories are stable under limits and colimits in the condensed categories which are presentable).

\begin{itemize}

\item[(1)] By  \Cref{PropCohomologyProfinites}, for all $S\in \Prof$ we have that $R\iHom(\bb{Z}[S], \bb{Z})=C(S,\bb{Z})$ is the space of locally constant functions. This implies that 
\[
\iHom(P,\bb{Z})=\bigoplus_{n\in \bb{N}} \bb{Z}.
\]
Then, the action of $1-\Shift^*$ maps a sequence $(a_0,a_1,\ldots)$ to $(a_0-a_1,a_1-a_1,\ldots)$, which clearly has by inverse 
\[
(b_0,b_1,b_2,\ldots)\mapsto (\sum_{i\geq 0} b_i, \sum_{i\geq 1} b_i,\ldots).
\]
since the sequences are eventually zero.

\item[(2) and (5)] The existence of the left adjoint follows from the adjoint functor theorem \cite[Corollary 5.5.2.9]{HigherTopos}.

\item[(3) and (8)]  It suffices to show that the kernel of the adjoints $(-)^{\sol}$ and $(-)^{L \sol}$ are tensor ideals in $\Solid$ and $\ob{D}(\CondAb)^{\sol}$ respectively. Let us just explain the proof for $(-)^{L\sol}$. Let $A\in \ob{D}(\CondAb)$ be such that $A^{L\sol}=0$ and let $M\in \ob{D}(\CondAb)$. To prove that $(M\otimes^L A)^{L\sol}=0$ it suffices to show that for all $B\in \ob{D}(\CondAb)^{\sol}$ we have 
\[
R\Hom(A\otimes^L M,B)=0,
\]
but we  have that 
\begin{equation}
\label{equationadjointtensor}
R\Hom(A \otimes^L M,B)=R\Hom(A,R\iHom(M,B)),
\end{equation}
and $R\iHom(M,B)$ is solid by (7) down below. This shows that  \eqref{equationadjointtensor} vanishes.

\item[(4)] Since $\bb{R}$ is an algebra and the functor $(-)^{L\sol}$ is symmetric monoidal, it suffices to show that 
\[
\pi_0(\bb{R}^{L\sol})=\bb{R}^{\sol}=0. 
\]
It suffices to show that the unit map $\bb{Z}\to \bb{R}^{\sol}$ is zero. For this, consider the null-sequence in $\bb{R}$
\[
(1,\frac{1}{2},\frac{1}{2},\frac{1}{4},\frac{1}{4},\frac{1}{4},\frac{1}{4},\cdots)
\]
defining a map $f:P\to \bb{R}$. By definition of the solidification, there is a unique map $g:P\to \bb{R}^{\sol}$ making the following diagram commutative
\[
\begin{tikzcd}
P \ar[r,"f"] \ar[d,"1-\Shift"] & \bb{R}  \ar[d] \\
P   \ar[r,"g"] & \bb{R}^{\sol}.
\end{tikzcd}
\]
Let $[0]:\bb{Z}\to P$ be the inclusion in the $0$-th component, then $g\circ [0]:\bb{Z}\to \bb{R}^{\sol}$ defines an element $x_0$ (virtually given by $1+\frac{1}{2}+\frac{1}{2}+\cdots$).  We claim that $x_0=1+x_0$, this would show that $1=0$ and that $\bb{R}^{\sol}=0$.

Consider the maps 
\[
F: \bb{Z}[\bb{N}]\to \bb{Z}[\bb{N}] : [n]\mapsto [2n+1]+[2n+2]
\]
\[
G:\bb{Z}[\bb{N}] \to \bb{Z}[\bb{N} ] : [n]\mapsto [2n+1].
\]
These maps naturally extend to endomorphisms of $P$. We claim that we have a commutative diagram 
\[
\begin{tikzcd}
P \ar[r, "F"] \ar[d,"1-\Shift"] \ar[d, "1-\Shift"] & P \ar[d,"1-\Shift"] \\ 
P \ar[r, "G"] & P,
\end{tikzcd}
\]
namely, we have 
\[
(1-\Shift)\circ F ([n])= (1-\Shift)([2n+1]+[2n+2])=[2n+1]-[2n+2]+[2n+2]-[2n+3]=[2n+1]-[2n+3]
\]
and 
\[
G\circ (1-\Shift) ([n]) = G([n]-[n+1])=[2n+1]- [2n+3]. 
\]
On the other hand, we have that $f\circ F=f$, namely it is the sequence
\[
((\frac{1}{2}+\frac{1}{2}), (\frac{1}{4}+\frac{1}{4}),(\frac{1}{4}+\frac{1}{4}),( \frac{1}{8}+\frac{1}{8}),(\frac{1}{8}+\frac{1}{8}),(\frac{1}{8}+\frac{1}{8}),(\frac{1}{8}+\frac{1}{8}),\cdots)=(1,\frac{1}{2},\frac{1}{2},\frac{1}{4},\frac{1}{4},\frac{1}{4},\frac{1}{4},\cdots).
\]
By uniqueness of the lift $g:P\to \bb{R}^{\sol}$, we must have $g\circ G=g$. Then, if $g$ represents the null sequence $(x_0,x_1,x_2,x_3,\cdots)$, we must have $x_n=x_{2n+1}$ for all $n\in \bb{Z}$. In particular, $x_0=x_1$, so that 
\[
0=x_0-x_1= 1,
\]
proving what we wanted.

\item[(6)] This follows from the fact that for all $C\in \ob{D}(\CondAb)$ we have 
\[
H^i(R\iHom(P,C))= \iHom(P,H^i(C)) \mbox{ for }i\in \bb{Z}
\]
since $P$ is internally projective. 

\item[(7)] Let $M\in \ob{D}(\CondAb)$ and $C\in \ob{D}(\CondAb)^{\sol}$, then the claim follows from the isomorphism
\[
R\iHom(P,R\iHom(M,C))= R\iHom(M,R\iHom(P,B)),
\] 
and the fact that $B$ is solid.

\end{itemize}

Now let us assume that properties (11) and (12) hold. 

\begin{itemize}

\item[(9)] The map $\Solid\to \CondAb$ induces a functor of derived categories $\ob{D}(\Solid)\to \ob{D}(\CondAb)$, by \cite[Proposition 1.3.3.7]{HigherAlgebra}, and since $P^{\sol}=\prod_{\bb{N}}\bb{Z}$ is a compact projective generator of $\Solid$,  it suffices to show that for $A\in \Solid$  we have 
\[
R\iHom(P^{\sol}, A)= \iHom(P^{\sol},A).
\]
But we know that $P^{\sol}=P^{L\sol}$, and by the left adjoints of (2) and (5) we have 
\[
R\iHom(P^{\sol},A)=R\iHom(P^{L\sol},A)=R\iHom(P,A)=\iHom(P,A)=\iHom(P^{\sol},A).
\]

\item[(10)] This follows from the fact that $\bb{Z}[S]^{L\sol}=\bb{Z}[S]^{\sol}$ sits in degree zero. Indeed, since both derived categories are right complete, it suffices to show that the restriction of $(-)^{L\sol}$ to  connective objects $\ob{D}_{\geq 0}(\CondAb)$ (i.e. objects supported in  non-negative homological degrees) is the left derived functor. This statement boils down to the fact that $(-)^{L\sol}:\ob{D}_{\geq 0}(\CondAb)\to \ob{D}_{\geq 0}(\Solid)$ is the left Kan extension of its restriction to the full subcategory of generators $\s{C}^{0}=\{\bb{Z}[S]\}_{S\in \Prof^{\light}}\subset \ob{D}_{\geq 0}(\CondAb)$ \footnote{Note that the full subcategory $\s{C}^0\subset \ob{D}_{\geq 0}(\CondAb)$ is not a full subcategory of  $\CondAb$ since the objects of $\s{C}^0$ are not necessarily projective}. In other words, that for $C\in \ob{D}_{\geq 0}(\CondAb)$ we have 
\[
C^{L\sol}= \varinjlim_{\bb{Z}[S]\in \s{C}^0/C} \bb{Z}[S]^{L\sol}= \varinjlim_{\bb{Z}[S]\in \s{C}^0/C} \bb{Z}[S]^{\sol}. 
\]

\item[(11)] Finally, to show that $\otimes^{L}_{\sol}$ is the left derived functor of $\otimes_{\sol}$, it sufices to show that there is a family of compact  projective generators $\s{C}^0\subset \Solid$ stable under the solid tensor product, such that for $A,B\in \s{C}^0$ we have $A\otimes^{L}_{\sol} B=A\otimes_{\sol} B$. Taking $\s{C}^0$ as the full subcategory spanned by $\bb{Z}_{\sol}[S]$ with $S$ light profinite we are done thanks to property (13).

\end{itemize}

\end{proof}

\begin{corollary}
\label{coroRdies}
Let $C$ be a real condensed vector space. Then $C^{L\sol}=0$. 
\end{corollary}
\begin{proof}
The solidification functor $(-)^{L\sol}$ is symmetric monoidal, in particular $\bb{R}^{L\sol}$ is an algebra and  $C^{L\sol}$ has a natural $\bb{R}^{L\sol}$ -module structure. But $\bb{R}^{L\sol}=0$, which implies that $C^{L\sol}=0$.
\end{proof}

We have proven most of  \Cref{TheoSolidAb}, it is left to show points (12)-(14) regarding the explicit description of the free objects $\bb{Z}_{\sol}[S]:=\bb{Z}[S]^{L\sol}$, their solid tensor products, and the flatness of $\prod_{\bb{N}} \bb{Z}$ in $\Solid$, we left those properties for the next sections.

\subsection{Computing measures in solid abelian groups}
\label{SubsecMeasuresSolid}

The objective in this section is to prove the following theorem:

\begin{theorem}
\label{TheoSolidMeasures}
Let $S=\varprojlim_n S_n$ be a light profinite set. Then the natural map of solid abelian complexes
\[
\bb{Z}[S]^{L\sol} \to \varprojlim_n \bb{Z}[S_n]
\]
is an equivalence. Furthermore, the following hold:
\begin{enumerate}

\item $\prod_{\bb{N}}\bb{Z}$ is a compact projective generator of $\Solid$. 

\item The natural map $\ob{D}(\Solid)\to \ob{D}(\CondAb)^{\sol}$ is an equivalence of $\infty$-categories.

\item The functor $(-)^{L\sol}$ is the left derived functor of $(-)^{\sol}$.

\end{enumerate}

\end{theorem}

By  \Cref{PropFirstPartTheoSolid}, we only need to prove the first assertion of the theorem, this will require some lemmas.  Recall that  $P=\bb{Z}[\bb{N}\cup\{\infty\}]/(\infty)$ is the solid abelian group parametrizing null-sequences.

First, we note that it suffices to compute the solidification of $P$ in order to compute the solidification of $\bb{Z}[S]$ for $S$ a light profinite set:

\begin{lemma}
\label{LemmaReductionSolidtoP}
Let $S$ be a light profinite set, there is a map $P\to \bb{Z}[S]$ that induces isomorphisms on solidifications 
\[
P^{L\sol}\xrightarrow{\sim} \bb{Z}[S]^{L\sol}.
\]
\end{lemma}
\begin{proof}
Let us write $S=\varprojlim_{n} S_n$ as a limit of finite sets with surjective transition maps and projections $\pi_n:S\to S_n$. We can find a sequence of compatible lifts $S_0\to S_1\to S_2\to \cdots \to S$ with $\iota_n:S_n\to S$. Enumerating $\bigcup_n \iota_{n}(S_n)\cong \bb{N}$ along the previous inclusions, we get an injection $\bb{N}\to S$. Then for $a\in \iota_n(S_n)\backslash \iota_{n-1}(S_{n-1}) \subset \bb{N}$  consider the element $\iota_{n}(a)-\iota_{n-1}(a)$. The sequence $(\iota_{n}(a)-\iota_{n-1}(a))_{n\in \bb{N}}$ converges to zero in $\bb{Z}[S]$ and defines an injective map $g:P\to \bb{Z}[S]$.  We claim that $g$ induces an isomorphism after solidification.  

We claim that we have a commutative diagram 
\begin{equation}
\label{eqSolidDiagramMeasures}
\begin{tikzcd}
P\otimes \bb{Z}[S] \ar[r,"F"] & \bb{Z}[S] \\ 
P \otimes \bb{Z}[S]  \ar[u, " (1-\Shift)\otimes \id_{S}"] \ar[r,"G"]& P   \ar[u]
\end{tikzcd}
\end{equation}
where the top horizontal arrow $F$ arises from a map $(\bb{N}\times \{\infty \})\times S\to \bb{Z}[S]$ that vanishes at $\infty \times S$. This map is given by the sequence of maps $\{n\}\times S\to \bb{Z}[S]$ given  by $\id_{\bb{Z}[S]}$ if $n=0$ and $\id_{\bb{Z}[S]}- \iota_{n-1}\circ \pi_{n-1}$  if $n\geq 1$, which vanish uniformly on $S$ at $\infty$. Then, to define the lower horizontal arrow $G$ we need to show that the composite $F\circ (1-\Shift)$ lands in $P$, but the composite corresponds to the map of condensed sets
\[
G:(\bb{N}\cup\{\infty\})\times S \to \bb{Z}[S]
\]
vanishing at $\infty\times S$, and given by $\iota_{n}\circ \pi_{n}- \iota_{n-1}\circ \pi_{n-1} : S\to \bb{Z}[S]$ on $\{n\}\times S$ (where we make the convention $\iota_{-1}\circ \pi_{-1}=0$). In particular, $G(\{n\}\times S)$ lands in $P$, and so it extends to a map $G: (\bb{N}\cup\{\infty\})\times S \to P$ that vanishes at $\{\infty\}\times S$, producing the desired factorization. 

Taking solidifications of \eqref{eqSolidDiagramMeasures}, we get a commutative diagram

\begin{equation}
\label{eqSolidDiagramMeasures2}
\begin{tikzcd}
(P\otimes \bb{Z}[S])^{L\sol} \ar[r,"F"] & \bb{Z}[S]^{L\sol}  \ar[l, bend right, dashed ]\\ 
(P \otimes \bb{Z}[S])^{L\sol}  \ar[u, "\wr"] \ar[r,"G"]& P^{L\sol}   \ar[u]
\end{tikzcd}
\end{equation}
where the left vertical arrow is an isomorphism, and the top horizontal arrow has a section induced from the map $\{0\}\times S \to P\otimes \bb{Z}[S]$. The previous shows that $\bb{Z}[S]^{L\sol}$ is a retract of $P^{L\sol}$ with idempotent morphism $r:P^{L\sol}\to P^{L\sol}$. To show that the map is an actual isomorphism we need to show that $r$ is the identity.    To prove this last claim, note that the diagram  \eqref{eqSolidDiagramMeasures} restricts to a diagram 
\[
\begin{tikzcd}
P\otimes P \ar[r,"F"] & P \\ 
P\otimes P \ar[u, "(1-\Shift)\otimes \id_P"] \ar[r]& P \ar[u, "\id_P"]
\end{tikzcd}
\]
via the inclusion $P\subset \bb{Z}[S]$. Indeed, the  map $F$ is given by the sequence of endomorphisms $\id_{\bb{Z}[S]}- \iota_{n-1}\circ \pi_{n-1}$ of $\bb{Z}[S]$, which restrict to the  endomorphisms $\id_P-\iota_{n-1}\circ \pi_{n-1}$ of $P$. Taking solidifications we get
\begin{equation}
\label{eqSolidDiagramMeasures3}
\begin{tikzcd}
(P\otimes P)^{L\sol}  \ar[r,"F"] & P^{L\sol}  \ar[l, bend right, dashed ] \\ 
(P\otimes P)^{L\sol} \ar[u, "\wr"] \ar[r]& P^{L\sol} \ar[u, "\id_P"],
\end{tikzcd}
\end{equation}
and the idempotent $r$ obtained from \eqref{eqSolidDiagramMeasures2} is the same as the idempotent obtained from \eqref{eqSolidDiagramMeasures3} which is the identity. 
\end{proof}

Now, we compute the solidification of $P$. We apply the same trick as in the proof of  \Cref{LemmaReductionSolidtoP} to replace $P$ by a simpler condensed abelian group. 

\begin{lemma}
\label{Lemma1Solid}
Let $\prod_{\bb{N}}^{\bnd} \bb{Z}= \bigcup_{n\in \bb{N}}\prod_{\bb{N}}  \bigg(\bb{Z}\cap [-n,n] \subset \prod_{\bb{N}} \bb{Z}\bigg)$ be the condensed abelian group of bounded sequences of integers. Consider the natural map $P\to \prod_{\bb{N}}^{\bnd} \bb{Z}$ induced by the null sequence $e_n\in \prod_{\bb{N}}^{\bnd} \bb{Z}$ with $e_n=(0,0,\cdots,0,1,0,\cdots)$, where $e_n$ has zero entries except for an occurence of $1$ in the $n$-th component. Then the natural map 
\[
P^{L\sol}\to (\prod_{\bb{N}}^{\bnd} \bb{Z})^{L\sol}
\]
is an isomorphism. 
\end{lemma}
\begin{proof}
We claim that there is a commutative square
\begin{equation}
\label{eqrefSolidEquiationBoundedNull}
\begin{tikzcd}
P\otimes \prod_{\bb{N}}^{\bnd} \bb{Z} \ar[r, "F"] & \prod_{\bb{N}}^{\bnd} \bb{Z} \\ 
P\otimes \prod_{\bb{N}}^{\bnd} \bb{Z}  \ar[r,"G"] \ar[u,"(1-\Shift)\otimes \id"]& P \ar[u]
\end{tikzcd}
\end{equation}
where the top horizontal arrow $F$ is given by the null-sequence of endomorphisms of $ \prod_{\bb{N}}^{\bnd} \bb{Z}$ given by the projection $\pi_{\geq n}$ in the $\geq n$-components. To prove the claim, we need to see that the map $G=F\circ (1-\Shift)$ lands in $P$, but it is given by the null-sequence of endomorphisms of $\prod_{\bb{N}}^{\bnd} \bb{Z}$ given by the projections $\pi_n=\pi_{\geq n}-\pi_{\geq n+1}$, whose target is in $P$.  Taking solidifications  of \eqref{eqrefSolidEquiationBoundedNull} we get a commutative diagram 
\[
\begin{tikzcd}
(P\otimes \prod_{\bb{N}}^{\bnd} \bb{Z})^{L\sol} \ar[r, "F"] & (\prod_{\bb{N}}^{\bnd} \bb{Z})^{L\sol} \ar[l, bend right, dashed ] \\ 
(P\otimes \prod_{\bb{N}}^{\bnd} \bb{Z})^{L\sol}  \ar[r,"G"] \ar[u,"\wr"]& P^{L\sol} \ar[u]
\end{tikzcd}
\]
such that the top horizontal arrow has a section given by the embedding in the $0$-th component of the tensor. Then, as in the proof of  \Cref{LemmaReductionSolidtoP}, one gets an idempotent endomorphism $r:P^{L\sol}\to P^{L\sol}$ whose retract is $\bb{Z}[S]^{L\sol}$, and to see that $r$ is the identity, it suffices to notice that   \eqref{eqrefSolidEquiationBoundedNull} restricts to a commutative diagram of the form \eqref{eqSolidDiagramMeasures3}, and then one applies the argument as in the proof of  \Cref{LemmaReductionSolidtoP}.
\end{proof}

\begin{lemma}
\label{LemmaSolidificationBoundedSolid}
The natural map ${\prod_{\bb{N}}}^{\bnd} \bb{Z} \to \prod_{\bb{N}} \bb{Z}$ induces an isomorphism in solidifications 
\[
(\prod_{\bb{N}}^{\bnd} \bb{Z})^{L\sol} =  \prod_{\bb{N}} \bb{Z}.
\]
\end{lemma}
\begin{proof}
Let $\prod_{\bb{N}}^{\bnd}\bb{R}= \bigcup_n \prod_{\bb{N}}\bigg( \bb{R}\cap [-n,n] \bigg)$ be the condensed real vector space. We have isomorphisms of condensed abelian groups
\[
\prod_{\bb{N}} \bb{Z}/\prod_{\bb{N}}^{\bnd} \bb{Z} = \prod_{\bb{N}} \bb{R} /\prod_{\bb{N}}^{\bnd}\bb{R}.
\]
Indeed, this follows from the fact that we have short exact sequences
\[
0\to \prod_{\bb{N}} \bb{Z}\to \prod_{\bb{N}} \bb{R} \to\prod_{\bb{N}} \bb{R}/\bb{Z} \to 0 
\]
and
\[
0\to \prod_{\bb{N}}^{\bnd} \bb{Z}\to \prod_{\bb{N}}^{\bnd}\bb{R}\to  \prod_{\bb{N}} \bb{R}/\bb{Z} \to 0.
\]
In particular, the quotient $\prod_{\bb{N}} \bb{Z}/\prod_{\bb{N}}^{\bnd} \bb{Z}$ can be endowed with an structure of condensed $\bb{R}$-vector space, and so its solidification vanishes by  \Cref{coroRdies}. This shows that 
\[
(\prod_{\bb{N}}^{\bnd} \bb{Z})^{L\sol}=(\prod_{\bb{N}} \bb{Z})^{L\sol}= \prod_{\bb{N}} \bb{Z}
\]
as wanted. 
\end{proof}

\begin{corollary}
\label{CoroDualitySolidMeasures}
Let $S=\varprojlim_n S_n$ be a light profinite set, then we have natural isomorphisms 
\[
\bb{Z}_{\sol}[S]=R\iHom(C(S,\bb{Z}),\bb{Z})
\]
and 
\[
C(S,\bb{Z})=R\iHom(\bb{Z}_{\sol}[S],\bb{Z}).
\]
\end{corollary}
\begin{proof}
The first isomorphism follows from the fact that $C(S,\bb{Z})=\varinjlim_n C(S_n,\bb{Z})$ and that $\bb{Z}_{\sol}[S]=\varprojlim_n \bb{Z}[S_n]$. The second isomorphism follows from  the left adjoint $(-)^{L\sol}$
\[
R\iHom(\bb{Z}_{\sol}[S],\bb{Z}) = R\iHom(\bb{Z}[S],\bb{Z})=C(S,\bb{Z}).
\] 
\end{proof}

\begin{corollary}
\label{CoroSolidProperties}
 \Cref{TheoSolidMeasures} holds. Moreover, for $I$ and $J$ countable sets, the natural map   $\prod_{I} \bb{Z} \otimes^{L}_{\sol} \prod_{J} \bb{Z} \xrightarrow{\sim} \prod_{I\times J} \bb{Z}$ is an isomorphism. In particular, $\otimes^{L}_{\sol}$ is the left derived functor of $\otimes_{\sol}$. 
\end{corollary}
\begin{proof}
The consequences (1)-(3) of the theorem were proven in  \Cref{PropFirstPartTheoSolid}. By  \Cref{LemmaReductionSolidtoP,Lemma1Solid,LemmaSolidificationBoundedSolid}, we know that $\bb{Z}[S]^{\sol} \simeq \prod_{I} \bb{Z}$ (for $I$ a countable set) abstractly as solid abelian groups. Following the explicit isomorphisms constructed in the lemmas, one can verify that the previous isomorphism actually identifies with the natural arrow
\begin{equation}
\label{eqComputationMeasuresSolid}
\bb{Z}[S]^{L\sol}\xrightarrow{\sim} \varprojlim_n \bb{Z}[S_n].
\end{equation}
More explicitly, this holds true for $P$ by the proof of  \Cref{Lemma1Solid,LemmaSolidificationBoundedSolid}. In particular, we have natural isomorphisms $R\iHom(\prod_{I}\bb{Z}, \bb{Z})=\bigoplus_{I} \bb{Z}$ and $R\iHom(\bigoplus_{I} \bb{Z},\bb{Z})=\prod_{I} \bb{Z}$ for $I$ countable. This shows that the objects $\bb{Z}[S]^{L\sol}$ are reflexive over $\bb{Z}$, and it suffices to prove that the map \eqref{eqComputationMeasuresSolid} becomes an isomorphism after taking duals.  This follows from the fact that
\[
\begin{gathered}
R\iHom(\bb{Z}[S]^{L\sol},\bb{Z}) = R\iHom(\bb{Z}[S],\bb{Z})=C(S,\bb{Z})\\
=\varinjlim_i C(S_i,\bb{Z}) = \varinjlim_i R\iHom(\bb{Z}[S_i], \bb{Z}) = R\iHom(\varprojlim_i\bb{Z}[S_i],\bb{Z}),
\end{gathered}
\]
where in the last equality we use that $\varprojlim_i\bb{Z}[S_i]$ is isomorphic to $\prod_{\bb{N}} \bb{Z}$ by  \Cref{PropFreeContinuousFunctions}.

On the other hand, we have an isomorphism $P\times P \xrightarrow{\sim } P$ given by taking an anti-diagonal enumeration of $\bb{N}\times \bb{N}$. This shows that 
\begin{equation}
\label{eqTensorSolid}
\prod_{\bb{N}} \bb{Z} \otimes^L_{\sol} \prod_{\bb{N}} \bb{Z} \cong (P\otimes P)^{L\sol} \cong P^{L\sol} \cong \prod_{\bb{N}} \bb{Z}.
\end{equation}
An explicit description of this enumeration shows that the isomorphism \eqref{eqTensorSolid} is given by the natural map 
\[
\prod_{\bb{N}} \bb{Z} \otimes^L_{\sol} \prod_{\bb{N}} \bb{Z} \xrightarrow{\sim} \prod_{\bb{N}\times \bb{N}} \bb{Z}.
\]
\end{proof}

A first interesting property of the solidification functor is that it computes singular cohomology of CW complexes. 

\begin{proposition}
\label{PropSingularCohomologyCW}
Let $X$ be a CW complex, then $\bb{Z}[X]^{L\sol}$ is equivalent to the complex of singular chains in $X$. 
\end{proposition}
\begin{proof}
Writing $X$ as a colimit of finite CW complexes it suffices to construct a natural quasi-isomorphism between $\bb{Z}[X]^{L\sol}$ and the chain complex of $X$, we can then assume $X$ to be compact. Let $S\to X$ be a surjection from a light profinite set with \v{C}ech nerve $S_{\bullet}\to X$. We have a resolution 
\[
\cdots \to \bb{Z}[S_2]\to \bb{Z}[S_1]\to \bb{Z}[S_0]\to \bb{Z}[X]\to 0
\]
proving that $\bb{Z}[X]^{L\sol}$ is given by the connective complex. 
\[
\cdots \to  \bb{Z}_{\sol}[S_2] \to \bb{Z}_{\sol}[S_1] \to \bb{Z}_{\sol}[S_0]\to 0.
\]
By  \Cref{CoroDualitySolidMeasures} the complex $\bb{Z}[X]^{L\sol}$ is reflexive, and to naturally identify it with singular chains it suffices to naturally identify its dual with singular cochains. But
\[
R\iHom(\bb{Z}[X]^{L\sol}, \bb{Z})=R\iHom(\bb{Z}[X], \bb{Z})= R\Gamma_{\cond}(X,\bb{Z})
\]
is the condensed cohomology of $X$, that we identified with sheaf cohomology on $X$ by  \Cref{PropCondensedAsSheafCoho}, and so with singular cochains. 
\end{proof}

\subsection{Flatness of $\prod_{\bb{N}}\bb{Z}$ and the structure of the category $\Solid$}
 In this section we prove the last part of  \Cref{TheoSolidAb} regarding the flatness of $\prod_{\bb{N}} \bb{Z}$ as solid abelian group. The proof strategy begins by first describing all the finitely presented solid abelian groups.

\begin{definition}
A solid abelian group is said finitely generated if it is a quotient of $\prod_{\bb{N}} \bb{Z}$. A solid abelian group is said finitely presented if it is a cokernel of a map $\prod_{\bb{N}} \bb{Z} \to \prod_{\bb{N}} \bb{Z}$, equivalently, if it is compact in the category of solid abelian groups. 
\end{definition}

\begin{theorem}
\label{TheoStructureTheoSolid}
The finitely presented objects  of $\Solid$ form an abelian category stable under kernels, cokernels and extensions in $\Solid$, such that $\Solid=\Ind(\Solid^{\finpres})$. Moreover,  any $M\in \Solid^{\finpres}$ has a resolution 
\[
0\to \prod_{\bb{N}} \bb{Z} \to \prod_{\bb{N}} \bb{Z} \to M\to 0.
\]
\end{theorem}

A first corollary is the flatness of $\prod_{\bb{N}} \bb{Z}$.

\begin{corollary}
\label{CoroflatnessProdZ}
The solid abelian group $\prod_{\bb{N}} \bb{Z}$ is flat for the solid tensor product. 
\end{corollary}
\begin{proof}
Since $\Solid=\varinjlim(\Solid^{\finpres})$, it suffices to show that for $M$ a finitely presented solid abelian group, $M\otimes^L_{\sol} \prod_{\bb{N}} \bb{Z}$ sits in degree $0$.  By  \Cref{TheoStructureTheoSolid} we have a resolution 
\[
0\to  \prod_{\bb{N}} \bb{Z} \to \prod_{\bb{N}} \bb{Z} \to M\to 0.
\]
Tensoring with $\prod_{\bb{N}} \bb{Z}$, and using  \Cref{CoroSolidProperties} we see that 
\[
M\otimes^L_{\sol} \prod_{\bb{N}} \bb{Z}=\prod_{\bb{N}} M
\]
which clearly sits in degree $0$.
\end{proof}

In order to prove \Cref{TheoStructureTheoSolid} we shall need the following lemma.

\begin{lemma}
\label{KeyLemmaFinitelyGenerated}
Any finitely generated submodule of $\prod_{\bb{N}} \bb{Z}$ is isomorphic to $\prod_{I} \bb{Z}$ with $I$ countable. 
\end{lemma}
\begin{proof}
 Let $M\subset \prod_{\bb{N}} \bb{Z}$ be a finitely generated subobject, then $M$ is the image of a map $f:\prod_{\bb{N} } \bb{Z} \to \prod_{\bb{N}} \bb{Z}$, which is the dual of a map 
\begin{equation}
\label{eqMapgdualtoProd}
g: \bigoplus_{\bb{N}} \bb{Z}\to \bigoplus_{\bb{N}} \bb{Z}. 
\end{equation}
If $M$ is discrete, then the map $f$ factors through a  map $\Z^{k}\to M$ with $k\in \N$, which implies that $M$ is finitely generated as discrete abelian group. Since $M$ is torsion free, it is finite free and hence isomorphic to $\prod_I \Z$ with $I$ finite. Hence, we can assume without loss of generality that $M$ is not discrete, in this case we will show that $M\cong \prod_{\N} \Z$.

We shall need the following claim: 

\begin{claim}
Let $N$ be a countable abelian group that embeds in a  direct product of $\bb{Z}$, then $N$ is free. 
\end{claim}
\begin{proof}[Proof of the claim]
Let us pick a basis $\{e_{n}\}_{n\in \bb{N}}$ of $\bb{Q}\otimes N$, and let $N_n= \langle  e_0,\ldots, e_n \rangle_{\bb{Q}} \cap N$. It suffices to show that each $N_n$ is finite free, namely, we have $N=\varinjlim_{n} N_n$ and the quotient $N_{n+1}/N_n$ is torsion free.  We can assume without loss of generality that $\{e_n\}_{n\in \bb{N}}\subset N$. Then, it suffices to prove that $M_n=N_n/ \langle e_1,\ldots, e_n \rangle_{\bb{Z}}$ is finite. Suppose it is not, then we can find elements $x_m\in M_n$ of exactly  $b_m$ torsion  for $m\in \bb{N}$, so that $b_m\to \infty$ as $m\to \infty$. Taking lifts $y_m\in N_n$ of $x_n$ this implies that $y_m=\sum_{i=0}^n \frac{c_{i,m}}{d_{i,m}} e_i$ with coefficients satisfying the following properties:

\begin{itemize}

\item $c_{i,m}=0$ or  $\mathrm{GCD}(c_{i,m}{d}_{i,m})=1$, 

\item $\mathrm{lcm}(d_{i,m})=b_m$.

\end{itemize} 
By hypothesis $N$ embeds into $\prod_I \bb{Z}$. Then, there is some projection $\prod_{I}\bb{Z}  \to \prod_{J\subset I} \bb{Z}$ with $J$ finite such that the image of the elements $\{e_1,\ldots, e_n\}$ are linearly independent, proving that for $m>>0$ the element $y_{m}$ cannot be mapped into $\prod_{i=0}^k \bb{Z}$ as $b_m\to \infty$ as $m\to \infty$, which is a contradiction.   This proves the claim. 
\end{proof}

We can decompose the map $g=j\circ h$ in \eqref{eqMapgdualtoProd} as a split surjection $h:\bigoplus_{\bb{N}} \bb{Z} \to M$ and an injection $j:M\to \bigoplus_{\bb{N}} \bb{Z}$.  We can then write short exact sequences
\[
0\to K \to \bigoplus_{\bb{N}} \bb{Z} \xrightarrow{h} M \to 0
\]
and
\[
0\to M \xrightarrow{j} \bigoplus_{\bb{N}} \bb{Z} \to Q\to 0
\]
with $M$ and $K$ free abelian groups.  Taking duals we get  exact sequences
\[
0\to M^{\vee} \to \prod_{\bb{N}} \bb{Z}\to K^{\vee} \to 0
\]
and
\[
0\to \iHom(Q,\bb{Z}) \to \prod_{\bb{N}} \bb{Z} \to M^{\vee} \to \underline{\Ext}^1(Q,\bb{Z}) \to 0.
\]
Then, the composite 
\[
\prod_{\bb{N}} \bb{Z} \xrightarrow{f} \prod_{\bb{N}} \bb{Z} \to K^{\vee}
\]
is zero and we can assume without loss of generality that $K=0$ and so $g$ is injective. Thus, we have an exact sequence
\begin{equation}
\label{equation1SolidStructure}
0\to \bigoplus_{\bb{N}} \bb{Z} \xrightarrow{g} \bigoplus_{\bb{N}} \bb{Z} \to Q \to 0.
\end{equation}
Consider the natural map 
\[
Q\to \prod_{\Hom(Q,\bb{Z})} \bb{Z}
\]
and let $\overline{Q}$ be its image. By the previous claim $\overline{Q}$ is a free abelian group, and so $Q\to \overline{Q}$ is a split surjection. Thus, by taking out the free direct summand,  we can assume without of generality that $\Hom(Q,\bb{Z})=0$. Then, one actually has that $\iHom(Q,\bb{Z})=0$, namely, the $S$-valued points of the $\iHom$ space are equal to $\Hom(Q,C(S,\bb{Z}))$ and $C(S,\bb{Z})$ is a free $\bb{Z}$-module. We deduce that the dual of \eqref{equation1SolidStructure} is the short exact sequence
\[
0\to \prod_{\bb{N}} \bb{Z} \to \prod_{\bb{N}} \bb{Z} \to \underline{\Ext}^1(Q,\bb{Z})\to 0, 
\]
getting that the image of $f$ is isomorphic to $\prod_{\bb{N}} \bb{Z}$ as wanted. 
\end{proof}

\begin{proof}[Proof of  \Cref{TheoStructureTheoSolid}]
By the proof of  \Cref{KeyLemmaFinitelyGenerated} any finitely presented module $M\in \Solid$ is of the form $M=\prod_{I} \bb{Z} \oplus \underline{\Ext}^1(Q,\bb{Z})$ with $I$ a countable set, and  $Q$ a countable abelian group such that $\Hom(Q,\bb{Z})=0$. By taking duals of a  free resolution 
\[
0\to \bigoplus_{\bb{N}} \bb{Z} \to \bigoplus_{\bb{N}} \bb{Z} \to Q\to 0, 
\]
we get a presentation 
\[
0\to \prod_{\bb{N}} \bb{Z} \to \prod_{\bb{N}} \bb{Z} \oplus \prod_{I} \bb{Z} \to M \to 0
\]
proving the second statement of the theorem. The stability of finitely presented solid modules under kernels, cockernels and extensions is then a standard fact for abelian categories for which finitely presented objects admit a resolution by compact projective generators (i.e. are pseudo-compact or  pseudo-coherent, cf. \cite[\href{https://stacks.math.columbia.edu/tag/064N}{Tag 064N}]{stacks-project} for the case of modules over rings).
\end{proof}

\subsection{Examples of solid tensor products}

We finish the discussion of solid abelian groups with some computations of solid tensor products that appear very often  in practice. 

\begin{example}[Power series ring]\label{ExampleTensor}
Let $\bb{Z}[[q]]$ be the ring of power series in one variable seen as a condensed ring. It is  a solid abelian group s since   $\bb{Z}[[q]]=\varprojlim_{n} \bb{Z}[q]/q^n$ is a limit of discrete modules. Indeed, if $\bb{Z}[\widehat{q}]=\bb{Z}[\bb{N}\cup\{\infty\}]/(\infty)$ is the algebra of null-sequences, see  \Cref{PropNullSequenceAlgebra}, we have $\bb{Z}[\widehat{q}]^{L\sol}= \bb{Z}[[q]]$.  \Cref{CoroSolidProperties} implies that 
\[
\bb{Z}[[q_1]]\otimes^L_{\sol} \bb{Z}[[q_2]]= \bb{Z}[[q_1,q_2]]. 
\]
On the other hand,  the morphism of algebras $\bb{Z}[q] \to \bb{Z}[[q]]$ is idempotent when considered as solid algebras, namely, 
\[
\bb{Z}[[q]]\otimes^L_{\bb{Z}[q]} \bb{Z}[[q]]=(\bb{Z}[[q_1]]\otimes^L_{\sol} \bb{Z}[[q_2]])\otimes_{\bb{Z}[q_1-q_2]}^L \bb{Z}= \bb{Z}[[q_1,q_2]]\otimes_{\bb{Z}[q_1-q_2]}^L \bb{Z} = \bb{Z}[[q_1,q_2]]/^{\bb{L}}(q_1-q_2) = \bb{Z}[[q]],
\]
where $\bb{Z}[[q_1,q_2]]/^{\bb{L}}(q_1-q_2)$ is the derived quotient, represented by a Koszul complex. 
\end{example}

\begin{example}[$p$-adic integers]
\label{ExampleZpTensor}
The $p$-adic integers $\bb{Z}_p=\varprojlim_n \bb{Z}/p^n$ is a solid abelian group being a limit of discrete abelian groups.  We have a short exact sequence of solid abelian groups
\[
0\to \bb{Z}[[X]] \xrightarrow{X-p} \bb{Z}[[X]] \to \bb{Z}_p\to 0,
\]
indeed, this is the limit of the short exact sequences
\[
0\to \bb{Z}[X]/X^n \xrightarrow{X-p} \bb{Z}[X]/X^n \to \bb{Z}/p^n\to 0. 
\] 
Thus, the tensor $\bb{Z}_p\otimes^{L}_{\sol} \bb{Z}[[Y]]$ is nothing but $\bb{Z}_p[[Y]]$.

On the other hand, the tensor product $\bb{Z}_p\otimes^L_{\sol} \bb{Z}_{\ell}$ is represented by the complex
\[
\bb{Z}_p[[X]] \xrightarrow{X-\ell} \bb{Z}_p[[X]],
\] 
if $\ell\neq p$ then $\bb{Z}_p\otimes^L_{\sol} \bb{Z}_{\ell}=0$ while if $\ell=p$ we get $\bb{Z}_p\otimes^L_{\sol} \bb{Z}_{p}= \bb{Z}_p$. In particular, $\bb{Z}_p$ is an idempotent $\bb{Z}$-algebra for the solid tensor product. In other words, being a $\bb{Z}_p$-module is not additional structure but a property for solid abelian groups!
\end{example}

\begin{example}[$I$-adically complete modules]
\label{TensorCompleteIadic}
Given a discrete ring $A$ and $I$ a finitely generated ideal, there is a notion of being derived $I$-adically complete (see \cite[Definition 2.12.3]{MannSix} and \cite[\href{https://stacks.math.columbia.edu/tag/091N}{Tag 091N}]{stacks-project}). When $I=(a)$ is generated by a single element, and $A\xrightarrow{a} A$ is the multiplication by $a$,  for an object $C$ in the derived category of (condensed) $A$-modules being $I$-adically complete is equivalent  to the vanishing of the limit  $R\varprojlim_{a} C =0$ given by multiplication along the complex $A \xrightarrow{a} A$.  If we write $J \to A$ for $A \xrightarrow{a} A$, we can think of $J$ as a generalized Cartier divisor, namely, an invertible $A$-module together with a map $J\to A$. We can define powers of $J$ by tensoring, obtaining generalized Cartier divisors $J^n\to A$. Then,  a $A$-modules $C$ is derived $I$-adically complete  if the natural map
\[
C\to R\varprojlim C/^{\bb{L}} J^n,
\]
where the quotient $C/^{\bb{L}} J^n$ is the pushout of $C$ along the map of derived rings $A\to A/^{\bb{L}} J^n$, where  $A/^{\bb{L}} J^n$ is the dg-algebra given by the Koszul complex $J\to A$.  

By \cite[Lemma 2.12.9]{MannSix}, if $A$ is a finitely generated $\bb{Z}$-algebra and $N,M$ are connective derived $I$-adically complete modules, then $N\otimes^L_{A,\sol} M$ is also derived $I$-adically complete (here the tensor product is the natural one attached for a commutative ring object in $\Solid$, equivalently, it is the solidification of the condensed tensor product over $A$). 
\end{example}

\begin{example}[Tensor product of $\bb{Q}_p$-Banach spaces]
Specializing  \Cref{TensorCompleteIadic} to Banach spaces we get the following computation: let $I$ and $J$ be two countable sets, then
\begin{equation}
\label{equationSolidTensorBanach}
\widehat{\bigoplus_I} \bb{Q}_p \otimes^L_{\bb{Q}_p,\sol} \widehat{\bigoplus_J} \bb{Q}_p = \widehat{\bigoplus_{I\times J}} \bb{Q}_p.
\end{equation}
To prove this, since $\widehat{\bigoplus_I} \bb{Q}_p  =(\widehat{\bigoplus_I} \bb{Z}_p) [\frac{1}{p}] $ it suffices to do the analogue computation for $\bb{Z}_p$. By  \Cref{ExampleZpTensor}, the ring $\bb{Z}_p$ is an idempotent solid $\bb{Z}$-algebra, and so the $\bb{Z}$-solid or $\bb{Z}_p$-solid tensor products are the same. Then,    \Cref{TensorCompleteIadic} implies that the solid tensor product
\[
\widehat{\bigoplus_I} \bb{Z}_p \otimes^L_{\sol} \widehat{\bigoplus_J} \bb{Z}_p
\]
is $p$-adically complete, and so it is equal to 
\[
R\varprojlim_{n} ( \bigoplus_I \bb{Z}/p^n \otimes^L_{\sol} \bigoplus_J \bb{Z}/p^n ) =R\varprojlim_n \bigoplus_{I\times J} \bb{Z}/p^n = \widehat{\bigoplus_{I\times J}} \bb{Z}_p. 
\]
For a more direct proof of this fact see \cite[Lemma 3.13]{RRLocallyAnalytic}.
\end{example}

\begin{example}[Tensor product Fr\'echet spaces]
\label{TensorFrechet}
A Fr\'echet $\bb{Q}_p$-vector space is by definition\footnote{A more classical definition is as a topological space with topology given by a countable familiy of seminorms. After taking completions for these seminorms, and passing to the associated condensed set one gets the definition described in the example.} a sequential limit $F=\varprojlim_n V_n$  of Banach spaces,  in particular they are naturally solid $\bb{Q}_p$-vector spaces. If $G= \varprojlim_{n} W_n$ is another Fr\'echet space then 
\[
F\otimes^L_{\sol} G =\varprojlim_{n} (V_n\otimes_{\sol} W_n)
\]
is the projective tensor product in classical functional analysis. In particular, we have that  for $I$ and $J$ countable sets we get
\[
\prod_{I} \bb{Q}_p \otimes^L_{\sol} \prod_J \bb{Q}_p = \prod_{I\times J} \bb{Q}_p.
\]
For a proof of this fact see for example \cite[Lemma 3.28]{RRLocallyAnalytic}. 
\end{example}

\section{Analytic rings}
\label{SectionAnalyticRings}

The building blocks of algebraic geometry are given by commutative rings. In analytic geometry the building blocks are the so called ``analytic rings''. The notion of analytic ring arises from the following desiderata:
\begin{itemize}

\item  An analytic ring $A$ should have an underlying ``topological'' or condensed ring $A^{\triangleright}$.

\item  An analytic ring $A$ should be endowed with a category of complete $A$-modules $\ob{D}(A)$, and with a completed tensor product $\otimes_A$. 
 
\end{itemize}

In the next section we introduce analytic rings and prove some of their most fundamental properties. We will see how the new light foundations of the theory help to construct new examples of analytic rings. Before that, let us give the general definition of light condensed objects in an $\infty$-category:

\begin{definition}\label{DefCondObjectCat}
Let $\s{C}$ be an $\infty$-category with countable limits. The category of light condensed $\s{C}$-objects, denoted by $\Cond(\s{C})^{\light}$, is the $\infty$-category $\widehat{\ob{Shv}}(\Prof^{\light}, \s{C})$ of $\s{C}$-valued hypersheaves on light profinite sets. Equivalently, it is the full subcategory of functors $T\colon \Prof^{\light,\op}\to \s{C}$ such that:

\begin{enumerate}

\item $T(\emptyset)=*$ where $*$ is a final object of $\s{C}$.

\item $T(S_1\sqcup S_2)=T(S_1)\times T(S_2)$

\item For an hypercover $S_{\bullet}\to S$ of light profinite sets, the natural map of anima 
\[
T(S)\to \ob{Tot}(T(S_{\bullet}))=\varprojlim_{[n]\in \Delta} T(S_{\bullet})
\]
is an equivalence. 

\end{enumerate}
\end{definition}

\begin{example}\label{ExamCondensedAnima}
The category of light condensed anima is defined as $\CondAni:=\Cond(\Ani)$. Similarly, if $\Ring$\footnote{Since we use sistematically higher categories it is more natural to use $\Ring$ for the category of animated rings and not for the classical category of commutative rings in sets. The change of notation from non derived to derived world will be recurrent in the rest of the notes.} denotes the $\infty$ -category of animated rings, the category of condensed animated rings is defined as $\CondRing:=\Cond(\Ring)$. 
\end{example}

\begin{remark}\label{RemarkCondensedvsAll}
Let $\mathcal{C}$ be a presentable $\infty$-category  where totalizations commute with $\aleph_1$-filtered colimits (eg. an $\aleph_1$-compactly generared presentable $\infty$-category). Then  $\Cond(\s{C})^{\light}\subset \Cond(\s{C})$ is the full subcategory of all condensed $\s{C}$-objects (see \cite[Definition 2.1.1]{MannSix}) such that for $S=\varprojlim_i S_i$ a profinite set writen as an $\aleph_1$-cofiltered limit of light profinite sets, the natural map 
\[
\varinjlim_{i} T(S_i)\to T(S)
\]   
is an equivalence. 
\end{remark}

\subsection{First definitions and properties}
\label{SubsecFirstDef}

We want to define building blocks for analytic geometry for which we can naturally attach a category of ``complete modules''. It turns out that in condensed mathematics a category of complete modules  for a condensed ring is additional datum; given a condensed ring $A$ there could be very different ways to complete condensed $A$-modules, and none of them should have a preference. Nonetheless, once a category of ``complete modules'' is fixed, being a complete module should be just a property.

On the other hand, derived algebraic geometry \cite{LurieDerivedAlgebraic,ToenDerived}  has shown that the correct framework to study geometric properties of algebraic varieties such as intersections is within higher category theory. In analytic geometry the requirement of higher category theory and higher algebra (taken in the form of \cite{HigherTopos,HigherAlgebra,LurieSpectralAlg}) is even more notorious: even open localizations of rigid or complex spaces are not flat. In particular,  the only way to obtain actual useful new descent results is by looking at the derived $\infty$-categories of modules.

 This desiderata for the notion of analytic ring is formalized in the following definition (see \cite[Definition 12.1 and  Proposition 12.20]{ClauseScholzeAnalyticGeometry} and \cite[Definition 2.3.1]{MannSix}). 

\begin{definition}[Analytic ring]
\label{DefinitionAnalyticRing}
An \textit{uncompleted analytic ring} is a pair $A=(A^{\triangleright},\ob{D}(A))$ consisting on a condensed animated ring $A^{\triangleright}$ and a full subcategory $\ob{D}(A)\subset \ob{D}(A^{\triangleright})$ of the $\infty$-category of condensed $A^{\triangleright}$-modules satisfying the following properties. 
\begin{enumerate}

\item $\ob{D}(A)$ is stable under limits and colimits in $\ob{D}(A^{\triangleright})$ (which is automatically presentable by \cite{ReflectionPrinciple}). 

\item $\ob{D}(A)$ is linear over $\ob{D}(\CondAb)$\footnote{This condition implies that $\ob{D}(A)$ is actually enriched in condensed abelian groups. It can be heuristically thought as a suitable "continuity" or "condensed" condition for $\ob{D}(A)$.}. More precisely for all $C\in \ob{D}(\CondAb)$ and $M\in \ob{D}(A)$  the object $R\iHom_{\bb{Z}}(C, M)$ is in $\ob{D}(A)$.

\item The left adjoint $F$ (which exists by  the adjoint functor theorem \cite[Corollary 5.5.2.9]{HigherTopos}) sends connective objects to connective objects. In particular, $\ob{D}(A)$ has a natural $t$-structure induced from $\ob{D}(A^{\triangleright})$ (see  \Cref{PropCompletionGroups}). 

\end{enumerate}

\begin{itemize}

\item We say that $A$ is an \textit{uncompleted analytic ring structure} of $A^{\triangleright}$.  Finally, we say that  $A$ is an \textit{analytic ring} if  in addition $A^{\triangleright}\in \ob{D}(A)$.  We often write $A\otimes_{A^{\triangleright}} -$ for the left adjoint $F$ (note the drop of derived notation). 

\item Given $T$ a condensed (animated) set we let $A[T]:= A\otimes_{A^{\triangleright}} A^{\triangleright}[T]$, where  $A^{\triangleright}[T]$ is the free $A^{\triangleright}$-module generated by $T$.

\item A morphism of analytic rings $f:A\to B$ is a morphism of animated condensed rings $f:A^{\triangleright} \to B^{\triangleright}$ such that the forgetful functor $f_* :\ob{D}(B^{\triangleright})\to \ob{D}(A^{\triangleright})$ sends $\ob{D}(B)$ to $\ob{D}(A)$.  

\item We let $\AnRing^{un}$ denote the $\infty$-category of uncompleted analytic rings. Let $\AnRing\subset \AnRing^{un}$ be the full subcategory of (completed) analytic rings. 
 
\end{itemize} 
\end{definition}

\begin{remark}
\label{RemarkCondition2}
Condition (2) of  \Cref{DefinitionAnalyticRing} is equivalent to the following:

\begin{itemize}
\item[(2')] For all $C\in \ob{D}(A^{\triangleright})$ and $M\in \ob{D}(A)$ then $R\iHom_{A^{\triangleright}}(C,M)$ is in $\ob{D}(A)$. 

\end{itemize}

Indeed, it suffices to check the condition (2') and (2) on generators of $\ob{D}(A^{\triangleright})$ and $\ob{D}(\CondAb)$ respectively. Then we can suppose without loss of generality that $C=A^{\triangleright}[S]$ or $C=\bb{Z}[S]$ for $S\in \Prof^{\light}$. In this case we have 
\[
R\iHom_{A^{\triangleright}}(A^{\triangleright}[S],M)= R\iHom_{\bb{Z}}(\bb{Z}[S],M). 
\]
\end{remark}

\begin{remark}
\label{RemarkPresentabilityCategories}
Recall that in the new foundations we work with light profinite sets, and so for a condensed animated ring $A^{\triangleright}$ the category $\ob{D}(A^{\triangleright})$ is presentable. By the reflection principle \cite{ReflectionPrinciple} the category $\ob{D}(A)$ of complete $A$-modules is also a presentable category. 
\end{remark}

\begin{example}
\label{ExampleTwoSolidRings}
So far we have   seen essentially only two examples of analytic rings. 

\begin{enumerate}

\item The initial analytic ring is $\bb{Z}^{\cond}:=(\bb{Z},\ob{D}(\CondAb))$, the ring of condensed integers. More generally, given $B^{\triangleright}$ a condensed animated ring, we let $B^{\triangleright, \cond}=(B^{\triangleright},\ob{D}(B^{\triangleright}))$ denote the trivial analytic ring structure on $B^{\triangleright}$. 

 \item A more ``complete'' analytic ring is $\bb{Z}_{\sol}= (\bb{Z},\ob{D}(\Solid))$, the ring of solid integers. Later in  \Cref{SubsecMoreSolidRings} we shall introduce more examples of analytic rings arising in solid geometry.

 \item  Other analytic rings are the liquid rings of \cite{ClauseScholzeAnalyticGeometry} and the gaseous ring of  \Cref{SubsecI:Examples}; these rings are global in the sense that they define analytic ring structures over the subring $\bb{Z}[\widehat{q}]\subset \bb{Z}[[q]]$ of null-sequences that specializes to analytic ring structures over all type of local fields (reals, $p$-adics, and modulo $p$).
 
 \item In   \Cref{SubsectionKillingAlgebras} we discuss a general way to construct analytic rings. This addresses a problem in the previous foundations of condensed mathematics, namely, the difficulty of constructing analytic rings.

\end{enumerate} 
 
\end{example}

Condensed rings embed fully faithfully into analytic rings via the trivial analytic ring structure.

\begin{lemma}
\label{PropTrivialAnRing}
The functor $F:\CondRing\to \AnRing^{un}$ mapping a condensed  condensed ring $A^{\triangleright}$ to $A^{\triangleright,\cond}$ is fully faithful.  Moreover, $F$ has a right adjoint mapping  an uncompleted analytic ring $B$ to its underlying condensed ring $B^{\triangleright}$.
\end{lemma}
\begin{proof}
By definition, given two  uncompleted analytic rings $A$ and $B$ the mapping space  $\Map_{\AnRing^{un}}(A,B)$ is the full subspace of $\Map_{\CondRing}(A^{\triangleright}, B^{\triangleright})$ such that the forgetful functor $\ob{D}(B^{\triangleright})\to \ob{D}(A^{\triangleright})$ sends  $B$-complete objects to $A$-complete objects. If $A$ has the trivial analytic ring structure this condition is tautological which yields  
\[
\Map_{\AnRing^{un}}(A^{\triangleright},B) = \Map_{\CondRing}(A^{\triangleright}, B^{\triangleright}),
\]
proving the fully-faithfulness and the adjunction. 
\end{proof}

The category of complete modules of an uncompleted analytic ring has a natural symmetric monoidal structure. 

\begin{prop}[{\cite[Proposition 12.4]{ClauseScholzeAnalyticGeometry} and \cite[Proposition 2.3.2]{MannSix} }]
\label{PropSymmetricMonoidalStructureDA}
The category $\ob{D}(A)$ has a unique symmetric monoidal structure $\otimes_A$ making $A\otimes_{A^{\triangleright}}-: \ob{D}(A^{\triangleright})\to \ob{D}(A)$ symmetric monoidal. Moreover, given $A\to B$ a morphism of analytic rings, the functor
\[
\ob{D}(A^{\triangleright})\xrightarrow{B^{\triangleright}\otimes_{A^{\triangleright}}-} \ob{D}(B^{\triangleright}) \xrightarrow{B \otimes_{B^{\triangleright}}} \ob{D}(B)
\] 
factors (uniquely) through a functor 
\[
\ob{D}(A^{\triangleright})\xrightarrow{A\otimes_{A^{\triangleright}}-} \ob{D}(A) \xrightarrow{B\otimes_A-} \ob{D}(B).
\]
The functor $B\otimes_A-$ is the left adjoint of the forgetful functor $\ob{D}(B)\to \ob{D}(A)$. 
\end{prop}
\begin{proof}
To show that $\ob{D}(A)$ has a natural symmetric monoidal structure such that $A\otimes_{A^{\triangleright}}$ is symmetric monoidal, it suffices to show that the kernel $K$ of the completion functor is a $\otimes$-ideal by \cite[Theorem I.3.6]{NikolausScholze}. Let $M\in \ob{D}(A^{\triangleright})$ be such that $A\otimes_{A^{\triangleright} } M =0$ and let $C\in \ob{D}(A^{\triangleright})$. Then, for $N\in \ob{D}(A)$, we have
\[
\begin{aligned}
R\iHom_{A^{\triangleright}}( A\otimes_{A^{\triangleright}}( C\otimes_{A^{\triangleright}} M),N) & = R\iHom_{A^{\triangleright}}(C\otimes_{A^{\triangleright}} M,N) \\
				& = R\iHom_{A^{\triangleright}}( M,R\iHom_{A^{\triangleright}}(C,N)) \\
				& = R\iHom_{A^{\triangleright}}(A\otimes_{A^{\triangleright} }M, R\iHom_{A^{\triangleright}}(C,N)) \\
				& = 0,
\end{aligned}
\]
where the first two equalities are the obvious adjunctions, and the third equality follows since  $R\iHom_{A^{\triangleright}}(C,N)$ is $A$-complete by (2) of  \Cref{DefinitionAnalyticRing} (cf.  \Cref{RemarkCondition2}).  The previous shows that $A\otimes_{A^{\triangleright}}( C\otimes_{A^{\triangleright}} M)=0$ as wanted.  

Now, in order to see that the composite 
\[
\ob{D}(A^{\triangleright})\xrightarrow{B^{\triangleright}\otimes_{A^{\triangleright}}-} \ob{D}(B^{\triangleright}) \xrightarrow{B \otimes_{B^{\triangleright}}} \ob{D}(B)
\]
factors through $\ob{D}(A)$, it suffices to see that it kills the kernel of $A\otimes_{A^{\triangleright}}$ (then it would be immediate that the resulting functor is symmetric monoidal). Let $M\in \ob{D}(A)$ be an object killed by $A$-completion and let $K\in \ob{D}(B)$, then 
\[
\begin{aligned}
R\iHom_{B^{\triangleright}}( B\otimes_{B^{\triangleright}}(B^{\triangleright}\otimes_{A^{\triangleright}} M), K) & = R\iHom_{B^{\triangleright}}(B^{\triangleright}\otimes_{A^{\triangleright}} M,K) \\
					& = R\iHom_{A^{\triangleright}}(M,K) \\
					& = R\iHom_{A^{\triangleright}}(A\otimes_{A^{\triangleright}} M,K)\\
					& = 0,
\end{aligned}
\]
where the first three equalities are adjunctions, and the last follows since $K$ is an $A$-complete module by definition of analytic ring. 
\end{proof}

\begin{remark}\label{RemDiscreteModules}
Let $A$ be an analytic ring. Then  $A[*]$ is a compact projective object in $\ob{D}(A)$ and has endormorphisms $\End_{A}(A[*])=A^{\triangleright} (*)$. This produces a symmetric monoidal fully faithful embedding of categories of modules 
\[
\ob{D}(A^{\triangleright}(*))\hookrightarrow \ob{D}(A)
\]
from the category of $A^{\triangleright}(*)$-modules in $\Z$-modules into complete $A$-modules. We call the essential image of this map the full subcategory of \textit{discrete $A$-modules}, and denote it as $\ob{D}^{\delta}(A)$. 

Note that an object in $M\in \ob{D}^{\delta}(A)$ is not necessarily discrete as a condensed set, eg. $A^{\triangleright}$ could have a non-discrete condensed structure. However,  the module $M$ is discrete relative to $A^{\triangleright}$ as it belongs to the smallest full subcategory of complete $A$-modules containing $A^{\triangleright}$ and stable under colimits. 
\end{remark}

Completion of modules for analytic rings can be detected at the level of cohomology groups. 

\begin{prop}[{\cite[Proposition 12.4]{ClauseScholzeAnalyticGeometry}}]
\label{PropCompletionGroups}
Let $A$ be an analytic ring. An object $M\in \ob{D}(A^{\triangleright})$ is $A$-complete if and only if $\pi_{i}(M)=H^{-i}(M)$ is $A$-complete for all $i\in \bb{Z}$. 
\end{prop}
\begin{proof}
Let us first show the statement for connective objects (i.e. concentrated in positive homological degrees). Let $M\in \ob{D}(A)_{\geq 0}$ and consider the fiber sequence
\[
\pi_{\geq 1} M \to M \to \pi_0 M. 
\]
Taking completions we get a fiber sequence
\[
A\otimes_{A^{\triangleright}} (\pi_{\geq 1} M) \to M \to A\otimes_{A^{\triangleright}}(\pi_0 M).
\]
Since completion preserves connective objects, taking $\geq 1$-truncations we get a map 
\[
A\otimes_{A^{\triangleright}} (\pi_{\geq 1} M) \to \pi_{\geq 1} M
\]
which exhibits $\pi_{\geq 1} M$ as a retract of $A\otimes_{A^{\triangleright}} (\pi_{\geq 1} M)$. Since $\ob{D}(A)$ is stable under colimits we deduce that $\pi_{\geq 1} M$ and so $\pi_0(M)$ are in $\ob{D}(A)$. An inductive argument shows that $\pi_i(M)$ is $A$-complete for all $i\geq 0$. Conversely, let $M\in \ob{D}_{\geq 0}(A^{\triangleright})$ be such that all its homotopy groups $\pi_i M$ are $A$-complete.  Then $M=\varprojlim_{n} \tau_{\leq n} M$ is the limit of its Postnikov tower. By induction, each truncation $ \tau_{\leq n} M$ is $A$-complete and then so is $M$ since $\ob{D}(A)$ is stable under limits. 

We now prove the general case. Let $M\in \ob{D}(A)$, then we can write 
\[
M=\varinjlim_{n} \tau_{\geq -n} M,
\] 
and by the connective case it suffices to show that each $\tau_{\geq -n} M$ is $A$-complete. Since $A$-completion preserves connective objects,  $\tau_{\geq -n} M$ is a retract of $A\otimes_{A^{\triangleright}}(\tau_{\geq -n} M)$, and so $A$-complete since $\ob{D}(A)$ is stable under colimits. Conversely, suppose that $M\in \ob{D}(A^{\triangleright})$ is such that $\pi_i(M)$ is $A$-complete for all $i\in \bb{Z}$. By the connective case we know that $\tau_{\geq -n} M$ is $A$-complete for all $n\in \bb{N}$. The proposition follows by writing $M=\varinjlim_{n} \tau_{\geq -n} M$. 
\end{proof}

Our next goal is to show that analytic rings admit small colimits. As a first approximation let us show that uncompleted analytic rings have small colimits. First we will recall  induced analytic structures \cite[Definition 2.3.13]{MannSix}.

\begin{lemma}[Induced analytic structure]
\label{LemmaInducedStructures}
Let $A$ be an uncomplete analytic ring and let $B^{\triangleright}$ be an animated $A^{\triangleright}$-algebra. Then there is a natural induced analytic structure $B_{A/}^{\triangleright}$ on $B^{\triangleright}$ such that $\ob{D}(B_{A/}^{\triangleright}) \subset \ob{D}(B^{\triangleright})$ is the full subcategory of $B^{\triangleright}$-modules whose underlying $A^{\triangleright}$-module is $A$-complete. The uncompleted analytic ring $B_{A/}^{\triangleright}$ is the pushout $A\otimes_{A^{\triangleright,\cond}} B^{\triangleright,\cond}$.
\end{lemma}
\begin{proof}
We want to see that $B_{A/}^{\triangleright}$ defines an  (uncompleted) analytic ring structure on $B^{\triangleright}$. Stability under limits and colimits is clear since the forgetful functor $\ob{D}(B^{\triangleright})\to \ob{D}(A^{\triangleright})$ commutes with limits and colimits.  Since we can describe $\ob{D}(B_{A/}^{\triangleright})=\ob{Mod}_B(\ob{D}(A))$ and $\ob{D}(B)=\Mod_{B}(\ob{D}(A^{\triangleright}))$ and the functor $A\otimes_{A^{\triangleright}}-$ is symmetric monoidal, it is formal that it sends $B^{\triangleright}$-modules to $B^{\triangleright}$-modules and that it defines a left adjoint for the forgetful functor $\ob{D}(B_{A/}^{\triangleright})\to \ob{D}(B^{\triangleright})$ proving that $B_{A/}^{\triangleright}\otimes_{B^{\triangleright}}-= A\otimes_{A^{\triangleright}}-$.

  Stability under $R\iHom(C,-)$ for $C\in \ob{D}(\CondAb)$ is obvious. It is also clear that the left adjoint  $B_{A/}^{\triangleright}\otimes_{B^{\triangleright}}-$ sends connective objects to connective objects. We have proven that $B_{A/}^{\triangleright}$ is an uncompleted analytic ring. 
 
Let us now check that $B_{A/}^{\triangleright}= A\otimes_{A^{\triangleright,\cond}} B^{\triangleright,\cond}$ as uncompleted analytic rings. Let $C$ be an uncomplete analytic ring.  Since $B^{\triangleright,\cond}$ and $A^{\triangleright,\cond}$ have the trivial analytic ring structure,  \Cref{PropTrivialAnRing} implies that a map  $B^{\triangleright, \cond} \to C$ is just given by a map of condensed animated rings  $B^{\triangleright}\to C^{\triangleright}$. Thus, it suffices to see that the following diagram  of mapping spaces is cartesian
  \begin{equation}
  \label{eqCartesianDiagramAnrings}
  \begin{tikzcd}
  \Map_{\AnRing^{uc}}(B_{A/}^{\triangleright}, C)  \ar[r] \ar[d] & \Map_{\CondRing}(B^{\triangleright}, C^{\triangleright}) \ar[d] \\
  \Map_{\AnRing^{uc}}(A,C) \ar[r] & \Map_{\CondRing}(A^{\triangleright}, C^{\triangleright}). 
  \end{tikzcd}
  \end{equation}
The bottom horizontal map of \eqref{eqCartesianDiagramAnrings} is an inclusion. The pullback $\s{C}$ of \eqref{eqCartesianDiagramAnrings} is the full subanima of $\Map_{\CondRing}(B^{\triangleright},C^{\triangleright})$ consisting on those maps $B^{\triangleright}\to C^{\triangleright}$  of $A^{\triangleright}$-algebras such that the forgetful  functor $\ob{D}(C^{\triangleright}) \to \ob{D}(B^{\triangleright})$ sends $C$-complete objects to $A$-complete modules. But this is by definition the mapping space $\Map_{\AnRing^{un}}(B_{A/}^{\triangleright}, C)$, proving what we wanted. 
\end{proof}

A second  important example of colimits of uncompleted analytic rings is obtained by taking intersections of analytic ring structures. 

\begin{lemma}
\label{LemmaIntersectionAnalyticStructures}
Let $A^{\triangleright}$ be a condensed animated ring and let $\{A_i\}_{i\in I}$ be a diagram of (uncompleted) analytic ring structures over $A^{\triangleright}$. Then the pair  $B=(A^{\triangleright}, \bigcap_i \ob{D}(A_i))$ is an  (uncompleted) analytic ring representing the colimit $\varinjlim_{i} A_i$ in the category $\AnRing^{(un)}_{A^{\triangleright}/}$ of  (uncompleted) analytic rings over $A^{\triangleright}$. 
\end{lemma}
\begin{proof}
Let $B$ denote the pair $(A^{\triangleright}, \bigcap_i \ob{D}(A_i))$ where the intersection takes place in $ \ob{D}(A^{\triangleright})$. Note that conditions  (1)-(3) of  \Cref{DefinitionAnalyticRing} are stable under intersection; conditions (2) and (3) are obvious once (1) is proven.  Stability under limits and colimits in (1) is clear. If $A^{\triangleright}$ is $A_i$-complete for all $i$, it is also $B$-complete proving that $B$ is an analytic ring if all the $A_i$ are so. 

It is left to show that $B$ is the colimit of the diagram $A_i$ in the category of (uncompleted) analytic rings over $A^{\triangleright}$. This follows from the fact that for any $C\in \AnRing^{un}$ the maps
\[
\Map_{\AnRing^{un}}(A_i, C) \to \Map_{\CondRing}(A^{\triangleright}, C^{\triangleright})
\]
are fully-faithful embeddings for all $i$, and then so is their  limit. Then, the limit   $\varprojlim_i \Map_{\AnRing^{un}}(A_i, C) $ over $\Map_{\CondRing}(A^{\triangleright}, C^{\triangleright})$  is the full-subanima of  $\Map_{\CondRing}(A^{\triangleright}, C^{\triangleright})$ whose connected components are those maps $A^{\triangleright}\to C^{\triangleright}$ such that the forgetful functor sends $C$-complete modules to $A_i$-complete modules for all $i$. This is exactly the mapping space $\Map_{\AnRing^{un}}(B,C)$ proving what we wanted. 
\end{proof}

We can finally prove the existence of colimits in uncomplete analytic rings. 

\begin{proposition}
\label{PropColimitsUncompletedAnalytic}
The category $\AnRing^{un}$ of uncompleted analytic rings has small colimits. More precisely, let $\{A_i\}_I$ be a diagram of uncompleted  analytic rings. Then $B=\varinjlim_i A_i$ is the uncompleted analytic ring with underlying ring $B^{\triangleright}=\varinjlim_i A_i^{\triangleright}$ and with category of complete modules $\ob{D}(B)\subset \ob{D}(B^{\triangleright})$ given by those $B^{\triangleright}$-modules $M$ whose restrictions to an $A_i^{\triangleright}$-module is $A_i$-complete for all $i$. 
\end{proposition} 
\begin{proof}
First, let us show that the pair $B=(B^{\triangleright}, \ob{D}(B))$ constructed in the statement of the  proposition is an analytic ring. This follows from the fact that $B$ can be written as the colimit
\[
B=\varinjlim_i B^{\triangleright}_{A_i/},
\]
of uncompleted analytic ring structures over  $B^{\triangleright}=\varinjlim_{i} A_i^{\triangleright}$ (\Cref{LemmaIntersectionAnalyticStructures}), where $B^{\triangleright}_{A_i/}$ is the induced analytic ring structure of   \Cref{LemmaInducedStructures}. 

Let us now consider the underlying diagram of condensed animated rings $\{A_i^{\triangleright}\}_i$. Let $C\in \AnRing^{un}$. By definition of the category of analytic rings the limit 
\begin{equation}
\label{eqLimitAnimaAnRings}
\varprojlim_{ i} \Map_{\AnRing^{un}}(A_i,C)
\end{equation}
is a full-subanima of the space 
\[
\varprojlim_i \Map_{\AnRing^{un}}(A_i^{\triangleright}, C^{\triangleright})= \Map(B^{\triangleright}, C^{\triangleright}). 
\]
Furthermore, it is the full subanima of connected components consisting on those maps $B^{\triangleright}\to C^{\triangleright}$ for which a complete $C$-module is $A_i$-complete for all $i$, equivalenty, for which a complete $C$-module is $B^{\triangleright}_{A_i/}$-complete for all $i$. This shows that  \eqref{eqLimitAnimaAnRings} is the  full anima  $\Map_{\AnRing^{un}}(B, C) \subset \Map(B^{\triangleright}, C^{\triangleright})$, proving that $B=\varinjlim_i A_i$ as wanted. 
\end{proof}

A first consequence of the previous lemma is the stability of analytic rings under sifted colimits in the category of uncompleted analytic rings. 

\begin{corollary}
\label{LemmaSiftedAnalytic}
The $\infty$-category $\AnRing$ of analytic rings is stable  under sifted colimits in $\AnRing^{un}$. Moreover, let $B=\varinjlim_{i\in I} A_i$  be a sifted colimit of uncompleted analytic rings. Then for $S\in \Prof^{\light}$ we have 
\begin{equation}\label{eqopjoaegasd}
B[S]=\varinjlim_{i\in I} A_i[S]
\end{equation}
\end{corollary}
\begin{proof}
It suffices to prove the second claim, namely, if that is the case then 
\[
B[*]=\varinjlim_{i\in I} A_i[*]=\varinjlim_{i\in I}A_i^{\triangleright}=B^{\triangleright}
\]
proving that $B^{\triangleright}$ is $B$-complete and so that $B$ is an analytic ring. 

Since $I$ is sifted, $B^{\triangleright}=\varinjlim_i A_i^{\triangleright}$ has a natural structure of condensed animated ring and we have a natural equivalence of categories $\ob{D}(B^{\triangleright})=\varprojlim_i \ob{D}(A_i^{\triangleright})$ with transition maps given by forgetful functors.    Thus by \Cref{PropColimitsUncompletedAnalytic} we have that
\[
\ob{D}(B)=\varprojlim_{i} \ob{D}(A_i) \subset \varprojlim_{i} \ob{D}(A_i^{\triangleright}) =\ob{D}(B^{\triangleright})
\]
In this way, we see that $B=\varinjlim_{i} B^{\triangleright}_{A_i/}$.  Thus, for $M$ a $B^{\triangleright}$-module we have that 
 \[
 \begin{aligned}
 B\otimes_{B^{\triangleright}}M & = \varinjlim_{i\in I} B^{\triangleright}_{A_i/} \otimes_{B^{\triangleright}} M \\ 
 					& = \varinjlim_{i} A_i\otimes_{A_{i}^{\triangleright}} M. 
 \end{aligned}
 \]
Applying the previous to $B^{\triangleright}[S]=\varinjlim_i A^{\triangleright}[S]$ we get
\[
B[S]= \varinjlim_{i} A_i\otimes_{A_i^{\triangleright}} (\varinjlim_{j} A_j^{\triangleright}[S]) = \varinjlim_{i} A_i \otimes_{A_i^{\triangleright}} A_i^{\triangleright}[S] = \varinjlim_i A_i[S] 
\]
where in the second equivalence we used the fact that $I$ is sifted. This proves the corollary. 
\end{proof}

In order to show that analytic rings admit arbitrary colimits we first need to discuss completions of analytic rings. 

\begin{theorem}[{\cite[Proposition 2.3.12]{MannSix}}]
\label{TheoCompletedAnRing}
The functor $\AnRing \to \AnRing^{uc}$ has a left adjoint $A\mapsto A^{\cong}$ called the "completion functor". We have $\ob{D}(A)= \ob{D}(A)^{\cong}$ and $A^{\cong,\triangleright}= A\otimes_{A^{\triangleright}} A^{\triangleright}$ is the $A$-completion of  $A^{\triangleright}$ (i.e. the unit in $\ob{D}(A)$).  In particular, $\AnRing$ admits small colimits. A diagram $\{A_i\}_{i}$ of analytic rings has colimit $B^{\cong}$ where  $B=\varinjlim_i A_i$ is the colimit in the category of uncompleted analytic rings. 
\end{theorem}
\begin{proof}
We will prove a weaker version of the theorem where ``animated ring'' gets replaced by ``commutative or $\bb{E}_{\infty}$-ring''. Indeed, the difficult part of the theorem is to show that the unit $A^{\cong,\triangleright}$ has a natural animated ring structure.  This will be handle in the next section. 

Let $B$ be an analytic ring and $A$ an uncomplete analytic ring. By definition, $\Map_{\AnRing^{un}}(A,B)$ is the full subanima of maps $\Map_{\CAlg(\ob{D}(\Cond))}(A^{\triangleright}, B^{\triangleright})$ of commutative condensed algebras  such that the forgetul functor 
\[
\ob{D}(B^{\triangleright}) \to \ob{D}(A^{\triangleright})
\]
sends $\ob{D}(B)$ to $\ob{D}(A)$. By \cite[Corollary 4.8.5.21]{HigherAlgebra}, the anima $\Map_{\CAlg(\ob{D}(\Cond))}(A^{\triangleright}, B^{\triangleright})$ is naturally equivalent to the anima of $\ob{D}(\CondAb)$-linear symmetric monoidal functors $\ob{D}(A^{\triangleright})\to \ob{D}(B^{\triangleright})$. Therefore, $\Map_{\AnRing^{un}}(A,B)$ gets identified with the full subcategory of symmetric monoidal functors as above that factor through 
\[
\begin{tikzcd}
\ob{D}(A^{\triangleright})  \ar[r, "B^{\triangleright}\otimes_{A^{\triangleright}}"] \ar[d," A\otimes_{A^{\triangleright}}"]& \ob{D}(B^{\triangleright})  \ar[d, "B\otimes_{B^{\triangleright}}"]\\
\ob{D}(A)  \ar[r]& \ob{D}(B).
\end{tikzcd}
\]
Since both $\ob{D}(A)$ and $\ob{D}(B)$ are localizations of $\ob{D}(A^{\triangleright})$ and $\ob{D}(B^{\triangleright})$ respectively, the space $\Map_{\AnRing^{un}}(A,B)$ is naturally equivalent to the space of $\ob{D}(\CondAb)$-linear symmetric monoidal functors $\ob{D}(A)\to \ob{D}(B)$, which is also  equivalent to $\Map_{\AnRing}(A^{\simeq},B)$, proving the desired adjunction. 

The last claim about the computation of the colimit of analytic rings follows directly from the existence of the left adjoint $(-)^{\simeq}$. 
\end{proof}

A consequence of the previous discussion is that the functor $A\mapsto \ob{D}(A)$ commutes with colimit in  presentable categories, in particular it satisfies categorical K\"unneth formula.

\begin{proposition}\label{PropColimitsAnRingAndMod}
The functor $\ob{D}\colon \AnRing\to \ob{CAlg}(\Pr^L_{\ob{D}(\Z^{\cond})})$ sending an analytic ring $A$ to its symmetric monoidal category of complete modules $\ob{D}(A)$ commutes with colimits. 
\end{proposition}
\begin{proof}
By \Cref{TheoCompletedAnRing} the category of analytic rings admits colimits. To show that $\ob{D}$ commutes with colimits it suffices to prove that it commutes with sifted colimits and with pushout of analytic rings.

 We first prove the claim for sifted colimis.  Let $\{A_i\}_i$ be a sifted diagram of analytic rings. Then, as $I$ is sifted, the colimit $\varinjlim_{i\in I} \ob{D}^*(A_i)$ in $ \ob{CAlg}(\Pr^L_{\ob{D}(\Z^{\cond})})$  (with transition maps given by pullbacks) agrees with the colimit in $ \Pr^L_{\ob{D}(\Z^{\cond})}$ and the underlying category can be computed as the limit in $\Pr^R_{\ob{D}(\Z^{\cond})}$ given by $\varprojlim_{i\in I}\ob{D}_*(A_i)$ where the transition maps are forgetful functors. We have a fully faithful embedding 
 \[
 \varprojlim_{i\in I}\ob{D}_*(A_i)\subset \varprojlim_{i\in I}\ob{D}_*(A_i^{\triangleright, \cond}),
 \]
 where $\varprojlim_{i\in I}\ob{D}_*(A_i^{\triangleright, \cond})$ also computes the colimit $\varinjlim_{i\in I}\ob{D}^*(A_i^{\triangleright, \cond})$ in $\ob{CAlg}(\Pr^L_{\ob{D}(\Z^{\cond})})$. We can write $\ob{D}(A_i^{\triangleright, \cond})=\ob{Mod}_{A_i^{\triangleright}}(\ob{D}(\Z^{\cond}))$ and   by \cite[Corollary 4.8.5.14]{HigherAlgebra} we have 
 \[
 \begin{aligned}
 \varinjlim_{i\in I}\ob{D}^*(A_i^{\triangleright, \cond}) & = \varinjlim_{i} \ob{Mod}_{A_i^{\triangleright}}(\ob{D}(\Z^{\cond})) \\ 
  & = \Mod_{\varinjlim_i A^{\triangleright}_i}  (\ob{D}(\Z^{\cond})) \\ 
  & = \ob{D}(\varinjlim_i A^{\triangleright,\cond}_i).
 \end{aligned}
 \]
 We deduce that $\varinjlim_{i} \ob{D}^*(A_i)= \varprojlim_{i} \ob{D}_*(A_i) \subset  \ob{D}(\varinjlim_i A^{\triangleright,\cond}_i)$  is the full subcategory of condensed $\varinjlim_i A^{\triangleright,\cond}_i$-modules that are $A_i$-complete for all $i\in I$, this is precisely $\ob{D}(\varinjlim_i A_i)$ proving what we wanted.

Next we prove the claim for pushouts, let $A$ be an analytic ring and $B$  and $C$ analytic $A$-algebras. We want to see that $\ob{D}(B\otimes_A C)= \ob{D}(B)\otimes_{\ob{D}(A)} \ob{D}(B)$ in $\ob{CAlg}(\Pr^L_{\ob{D}(\Z^{\cond})})$ (that is, the categorical K\"unneth formula).   By writing $A\to B$ as the composite $A\to B_{A/} \to B$, it suffices to show the compatibility for pushouts in the case where either $A\to B$ has the induced structure, or is given by an analytic ring structure. In the first situation, i.e. when $B=B_{A/}$, we have by definition $\ob{D}(B_{A/})=\ob{Mod}_{B}(\ob{D}(A))$ and the K\"unneth formula then follows from \cite[Theorem 4.8.46]{HigherAlgebra}.  Let us then assume that $A\to B$ arises from an analytic ring structure, in other words, that the forgetful functor $\ob{D}(B)\to \ob{D}(A)$ is fully faithful. Then $\ob{D}(A)\to \ob{D}(B)$ is a localization, and by  \cite[Proposition A.5]{NikolausScholze}, for any map of symmetric monoidal categories $\ob{D}(A)\to \s{C}$, one has a fully faithful embedding  
\[
\ob{Fun}^{\otimes}(\ob{D}(B), \s{C}) \subset \ob{Fun}^{\otimes}(\ob{D}(A), \s{C})
\]
with essential image consisting on those symmetric monoidal maps $\ob{D}(A)\to \s{C}$ that kills the tensor ideal given by the kernel $K$ of $\ob{D}(A)\to \ob{D}(B)$. Thus, we have a pullback square 
\[
\begin{tikzcd}
\ob{Fun}^{\otimes}(\ob{D}(B)\otimes_{\ob{D}(A)}\ob{D}(C), \s{C})  \ar[r] \ar[d]& \ob{Fun}^{\otimes}(\ob{D}(C), \s{C})  \ar[d] \\ 
\ob{Fun}^{\otimes}(\ob{D}(B), \s{C}) \ar[r] & \ob{Fun}^{\otimes}(\ob{D}(A), \s{C})
\end{tikzcd}
\]
and the map $\ob{Fun}^{\otimes}(\ob{D}(B)\otimes_{\ob{D}(A)}\ob{D}(C), \s{C})  \subset  \ob{Fun}^{\otimes}(\ob{D}(C), \s{C})$ is fully faithful with essential image given by those symmetric monoidal functors $\ob{D}(C)\to \s{C}$ that kill the tensor ideal generated by $K$ along $\ob{D}(A)\to \ob{D}(C)$. But then the map $\ob{D}(C)\to \ob{D}(B)\otimes_{\ob{D}(A)}\ob{D}(C)$ is a localization and its fully faithful right adjoint $f_*\colon \ob{D}(B)\otimes_{\ob{D}(A)}\ob{D}(C)\hookrightarrow \ob{D}(C)$ has by essential image those objects $M\in \ob{D}(C)$ such that for all $N\in K=\ker(\ob{D}(A)\to \ob{D}(B))$ we have $R\Hom_{C}(C\otimes_A N,M)= R\Hom_A(N , M)=0$. This implies that we have a pullback square
\[
\begin{tikzcd}
\ob{D}(B)\otimes_{\ob{D}(A)}\ob{D}(C)  \ar[r] \ar[d] & \ob{D}(C) \ar[d] \\ 
\ob{D}(B) \ar[r] & \ob{D}(A)
\end{tikzcd}
\] 
where the right vertical and lower horizontal maps are forgetful functors. This shows, thanks to \Cref{TheoCompletedAnRing} that  $\ob{D}(B)\otimes_{\ob{D}(A)}\ob{D}(C) =\ob{D}(B\otimes_A C)$ proving what we wanted.
\end{proof}

\subsection{Completions of analytic rings}
\label{SubsectionCompletions}

In this section we will complete the proof of Theorem \Cref{TheoCompletedAnRing}. For this we need to recall how animated rings are constructed out of connective modules. 

\begin{definition}
\label{DefinitionAnimation}
Let $\s{C}$ be a (presentable) compactly projective generated $1$-category. Let $\s{C}^0\subset \s{C}$ be the full subcategory of compact projective objects. The \textit{animation} of $\s{C}$ (or its non-abelian derived category) is defined as the sifted ind-completion of $\s{C}^0$: $\Ani(\s{C}):= \sInd(\s{C}^0)$ (also denote as $\n{P}_{\Sigma}(\s{C}^0)$ in  \cite[\S 5.5.8]{HigherTopos}). More precisely, it is the full subcategory  
\[
\sInd(\n{C}^0)\subset \Fun(\s{C}^{0,\op}, \Ani)
\]
of accessible presheaves $F$ preserving finite products (i.e. $F(X\sqcup Y)= F(X)\times F(Y)$ for $X,Y\in \s{C}^{0}$). 

\end{definition}

\begin{example}
\label{ExampleAnimatedCategories}
Standard examples of animation are the following:
\begin{enumerate}

\item If $\s{C}=\Sets$ is the category of  sets, then $\s{C}^0$ is the category of finite sets and $\Ani(\s{C})=\Ani$ is the category of anima or of ``spaces'' (this is a tautology from the definition). 

\item If $\s{C}=\Ab$ is the category of abelian groups, then $\s{C}^0$ is the category of free abelian groups and $\Ani(\s{C})$ is the category of animated abelian groups (also known in the literature as \textit{simplicial abelian groups}). Thanks to the Dold-Kan-correspondence \cite[Theorem 1.2.3.7]{HigherAlgebra}, it is also equivalent to the category $\ob{D}_{\geq 0}(\bb{Z})$ of connective objects in the $\infty$-derived category of abelian groups. 

\item If $\s{C}=\Ring^{\ob{classic}}$ is the category of static/classical commutative rings in sets, then $\s{C}^0$ is the category of retracts of polynomial rings in finitely many variables and $\Ani(\s{C}^0)$ is the category $\Ring$ of animated commutative rings (also known as the category of \textit{simplicial commutative rings} in the literature).

\end{enumerate}

\end{example}

\begin{definition}[Symmetric functors]
\label{DefinitionSymmetricFunctor}
Consider $\ob{D}_{\geq 0}(\bb{Z})$ the infinity category of animated abelian groups. The symmetric power functors
\[
\Sym^n: \ob{D}_{\geq 0} (\bb{Z}) \to \ob{D}_{\geq 0} (\bb{Z})
\]
are defined as the left derived functors of the usual symmetric power functors in static rings and abelian groups. More explicitly, it is the unique functor preserving sifted colimits and mapping a finite free abelian group $F$ to its symmetric power $\Sym^n  F$. 
\end{definition}

The importance of the symmetric functors for us is that they appear in the monad defining animated rings. 

\begin{proposition}
\label{PropMonadAnimatedRings}
Let $\Ring$ be the $\infty$-category of animated commutative rings. Let $\ob{D}_{\geq 0}(\bb{Z})$ be the $\infty$-category of connective abelian group.  Then the forgetful functor
\[
G:\Ring\to \ob{D}_{\geq 0}(\bb{Z})
\]
has a left adjoint given by the left derived functor of the functor $\Ab\to \Ring$ mapping an abelian group $M$ to its symmetric algebra $\Sym^{\bullet} M$.  Furthermore, the previous adjunction is monadic. 
\end{proposition}
\begin{proof}
  The forgetful functor $F:\Ring^{\ob{classic}}\to \Ab$ has by left adjoint the symmetric power functor $\Sym^{\bullet}:\Ab\to \Ring^{\ob{classic}}$.  Let $\Ab^0\subset \Ab$ and $\Ring^0\subset \Ring^{\ob{classic}}$ denote the full subcategories of compact projective objects, namely, $\Ab^0$ is the category of finite free abelian groups and $\Ring^0$ is the category of (retracts of) polynomial algebras of finite type. The symmetric power functor  $\Sym^{\bullet}$ restricts to a coproduct preserving functor $\Sym^{\bullet}: \Ab^0\to \Ring^0$. We can then form the sifted ind categories $\sInd$ obtaining the left derived functor 
\begin{equation}
\label{eqSymmetricAnimation}
\Sym^{\bullet}\colon  \ob{D}_{\geq 0}(\bb{Z})\cong \sInd(\Ab^0)\to \sInd(\Ring^0)=\Ring.
\end{equation}
By construction $\Sym^{\bullet}$, preserves coproducts when restricted to $\Ab^0$, namely, if $F_1$ and $F_2$ are finite free abelian groups then $\Sym^{\bullet}(F_1\oplus F_2)= \Sym^{\bullet}F_1 \otimes \Sym^{\bullet} F_2$. Then,   \cite[Proposition 5.5.8.15]{HigherTopos} (3) implies that \eqref{eqSymmetricAnimation} preserves colimits. By the adjoint functor theorem \cite[Corollary 5.5.2.9]{HigherTopos} the functor $\Sym^{\bullet}$ has a right adjoint $G$. Note that by uniqueness of the adjunction, $G$ restricted to $\Ring^{\ob{classic}}\subset \Ring$ is the forgetful functor $G:\Ring^{\ob{classic}}\to \Ab \subset  \ob{D}_{\geq 0}(\bb{Z})$. Then, to see that $G$ is the ``forgetful functor'' on the category $\Ring$ it will suffice to show that it commutes with sifted colimits.  This follows from the fact that $\Sym^{\bullet}$ sends compact projective objects to compact projective objects: given a sifted diagram $\{A_i\}_{i\in I}$ in $\Ring$ and $F\in \Ab^0$ we have
\[
\begin{aligned}
\Map_{\ob{D}_{\geq 0}(\bb{Z})}(F, G (\varinjlim_{i} A_i)) & = \Map_{\Ring}(\Sym^{\bullet} F,  (\varinjlim_{i} A_i)) \\ 
															& = \varinjlim_{i}  \Map_{\Ring}(\Sym^{\bullet} F,  (A_i)) \\ 
															& = \varinjlim_{i} \Map_{\ob{D}_{\geq 0}(\bb{Z})}(F, G A_i)\\
															& = \Map_{\ob{D}_{\geq 0}(\bb{Z})}(F, \varinjlim_{i} G A_i)
\end{aligned} 
\]
where the first equivalence is the adjunction, the second follows since  $\Sym^{\bullet} F$ is compact projective in $\AniRing$, the third is another adjunction, and the last follows since $F$ is compact projective in $\ob{D}_{\geq 0}(\bb{Z})$.  This proves that the natural map $\varinjlim_{i} G A_i\to G(\varinjlim_i A_i)$ is an equivalence. 

 Finally, to show that  the adjunction is monadic, by the monadicity  theorem \cite[Theorem 4.7.3.5]{HigherAlgebra} it suffices to see that $G$ is conservative; this is obvious since $\Sym^{\bullet}$ sends $\Ab^0$ to a set of generators of $\Ring$. 
\end{proof}

\begin{remark}
\label{RemarkAdjunctionPassesToAnyinftyTopos}
%
%
Let $\s{C}$ be an $\infty$-category with finite limits. The adjunction $G\colon \Ring \to \ob{D}_{\geq 0}(\bb{Z})\colon \Sym^{\bullet}$ extends to an adjunction at the level of presheaves on $\s{C}$
\begin{equation}
\label{eqSymmetricPowerDiagram}
G\colon \PShv(\s{C}, \Ring)\to \PShv(\s{C}, \ob{D}_{\geq 0}(\bb{Z})): \Sym^{\bullet}.
\end{equation}
Suppose that $\s{C}$ has in addition a Grothendieck topology $\n{T}$, and for a presentable $\infty$-category $\ob{D}$ let $\widehat{\Shv}(\s{C}, \ob{D})$ denote the full subcategory of $\ob{D}$-valued hypersheaves of $\s{C}$. Then the adjunction \eqref{eqSymmetricPowerDiagram} restricts to an adjuction 
\[
\widehat{G}\colon \widehat{\Shv}(\s{C}, \Ring) \to \widehat{\Shv}(\s{C}, \ob{D}_{\geq 0}(\bb{Z})): \Sym^{\bullet}.
\]
Indeed, since $G$  preserves limits it maps the full subcategory $\widehat{\Shv}(\s{C}, \Ring)\subset \PShv(\s{C}, \Ring)$ to $\widehat{\Shv}(\s{C}, \ob{D}_{\geq 0}(\bb{Z}))\subset \PShv(\s{C},\ob{D}_{\geq 0}(\bb{Z}))$. On the other hand, the inclusion of hypersheaves has by left adjoint the hypercompletion functor
\[
(-)^{\wedge}:\PShv(\s{C},\ob{D})\to \widehat{\Shv}(\s{C}, \ob{D}).
\]
Thus, the forgetful functor 
\[
\widehat{\Shv}(\s{C},\Ring)\to \PShv(\s{C},\ob{D}_{\geq 0}(\bb{Z}))
\]
has by left adjoint the hypercompletion of the symmetric functor, namely, $(\Sym^{\bullet})^{\wedge}$. This restricts to an adjunction 
\[
\widehat{G}\colon \widehat{\Shv}(\s{C}, \Ring) \to \widehat{\Shv}(\s{C}, \ob{D}_{\geq 0}(\bb{Z})): (\Sym^{\bullet})^{\wedge}.
\]
It is clear that $\widehat{G}$ is conservative. Moreover, sifted colimits of objects in $ \widehat{\Shv}(\s{C}, \Ring)$ are   hypersheafifications of sifted colimits in presheaves. This shows that $\widehat{G}$ also commutes with sifted colimits, and so it is monadic. 

Applying the previous construction to $\s{C}=\Prof^{\light}$ endowed with its natural topology (and dropping further notation in the hypercompletion functor), we get the monadic adjunction
\[
G\colon \CondRing \to  \ob{D}_{\geq 0}(\CondAb): \Sym^{\bullet} 
\]
between condensed animated rings and condensed animated abelian groups. 
\end{remark}

After the previous preparations, we can now state the key proposition regarding the completion of analytic rings.

\begin{prop}[{\cite[Proposition 12.26]{ClauseScholzeAnalyticGeometry}}]
\label{PropCompletionSymmetricAlgebras}
Let $A$ be an uncompleted analytic ring. Let $\Ring_{A}\subset \Ring_{A^{\triangleright}/}$ be the full subcategory of condensed animated $A^{\triangleright}$-algebras whose underlying module is $A$-complete. Consider the adjunction 
\begin{equation}
\label{eqAdjunctionSym}
\Sym^{\bullet}_{A^{\triangleright}} \colon \ob{D}_{\geq 0}(A^{\triangleright}) \to \Ring_{A^{\triangleright}} \colon G.
\end{equation}
Then, for any map $N\to M$ of $A^{\triangleright}$-modules which induces an equivalence after $A$-completion, the natural map 
\[
A\otimes_{A^{\triangleright}} \Sym^{\bullet}_{A^{\triangleright}} N \to A\otimes_{A^{\triangleright}} \Sym^{\bullet}_{A^{\triangleright}} 
\]
is also an equivalence. In particular, the monadic adjunction \eqref{eqAdjunctionSym} localizes to a monadic adjuction 
\[
\Sym^{\bullet}_A: \ob{D}_{\geq 0}(A) \to \Ring_{A} \colon G,
\]
where $\Sym^{\bullet}_A = A\otimes_{A^{\triangleright}} \Sym^{\bullet}_{A^{\triangleright}}M$.
\end{prop}
\begin{proof}
This is \cite[Lemma 12.27]{ClauseScholzeAnalyticGeometry}; its proof consists in studying the Goodwillie derivatives of the polynomial functors $\Sym^i$ and reduce the statement to the fact that for all prime $p$ the Frobenius $\phi: A\to A/p$ is a morphism of analytic rings. This last claim will be proven in  \Cref{SubsectionFrob}. 

We give the details in how the symmetric power monad $\Sym^{\bullet}_{A^{\triangleright}}$ descends to a monad $\Sym^{\bullet}_A$ on $\ob{D}_{\geq 0}(A)$. Let $\s{C}$ be an $\infty$-category, $T\in \ob{Alg}(\ob{End}(\s{C}))$ a monad in $\s{C}$ and $\iota\colon \s{C}_0\subset \s{C}$ a full subcategory for which the inclusion admits a left adjoint $F\colon \s{C}\to \s{C}_0$. Suppose that for all map $Y\to X$ in $\s{C}$ that is an $F$-isomorphism the natural map $T(Y)\to T(X)$ also induces an $F$-isomorphism. Then we claim that $T_F= F\circ T \circ F=F\circ T$ has a natural structure of monad in $\mathcal{C}_0$ such that $\ob{LMod}_{T_F}(\s{C}_0)$ is the full subcategory of $\ob{LMod}_T(\s{C})$ consisting on those objects $X$ such that the underlying object in $\s{C}$ belongs to $\s{C}_0$. 

Let $U=\iota\circ F$ be the monad associated to the localization $F$. Then $U$ is an idempotent algebra in  the monoidal category $\s{D}:=\ob{End}(\s{C})$. Let $\ob{BMod}_U(\s{D})\subset \s{D}$ be the full subcategory of endomorphisms $G$ of $\s{C}$ such that the natural maps $G\to UG$ and $G\to GU$ are isomorphisms (this agrees with the category of $U$-bimodules in $\s{D}$ as $U$ is an idempotent algebra). Then  $\ob{BMod}_U(\s{D})$ is a monoidal category with unit $U$, the functor $\ob{BMod}_U(\s{D})\to \s{D}$ is lax monoidal and preserves tensor products (hence, the only failure of monoidality is that this functor does not preserve the unit).  This inclusion admits an oplax left adjoint $\s{D}\to \ob{BMod}_U(\s{D})$ given by sending $X\mapsto UXU$. Notice that, since $\s{C}_0$ is a localization of $\s{C}$, we have a natural equivalence of monoidal categories 
\begin{equation}\label{eqloojoa93aefwe}
\ob{BMod}_U(\s{D})= \ob{End}(\s{C}_0).
\end{equation}

Let $\s{D}_0\subset \s{D}$ be the full subcategory of those functors $X\colon \s{C}\to \s{C}$ that send $F$-isomorphisms to $F$-isomorphisms, equivalently, the category of endofunctors $X$ for which the natural map $UX\to UXU$ is an isomorphism.  We claim that $\s{D}^0$ is a full monoidal subcategory of $\s{D}$, namely, it clearly contains the unit $1$, and if $X,Y\in \s{D}^0$ send $F$-isomorphisms to $F$-isomorphisms then so does its composite $XY$. 
Thus, we have a composite of functors
\[
\ob{BMod}_U(\s{D})\subset \s{D}_0\subset \s{D}
\]
where the first functor is lax monoidal and preserves tensors, and the second functor is monoidal. We deduce that the inclusion $\ob{BMod}_U(\s{D})\subset \s{D}_0$ admits a  monoidal left adjoint sending $X\mapsto UXU=UX$. Indeed, the inclusion admits an apriori oplax monoidal left adjoint given by $X\mapsto UXU=UX$, but the condition $UX\xrightarrow{\sim} UXU$ for the objects in $\s{D}_0$ guarantees that the oplax monoidal structure is actually monoidal. Now, by hypothesis the monad $T$ belongs to $\s{D}_0$, this implies that $T_U:=UT$ has a natural structure of algebra in $\ob{BMod}_U(\s{D})$, which by \eqref{eqloojoa93aefwe} produces the desired monad $T_F$ on $\s{C}_0$. It is left to show that $\ob{LMod}_{T_F}(\s{C}_0)$ is the full subcategory of $\ob{LMod}_T(\s{C})$ consisting on those objects $M$ such that the underlying object in $\s{C}$ belongs to $\s{C}_0$. For that, notice that since $T_U=UT=UTU$, we have that $\ob{LMod}_{T_U}(\s{C})=\ob{LMod}_{T_F}(\s{C}_0)$. Hence, it suffices to see that the natural map of algebras $T\to T_U$ gives rise to a fully faithful functor $\ob{LMod}_{T_U}(\s{C})\hookrightarrow \ob{LMod}_{T}(\s{C})$ with essential image those $M$ such that $U(M)=M$. But this last property follows formally from the  isomorphism $T_U\xrightarrow{1\otimes \id}T_U\otimes_T T_U$ arising from $T_U\otimes_T T_U = (UT)\otimes_{T} T_U = UT_U=T_U$, namely, the previous idempotency implies that the forgetful functor $\ob{LMod}_{T_U}(\s{C})\to \ob{LMod}_{T}(\s{C})$ is fully faithful and  admits a left adjoint given by $T_U\otimes_T(-)=U(-)$. 
\end{proof}

\begin{corollary}
Let $A$ be an uncompleted analytic ring, then  the completion $A^{\cong}$ of $A$ as $\bb{E}_{\infty}$-ring has a natural structure of analytic ring making $A\to A^{\cong}$ a morphism of analytic rings. In other words, $A^{\cong,\triangleright}=A[*]$ has a natural structure of condensed animated ring defined by the completed symmetric powers of  \Cref{PropCompletionSymmetricAlgebras}. In particular, the formation of $A\mapsto A^{\cong}$  is the left adjoint of the natural inclusion $\AnRing\to \AnRing^{uc}$. 
\end{corollary}
\begin{proof}
This follows from  \Cref{PropCompletionSymmetricAlgebras}  and the monadic adjunction  of  \Cref{PropMonadAnimatedRings}, see  \Cref{RemarkAdjunctionPassesToAnyinftyTopos}. 
\end{proof}

\subsection{Frobenius}
\label{SubsectionFrob}


In the proof of   \Cref{PropCompletionSymmetricAlgebras} we used the fact that Frobenius induces a morphism of analytic rings. The goal of this section is to prove this fact (\Cref{TheoFrobIso}). The key step is  \Cref{LemmaTateConstruction} comparing the Tate constructions of free modules on light profinite sets with $C_p$-action.

\begin{lemma}[{\cite[Assumption 12.25]{ClauseScholzeAnalyticGeometry}}] \label{LemmaTateConstruction}
Let  $A$ be an analytic ring. Let $S$ be a light profinite set endowed with a $C_p$-action and let $S_0=S^{C_p}$ be the fixed points. Then the natural map 
\[
A[S_0]^{tC_p}\to A[S]^{tC_p}
\]
is an equivalence, where $(-)^{tC_p}$ is the Tate construction. 
\end{lemma}
\begin{proof}
	Recall  the Tate construction for spectra: let $X\in \mathrm{Sp}$ and let $C_p$ be the cyclic group on $p$-elements. Suppose that we have an homotopic action of $C_p$ on $X$, then there is a norm map $\mathrm{Nm}: X_{C_p}\to X^{C_p}$ from the homotopic co-invariants to the invariants. The Tate construction is defined as the cofiber 
	\[
	X^{tC_p}:= \mathrm{cofib}(X_{C_p}\to X^{C_p}).
	\]
	
	Now let $S$ be a light profinite set endowed with a $C_p$-action  and let $S_0= S^{C_p}$ be its fixed points. For a light locally profinite set $U$ with compactification $U \subset T$ and boundary $\partial T$, let $A[\overline{U}]:= A[T]/A[\partial T]$ be the  module of $A$-measures on $U$ vanishing at $\infty$. Set $U= S\backslash S_0$. The module $A[\overline{U}]$ is independent of the compactification since we have a pushout diagram
	\[
	\begin{tikzcd}
	\partial T \ar[r] \ar[d] &  T \ar[d] \\ 
	\{\infty\} \ar[r] &  U\cup\{\infty\}.
	\end{tikzcd}
	\]
	
	It suffices to show that $A[\overline{U}]^{tC_p}=0$. By  \Cref{PropOpenLightDisjointUnion} we can write $U=\bigsqcup_{n} S'_n$ as a countable disjoint union of light profinite sets. The action of $C_p$ on $U$ is then totally discontinuous not having any fixed point. Then, $U/C_p=\bigsqcup_{n}  S_{n}^{''}$ is a countable disjoint union of light profinite sets, and by taking pullbacks of such decomposition along the map $U\to U/C_p$, we can write $U\cong \bigsqcup_{n} C_p\times S_{n}^{''}$. Therefore, if $ S''$ is a compactification of $U/C_p$, we see that $C_p\times S''$ is a compactification of $U$. This shows that 
	\[
	A[\overline{U}]^{tC_p}= \mathrm{cofib}(A[C_p\times \partial S'']^{tC_p} \to A[C_p\times S'']^{tC_p}).
	\]
	But for any condensed anima $T$ we have $A[C_p\times T]= A[C_p]\otimes A[T]$ and so it has vanishing Tate cohomology $A[C_p\times T]^{tC_p}=0$ (as $A[C_p]$  is the induced representation of $C_p$), proving what we wanted. 
\end{proof}

\begin{theorem}[{\cite[Proposition 12.24]{ClauseScholzeAnalyticGeometry}}]\label{TheoFrobIso}
Let $A$ be an analytic ring and let $A/p=A\otimes_{\bb{Z}} \bb{F}_p$. Then the Frobenius map $\phi: A^{\triangleright}\to A^{\triangleright}/p$ is a map of analytic rings $\phi: A\to A/p$.
\end{theorem}
\begin{proof}
This is proven in \textit{loc. cit.} where the only condition needed is Assumption 12.25 which always holds  in tge light set up thanks to  \Cref{LemmaTateConstruction}. 
\end{proof}

\subsection{Invariance of analytic ring structures}
\label{SubSectionInvarianceAnalyticRings}

It is useful for constructions of analytic rings to compare analytic ring structures between morphisms of condensed animated rings. In this section we shall prove that analytic ring structures are "formally \'etale" in the sense that they are invariant under nilpotent thickenings and higher animated structures. We follow \cite[Lecture XII Appendix 1]{ClauseScholzeAnalyticGeometry}. The first result in this direction if the following theorem that encodes the datum of an analytic ring structure in terms of its underlying  abelian category. 

\begin{theorem}\label{TheoAnRingStructure}
Let $A^{\triangleright}$ be a condensed animated ring. Then the set of (uncompleted) analytic ring structures $A$ over $A^{\triangleright}$ is in bijection with full subcategories $\s{C}$ of the abelian category $\Mod(\pi_0(A^{\triangleright}))$ satisfying the following properties:
\begin{enumerate}
\item $\s{C}$ is stable under all limits, colimits and extensions in $\Mod(\pi_0(A^{\triangleright}))$.

\item $\s{C}$ is presentable. 

\item $\s{C}$ is stable under arbitrary higher direct products $\prod^{(n)}_I$ (these are zero if $I$ is countable). 

\item For all $S\in \Prof^{\light}$ and $C\in \s{C}$ the Ext modules $\underline{\Ext}^{i}_{\bb{Z}}(\bb{Z}[S], C)$ are in $\s{C}$. 

\end{enumerate}

More precisely, given $A$ an analytic ring structure of $A^{\triangleright}$, the category $\s{C}=\ob{D}(A)\cap \ob{D}^{\heartsuit}(A^{\triangleright})$ satisfies the conditions (1)-(4) above. Conversely, given a subcategory $\s{C}$ as above then the category  $\ob{D}\subset \ob{D}(A^{\triangleright})$ consisting on those complexes $C$  with cohomology groups in $\s{C}$ defines an analytic ring structure on $A^{\triangleright}$. 

\end{theorem}

In order to prove the theorem let us first show a bijection for localizations with weaker conditions. 

\begin{prop}[{\cite[Proposition 12.19]{ClauseScholzeAnalyticGeometry}}]\label{PropCategoriesStableLimitsColim}
Let $A^{\triangleright}$ be a condensed animated ring. The collection of full sub $\infty$-categories 	$\ob{D} \subset \ob{D}_{\geq 0}(A^{\triangleright})$  stable under limits and colimits 	is in natural bijection with the collection of all full subcategories $\s{C}\subset \Mod(\pi_0(A^{\triangleright}))=\ob{D}^{\heartsuit}(A^{\triangleright})$ stable under limits, colimits, extensions and higher derived products, via sending $\ob{D}$ to the intersection with $\Mod(\pi_0(A^{\triangleright}))$, and $\s{C}$ to the full subcategory $\ob{D}$ of all $C\in \ob{D}_{\geq 0}(A^{\triangleright})$ such that $\pi_i(C)\in \s{C}$ for all $i\geq 0$. 
\end{prop}
\begin{proof}
Let $\ob{D}\subset \ob{D}_{\geq 0}(A^{\triangleright})$ be a full subcategory stable under limits and colimits. Define $\s{C}=\ob{D}\cap \ob{D}^{\heartsuit}(A^{\triangleright})$. Given $C\in \ob{D}$ the functor $\tau_{\geq 1} C$ is the suspension of the loops of $C$, and so it is in $\ob{D}$. This shows that $\pi_0(C)[0]\in \s{C}$ being the cofiber of $\tau_{\geq 1} C [1]\to C$. Then, $\pi_i(C)[0]\in \s{C}$ for all $i\geq 0$. Since $\ob{D}$ is stable under finite limits and colimits, this shows that $\s{C}$ is stable under finite limits, finite colimits and extensions. Since arbitrary direct sums are exact and $\ob{D}$ has all colimits, then $\s{C}$ has arbitrary direct sums and it is stable under all colimits. Finally, given a family of objects $X_i$ in $\s{C}$, the homotopy product $\prod_{i} (X_i[n])$ is in $\ob{D}$ for all $n\in \bb{N}$ as $\ob{D}$ is stable under all limits. Taking homotopy groups we see that the higher products $\prod^{(n)}_{i} X_i$ are in $\s{C}$ for all $n\in \bb{N}$. In particular, $\s{C}$ has arbitrary products and so it is stable under all limits. 

Conversely, let $\s{C}\subset \ob{D}^{\heartsuit}(A^{\triangleright})$ be a full subcategory stable under all limits, colimits, extensions and arbitrary higher products. Let $\ob{D}\subset \ob{D}_{\geq 0}(A^{\triangleright})$ be the full subcategory consisting on those objects $C$ with homotopy groups in $\s{C}$. Stability under finite limits and extensions in $\s{C}$ shows that $\ob{D}$ is stable under fibers and cofibers. It is also clear that $\ob{D}$ is stable under Postnikov limits.  Since arbitrary direct sums are exact then  $\ob{D}$ is stable under direct sums and so under all colimits. Stability under arbitrary higher products in $\s{C}$ implies that $\ob{D}$ is stable under arbitrary homotopy products, and so it is stable under all limits. 
\end{proof}


\begin{remark}\label{RemarkSubtleDifference}
There is a minor difference between the statement of  \Cref{PropCategoriesStableLimitsColim} and \cite[Proposition 12.19]{ClauseScholzeAnalyticGeometry}, namely in the former we ask for the stability of higher derived products. In the classical condensed framework arbitrary products are exact thanks to the existence of extremally totally disconnected spaces. However, in the light set up a priori only countable products are exact, and there could be higher derived functors for sets with non countable cardinality. 
\end{remark}

\begin{proof}[Proof of  \Cref{TheoAnRingStructure}]
Given an analytic ring structure $A$ on $A^{\triangleright}$, the category $\s{C}=\ob{D}(A)\cap \ob{D}^{\heartsuit}(A^{\triangleright})$ satisfies (1)-(3) of  \Cref{TheoAnRingStructure} thanks to  \Cref{PropCategoriesStableLimitsColim}. Moreover, condition (2) of \Cref{DefinitionAnalyticRing} and  \Cref{PropCompletionGroups} imply that $\s{C}$ is also stable under internal Ext functors. On the other hand,  \Cref{PropCompletionGroups} also says that $\ob{D}(A)\subset \ob{D}$ is the full subcategory consisting on complexes whose cohomology groups are in $\s{C}$. 

Conversely, let $\s{C}\subset \Mod(\pi_0(A^{\triangleright}))$ be a full subcategory as in the statement of the theorem and let $\ob{D}_{\geq 0}\subset \ob{D}_{\geq 0}(A^{\triangleright})$ be the full subcategory of objects whose cohomology groups are in $\s{C}$.  \Cref{PropCategoriesStableLimitsColim} implies that $\ob{D}_{\geq 0}$ is stable under all limits and colimits, it is also presentable by the reflexion principle \cite{ReflectionPrinciple}. The same holds true for its stabilization $\ob{D}\subset \ob{D}(A^{\triangleright})$ consisting on all the complexes whose cohomology groups are in $\s{C}$. Since $\ob{D}$ is stable under all limits and colimits, we have the left adjoint for the inclusion $L:\ob{D}(A^{\triangleright})\to \ob{D}$. Moreover, this left adjoint preserves connective objects since $\ob{D}_{\geq 0}$ is also stable under all limits and colimits in $\ob{D}_{\geq 0}(A^{\triangleright})$. This proves  conditions (1) and (3) of  \Cref{DefinitionAnalyticRing}. It is left to show that $\ob{D}$ is stable under $R\iHom_{\bb{Z}}(\bb{Z}[S],-)$ for $S\in \Prof^{\light}$. Let $M\in \ob{D}$. By writing $M=\varprojlim_{n} \tau_{\geq n} M$ as limit of its Postnikov tower we can assume that $M\in \ob{D}_{\leq 0}$ is co-connective. Then, there is a convergent expectral sequence with second page
\[
E^{p,q}_{2}= \underline{\Ext}^{p}(\bb{Z}[S], \pi_{-q}(M))\Rightarrow \pi_{-p-q}( R\iHom_{\bb{Z}}(\bb{Z}[S],M)).
\]
Since all the objects in the $E_2$-page of the spectral sequence are in $\s{C}$ by hypothesis, and since $\s{C}$ is stable under limits, colimits and extensions, one deduces that the cohomology groups of $R\iHom(\bb{Z}[S],M)$ are in $\s{C}$. We deduce that $\ob{D}$ satisfies condition (2) of  \Cref{DefinitionAnalyticRing}, and so it defines an analytic ring structure on $A^{\triangleright}$. 
\end{proof}

The first application of  \Cref{TheoAnRingStructure} is the homotopy invariance of the analytic ring structures. 

\begin{cor}[{\cite[Proposition 12.21]{ClauseScholzeAnalyticGeometry}}]\label{PropInvarianceHomotopy}
Let $A^{\triangleright}\to B^{\triangleright}$ be a map of animated condensed rings such that $\pi_0(A^{\triangleright})\to \pi_0(B^{\triangleright})$ is an isomorphism.  There is a bijection between (uncompleted) analytic ring structures of $A^{\triangleright}$ and $B^{\triangleright}$ given by mapping $A$ to $B^{\triangleright}_{A/}$. 
\end{cor}
\begin{proof}
By  \Cref{TheoAnRingStructure} analytic ring structures on $A^{\triangleright}$ are in bijection with suitable localizations of the abelian category $\Mod(\pi_0(A^{\triangleright}))$. This proves the corollary. 
\end{proof}

Another application of  \Cref{TheoAnRingStructure} is the invariance of analytic ring structures under nilpotent thickenings.

\begin{prop}[{\cite[Proposition 12.23]{ClauseScholzeAnalyticGeometry}}]\label{PropInvarianceNil}
Let $A^{\triangleright}\to B^{\triangleright}$ be a map of condensed animated rings such that the kernel $I$ of $\pi_0(A^{\triangleright})\to \pi_0(B^{\triangleright})$ is nilpotent. Then there is a natural bijection of uncompleted analytic ring structures on $A^{\triangleright}$ and $B^{\triangleright}$ mapping an analytic ring structure $A$ of $A^{\triangleright}$ to the induced analytic ring structure $B^{\triangleright}_{A/}$.  
\end{prop}
\begin{proof}
By  \Cref{PropInvarianceHomotopy} we can assume that $A^{\triangleright}$ and $B^{\triangleright}$ are static rings.  By induction, we can even assume that $I^2=0$.  

Let $B$ be an analytic ring structure on $B^{\triangleright}$ corresponding to some category $\n{C}_B$.  Let $\s{C}\subset \ob{D}^{\heartsuit}(A^{\triangleright})$ be the full subcategory of objects $M$ such that $IM$ and $M/IM$ are in $\ob{D}(B)$. The category $\s{C}$ is clearly presentable.  We claim that it is stable under limits, colimits, extensions,  arbitrary higher direct products and internal Ext groups from condensed abelian groups. It is clear from the definition that $\s{C}$ is stable under kernels, cokernels and extensions and that it contains $\s{C}_B$.  It also contains arbitrary direct sums as they are exact. To see that it contains arbitrary higher products consider  a family of objects $\{M_i\}_{i\in I}$    in $\s{C}$. We  have a fiber sequence of homotopy products
\[
{\prod_{i}}^{h} IM_i \to {\prod_i}^h M_i \to {\prod_i}^h M_i/IM_i.
\]
Taking the long exact complex we see that the higher product $\prod^{(n)} M_i$ are in $\s{C}$, namely, by  \Cref{TheoAnRingStructure} we know that arbitrary higher products of objects in $\s{C}_B$ stay in $\s{C}_B$.  Finally, stability under internal Ext of $\s{C}_B$ with condensed abelian groups follows by the long exact sequence induced by the short exact sequence $0\to IM\to M\to M/IM\to 0$.
\end{proof}

\subsection{Morphisms of analytic rings}

Let $A$ and $B$ be analytic rings and let $f:A^{\triangleright}\to B^{\triangleright}$ be a morphism of condensed animated rings. We want to have a criterion for the map $f$ to be a morphism of analytic rings $f:A\to B$.  The category $\ob{D}_{\geq 0}(B)$ is generated by the objects $B[S]$ for $S$ a light profinite set. Then, $f$ is a morphism of analytic rings if and only if $B[S]$ is $A$-complete for all $S$-light profinite. Suppose that instead we are given with functorial maps $A[S]\to B[S]$ linear over $A^{\triangleright}\to B^{\triangleright}$ commuting with the map from $S\in \Prof^{\light}$. Then this datum produces a map of analytic rings under a mild condition: 

\begin{prop}[{\cite[Proposition 7.14]{ClausenScholzeCondensed2019}}]\label{propMorphismAnRing}
Keep the previous notation. Suppose that for all $S\in \Prof^{\light}$ with a map $S\to A^{\triangleright}$, inducing $A[S]\to A[*]$ in $\ob{D}(A^{\triangleright})$, and from the composite $S\to A^{\triangleright}\to B^{\triangleright}$, a unique map $B[S]\to B[*]$ in $\ob{D}(B^{\triangleright})$, the diagram 
\[
\begin{tikzcd}
\pi_0(A[S]) \ar[r] \ar[d] & \pi_0(A[*]) \ar[d] \\ 
\pi_0(B[S]) \ar[r] & \pi_0(B[*])
\end{tikzcd}
\]
commutes.  Then $f: A^{\triangleright} \to B^{\triangleright}$ is a morphism of analytic rings $f:A\to B$. 
\end{prop}
\begin{proof}
Let $\s{C}_A$ and $\s{C}_B$ be the hearts of the categories of complete $A$ and $B$-modules respectively.  By  \Cref{TheoAnRingStructure} it suffices to show that objects in $\s{C}_{B}$ are in $\s{C}_A$ when seen as $A^{\triangleright}$-modules. Since the objects $\pi_0(B[S])$ are generators of $\s{C}_B$ it suffices to prove that they are in $\s{C}_{A}$. This reduces the question to the abelian situation, where the proof of \cite[Proposition 7.14]{ClausenScholzeCondensed2019} applies. 
\end{proof}

\subsection{Localizing by killing algebras}
\label{SubsectionKillingAlgebras}

In the ``old'' foundations of condensed mathematics the construction of analytic rings was a big challenge. The construction of the solid integers required a full understanding of extension groups of locally compact abelian groups, and the construction of the liquid rings involved a lot of non-locally convex functional analysis. In the new framework of light condensed mathematics it is much easier to construct analytic rings out from the internally compact projective object $P$ of null sequences. This simplifies the construction of solid rings, and gives a natural construction of gaseous rings motivated from the Tate curve. A disclaimer: the construction of the liquid rings remains as difficult as before and a priori the light theory does not help to simplify the key computations. Nevertheless, we can now construct localizations of categories of modules in a much more general way as we shall explain down below.

Let $\s{C}$ be a presentably symmetric monoidal stable $\infty$-category. Let $A\in \s{C}$ be an object endowed with the following two maps
\begin{enumerate}
\item $m\colon A\otimes A\to A$ and

\item $\mu\colon 1\to A$

\end{enumerate}
such that  the composite $A \xrightarrow{ \mu\otimes \id_A} A\otimes A \xrightarrow{m} A$ is the identity of $A$.  We let $\s{E}\subset \s{C}$ be the full subcategory of objects $M$ such that $\iHom(A,M)=0$.  It is clear that $\s{E}$ is presentable by the reflection principle \cite{ReflectionPrinciple},   it is stable under limits in $\s{C}$ and under internal Hom. 

Our goal is to construct explicitly the localization functor $\s{C}\to \s{E}$.  Let $C= \fib(1\to A)$ be the fiber. Define $F: \s{C}\to \s{C}$ to be the functor $\iHom(C,-)$. Since we have a map $C\to 1$, there is a natural transformation of functors $\id_{\s{C}}\to F$. 

\begin{lemma}\label{LemmaLocalizationF}
Let $X\in \s{C}$ and $M\in \s{E}$. Then $\iHom(F(X),M)\to \iHom(X,M)$ is an equivalence. 
\end{lemma}
\begin{proof}
By unraveling the constructions, it suffices to show that  
\begin{equation}
\label{equationVanishingF}
\iHom(\iHom(A,X),M)=0.
\end{equation}
We claim that $\iHom(A,X)$ is a retract of $A\otimes \iHom(A,X)$. Suppose the claim holds, then we get that 
\[
\iHom(A\otimes \iHom(A,X), M) =\iHom(\iHom(A,X), \iHom(A,M))=0,
\]
which implies the vanishing of \eqref{equationVanishingF}.  Let us now prove the claim. The multiplication map $m:A\otimes A\to A$ induces a map 
\[
\iHom(A,X)\to \iHom(A\otimes A,X)
\]
which is adjoint to a map 
\[
A\otimes \iHom(A,X)\to \iHom(A,X).
\]
On the other hand, the unit map  $\mu:1\to A$ induces a map $ \iHom(A,X)\to A\otimes \iHom(A,X)$. Then a diagram chasing  shows that the composite 
\begin{equation}
\label{equationMapsComposite1}
\iHom(A,X) \to A\otimes \iHom(A,X) \to \iHom(A,X)
\end{equation}
is the identity map, proving the claim.   Indeed, the diagram \eqref{equationMapsComposite1} is adjoint  to a diagram 
\[
f\colon A\otimes \iHom(A,X)\xrightarrow{\mu \otimes \id_A} A\otimes A\otimes \iHom(A,X) \xrightarrow{g} X
\]
where $g$ is the composite
\[
g\colon A\otimes A\otimes\iHom(A,X) \xrightarrow{m^*} A\otimes A\otimes \iHom(A \otimes A, X) \xrightarrow{\ev_{A\otimes A}} X.
\] Then, we have a commutative square 
\[
\begin{tikzcd}
A\otimes A \otimes  \iHom(A,X)  \ar[r, "m^*"] \ar[d, "m\otimes \id"]& A\otimes A\otimes \iHom(A\otimes A, X) \ar[d, "\ev_{A\otimes A}"] \\
A\otimes \iHom(A,X) \ar[r, "\ev_{A}"] & X
\end{tikzcd}
\]
Therefore, $f$ is also the composite 
\[
A\otimes \iHom(A,X)\xrightarrow{\mu\otimes \id_A} A\otimes A \iHom(A,X) \xrightarrow{m\otimes \id} A\otimes\iHom(A,X)\xrightarrow{ev_A} X 
\]
which is the same as the evaluation map  $\ev_A$ since $m\circ (\mu\otimes \id_A)=\id_A$. Taking adjoints, one deduces that the composite  \eqref{equationMapsComposite1} is  the identity. 
\end{proof}

Let $n\in \bb{N}$ and let $F^n: \s{C}\to \s{C}$ be the $n$-th iteration of the functor $F$. The natural transformation $\id_{\s{C}}\to F$ produces a sequential diagram of natural transformations 
\[
\id_{\s{C}}\to F\to F^2\to \cdots \to F^n\to \cdots.
\]
We let $F^{\infty}=\varinjlim_n F^n$.  \Cref{LemmaLocalizationF} shows that for all $n\in [0,\infty]$, $X\in \s{C}$ and $M\in \s{E}$ the natural map 
\begin{equation}
\label{eqIsoFIteration}
\iHom(F^n(X),M)\xrightarrow{\sim} \iHom(X,M)
\end{equation}
is an equivalence. We want to impose some conditions on  $F$ for $F^{\infty}$ to be a left adjoint of the inclusion. 

\begin{prop}\label{PropStabilityF}
Suppose that one of the following conditions hold:
\begin{enumerate}
\item  The sequential colimit $ F^{\infty}(X)=\varinjlim_n F^n(X)$ stabilizes for all $X$ (eg. if $A$ is idempotent).

\item $\iHom(A, -):\s{C}\to \s{C}$ commutes with sequential colimits (eg. if $A$ is internally compact in $\s{C}$). 

\end{enumerate}
Then $F^{\infty}: \s{C}\to \s{C}$ lands in $\s{E}$ and is the left adjoint of the inclusion $\s{E}\subset \s{C}$. 
\end{prop}
\begin{proof}
By \eqref{eqIsoFIteration} it suffices to show that $F^{\infty}$ lands in $\s{E}$. Conditions (1) and (2) imply that for all $X\in \s{C}$ the natural map
\[
\varinjlim_{n} \iHom(A, F^n(X)) \to \iHom(A, F^{\infty}(X))
\]
is an equivalence. Note that we have a commutative diagram whose rows are fiber sequences
\[
\begin{tikzcd}
\iHom(A, F^{n+1}(X)) \ar[r] & F^{n+1}(X)  \ar[r] & F^{n+2}(X)  \\ 
\iHom(A, F^n(X)) \ar[u] \ar[r] & F^{n}(X) \ar[r] \ar[u] & F^{n+1}(X) \ar[u] \ar[ul,"\id", dashed].
\end{tikzcd}
\]
Then, taking colimits as $n\to \infty$ in the columns, we obtain a fiber sequence
\[
\varinjlim_n \iHom(A,F^n(X)) \to \varinjlim_n F^{n}(X) \xrightarrow{\sim} \varinjlim_{n} F^{n+1}(X),
\]
where the right arrow is an equivalence. This proves that $\iHom(A,F^{\infty}(X))=\varinjlim_n \iHom(A,F^n(X))=0$ as wanted. 
\end{proof}

\begin{example}
Some classical localizations in commutative algebra are examples of those in  \Cref{PropStabilityF}.

\begin{enumerate}

\item Let $R$ be an animated ring, $\s{C}=\ob{D}(R)$ and $P= R/^{\bb{L}}f$.   Then the category $\s{E}\subset \s{C}$ of objects $M$ such that $R\iHom_{R}(P,M)=0$ is precisely $\s{E}=\ob{D}(R[1/f])$. Indeed, an explicit computation shows that $F^{\infty}(M)=\varinjlim_{\times f} M =M[1/f]$.

\item Let us keep $R$ and $\s{C}$ as before and take $P=R[1/f]$. Then $P$ is an idempotent algebra and the category $\s{E}\subset \s{C}$ of $R$-modules $M$ such that $R\iHom_R(P,M)=0$ is precisely the category of $f$-adicaly complete modules. The functor $F^{\infty}$ stabilizes for $n=1$ and $F(M)=\varprojlim_n M/^{\bb{L}} f^n$ is the $f$-adic completion functor.

\end{enumerate}

\end{example}

\begin{example}\label{ExamplesLocalization}

Let $P=\bb{Z}[\bb{N}\cup\{\infty\}]/(\infty)$ be the free condensed abelian group of null sequences. By  \Cref{PropNullSequenceAlgebra} it has a natural algebra structure making $\bb{Z}[q]\to P$ a morphism of algebras, where $q$ is mapped to $[1]$. We will write $P=\bb{Z}[\widehat{q}]$. 

\begin{enumerate}

\item The multiplication by $q$ in $\bb{Z}[\widehat{q}]$ corresponds to the shift map $\Shift: P\to P$. Then, the category of solid abelian groups is precisely the category of those condensed abelian groups $M$ such that 
\[
R\iHom( \bb{Z}[\widehat{q}]/(1-q), M)=0. 
\]
The object $P$ is internally compact projective, then the previous localization lands in the case (2) of  \Cref{PropStabilityF}.

\item  Let $\bb{Z}_{\sol}$ be the ring of solid integers. We know that $\bb{Z}_{\sol}\otimes_{\bb{Z}} \bb{Z}[\widehat{q}]= \bb{Z}[[q]]$ is a power series ring in the variable $q$. We can construct additional solid structures arising from polynomial algebras as follows: we define the category of solid $\bb{Z}[T]_{\sol}$-modules, denoted by $\Mod(\bb{Z}[T]_{\sol})$ to be the full subcategory of $\bb{Z}$-solid $\bb{Z}[T]$-modules $M$ such that 
\[
R\iHom_{\bb{Z}[T]}( \bb{Z}[[q]] [T] /(1-Tq), M)=0.
\]
Heuristically, we are asking for a null sequence $(b_n)_{n\in \bb{N}}$ to be $T$-summable, i.e. for $\sum_n b_n T^n$ to converge.  Note that $\bb{Z}[[q]][T]/(1-qT)= \bb{Z}((T^{-1}))$ is the ring of Laurent power series in $T^{-1}$. By 
 \Cref{ExampleTensor} the algebra $\bb{Z}((T^{-1}))$ is  idempotent over $\bb{Z}[T]$. Then the previous localization  lands in both conditions (1) and (2) of 
 \Cref{PropStabilityF}.

\item The new kind of analytic rings that can be constructed abstractly using  \Cref{PropStabilityF} are the gaseous rings.  Let $A^{\triangleright}= \bb{Z}[\widehat{q}][q^{-1}]$  and consider the algebra $A^{\triangleright}\otimes_{\bb{Z}} P$. Let $T$ denote the variable of $P$. Then the gaseous structure over $A^{\triangleright}$ is the localization with respect to the algebra $A\otimes_{\bb{Z}} P /(1-qT)$. In other words, an object $M\in \ob{D}(A^{\triangleright})$ is gaseous if 
\[
R\iHom_{A^{\triangleright}}(A\otimes_{\bb{Z}} P /(1-qT), M)=0.
\]

\item More generally, given an analytic ring $A$ consider $P_A=A\otimes_{\bb{Z}} \bb{Z}[\widehat{q}]$. Then, for any $P_A$-algebra $R$ which is a perfect $P_A$-module one can consider the localization $\s{E}\subset \ob{D}(A)$ consisting on the objects $M$ such that $R\iHom(R,M)=0$. The category $\s{E}$ satisfies conditions (1) and (2) of  \Cref{DefinitionAnalyticRing}. The only constrain to define an analytic ring structure for $A^{\triangleright}$ is the connectivity condition (3) which is actually easy to check in practice. Nonetheless, this solves the problem of constructing several examples of analytic rings by a systematic procedure (after verifying condition (3) for connectivity).

Nonetheless, the difficulty in the construction of analytic rings is still preserved, and now is reflected in the computation of the free completed modules on light profinite sets. The advantage is that now the problem is not about a definition or construction but a computation, which is conceptually much more clean.

\end{enumerate}
\end{example}

\section{Solid analytic rings}
\label{SubsecMoreSolidRings}

In this section we give  examples of analytic rings arising from non-archimedean geometry. In  \Cref{SectionLightSolid} we constructed the analytic ring $\bb{Z}_{\sol}$ of solid integers, our first objective is to generalize this construction to finite type algebras over $\bb{Z}$ and then to arbitrary discrete animated rings. Then, we discuss the categories of solid quasi-coherent sheaves for schemes and discrete adic spaces.

\subsection{Smashing spectrum}

Before giving examples of solid analytic rings, let us discuss some general constructions in (stable) symmetric monoidal $\infty$-categories that will be crucial in the definition of analytic stacks. These are natural categorifications of open and closed immersions of spaces from the point of view of  six functor formalisms.  We shall follow \cite[Lecture V]{CondensedComplex}, see also \cite{aoki2023sheavesspectrum,aoki2024smashingspectrumsheaves}.   \textbf{Throughout this   subsection we shall work fully derived, namely, tensor products and Hom spaces will  be assume to be derived unless otherwise specified.}

\subsubsection{Topological six functors}

As motivation let us recall some basic facts about sheaves on topological spaces. Let $X$ be a topological space, and let $D(X,\bb{Z})$ be the derived category of abelian sheaves on $X$.  Given $U\subset X$ an open subspace and $Z=X\backslash U$ the closed complement we have different functors relating the categories $D(U,\bb{Z})$, $D(Z,\bb{Z})$ and $D(X,\bb{Z})$. More precisely, we have 

\begin{itemize}

\item  Pullback functors 
\[
\begin{aligned}
\iota^{*}:D(X,\bb{Z})\to D(Z,\bb{Z}) \\
j^*: D(X,\bb{Z})\to D(U,\bb{Z}). 
\end{aligned}
\]
\item  The functor $\iota^*$ has a right adjoint given by a pushforward or extension by $0$ functor
\[
\iota_*: D(Z,\bb{Z})\to D(X,\bb{Z}).
\]

\item $\iota_*$  itself has also a right adjoint given by sections supported at $Z$
\[
\iota^!:D(X,\bb{Z})\to D(Z,\bb{Z}).
\]

\end{itemize} 

The functor $\iota_*$  satisfies the projection formula, namely,  the following arrow is an equivalence 
\begin{equation}
\label{eqClosedImmersion}
 \iota_*N\otimes_{\bb{Z}_X} M  \xrightarrow{\sim} \iota_*(N\otimes_{\bb{Z}_Z} \iota^*M).
\end{equation}
 More explicitly, $\iota_* \bb{Z}_Z$ is the  locally constant sheaf $\bb{Z}$ supported on $Z$, and we have 
\[
\iota_*\iota^* =\iota_* \bb{Z}_Z\otimes_{\bb{Z}_X} -.
\]
The functor $\iota^!$ is then described as  
\[
\iota_* \iota^! = \iHom_{\bb{Z}_X}(\iota_{*} \bb{Z}_Z, -). 
\]

 On the other hand, the functor $j^*$ has both left and right adjoints.

\begin{itemize}
 
\item  The right adjoint $j_*:D(U,\bb{Z})\to D(X,\bb{Z})$ is the natural pushforward functor.  

\item The left adjoint $j_!: D(U,\bb{Z})\to D(X,\bb{Z})$ is the natural extension by $0$ functor.

\end{itemize}

  Furthermore, the functor $j_!$ satisfies the projection formula, namely, the natural map 
\begin{equation}
\label{eqOpenImersion}
j_!(N\otimes_{\bb{Z}_U}M)\xrightarrow{\sim} j_!N \otimes_{\bb{Z}_X} M
\end{equation}
is an equivalence.

It turns out that  the functors $\iota_*$, $j_*$ and $j_!$ are fully faithful; this is a consequence  of the fact that the map $\bb{Z}_X\to \iota_* \bb{Z}_Z$ is idempotent in $D(X,\bb{Z})$,  namely, that the natural map $\iota_* \bb{Z}_{Z}\xrightarrow{\id\otimes 1} \iota_* \bb{Z}_Z\otimes^L_{\bb{Z}_X} \iota_* \bb{Z}_Z$ is an equivalence. Hence,  one has natural fiber sequences in $D(X,\bb{Z})$ involving all four functors:
\[
\begin{aligned}
\iota_*\iota^!\to \id_{X}\to j_*j^*\\
 j_!j^* \to \id_{X} \to  \iota_* \iota^* .
\end{aligned}
\]
The previous fiber sequences give rise to ``Verdier exact sequences'' of derived categories
\begin{equation}
\label{eqLocalizationDerivedtodSps}
\begin{tikzcd}
D(Z,\bb{Z}) \ar[r,"\iota_*"] & D(X,\bb{Z}) \ar[r,"j^*"] \ar[l,"\iota^!"',shift left =4] \ar[l,"\iota^*"',shift right=4] &  D(U,\bb{Z}). \ar[l, "j_*"',shift left=4] \ar[l, "j_!"',shift right =4] \\ 
\end{tikzcd}
\end{equation}

One would like to generalize the localization sequences of \eqref{eqLocalizationDerivedtodSps} altogether with the projection formulas \eqref{eqClosedImmersion} and  \eqref{eqOpenImersion} in order to talk about abstract open and closed immersions of symmetric monoidal categories.  This idea is realized thanks to the smashing spectrum. 

\subsubsection{Smashing spectrum and idempotent algebras}

For a  general notion of open and closed immersions the key objects are idempotent algebras.

\begin{definition}
  Let $\s{C}$ be a  presentably symmetric monoidal stable $\infty$-category with unit $1$. In particular, $\s{C}$ is closed, i.e., it has an internal $\iHom$.

  An idempotent algebra in $\s{C}$ is a map 
  \[
  \mu:1\to A
  \]
  in $\s{C}$ such that the natural map 
  \[
  A\xrightarrow{\id_{A}\otimes \mu} A\otimes A
  \]
  is an equivalence.  A morphism $A\to B$ of idempotent algebras  is a map in $\s{C}$ preserving the unit.  We let $\Sm(\s{C})$ be the opposite of the category of idempotent algebras in $\s{C}$ and call it the \textit{smashing spectrum of $\s{C}$}.
  \end{definition}

A priori an idempotent algebra is not required to have any  algebra structure. However it will be always endowed with a natural commutative  algebra structure arising from the natural equivalences
\[
A\cong A\otimes A\cong A\otimes A\otimes A\cong \cdots 
\]  
together with all the higher coherences, see \cite[Proposition 4.8.2.9]{HigherAlgebra}. 

Given $A$ an idempotent algebra in $\s{C}$ the natural map 
\[
\Mod_A(\s{C}) \to \s{C}
\]
is fully faithful, namely, $M\in \s{C}$ is an $A$-module if and only if $M\xrightarrow{\mu\otimes \id_M} A\otimes M$ is an equivalence. In other words, being an $A$-module for an object in $\s{C}$ is a property and not additional structure. One also has that $M$ is an $A$-module if and only if $\iHom(A,M)\to M$ is an isomorphism. 
 
Finally, given $A$ and $B$ idempotent algebras, the mapping space $\Map_{1/}(A,B)$ of idempotent algebras is either contractible or empty. This shows that $\Sm(\s{C})$ is  a poset. Actually,  the category $\Sm(\s{C})$ has the structure of a  \textit{locale}, i.e. it behaves as the poset of  closed subspaces of a topological space:  
  
  \begin{prop}[{\cite[Proposition 5.3]{CondensedComplex}, \cite[Theorem 3.8]{aoki2023sheavesspectrum}}]
  The poset $\Sm(\s{C})$ is a locale whose closed subspaces $Z\subset \Sm(\s{C})$ correspond to idempotent algebras $A$ of $\s{C}$, so that
  \begin{enumerate}
\item  $Z\cap Z'$ corresponds to $A\otimes A'$; 

\item $Z\subset Z'$ if and only if $A\otimes A'=A$;

\item $Z\cup Z'$ corresponds to $\fib(A\oplus A'\to A\otimes A')$, the unit is given by $1\xrightarrow{\mu_{A}\oplus (-\mu_{A'})} A\oplus A'$. Equivalently, $Z\cup Z'$ corresponds to the pullback $A\times_{A\otimes A'} A'$. 

\item  $\bigcap_i A_i$ corresponds to $\varinjlim_i A_i$.

\end{enumerate}   
  \end{prop}
 
 We will often let $Z\subset \Sm(\s{C})$ denote a closed subspace of the local of $\s{C}$ and let $A(Z)$ be the attached idempotent algebra.

\subsubsection{Open and closed immersions of symmetric monoidal categories} \label{sss:OpenClosedSymmetricMonoidal}
 
Let $\ob{CAlg}(\Pr^L_{\st})$ be the $\infty$-category of presentably symmetric monoidal stable $\infty$-categories. We will take inspiration from the six functors of topological spaces to define open and closed immersions.  For geometric reasons, we will work with the opposite category $\ob{Affine}(\Pr^L_{\st}):=\ob{CAlg}(\Pr^L_{\st})^{\op}$, and given $\s{C}\in \ob{CAlg}(\Pr^L_{\st})$ we denote $\Spec \s{C}\in \ob{Affine}(\Pr^L_{\st})$ for its object in the opposite category.  Given an arrow $f:\Spec \s{C}\to \Spec \s{E}$ in $\ob{Affine}(\Pr^L_{\st})$ we shall write $f^*:\s{E}\to \s{C}$ for the corresponding map in $\ob{CAlg}(\Pr^L_{\st})$.

 \begin{definition}
Let $f: \Spec \s{C}\to \Spec \s{E}$ be a map in $\ob{Affine}(\Pr^L_{\st})$.  
\begin{enumerate}

\item We say that $f$ is a closed immersion if $f^*$ has a colimit preserving fully faithful right adjoint $f_*: \s{C}\to \s{E}$ such that for $M\in \s{C}$ and $N\in \s{E}$ the natural map 
\[
f_*M\otimes N \to f_*(M\otimes f^*N)
\]
is an equivalence. 

\item We say that $f$ is an open immersion if $f^*$ has a fully faithful left adjoint $f_!: \s{C}\to \s{E}$ such that for $M\in \s{C}$ and $N\in \s{E}$ the natural map 
\[
f_!(M\otimes f^*N)\to f_!M \otimes N
\]
is an equivalence. 
\end{enumerate}
 \end{definition}

In a few words, a morphism $f: \Spec \s{C}\to  \Spec \s{E}$ in $\ob{Affine}(\Pr^L_{\st})$ is a closed immersion if and only if $f_*$ is colimit preserving, fully faithful and satisfies the projection formula (so we shall have $f_!=f_*$). Similarly, $f$ is an open immersion if and only if and only if  $f^*$ has a left adjoint $f_!$ which is fully faithful and satisfies projection formula (so we  have $f^*=f^!$). The following proposition characterizes closed and open immersions in terms of the smashing spectrum.

\begin{prop}[{\cite[Proposition 6.5]{CondensedComplex}}]
  Let $f: \Spec \s{C}\to  \Spec \s{E}$ be a morphism in $\ob{Affine}(\Pr^L_{\st})$.
  
 \begin{enumerate}
 
\item $f$ is a closed immersion if and only if there is a (necessarily unique) idempotent algebra $A\in \s{E}$ such that $1_{\s{C}}\to f^* A$ is an isomorphism,  and the induced natural  map $\Mod_A(\s{E})\to \s{C}$ is an equivalence.

\item  $f$ is an open immersion if and only if there is a (necessarily unique) idempotent algebra $A$ such that $f^* A=0$, and the induced natural map $\s{E}/\Mod_A(\s{E})\to \s{C}$ is an equivalence.

 \end{enumerate}

\end{prop}

\begin{remark}\label{rem:ComputationsSixFuncOpenClosed}
Let $\s{C}\in \ob{CAlg}(\Pr^L_{\st})$ and let $A\in\s{C}$ be an idempotent algebra with associated closed subspace $Z\subset \Sm(\s{C})$ and open complement $U$. Let $\s{C}(Z)=\Mod_A(\s{C})$ and $\s{C}(U)=\s{C}/\s{C}(U)$ denote the closed and open localizations associated to $Z$ and $U$ respectively.   For future reference we shall write explicitly the six functors for open and closed immersions.

\begin{enumerate}

\item The pullback map $\iota^*: \s{C}\to \s{C}(Z)$ is given by the base change $\iota^*M=A\otimes M$ for $M\in \s{C}$. The pushforward $\iota_*:\s{C}(Z)\to \s{C}$ is the forgetful functor and the $\iota^!$ is determined by 
\[
\iota_* \iota^! M = \iHom_{\s{C}}(A,M)
\]
for $M\in \s{C}$.
\item We let $j^*: \s{C}\to \s{C}(U)$ denote the pullback or localization map. The functors $j_!$ and $j_*$ are determined by 
\[
j_!j^*M = \fib(1\to A)\otimes M
\]
and
\[
j_*j^* = \iHom_{\s{C}}(\fib(1\to A),M)
\]
for $M\in \s{C}$.
\end{enumerate} 
\end{remark}

Closed and open immersions behave  topologically as expected along pullbacks in $\ob{Affine}(\Pr^L_{\st})$, inclusions, unions and intersections (see   \cite[Lemma 6.4 and Corollary 6.6]{CondensedComplex}). Furthermore, one gets ``for free'' the descent of the underlying symmetric monoidal categories along the topology of the locales.

\begin{theorem}[{\cite[Theorem 6.7]{CondensedComplex}}]\label{theo:openCoversSym}
  \begin{enumerate}
  
\item There is a Grothendieck topology in $\ob{Affine}(\Pr^L_{\st})$ where the coverings of $\s{C}$ are given by open localizations of $\s{C}$  whose corresponding open subsets cover $\Sm(\s{C})$. 
   
\item The identity functor $\ob{Affine}(\Pr^L_{\st})^{\op}\to \ob{CAlg}(\Pr^L_{\st})$ is a sheaf with respect to this Grothendieck topology.

\item   The posets of open (resp. closed) immersions also satisfies descent for this Grothendieck topology.
   
  \end{enumerate}
\end{theorem}

\begin{remark}\label{rem:ClosedVariant}
There is also a variant of  \Cref{theo:openCoversSym} where the covers are given by finitely many closed localizations whose closed subspaces cover $\Sm(\s{C})$. Later we will see that a more general Grothendieck topology, called the $!$-topology, has open and closed covers of symmetric monoidal categories as  particular covers.
\end{remark}

\subsubsection{Stone duality} \label{ss:StoneDuality}

We finish this section with a brief discussion of a technical but useful duality between topological spaces and distributive lattices called \textit{Stone duality}, see  \cite{StoneDuality} and \cite{zbMATH07782961}.

\begin{definition}[{\cite[Section 1.2]{zbMATH07782961}}]\label{Def:DistributiveLattice}
A \textit{lattice} $L$ is a poset for which every finite set has a supremum and an infimum. More elementary, a lattice is a tuple $(L, \vee, \wedge,\perp, \top )$ satisfying the following axioms: 
\begin{enumerate}

\item the operations $\vee$ and $\wedge$ are commutative; 

\item the operations $\vee$ and $\wedge$ are associative;

\item the operations $\vee$ and $\wedge$ are idempotent; 

\item the absorption law $a\wedge (a\vee b) = a$ and $a\vee (a\wedge b)=a$ holds for $a,b\in L$.

\item the \textit{bottom} element $\perp$ is a unit  for $\vee$ and the \textit{top} element $\top$ is a unit for $\wedge$. 

\end{enumerate}
A map $f\colon L\to M$ between lattices is a \textit{lattice} homomorphisms if if preserves the binary operations $\vee, \wedge$ and the objects $\perp, \top$, equivalently, if $f$ is a map of posets that preserve infimum and supremum elements of finite subsets. 
\end{definition}

\begin{remark}\label{RemarkAnalogy}
Following the notation of \Cref{Def:DistributiveLattice}, if $L$ is a poset for which every finite set has a supremum, one defines the operations:
\begin{enumerate}

\item  $a\vee b= \inf_{c\in L}(c\geq a \mbox{ and } c\geq b)$, that is $a\vee b$ is the supremum of $\{a,b\}$.

\item $a\vee b = \sup_{c\in L} (c\leq a \mbox{ and } c\leq b)$, that is $a\wedge b$ is the infimum of $\{a,b\}$.

\item The element $\perp$ is the supremum of  $\emptyset$, and the element $\top$ is the infimum of  $\emptyset$.

\end{enumerate}

Conversely, given a tuple $(L,\vee,\wedge,\perp,\top)$ as before, one defines the poset structure on $L$ by declaring $a\leq b$ if and only if $a\wedge b =a$, or equivalently by the absorption law, that $a\wedge b=b$. The element $\perp$ is then the minimal element of $L$ and $\top$ is the maximal element.  
\end{remark}

\begin{example}\label{ExampleTopSpaces}
Let $X$ be a topological space, attached to $X$ there are two natural lattices, namely, the lattices of open and closed subspaces $\ob{Op}(X)$ and $\ob{Cl}(X)$ respectively. Indeed, both $\ob{Op}(X)$ and $\ob{Cl}(X)$ are posets ordered by the inclusion where $\vee$ and $\wedge$ are given by intersection and union. The top and the bottom elements   are $X$ and $\emptyset$ respectively. There is a natural antiequivalence of lattices $\ob{Op}(X)^{\op}\cong \ob{Cl}(X)$ given by sending an open subspace $U$ to its complement $Z=X\backslash U$.  A continuous map $f\colon Y\to X$ of topological spaces give rise to homorphisms of lattices $f^*\colon \ob{Op}(X)\to \ob{Op}(Y)$ and $f^*\colon \ob{Cl}(X)\to \ob{Cl}(Y)$. 
\end{example}

\begin{definition}\label{def:DistributiveLattice}
A lattice $L$ is said \textit{distributive} if for all $a,b,c\in L$ one has $a\wedge(b\vee c)=(a\wedge b) \vee (a \wedge c)$, or  equivalently, if $a\vee (b\wedge c)= (a\vee b) \wedge (a\vee c)$.
\end{definition}

\begin{example}\label{ExampleDistributive}
As a prototypical example of a distributive lattice we have the poset of  of open or closed subspaces of a topological space $X$.  If in addition $X$ is spectral, see  \cite[\href{https://stacks.math.columbia.edu/tag/08YF}{Section 08YF}]{stacks-project}, then the poset $\ob{Op}(X)^{\ob{qc}}$ of quasi-compact open subspaces of $X$ is stable under finite unions and finite intersections, and therefore it forms a distributive lattice. 
\end{example}

Stone duality is then a duality between spectral spaces and distributive lattices. 

\begin{theorem}[{\cite[Theorem 6.4]{zbMATH07782961}}]\label{Theo:StoneDualityLattices}
Let $\cat{Spec}$ denote the category of spectral spaces with spectral maps and let $\cat{Lat}^{\ob{dis}}$ denote the category of distributive lattices. The functor $\cat{Spec}^{\op}\to \cat{Lat}^{\ob{dis}}$ sending a spectral space $X$ to the poset $\ob{Op}(X)^{\ob{qc}}$ of quasi-compact open subspaces is an equivalence of categories.
\end{theorem}
\begin{proof}[Sketch of the proof]
We just mention the construction of the inverse of $\ob{Op}(-)^{\ob{qc}}$ of \cite[Definition 6.3]{zbMATH07782961}, we refer to \textit{loc. cit.} for more details in the equivalence. Given a distributive lattice $L$ one constructs its Stone dual $\ob{St}(L)$ to be the topological space with underlying set  $\Hom_{\ob{Latt}}(L, [1])$ where $[1]=\{0<1\}$.  The topology on $\ob{St}(L)$ is defined as follows: given $a\in L$ one defines the subspace $U_a\subset \ob{St}(L)$ to be the set of maps of lattices $x\colon L\to [1]$ such that $x(a)=1$. Then, one checks that: 
\begin{enumerate}

\item If $a<b$ then $U_{a}\subset U_b$.

\item For $a,b\in L$ one has $U_{a\wedge b}=U_a\cap U_b $ and $U_{a\vee b}= U_a\cup U_b$. 

\item $U_{\perp}=\emptyset$ and $U_{\top}= \ob{St}(L)$. 

\end{enumerate}
The topology of $\ob{St}(L)$ is then the topology generated by the subspaces $U_a$ with $a\in L$.  
\end{proof}

\begin{example}\label{ExampleStoneProfinite}
Stone duality specializes to profinite sets as follows. Let $L$ be a lattice, an element $a\in L$ is called \textit{clopen} if there is $a^c\in L$ such that $a\wedge a^c=\perp$ and $a\vee a^c=\top$, we call $a^c$ the \textit{complement of $a$}, it is an easy excercise to show that if $a$ admits a complement then it is unique. We say that a lattice $L$ is \textit{Boolean} if all element $a\in L$ admits a complement. Let $L$ be a Boolean lattice, note that $L$ has a natural structure of a Boolean ring with product $ab:= a\wedge b$ and sum $a+b=(a\vee b)\wedge (a\wedge b)^c$, the multiplicative unit of $L$ is the element $\top$, and the zero element is $\perp$. The functor sending a Boolean lattice $(L,\vee,\wedge,\perp,\top)$ to the Boolean ring $(L,+,\cdot)$ induces an equivalence of categories 
\[
\cat{Lat}^{\ob{Boolean}}\cong \cat{Ring}^{\ob{Boolean}}
\]
between Boolean lattices and Boolean rings. Then, the Stone duality of \Cref{Theo:StoneDualityLattices} specializes to the equivalence of categories
\[
\cat{Prof}^{\op}\cong \cat{Ring}^{\ob{Boolean}}
\]
sending a profinite set $S$ to the Boolean ring $C(S,\F_2)$ of locally constant functions from $S$ to $\F_2$, the inverse of this functor sends a Boolean ring $R$ to the topological space $\Spec R$. 
\end{example}

\begin{example}\label{Example:keyStoneDualityLocale}
Let us give the prototypical example where \Cref{Theo:StoneDualityLattices} will be applied. Let $\s{C}$ be a presentable symmetric monoidal stable $\infty$-category, and  let $\ob{Idem}(\s{C})$ be the poset of idempotent algebras of $\s{C}$ such that $A\leq A'$ if $A\to A'$. By definition, we have $\ob{Sm}(\s{C})=\ob{Idem}(\s{C})^{\op}$ and we think  of $\ob{Sm}(\s{C})$ as the poset of \textit{open localizations} of $\s{C}$. Let $L\subset \ob{Sm}(\s{C})$ be a subspace of the smashing localizations stable under finite unions, finite intersections, and containing the idempotent algebras $0$ and $A$. Then $L$ is a distributive sublattice in $\ob{Sm}(\s{C})$, and by Stone duality there is a natural map of locales
\[
f\colon \ob{Sm}(\s{C})\to \ob{St}(L)
\]
uniquely determined by the fact that for $a\in L$ with associated open subspace $U_a\subset\ob{St}(L)$ one has $f^{-1}(U_a)=a\in \ob{Sm}(\s{C})$.  
\end{example}

\subsection{The ring $\bb{Z}[T]_{\sol}$}

Let $P=\bb{Z}[\bb{N}\cup\{\infty\}]/(\infty)$ be the condensed abelian group of null sequences and let $\bb{Z}[\hat{q}]$ be $P$ considered as an algebra.  The multiplication by $q=[1]$ in $\bb{Z}[\hat{q}]$ corresponds to the shift map. Since $\bb{Z}[\hat{q}]\subset \bb{Z}[[q]]$, this ring is an integral domain and  we have a short exact sequence
\[
0\to \bb{Z}[\hat{q}]\xrightarrow{1-q} \bb{Z}[\hat{q}]\to \bb{Z}[\hat{q}]/(1-q)\to 0.
\]

 By definition, $\Solid$ is the full subcategory of condensed abelian groups $M$ such that 
\[
R\iHom_{\bb{Z}}(\bb{Z}[\hat{q}]/(1-q),M)=0.
\]
This localization process fits in the general framework of  \Cref{PropStabilityF} (2) since $P$ is internally compact projective in condensed abelian groups.   From  \Cref{TheoSolidAb} (12) we know that 
\[
\bb{Z}[\hat{q}]^{L\sol}= \bb{Z}[[q]]. 
\]

Let us consider the polynomial algebra in one variable  $\bb{Z}[T]$ seen as a solid abelian group. We let $A=(\bb{Z}[T],\bb{Z})_{\sol}$ denote the induced analytic structure $\bb{Z}[T]_{\bb{Z}_{\sol}/}$. Then 
\[
A\otimes_\bb{Z} P= (\bb{Z}[T],\bb{Z})_{\sol}\otimes_{\bb{Z}} \bb{Z}[\hat{q}]= \bb{Z}[[q]][T]
\]
is a polynomial algebra over the power series ring in the variable $q$. Then, we could solidify  the variable $T$ by asking that a null-sequence $(m_n)$ in an $A$-module $M$ is ``$T$-summable'', i.e. that $\sum_{n} m_n T^n$ converges (uniquely and functorially) in $M$. This leads to the following definition.
\begin{definition}\label{SolidZT}
An object $M\in \Mod(A)$ (resp. in $\ob{D}(A)$) is said \textit{$T$-solid} if the natural map
\[
R\iHom_{A}(A\otimes_{\bb{Z}} P,M ) \xrightarrow{1-T\Shift^*} R\iHom_{A}(A\otimes_{\bb{Z}} P, M)
\]
is an equivalence. We let $\Mod(\bb{Z}[T]_{\sol})$ (resp. $\ob{D}(\bb{Z}[T]_{\sol})$) be the full subcategory of \textit{$T$-solid modules} (in the abelian or derived category respectively).  
\end{definition}

Since $A\otimes_{\bb{Z}} P$ is a compact projective $A$-algebra,  \Cref{PropStabilityF} (2) shows that $\ob{D}(\bb{Z}[T]_{\sol})$ is  an (uncompleted) analytic ring structure on $A^{\triangleright}$;  the only condition to verify is the right $t$-exactness of the localization but this is automatic being the left adjoint of the inclusion which is $t$-exact. We will do better, and we will actually compute the free solid modules $\bb{Z}[T]_{\sol}[S]$ for $S\in \Prof^{\light}$.  

To begin with, note that the fiber  sequence $A\otimes_{\bb{Z}} P \xrightarrow{1-T\Shift} A\otimes_{\bb{Z}} P \to Q $ is actually exact and is equivalent to the short exact sequence
\begin{equation}
\label{eqsequenceZLaurent}
0\to \bb{Z}[[q]][T]\xrightarrow{1-Tq} \bb{Z}[[q]][T]\to \bb{Z}((T^{-1}))\to 0
\end{equation}
where $\bb{Z}((T^{-1}))=\bb{Z}[[T^{-1}]][T]$ is the algebra of Laurent power series in the variable $T^{-1}$. In particular, since $\bb{Z}[[X]]$  is an idempotent $(\bb{Z}[X],\bb{Z})_{\sol}$-algebra, one deduces that $\bb{Z}((T^{-1}))$ is an idempotent $(\bb{Z}[T],\bb{Z})_{\square}$-algebra. We get the following proposition.

\begin{prop}\label{PropSolidZTidempotentComplement}
The category $\ob{D}(\bb{Z}[T]_{\sol})$ is the localization of $\ob{D}((\bb{Z}[T],\bb{Z})_{\sol})$ with respect to the idempotent algebra $\bb{Z}((T^{-1}))$. More precisely, we have a semi-orthogonal decomposition 
\[
\ob{D}(\bb{Z}((T^{-1}))_{\bb{Z}_{\sol}/})\xrightarrow{\iota_*} \ob{D}((\bb{Z}[T],\bb{Z})_{\sol})\xrightarrow{j^*} \ob{D}(\bb{Z}[T]_{\sol})
\]
where $\iota_*$ is the natural inclusion and $j^*$ is the localization functor. We let $\iota^*$ be the base change $\bb{Z}((T^{-1}))\otimes_A -$ and let $j_*$ be the right adjoint of $j^*$.
\end{prop}

From the general non-sense of smashing localizations in presentably symmetric monoidal stable $\infty$-categories (see \cite[Lecture V]{CondensedComplex} and  \Cref{rem:ComputationsSixFuncOpenClosed}), we can explicitly compute the functor $j_*$:
\[
j_*j^* M = R\iHom_{A}( \mathrm{fib}(\bb{Z}[T]\to \bb{Z}((T^{-1}))), M),
\]
in particular the functor $j_*j^*$ commutes with limits and colimits.

Using the resolution \eqref{eqsequenceZLaurent} of $\bb{Z}((T^{-1}))$ one can easily compute that $R\iHom_A(\bb{Z}((T^{-1})),\bb{Z}[T])=0$ so that $j_*j^* \bb{Z}[T]= \bb{Z}[T]$. Indeed, we can write 
\[
R\iHom_{A}(\bb{Z}[[q]][T], \bb{Z}[T]) = \bb{Z}[T]((q))/ q\bb{Z}[T][[q]]
\]
as $\bb{Z}[q,T]$-module, and multiplication by $1-Tq$ is invertible on $\bb{Z}[T][[q]]$.

With this computation we can show the following theorem (for more details see \cite[Lecture VIII]{ClausenScholzeCondensed2019}).

\begin{theorem}\label{TheoSolidZT}
The full subcategory $\ob{D}(\bb{Z}[T]_{\sol})\subset \ob{D}((\bb{Z}[T],\bb{Z})_{\sol})$ defines an analytic ring structure on $\bb{Z}[T]$. For $S=\varprojlim_n S_n$ a light profinite set written as limit of finite sets, the free $\bb{Z}[T]_{\sol}$-module generated by $S$ is given by 
\[
\bb{Z}[T]_{\sol}[S]=\varprojlim_{n\in \bb{N}} \bb{Z}[T][S_n].
\]
\end{theorem}
\begin{proof}
Let us write $A=(\bb{Z}[T],\bb{Z})_{\sol}$. Conditions (1) and (2) of  \Cref{DefinitionAnalyticRing} follow immediately from   \Cref{PropStabilityF}   and the fact that the object $A\otimes_{\bb{Z}} P /(1-Tq)$ is compact in $\ob{D}((\bb{Z}[T],\bb{Z})_{\sol})$. It is left to show condition (3), this one follows from the computation of $\bb{Z}[T]_{\sol}[S]$. Recall that, if $S$ is infinite,  $\bb{Z}_{\sol}[S]\cong \prod_{\bb{N}}\bb{Z}$. Since $\ob{D}(\bb{Z}[T]_{\sol})$ is the localization with respect to the objects in $\ob{D}(\bb{Z}((T^{-1})))$, it suffices to show that the cofiber $Q$ of the map $(\prod_{\bb{N}} \bb{Z})[T]\to \prod_{\bb{N}} (\bb{Z}[T]) $ is a $\bb{Z}((T^{-1}))$-module. But one has that
\[
\prod_{\bb{N}}(\bb{Z}[T])/((\prod_{\bb{N}} \bb{Z})[T]) = \prod_{\bb{N}} (\bb{Z}((T^{-1})))/((\prod_{\bb{N}} \bb{Z})\otimes_{\bb{Z}_{\sol}} \bb{Z}((T^{-1}))).
\]

We deduce that 
\[
j_*j^* Q = R\iHom_{\bb{Z}[T]}( \mathrm{fib}(\bb{Z}[T]\to \bb{Z}((T^{-1}))), Q) = 0
\]
as $Q$ is a $\bb{Z}((T^{-1}))$-module. We get that 
\[
\begin{aligned}
\bb{Z}[T]_{\sol}[S] & = j_*j^* (\bb{Z}_{\sol}[S] [T]) \\
					 & = j_*j^* (\varprojlim_n \bb{Z}[T][S_n]) \\ 
					 & = \varprojlim_n (j_* j^* \bb{Z}[T][S_n]) \\ 
					 & = \varprojlim_n \bb{Z}[T][S_n]\\
\end{aligned}
\]
finishing the proof of the theorem. 
\end{proof}

\subsection{Solid rings of finite type algebras}

Let us now generalize the construction of the ring $\bb{Z}[T]_{\sol}$ to finite type algebras over $\bb{Z}$.

\begin{definition}\label{DefinitionRsolidfinitetype}
\begin{enumerate}
\item Let $R$ be a solid animated algebra and let $r\in \pi_0(R)$. We say that a solid $R$-module $M$ is \textit{$r$-solid} if for any  lift $\bb{Z}[T]\to R$ of $r$, the restriction of $M$ to a $\bb{Z}[T]$-module is $T$-solid (by \Cref{PropInvarianceHomotopy} this condition is independent of the lift of $r$ to a map $\Z[T]\to R$).

\item Let $R$ be a finite type algebra over $\bb{Z}$. We let $R_{\sol}$ be the analytic ring structure on $R$ making an $R$-module $M$ complete if and only if $M$ is $r$-solid for all $r\in R$. 
\end{enumerate}
\end{definition}

\begin{theorem}\label{theoSolidR}
Let $R$ be a finite type algebra over $\bb{Z}$. Then for $S=\varprojlim_n S_n$ a light profinite set the natural map 
\[
R_{\sol}[S] \to \varprojlim_n R[S_n]
\]
is an isomorphism. In particular, $\ob{D}(R_{\sol})$ is the derived category of its heart and $\otimes^L_{R_{\sol}}$ is the left derived functor of $\otimes_{R_{\sol}}$. Moreover, $\prod_{\N} R$ is flat for the $R_{\sol}$-tensor product. 
\end{theorem}

Our first task is to compute the free solid generators $R_{\sol}[S]$ for $S\in \Prof^{\light}$. Note that any algebra of finite type is a quotient of a polynomial algebra, let us then start with those:

\begin{proposition}\label{PropSolidPolynomials}
Let $T_1,\ldots, T_n$ be a set of variables, then the natural map 
\begin{equation}
\label{eqMapAnayticRingPolySolid}
\bb{Z}[T_1]_{\sol}\otimes_{\bb{Z}_{\sol}} \cdots \otimes_{\bb{Z}_{\sol}} \bb{Z}[T_n]_{\sol}\to \bb{Z}[T_1,\ldots, T_n]_{\sol}
\end{equation}
is an isomorphism. In other words, a $\Z$-solid $\Z[T_1,\ldots, T_n]$-module is $\Z[T_1,\ldots, T_n]$-solid if and only if it is $T_1,\ldots, T_n$-solid. Furthermore, for $S=\varprojlim_n S_n$ a light profinite set written as a countable limit of finite sets we have 
\[
\bb{Z}[T_1,\ldots, T_n]_{\sol}[S]=\varprojlim_{n} \bb{Z}[T_1,\ldots, T_n]_{\sol}[S_n]. 
\]
\end{proposition}
\begin{proof}
Let $A= \bb{Z}[T_1]_{\sol}\otimes_{\bb{Z}_{\sol}} \cdots \otimes_{\bb{Z}_{\sol}} \bb{Z}[T_n]_{\sol}$, it is the analytic ring structure on $\bb{Z}[T_1,\ldots, T_n]$ making a module $A$-complete if and only it it is $T_i$-solid for all $i=1,\ldots, n$. We clearly have that $\ob{D}(\bb{Z}[T_1,\ldots, T_n]_{\sol}) \subset \ob{D}(A)\subset \ob{D}(\bb{Z}[T_1,\ldots, T_n]^{\cond})$. We need to show the opposite inclusion $\ob{D}(A)\subset \ob{D}(\bb{Z}[T_1,\ldots, T_n]_{\sol})$, i.e. that if a solid  $\bb{Z}[T_1,\ldots, T_n]$-module is  $T_i$-solid for all $i=1,\ldots, n$, then it is $p(T)$-solid for all $p(T)\in \bb{Z}[T_1,\ldots, T_n]$. 

As a first step, let us compute the compact projective generators of the ring $A$. We show by induction on the number of variables that  for $S=\varprojlim_k S_k$ profinite
\[
A[S]=\varprojlim_n \bb{Z}[T_1,\ldots, T_n] [S_k],
\]
the case $n=1$ being  \Cref{TheoSolidZT}. Suppose that the claim follows for $n$ and consider $B=\bb{Z}[T_1,\ldots, T_n,T_{n+1}]_{A/}$ the induced analytic rings structure. Let $C=\bb{Z}[T_1]_{\sol}\otimes_{\bb{Z}_{\sol}} \cdots \otimes_{\bb{Z}_{\sol}} \bb{Z}[T_{n+1}]_{\sol}$. By definition $\ob{D}(C)\subset \ob{D}(B)$ is the full subcategory of objects that are $\bb{Z}[T_{n+1}]$-solid (since they are already $\bb{Z}[T_i]$-solid for all $i=1,\ldots,n$). By  \Cref{PropSolidZTidempotentComplement}  an object $M\in \ob{D}(A)$ is $\bb{Z}[T_{n+1}]$-solid if and only if 
\[
R\iHom_{\bb{Z}[T_{n+1}]}( \bb{Z}((T_{n+1}^{-1})),M)=0.
\]
By taking base change along $(\bb{Z}[T_{n+1}],\bb{Z})_{\sol} \to B$ this is equivalent to   
\[
R\iHom_{\bb{Z}[T_1,\ldots, T_{n+1}]}( B\otimes_{\bb{Z}[T_{n+1}]} \bb{Z}((T_{n+1}^{-1})), M)=0.
\]
Recall that we have the resolution
\[
0\to \bb{Z}[[q]][T_{n+1}]\xrightarrow{1-qT_{n+1}} \bb{Z}[[q]][T_{n+1}]\to \bb{Z}((T^{-1}_{n+1}))\to 0.
\]
Then, by induction, we have that 
\[
B\otimes_{\bb{Z}[T_{n+1}]} \bb{Z}((T_{n+1}^{-1}))= \bb{Z}[T_1,\ldots, T_n]((T_{n+1}^{-1})),
\]
which is an idempotent $B$-algebra. Then, the same argument of  \Cref{TheoSolidZT} will show that 
\[
C_{\sol}[S]= \varprojlim_n \bb{Z}[T_1,\ldots, T_{n+1}] [S_k]
\]
as wanted. 

Now, note that any discrete $\bb{Z}[T_1,\ldots, T_n]$-module is immediately $\bb{Z}[T_1,\ldots, T_n]_{\sol}$-complete since it is $a$-solid for any  $a\in \bb{Z}[T_1,\ldots, T_n]$. This shows that $A[S]\cong \prod_{I} \bb{Z}[T_1,\ldots, T_n]$ is $\bb{Z}[T_1,\ldots, T_n]_{\sol}$-complete (since complete modules are stable under products), and so that any complete $A$-module is $\bb{Z}[T_1,\ldots, T_n]_{\sol}$-complete (being stable under colimits). One deduces that 
\[
\ob{D}(A)\subset \ob{D}(\bb{Z}[T_1,\ldots ,T_n]_{\sol})\subset \ob{D}(\bb{Z}[T_1,\ldots, T_n]^{\cond})
\] 
and so that we have the equivalence
\[
\ob{D}(A)=\ob{D}(\bb{Z}[T_1,\ldots, T_n]_{\sol}),
\]
proving that the map \eqref{eqMapAnayticRingPolySolid} is  an equivalence. 
\end{proof}

\begin{corollary}\label{CorollaryFreeProfiniteR}
Let $R$ be a finite type algebra and let $\bb{Z}[T_1,\ldots, T_n]\to R$ be a surjection. Then 
\[
R_{\sol}= R_{\bb{Z}[T_1,\ldots, T_n]_{\sol}/}
\]
has the induced analytic structure. Moreover, for $S=\varprojlim_k S_k$ a light profinite set we have 
\[
R_{\sol}[S]=\varprojlim_k R[S_n]. 
\]
\end{corollary}
\begin{proof}
By definition, a $M$-module is $\bb{Z}[T_1,\ldots, T_n]_{\sol}$-complete if it is $a$-complete for all $a\in \bb{Z}[T_1,\ldots, T_n]$. By definition of $R_{\sol}$ this shows that it has the induced analytic structure. In particular, 
\[
R_{\sol}[S]=R\otimes_{\bb{Z}[T_1,\ldots, T_n]_{\sol}}^L \bb{Z}[T_1,\ldots, T_n]_{\sol}[S].
\]
We will prove a more general fact: let $M$ be a finite type $\bb{Z}[T_1,\ldots, T_n]$-module, then 
\[
M\otimes^L_{\bb{Z}[T_1,\ldots, T_n]_{\sol}} \prod_k (\bb{Z}[T_1,\ldots, T_n]_{\sol}) = \prod_k M. 
\]

Indeed, consider a finite projective resolution $P_{\bullet}\to M$ where all the terms $P_n$ are finitely many copies of $\bb{Z}[T_1,\ldots, T_n]$. Then $M\otimes^L_{\bb{Z}[T_1,\ldots, T_n]_{\sol}} \prod_k (\bb{Z}[T_1,\ldots, T_n]_{\sol})$ is equivalent to the complex $P_{\bullet} \otimes^L_{\bb{Z}[T_1,\ldots, T_n]_{\sol}} \prod_k (\bb{Z}[T_1,\ldots, T_n]_{\sol})$. Since each term $P_n$ is a finite free module we actually have 
\[
P_{\bullet} \otimes^L_{\bb{Z}[T_1,\ldots, T_n]_{\sol}} \prod_k (\bb{Z}[T_1,\ldots, T_n]_{\sol}) = \prod_{k} P_{\bullet}.
\]
Since countable products are exact we have an equivalence
\[
\prod_{k} P_{\bullet} \xrightarrow{\sim} \prod_k M.
\]
\end{proof}

\begin{proof}[Proof of  \Cref{theoSolidR}]
The claim about the free objects on profinite sets is  \Cref{CorollaryFreeProfiniteR}. The fact that $\ob{D}(R_{\sol})$ is the derived category of its heart and that $\otimes^L_{R_{\sol}}$ is the left derived functor of $\otimes_{R_{\sol}}$ follows the same argument of the analogue statements in  \Cref{TheoSolidAb}. The final statement about flatness of $\prod_n R$ will be proven in  \Cref{propFlatnesprodR}.
\end{proof}

It is left to prove flatness of $\prod_n R$, this will follow a similar argument as the one for $\bb{Z}$. 

\begin{definition}\label{DefinitionFinitelyPresentedstuff}
Let $R$ be finite type algebra.  An $R_{\sol}$-module is \textit{finitely generated} if it is a quotient of $\prod_{I} R$ for $I$ a countable set.  An $R_{\sol}$-module is of finite presentation if it is a cokernel of a map $\prod_{I} R\to \prod_{J} R$ with $I$ and $J$ countable sets. 
\end{definition}

We want to have a good understanding of finitely presented $R_{\sol}$-modules. Under some topological hypothesis these are easier to describe:

\begin{lemma}\label{LemmaquasiSepCoherentR}
Let  $R$ be a finite type algebra. If $M\in \Mod_{R_{\sol}}$ is quasi-separated, then the following are equivalent. 
\begin{itemize}
\item[i.] $M$ is of finite presentation.

\item[ii.] $M$ is finitely generated. 

\item[iii.] $M=\varprojlim M_n$ is a limit of finitely generated discrete $R$-modules  with surjective transition maps. 

\end{itemize}
\end{lemma}

\begin{proof}
Let $(M_n)_n$ be a projective system  as in (iii).  Since $R$ is noetherian we can find  a right exact resolution 
\[
R^{s_n}\to R^{k_n}\to M_n\to 0
\]
for each $n$. By constructing the resolution step by step, we can construct liftings 
\[
\begin{tikzcd}
R^{s_{n+1}}  \ar[d] \ar[r]& R^{k_{n+1}}  \ar[d]\ar[r]& M_{n+1} \ar[d] \ar[r] & 0 \\
R^{s_{n}} \ar[r] & R^{k_{n}} \ar[r] & M_{n} \ar[r] & 0
\end{tikzcd}
\]
such that all the vertical maps are surjective. Taking limits we find a right exact sequence
\[
\prod_{\bb{N}} R \to \prod_{\bb{N}} R\to M\to 0
\]
proving (iii)$ \Rightarrow$ (i). 

It is clear that (i) $\Rightarrow$ (ii). It is left to show that (ii) implies (iii). We have a surjection $\prod_{I} R\to M\to 0$ with kernel $K$. We can assume without loss of generality that $I=\N$. The space $\prod_{\N} R$ arises from a metrizable topological space, and since $M$ is quasi-separated the space $K\subset \prod_{\N} R$ is closed. Writing $\prod_{\N} R = \varprojlim_{n} R^n$ as the limit of finite free $R$-modules, we see that $K=\varprojlim_n K_n$ where $K_n$ is the image of $K$ in the projection $\prod_{\N} R\to R^n$. This shows that 
\[
M=\varprojlim_n R^n/K_n
\]
proving what we wanted. 
\end{proof}

\begin{prop}\label{propFlatnesprodR}
Let $R$ be a finite type $\Z$-algebra.  Let $\Mod(R_{\sol})$ be the abelian category of $R_{\sol}$-modules and let $\Mod(R_{\sol})^{fp}\subset \Mod(R_{\sol})$ be the full subcategory of finitely presented solid modules. The following hold:

\begin{enumerate}

\item We have $\Mod(R_{\sol})=\Ind(\Mod(R_{\sol})^{fp})$.

\item The category $\Mod(R_{\sol})^{fp}$ of finitely presented modules is an abelian category stable under all kernels, cokernels and extensions. 

\item Any finitely presented module $M$ is pseudo-compact, namely, it has a resolution of the form $P_{\bullet}\to M$ with $P_{\bullet}$ a free $R_{\sol}$-module on a light profinite set $S$. 

\item The $R_{\sol}$-module $\prod_{\N} R$ is flat. 

\end{enumerate}
\end{prop}
\begin{proof}

\begin{enumerate}

\item The category $\Mod_{R_{\sol}}$ has a compact projective generator given by $\prod_{\bb{N}} R$. This formally shows that finitely presented modules are the compact objects in $\Mod_{R_{\sol}}$ and the description as inductive category of (1). 

\item By the standard arguments in commutative algebra it suffices to show that any finitely generated module  $M$ of $\prod_{\bb{N}} R$ is actually finitely presented. But then $M$ is quasi-separated and   \Cref{LemmaquasiSepCoherentR} implies that $M$ is finitely presented.

\item Let $M$ be a finitely presented module and consider a right exact sequence
\[
\prod_{\bb{N}} R \xrightarrow{f}  \prod_{\bb{N}} R \to M\to 0.
\]
By part (2) the kernel $K=\ker f$ is a finitely presented module. An inductive argument allow us to construct the resolution  $P_{\bullet}\to M$ of $M$ where each term is isomorphic to $\prod_{\bb{N}} R$.

\item It suffices to show that for any finitely presented module $M$ we have $M\otimes^L_{R_{\sol}} \prod_{\bb{N}} R= \prod_{\bb{N}} M$ as then the tensor product will be exact. This follows from the fact that $P\otimes_{R_{\sol}}^L\prod_{\bb{N}} R= \prod_{\bb{N}} P$ for $P\cong \prod_{\bb{N}} R$ and the resolution of $M$ of part (3). 
\end{enumerate}
\end{proof}

\begin{remark}\label{RemarkTensorPseudoCoherent}
Let $R$ be a finite type $\bb{Z}$-algebra.  Then for any pseudo-compact solid $R_{\square}$-module we have $M\otimes_{R_{\square}}^L \prod_{\bb{N}} R=\prod_{\bb{N}} M$, namely, this follows from the computation of the solid tensor product of countable products of $R$, and the fact that  $M$ has a projective resolution  whose terms are given $\prod_{\bb{N}} R$

As a special consequence, if $R$ is an animated algebra of almost finite presentation, that is  with $\pi_0(R)$ a finitely generated $\bb{Z}$-algebra and  $\pi_i(R)$ a finite $\pi_0(R)$-module for all $i$, then the induced analytic ring structure $R_{\square}:=R_{\pi_0(R)_{\square}/}$ from   \Cref{PropInvarianceHomotopy} is such that 
\[
R_{\square}[S]= \varprojlim_n R[S_n],
\]
so that $\prod_{I} R\otimes_{R_{\square}} \prod_J R = \prod_{I\times J} R$ for countable sets $I$ and $J$. 
\end{remark}

The solid ring structures are also independent of integral extensions. 

\begin{corollary}\label{cor:SolidIntegralInvariance}
Let $R\to A$ be an integral map of finitely generated algebras. Then the natural map of analytic rings $A_{R_{\square}/}\to A_{\square}$ is an isomorphism.
\end{corollary}
\begin{proof}
By hypothesis both $R$ and $A$ are finitely generated $\bb{Z}$-algebras with $A$ integral over $R$, so that $A$ is a finite $R$-module. By  \Cref{RemarkTensorPseudoCoherent} we find that $A\otimes_{R_{\square}} \prod_{\bb{N}} R = \prod_{\bb{N}} A$, which implies the corollary. 
\end{proof}

Finally, the following lemma is useful to understand tensor products of solid modules.

\begin{lemma}\label{LemTensorSumsProd}
Let $R$ be an algebra of  finite type   over $\Z$. In the following all the tensor products are derived. 

\begin{enumerate}

\item Let $P$ and $Q$ be pseudo-compact $R$-modules, then $P\otimes_{R_{\sol}} Q$ is pseudo-compact.

\item  let $(P_{n})_{n\in \N}$ be a countable family of connective pseudo-compact $R$ modules. Then $\prod_{n} P_n$ is pseudo-compact. Moreover, if $Q$ is another pseudo-compact object then $Q\otimes_{R_{\sol}} \prod_n P_n=\prod_n (Q\otimes_{R_{\sol}} P_n)$. 

\item Let $(P_n)_{n\in N}$ and $(Q_n)_{n}$ be sequential projective diagrams of connective pseudo-compact $R_{\sol}$-modules, then $\varprojlim_{n} P_n$ and $\varprojlim_n Q_n$ are pseudo-compact, $-1$-connective and  
\[
\bigg( \varprojlim_{n} P_n \bigg)\otimes_{R_{\sol}} \bigg( \varprojlim_{n} Q_n \bigg)=\varprojlim_n  P_n\otimes_{R_{\sol}} Q_n.
\]

\item Let $(P_{n,m})_{n,m\in \N}$ be a family of connective pseudo-compact $R$-modules and let $Q$ be a pseudo-compact $R$-module. Then the natural map of solid $R$-modules\footnote{We thank Yutaro Mikami for corrections in a previous version of this lemma.} 
\[
\bigg(\prod_{n} \bigoplus_{m} P_{n,m} \bigg)\otimes_{R_{\sol}} Q \xrightarrow{\sim}  \prod_{n} \bigoplus_{m}  \big( P_{n,m}\otimes_{R_{\sol}} Q \big)
\]
is an equivalence. 

\end{enumerate}

\end{lemma}
\begin{proof}

Parts (1) and (2) follow from \Cref{RemarkTensorPseudoCoherent} as any pseudo-compact object has a projective resolution with terms given by countable products of $R$, we left the details to the reader. Part (3) follows from (1) and (2) after writing the limit as the fiber 
\[
\varprojlim_n P_n\to \prod_{n\in \N} P_n \to  \prod_{n\in \N} P_n.
\]

We now prove part (4). By the AB6 property we can write 
\[
\bigg(\prod_{n} \bigoplus_{m} P_{n,m} \bigg)=\varinjlim_{f\colon \N\to \N} \prod_{n} \bigoplus_{m\leq f(n)} P_{n,m}
\]
where the maps $f\colon \N\to \N$ have the point-wise partial order. We see that 
\[
\begin{aligned}
\bigg(\prod_{n} \bigoplus_{m} P_{n,m} \bigg)\otimes_{R_{\sol}} Q & = \varinjlim_{f\colon \N\to \N}  \bigg( \prod_{n} \bigoplus_{m\leq f(n)} P_{n,m} \bigg) \otimes_{R_{\sol}} Q  \\ 
& =   \varinjlim_{f\colon \N\to \N}  \bigg( \prod_{n} \bigoplus_{m\leq f(n)} P_{n,m} \otimes_{R_{\sol}} Q\bigg)  \\ 
& = \prod_{n} \bigoplus_{m} \big( P_{n,m}\otimes_{R_{\sol}}Q\big)
\end{aligned}
\]
where in  the first and  last equivalences we use the AB6 property, and in the second equivalence we use part (2) of the lemma.  
\end{proof}

\begin{remark}\label{RemarkGeneralizationAnimatedRing}
\Cref{LemTensorSumsProd} generalizes to animated rings of almost finite presentation in the sense of \Cref{RemarkTensorPseudoCoherent}. Indeed, if $R$ is an animated ring of almost finite presentation, a $R_{\sol}$-module $M$ is said pseudo-compact if it satisfies one of the following equivalent conditions (see \cite[Definition 9.9 and Proposition 9.10]{CondensedComplex} for a related discussion):
\begin{enumerate}
\item For a filtered colimit $N=\varinjlim_i N_i$ of uniformly cohomologically  left bounded  $R_{\sol}$-modules the natural map 
\begin{equation}\label{eqihniahenoawfraq}
\varinjlim_i R\iHom_R(M, N_i)\to R\iHom_R(M, \varinjlim_i N_i)
\end{equation}
is an isomorphism. 

\item $M$ is $(-k)$-connective for some $k\in \N$ and $\tau_{\leq n} (M[k])$ is a $n$-truncated compact  $\tau_{\leq n} R_{\sol}$-module for all $n\in \N$. 

\item  There is some $k\in \Z$ such that  $M$ can be written as $M=\varinjlim_{n\in \Z_{\geq k}} M_n$ with $M_n$ a compact  $R_{\sol}$-module with homology supported in $[k,\infty)$ such that $\ob{cofib}(M_n\to M)$ is $(n+1)$-connective and for all $n\geq k$ there is a fiber sequence $M_{n}\to M_{n+1}\to R_{\sol}[S_{n+1}][n+1]$ with $S_{n+1}$ a light profinite set.   

\end{enumerate}

We give a brief sketch of the equivalence between (1)-(3) above. For $(1)\Rightarrow (2)$, the condition on the Hom spaces imply that $M$ must be cohomologically right bounded, i.e. $(-k)$-connective for some $k\in \N$. Indeed, by right completeness of $\ob{D}(R_{\sol})$ we have that $\tau_{\leq 0} M=\varinjlim_{n\in \N} \tau_{[-n,0]} M$, so the natural map $M\to \tau_{\leq 0} M $ factors through one of the truncations $M\to \tau_{[-k,0]}M$ which imply that $M$ must be $(-k)$-connective.  If $M$ is $(-k)$-connective,  it is also clear that condition (1) imply that $\tau_{\leq n} (M[k])$ is  compact as $n$-truncated connective module over $\tau_{\leq n}R_{\sol}$. 
 For $(2)\Rightarrow (3)$, we can assume without loss of generality that $M$ is connective, since $\pi_0(M)$ is a compact static $\pi_0(R_{\sol})$-module, there is a light profinite set $S_0$ and an epimorphism $R_{\sol}[S_0]\to M$; we set $M_{0}:=R_{\sol}[S_0]$ and $Q_1:=\ob{cofib}(M_0\to M)$. Then $Q_1$ also satisfies the condition (2) and it is $1$-connective, we can then find a light profinite set $S_1$ and a map $R_{\sol}[S_1][1]\to Q_1$ such that $Q_2:=\ob{cofib}(R_{\sol}[S_1][1]\to Q_1)$ is $2$-connective. We then define $M_1:=\ob{fib}(M\to Q_2)$, and one immediately checks from the definition a fiber sequence $M_0\to M_1\to R_{\sol}[S_1][1]$. By reapeating this construction inductively, one deduces the presentation of $M$ as a colimit in (3). Finally, for $(3)\Rightarrow (1)$, if $N=\varinjlim_i N_i$ is a filtered colimit of $0$-truncated objects in $\ob{D}(R_{\sol})$ and $M=\varinjlim_n M_n$ is as in (3) with all $M_n$ connective modules (which we can assume without loss of generality by shifting), one has that the Hom spaces \eqref{eqihniahenoawfraq} are $0$-truncated, and that their $\pi_{-m}$ homotopy group for $m\in \M$ only depends on $\tau_{\leq m}M=\tau_{\leq m} M_m$ which is an $m$-truncated compact $\tau_{\leq m} R_{\sol}$-module; this yields the equivalence of \eqref{eqihniahenoawfraq} proving what we wanted.

Going back to our original discussion, one can use the explicit presentation of a pseudo-compact $R_{\sol}$-module of (3) above to prove the statements (1)-(4) of  \Cref{LemTensorSumsProd}, we leave the details to the interested reader. 
\end{remark}

We deduce the following proposition

\begin{proposition}\label{PropCompletionSolidTensorProduct}
Let $R$ be an algebra of finite type over $\Z$ and let   $I\subset R$ be a finitely generated ideal of $R$. In the following all tensor products are derived. 

\begin{enumerate}

\item   Let $\widehat{R}_I$ be the derived $I$-adic completion of $R$ (note that as $R$ is noetherian this agrees with the classical $I$-adic completion). Then $\widehat{R}_I$ is pseudo-compact as solid $R$-module and $\widehat{R}_I\otimes_{R_{\sol}}\prod_{\N} R= \prod_{\N} \widehat{R}_{I}$ (for the definition of derived complete modules see \cite[Definition 2.12.3]{MannSix}).

\item Let $J\subset R$ be another finitely generated ideal. Then $ \widehat{R}_I\otimes_{R_{\sol}} \widehat{R}_J = \widehat{R}_{I+J}$. In particular, $\widehat{R}_I$ is an idempotent  $R_{\sol}$-algebra. 

\item Let $N$ and $M$ be connective  derived $I$-adically complete $R_{\sol}$-modules. Then $N\otimes_{R_{\sol}} M$ is derived $I$-adically complete.

\end{enumerate}

\end{proposition}
\begin{proof}
For part (1), let $I=(f_1,\ldots,f_k)$, as $R$ is noetherian  we can write 
\[
\widehat{R}_{I} = \varprojlim_{n} R/^{\bb{L}}(f_1^n,\ldots, f_k^n)
\]
where the derived quotient $R/^{\bb{L}}(f_1^n,\ldots, f_k^n)$ is defined to be the tensor product of animated rings $R\otimes_{\Z[T_1,\ldots,T_k]} \Z$ where $T_i\mapsto f_i^n$ in $R$ and $T_i\mapsto 0$ in $\Z$. By \Cref{LemTensorSumsProd} we deduce that $\widehat{R}_{I}$ is pseudo-compact as solid $R$-module and that 
\[
\widehat{R}_{I}\otimes_{R_{\sol}} \prod_{\N} R = \varprojlim_{n} \prod_{\N} R/^{\bb{L}}  (f_1^n,\ldots, f_k^n)= \prod_{\N}\widehat{R}_{I}. 
\]

For Part (2), write $J=(g_1,\ldots, g_s)$, then it  also follows formally from \Cref{LemTensorSumsProd} that 
\[
\widehat{R}_I\otimes_{R_{\sol}} \widehat{R}_J = \varprojlim_{n} \bigg( R/^{\bb{L}}(f_1^n,\ldots, f_k^n)\otimes_{R} R/^{\bb{L}}(g_1^n,\ldots, g_s^n)  \bigg) = \widehat{R}_{I+J}. 
\]

Finally, for part (3), by induction we can assume without loss of generality that $I=(f)$ is generated by a single element.    both $N$ and $M$ admits a resolution by a complex of objects of the form $\bigoplus_{K} \prod_{\N} R$ with $K$ an arbitrary large index set. Thus, to prove the claim we can assume without loss of generality that $N$ and $M$ are of the form $\widehat{\bigoplus}_{K} \prod_{\N} \widehat{R}_I$ where the completed direct sum is with respect to the $I$-adic topology. Since taking derived $I$-completions commutes with $\aleph_1$-filtered colimits, we can assume without loss of generality that $K$ is countable, and even that $K=\N$.

 We have then a fiber sequence 
\[
 \bigg(\bigoplus_{\N} \prod_{\N} R \bigg) [[T]] \xrightarrow{T-f} \bigg(\bigoplus_{\N} \prod_{\N} R \bigg) [[T]]  \to \widehat{\bigoplus}_{\N}  \prod_{\N} \widehat{R}_{I}. 
\]
From  \Cref{LemTensorSumsProd} (4) and the idempotency of $\widehat{R}_I$ for the tensor product $\otimes_{R_{\sol}}$,  we deduce that $\widehat{\bigoplus}_{\N}  \prod_{\N} \widehat{R}_{I} = \widehat{\bigoplus}_{\N} \widehat{R}_I \otimes_{R_{\sol}} \prod_{\N} \widehat{R}_I$. Hence, we are reduced to show that the natural map 
\[
\widehat{\bigoplus}_{\N} \widehat{R}_{I}\otimes_{R_{\sol}} \widehat{\bigoplus}_{\N} \widehat{R}_{I}\to \widehat{\bigoplus}_{\N\times \N} \widehat{R}_I
\]
is an isomorphism. Now write $\Z[T]\to \Z[T,X_1,\ldots, X_k]\to R$ where $\Z[T,X_1,\ldots, X_k]\to R$ is surjective and $T\mapsto f\in R$. We claim that 
\begin{equation}\label{eqooajwpbhaiwsrea}
\widehat{\bigoplus}_{\N} \Z[[T]]\otimes_{\Z[T]_{\sol}} R_{\sol}= \widehat{\bigoplus}_{\N} \widehat{R}_I,
\end{equation}
provided the claim, we can assume without loss of generality that $R=\Z[T]$ and $I=(T)$.

Let us prove the claim first, note that 
\[
\widehat{\bigoplus}_{\N} \Z[[T]]\otimes_{\Z[T]_{\sol}} R_{\sol} = \big( \widehat{\bigoplus}_{\N} \Z[[T]]\otimes_{\Z[T]_{\sol}} \Z[T,X_1,\ldots, X_n]_{\sol} \big) \otimes_{\Z[T,X_1,\ldots, X_n]_{\sol}} R_{\sol}.
\]
By \Cref{LemmaBaseChangePolynomial} the base change along $\Z[T]_{\sol}\to \Z[T,X_1,\ldots, X_n]_{\sol}$ commutes with limits. This implies that it sends $T$-adically complete objects to $T$-adically complete objects, and we see immediately that 
\[
 \widehat{\bigoplus}_{\N} \Z[[T]]\otimes_{\Z[T]_{\sol}} \Z[T,X_1,\ldots, X_n]_{\sol}  = \widehat{\bigoplus}_{\N} \Z[X_1,\ldots, X_n][[T]]. 
\]
Since $R$ is pseudo-coherent as $\Z[T,X_1,\ldots, X_n]$-module, it is pseudo-compact as solid module and the presentation 
\[
0\to \big(\bigoplus_{\N} \Z[T,X_1,\ldots, X_n]\big)[[Y]] \xrightarrow{Y-T}\big(\bigoplus_{\N} \Z[T,X_1,\ldots, X_n]\big)[[Y]] \to \widehat{\bigoplus}_{\N} \Z[X_1,\ldots, X_n][[T]]\to 0
\]
together with \Cref{LemTensorSumsProd} (4) yields the equivalence \eqref{eqooajwpbhaiwsrea}.

By the previous discussion, we are reduced to show that the natural map 
\begin{equation}\label{eqniehnlakwawr}
\widehat{\bigoplus}_{\N} \Z[[T]]\otimes_{\Z[T]_{\sol}} \widehat{\bigoplus}_{\N} \Z[[T]]\xrightarrow{\sim} \widehat{\bigoplus}_{\N\times \N} \Z[[T]]
\end{equation}
is an equivalence (we will see that the previous holds true also for the tensor product of the analytic ring $(\Z[T],\Z)_{\sol}$), this requires an explicit computation as follows.  We can write  as condensed $\Z[T]$-modules
\[
\widehat{\bigoplus}_{\N} \Z[[T]] = \varinjlim_{\substack{f\colon \N\to \N \\ f(n)\to \infty \\ n\to \infty}} \prod_{n\in \N} T^{f(n)}\Z[[T]].
\]
Thus, the left hand side term of \eqref{eqniehnlakwawr} is equivalent to 
\[
\widehat{\bigoplus}_{\N} \Z[[T]]\otimes_{\Z[T]_{\sol}} \widehat{\bigoplus}_{\N} \Z[[T]] = \varinjlim_{_{\substack{f,g\colon \N\to \N \\ f(n),g(m)\to \infty} \\ n,m\to \infty}} \prod_{(n,m)\in \N\times \N} T^{f(n)+ g(m)} \Z[[T]].
\]
On the other hand, the right term of \eqref{eqniehnlakwawr} is equivalent to 
\[
 \widehat{\bigoplus}_{\N} \Z[[T]]\xrightarrow{\sim} \widehat{\bigoplus}_{\N\times \N} \Z[[T]]=\varinjlim_{\substack{h\colon \N\times \N\to \N \\ h(n,m)\to \infty \\ 
 \sup(n,m)\to \infty}} \prod_{(n,m)\in \N\times \N} T^{h(n,m)} \Z[[T]]. 
\]
Thus, to prove the equivalence we have to show that given $h\colon \N\times \N\to \N$ with $h(n,m)\to \infty$ as $\sup(n,m)\to \infty$, there are $f,g\colon \N\to \N$ with $f(n),g(m)\to \infty$ as $n,m\to \infty$ such that $f(n)+g(m)\leq h(n,m)$ for all $(n,m)\in \N\times \N$, for this  one can take $f(n)=\lfloor \frac{1}{2}\inf_m h(n,m)\rfloor$ and $g(m)=\lfloor \frac{1}{2}\inf_n h(n,m)\rfloor$. 
\end{proof}

\begin{lemma}\label{LemmaBaseChangePolynomial}
Let $R$ be a finite type $\Z$-algebra and consider the morphism of rings $R\to R[T]$. Then the base change $R[T]_{\sol}\otimes_{R_{\sol}}-\colon \ob{D}(R_{\sol})\to \ob{D}(R[T])_{\sol}$ is naturally isomorphic to the functor $R\iHom_{R}(R((T^{-1}))/TR[T], -)$, in particular it commutes with limits. 
\end{lemma}
\begin{proof}
Let us denote $F:=R\iHom_{R}(R((T^{-1}))/TR[T], -)$ Let $R((T^{-1}))/TR[T]\to R$ be the projection onto the zeroth coefficient of the Laurent series expansion. For $M\in \ob{D}(R_{\sol})$ this produces a natural map $M\to F(M)$, since the module $F(M)$ has a natural structure of $R[T]$-module, this natural transformation naturally extends to $R[T]\otimes_{R} M\to F(M)$. We claim that $F(M)$ is $R[T]_{\sol}$-complete, this would yield a natural map $R[T]_{\sol}\otimes_{R_{\sol}}M\to F(M)$.  To prove the claim, note that $R((T^{-1}))/R[T]$ is a compact $R_{\sol}$-module, so that $F$ commutes with both limits and colimits. Thus, to see that it sends $R_{\sol}$-complete modules to $R[T]_{\sol}$-complete modules it suffices to prove it for $R$. We have that 
\begin{equation}\label{eqkhi93hlaweklhaetg}
R\iHom_{R}(R((T^{-1}))/TR[T],R) = R[T]
\end{equation}
as $R\iHom_{R}(R[T], R)= R((T^{-1}))/ TR[T]$ and compact $R_{\sol}$-modules are reflexive. Now, consider the natural transformation 
\[
R[T]_{\sol}\otimes_{R_{\sol}} M\to F(M),
\]
to see that it is an equivalence, it suffices to prove it in a family of compact projective generators, that is on $\prod_{\N} R$. But then both $R[T]_{\sol}\otimes_{R_{\sol}}-$ and $F(-)$ commutes with countable products of copies of $R$, whence it suffices to prove it for $M=R$ which reduces to \eqref{eqkhi93hlaweklhaetg}. 
\end{proof}

\begin{remark}\label{RemCompleteSolidTensor}
The solid tensor product preserves derived $I$-adically complete connective modules in greater generality, see \cite[Proposition 2.12.10]{MannSix}. In particular, \Cref{PropCompletionSolidTensorProduct} holds for animated rings almost of finite presentation over $\Z$ as in \Cref{RemarkTensorPseudoCoherent}, and also works for induced analytic ring structures. 
\end{remark}

\subsection{Schemes as analytic stacks}\label{ss:SchemesAsstacks}

With the introduction of the rings $R_{\square}$ for $R$ a finitely generated $\bb{Z}$-algebra, we are in shape to talk about two different realizations of schemes as analytic stacks (these will be introduced later in the notes).

\subsubsection{Classical approach}

We first need to see commutative rings as analytic rings. Let $\Ring$ be the $\infty$-category of (discrete) animated rings, given $R$ a discrete animated ring we write $\ob{D}(R)=\ob{Mod}_R(\ob{D}(\Z))$ for its category of modules in the derived category of discrete abelian groups.  

\begin{proposition}
 There is a fully faithful embedding 
\[
(-)^{\cond}: \Ring\to \AnRing
\]
mapping a ring $R$ to the analytic ring $R^{\cond}= (R, \ob{D}(R^{\cond}))$ to the trivial analytic ring structure on $R$, i.e. the analytic ring structure whose complete modules are all condensed $R$-modules. 
\end{proposition}
\begin{proof}
This follows from the natural fully faithful embedding of animated rings into animated condensed rings as discrete rings.  
\end{proof}

\begin{remark}\label{RemTRivialKunneth}
Let $R$ be a discrete animated ring. By definition we have that $\ob{D}(R^{\cond})=\Mod_{R}(\ob{D}(\Z^{\cond}))$. Thus, by \cite[Theorem 4.8.4.6]{HigherAlgebra} we have $\ob{D}(R^{\cond})=\ob{D}(R)\otimes_{\ob{D}(\Z)}\ob{D}(\Z^{\cond})$. By \Cref{RemDiscreteModules} $\ob{D}(R)$ is equivalent to the full subcategory $\ob{D}^{\delta}(\Z^{\cond})$ of discrete $R^{\cond}$-modules. Notice that as $R$ is also discrete as condensed set, all the objects of  $\ob{D}^{\delta}(\Z^{\cond})$ are discrete. 
\end{remark}

Let $R$ be an animated commutative ring and let $\Spec R$ be its Zariski spectrum defined as the Zariski  spectrum of $\pi_0(R)$. Let $\ob{D}(R)$ be the $\infty$-derived category of $R$-modules.  The  Zariski topology of $\Spec R$ has a basis of open affine schemes given by spectra of the form $\Spec R[f^{-1}]$ for $f\in R$ (were by definition $R[f^{-1}]=\varinjlim_{\times f} R=R\otimes_{\Z[T]} \Z[T^{\pm 1}]$ is given by the colimit of multiplication by $f$). We can rephrase the classical Zariski descent of $R$-modules in the language of  \Cref{theo:openCoversSym} and  \Cref{rem:ClosedVariant}. 

\begin{proposition}\label{Prop:Zariski Descent}
Let $\{U_i=\Spec R_i \}_{i=1}^n$ be a finite affine Zariski cover of $\Spec R$. Then the morphisms of symmetric monoidal categories $\{f_i^*:\ob{D}(R) \to \ob{D}(R_i)\}_{i=1}^n$ form a closed cover of $\ob{D}(R)$. In particular, we have descent  of quasi-coherent sheaves on $\Spec R$ for the Zariski topology. 
\end{proposition}
\begin{proof}
We prove the proposition in two steps. 

\textit{Step 1.} Let $f_1,\ldots, f_n\in R$ be elements generating the unit ideal and suppose that the cover is of the form $U_i=\Spec R[f_i^{-1}]$. It is clear that the $R$-algebras $R[f_i^{-1}]$ are idempotent so that they define closed subspaces in the locale $\ob{Sm}(\ob{D}(R))$. We want to see that they cover the locale, but this is equivalent to asking that there is an equivalence of complexes
\[
R \xrightarrow{\sim} [ \bigoplus_{i=1}^n R[\frac{1}{f_i}] \to \bigoplus_{i<j} R[\frac{1}{f_if_j}]\to \cdots \to R[\frac{1}{f_1\cdots f_n}]] 
\]
which amounts to Zariski descent for the underlying ring.

\textit{Step 2.} Now let $\{U_i\}_{i=1}^n$ be an arbitrary open cover of $R$. By Zariski descent we know that there is a natural equivalence of complexes
\[
R\xrightarrow{\sim} [ \bigoplus_{i=1}^n R_i \to \bigoplus_{i<j} R_i\otimes_R R_j \to \cdots \to \bigotimes_{i} R_i] 
\]
so the only thing to show is that each $R$-algebra $R_i$ is idempotent. For any open affine subspace $U=\Spec R'$ of $R$ there is a Zariski cover of $U$ of the form $R[f_i^{-1}]$ for suitable $f_i\in R$. It is clear that $R[f_{i}^{-1}]= R'[f_i^{-1}]$ so that we have an equivalence 
\[
R' \xrightarrow{\sim} [ \bigoplus_{i=1}^n R[\frac{1}{f_i}] \to \bigoplus_{i<j} R[\frac{1}{f_if_j}]\to \cdots \to R[\frac{1}{f_1\cdots f_n}]],
\]
but the right hand side is the idempotent algebra in $\ob{D}(R)$ corresponding to the union of the closed subspaces of the locale $\Sm(\ob{D}(R))$ associated to the algebras $R[f^{-1}_i]$.  This shows that $R'$ is idempotent which finishes the proof. 
\end{proof}

An immediate corollary of  \Cref{Prop:Zariski Descent} is the construction of quasi-coherent sheaves for schemes. 

\begin{corollary}\label{Coro:ZariskiDescent}
Let $X$ be a scheme with structural sheaf $\s{O}_X$ and let $|X|^{\op}$ be the topological space with underlying set $|X|$ and the coarsest topology given   by declaring closed subspaces the subsets of the form $|U|\subset X$ with $U$ an open Zariski subspace. With this topology $|X|^{\op}$ is a spectral space.

 Then the functor that maps an open affine Zariski subspace $U\subset X$ to $\ob{D}(\s{O}_X(U))$ is a sheaf. More precisely, let $\ob{D}(X)=\varprojlim_{U\subset X} \ob{D}(\s{O}_X(U))$ be the category of quasi-coherent sheaves on $X$, where $U$ runs over the poset of open affine subspaces of $X$. Then there is a unique natural surjective morphism of locales
\[
F:\Sm(\ob{D}(X))\to |X|^{\op}
\]
such that for $U\subset X$ open affinoid  we have $\ob{D}(F^{-1}(U))=\ob{D}(\s{O}_X(U))$. 
\end{corollary}
\begin{proof}
This is a consequence of  \Cref{Prop:Zariski Descent}, the only thing to verify is that the map is surjective. For the last claim, let $U\subset X$ be a Zariski open subset and let $A_{U}$ be the idempotent algebra associated to $U$. Then we have that 
\[
U=\{x\in X : k(x)\otimes_{\s{O}_X} A_U \neq 0 \}.
\]
This shows that if for two Zariski open sets one has $F^{-1}(U)=F^{-1}(U')$, i.e. $A_U=A_{U'}$, then $U=U'$ proving surjectivity. 
\end{proof}


\begin{remark}
In  \Cref{Coro:ZariskiDescent} it is not relevant that we have used the classical category of quasi-coherent sheaves. Indeed, this follows from the K\"unneth formula of \Cref{RemTRivialKunneth}. Moreover, Zariski descent generalizes in the following way: let $\Pr^L_{\ob{D}(R)}$ be the $(\infty,2)$-category of $\ob{D}(R)$-linear presentable categories. Given a morphism of rings $R\to S$ there is a natural base change functor $\Pr^L_{\ob{D}(R)}\to \Pr^L_{\ob{D}(S)}$ given by Lurie's tensor product $M\mapsto M\otimes_{\ob{D}(R)} \ob{D}(S)$. Then the functor mapping an open affine Zariski subspace $U=\Spec R'\subset \Spec R$ to $\Pr^L_{\ob{D}(R')}$ is a sheaf for the Zariski topology. This follows, for example, from \cite[Proposition 3.45]{MathewDescent} as any Zariski cover gives rise to a \textit{descendable morphism} of algebras. Descendable maps of rings will be prototypical examples of \textit{$!$-covers} of affinoid analytic stacks, see \Cref{ss:AnStacksConstruction}.
\end{remark}

An anti-intuitive phenomena is happening in the classical Zariski descent of quasi-coherent sheaves, namely, open Zariski subspaces of  $\Spec R$ are giving rise to closed subspaces of the locale $\ob{Sm}(\ob{D}(R))$! This explains why the classical theory of quasi-coherent sheaves on schemes do not have (apparently!) a well defined theory of cohomology with compact support outside the proper case. A way to overcome this discrepancy is to use a different realization of schemes into analytic stacks, or equivalently, a different embedding of the category of commutative rings into the category of analytic rings. A the end, it will be more convenient to work with a generalization of both theories of quasi-coherent sheaves; this is captured in the theory of discrete adic spaces (see \cite[Lectures IX and X]{ClausenScholzeCondensed2019}). But before that, let us finish the topological study of the locale defined by classical algebraic geometry.

\subsubsection{Formal completions}
As we previously observed, an open Zariski map  $U=\Spec R' \to \Spec R$  gives rise to a closed subspace of the locale $\Sm(\ob{D}(R))$. A natural question is to give a description of the open complements, they are given by formal schemes. 

\begin{proposition}
Let $X$ be a scheme and let $U\subset X$ be a qcqs open Zariski subspace with closed complement $Z\subset X$. Let $X^{\wedge_Z}$ be the (derived) formal completion of $X$ along $Z$.     Let $\ob{D}(X)$ be the $\infty$-category of derived quasi-coherent sheaves on $X$ and let $F: \Sm(\ob{D}(X))\to |X|^{\op}$ be the map of locales provided by  \Cref{Coro:ZariskiDescent}.  Then the open subspace $ F^{-1}(Z)$ has  as underlying category the  category of (derived)  complete quasi-coherent sheaves on $X^{\wedge Z}$. 
\end{proposition}
\begin{proof}
We will explain the simplest case when $X=\Spec R$ is affine and $U= \Spec R[f^{-1}]$ is given by inverting $f\in R$, the general case follows by Zariski descent and induction in the number of generators of the ideal $I$ defining $Z$.

By definition, $\ob{D}(F^{-1}(Z))$ is the Verdier quotient $j^*: \ob{D}(R)\to \ob{D}(R)/\ob{D}(R[f^{-1}])$. Moreover, the functor $j_*: \ob{D}(F^{-1}(Z))\to \ob{D}(R)$ is fully faithful and has by essential image the elements $M\in \ob{D}(R)$ such that  the natural map 
\[
M\to R\Hom_R(\fib(R\to R[f^{-1}]), M)
\]
is an equivalence (equivalently those $M$ such that $R\Hom_{R}(R[f^{-1}],M)=0$). But we can write $R[f^{-1}]= \varinjlim_{\times f} R$ as the colimit of multiplication by $f$ on $R$, and so we find that 
\[
R\Hom_R(\fib(R\to R[f^{-1}]), M) = R\varprojlim_n (M\otimes_{R}^L R/^{\bb{L}} f^n)
\]
where $ R/^{\bb{L}} f^n$ is the derived quotient of $R$ by $f^n$, represented by the Koszul complex $[R\xrightarrow{f^n} R]$. Thus, by definition, the essential image of $f$ consists on those modules which are derived $f$-adically complete.  

Now let $I=(f) \subset R$ be the ideal generated by $f$. To finish the proof of the proposition it suffices to see that for any element $g\in \mathrm{Rad}(I)$ in the radical of $I$,  a derived $f$-adically complete module is also derived $g$-adically complete.  Indeed, this will show that the category $\ob{D}(F^{-1}(Z))$ only depends on the formal completion of $\Spec R$ along $Z$ and that by definition it consists in the  complete modules of $X^{\wedge_Z}$. To see this, we note that $\Spec R[g^{-1}] \subset \Spec R[f^{-1}]$, i.e. that we have a map of idempotent $R$-algebras $R[f^{-1}]\to R[g^{-1}]$. But then 
\[
\begin{aligned}
R\Hom_R(R[g^{-1}],M) & = R\Hom_R(R[g^{-1}]\otimes_{R}^L  R[f^{-1}], M ) \\ 
& = R\Hom_R(R[g^{-1}], R\Hom_R(R[f^{-1}],M)) \\ 
& = 0
\end{aligned}
\]
proving what we wanted. 
\end{proof}

To conclude,  in classical algebraic geometry we are placed in the strange situation where the  qcqs open Zariski maps give rise to closed maps at the level of derived categories, while formal completions along finitely generated Zariski closed subschemes give rise to open immersions. To solve this ``paradox'' (which is actually not a bug but a feature of the theory!) one needs to consider  solid  quasi-coherent sheaves instead.

\subsubsection{Solid approach}

Let us now discuss the solid approach to the theory of quasi-coherent sheaves for schemes. We need a definition:

\begin{definition}\label{def:SolidRingGeneralR}
Let $R$ be a commutative ring, we define the analytic ring  $R_{\square}$  to be the colimit 
\[
R_{\square}=\varinjlim_{B\subset R} B_{\square}
\]
where $B$ runs over all the finitely generated $\bb{Z}$-subalgebras, and $B_{\square}$ is the analytic ring constructed in   \Cref{theoSolidR}. Equivalently, $R_{\square}$ is the analytic ring structure on $R$ where a condensed $R$-module is $R_{\square}$-complete if and only if it is $r$-solid for all $r\in R$.   

In general, for $R$ an animated commutative ring we let $R_{\square}$ be the analytic ring structure on $R$ induced from $\pi_0(R)_{\square}$ via   \Cref{PropInvarianceHomotopy}. Equivalently, an $R$-module $M$ is $R_{\sol}$-complete if and only if all its cohomology groups $H^*(M)$ are $\pi_0(R)_{\sol}$-complete. 
\end{definition}

\begin{remark}
By  \Cref{theoSolidR}, for $B$ a finitely generated $\bb{Z}$-algebra $B_{\sol}$ is an analytic ring structure of $B$, i.e. $B$ is $B_{\sol}$-complete. It follows formally that any discrete  (derived) $B$-module is also $B_{\sol}$-complete. This shows that, for any animated ring $R$,  the ring itself is $R_{\sol}$-complete and so $R_{\sol}$ is an analytic ring structure on $R$. 
 \end{remark}
 
\begin{remark}\label{RemarkGeneratorsSolid}
Let $R$ be an animated ring and let $S\subset \pi_0(R)$  be a set of generators as $\bb{Z}$-algebra.  By  \Cref{CorollaryFreeProfiniteR} the ring $R_{\sol}$ is also obtained as the analytic ring structure on $R$ making a $\Z$-solid $R$-module complete if and only if it is $s$-solid for all $s\in S$.
\end{remark}

\begin{remark}
The reader might ask why not to define $B_{\sol}$ to be the analytic ring whose values at a profinite set $S=\varprojlim_n S_n$ are given by $B_{\sol}[S]=\varprojlim_n B[S_n]$; this would define a category that we might want to call \textit{ultra solid modules}. The answer is that it does not give rise to an analytic ring structure on $B$  in general\footnote{The problem being that in the correct definition of the ultra solid modules attached to $B$, the forgetful functor $\ob{D}(B_{\sol})\to \ob{D}(\underline{B})$ from ultra-solid to condensed $B$-modules is not fully faithful in general.}. So far we have seen that this is a sensible definition when $B$ is a finitely generated $\bb{Z}$-algebra (or an animated ring of almost finite presentation over $\Z$). It turns out that this also works well when $B$ is \textit{essenitally of finite type}, i.e. when it is a Zariski localization (i.e. colimit of open Zariski localizations) of a finitely generated $\bb{Z}$-algebra. For example, when $B=\bb{Q}$ one obtains the \textit{ultra solid} rational numbers and the category $\ob{D}(\bb{Q}_{\sol}) \subset \ob{D}(\bb{Z}_{\sol})$ is the open localization which is complement to the idempotent $\bb{Z}_{\sol}$-algebra $\widehat{\bb{Z}}= \prod_p \bb{Z}_p$. 

In general, by going beyond the notion of analytic ring, one can define categories of ultra solid modules over arbitrary animated rings and $\bb{E}_{\infty}$-algebras, see \cite{aparicio2024ultrasolidhomotopicalalgebra} for the definition in the case of fields, and \cite{brav2025beilinsonparshinadelessolidalgebraic} for ultra solid rings constructed for rings essentially of finite type over $\Z$.
 \end{remark}

Let $R\to S$ be a map of animated  rings, by definition of the analytic ring structure  of solid rings there is a map of analytic rings $R_{\sol}\to S_{\sol}$. Indeed, a complete $S_{\sol}$-module is an $S$-module $M$ which is solid for all $s\in S$,  this means concretely that for all map $\bb{Z}[T]\to R$ the module $M$ is $\bb{Z}[T]_{\sol}$-solid. Hence, it is easy to see that its restriction to an $R$-module is also $R_{\sol}$-complete.  We get a functor 
\[
(-)_{\sol}:\Ring\to \AnRing
\]
from animated commutative rings to analytic rings.

 One of the first properties we have to check is the compatibility with colimits of the functor $(-)_{\sol}$.

\begin{proposition}\label{PropositionColimitPreservingSolidRing}
The functor $(-)_{\sol}:\Ring\to \AnRing$  is fully faithful and commutes with colimits. 
\end{proposition}
\begin{proof}

Let $R,S$ be animated rings. We want to see that the natural fully-faithful map of anima
\[
(-)^{\triangleright}: \Map_{\AnRing}(R_{\sol},S_{\sol}) \to \Map_{\Ring}(R,S)
\]
is essentially surjective. For this, recall that we have a map $\AnRing\to \CondRing$ mapping an analytic ring to its underlying condensed ring. The claim follows from the fact that the composite 
\[
\Ring\xrightarrow{(-)_{\sol}} \AnRing \xrightarrow{(-)^{\triangleright}} \CondRing
\]
is fully faithful.

Next, we show the second statement. Let $\{R_i\}_{i\in I}$ be a diagram of animated commutative rings and let $\{R_{i,\sol}\}_{i\in I}$ be its associated diagram of solid rings. By  \Cref{PropColimitsUncompletedAnalytic} the colimit $\varinjlim_{i} R_{i,\sol}$ is the analytic ring given by the completion of the analytic ring structure on $S:=\varinjlim_{i} R_i$ where an $S$-module $M$ is complete if and only if it is $R_{i,\sol}$-complete for all $i$.  But this is equivalent for $M$ to be $r$-solid for all  $r\in R_i$ and all $i\in I$, which is also equivalent to be $s$-solid for all $s\in S$. This shows that $\varinjlim_i R_{i,\sol}= S_{\sol}$ as wanted.   
\end{proof}

Our next step to construct solid analytic stacks (to be defined properly later in \Cref{s:AnStacks}) attached to schemes. This will turn out to be equivalent to  proving Zariski descent for the functor  $\Spec R\mapsto \ob{D}(R_{\sol})$ at the level of affine schemes.  We first prove this for coverings given by basic open Zariski subspaces. We need some technical lemmas. 

\begin{lemma}\label{LemmaOpenBaseChange}
Let $f:A\to B$ be a map of analytic rings such that the induced map of derived categories $f^*:\ob{D}(A)\to \ob{D}(B)$ is an open immersion of locales. Then for all analytic ring $C$ and any map $A\to C$ the base change $C\to C\otimes_A B$ also induces an open immersion at the level of locales. 
\end{lemma}
\begin{proof}
This follows formally from the K\"unneth formula of \Cref{PropColimitsAnRingAndMod}, we also give a self contained proof. By hypothesis there is $D\in \ob{D}(A)$ an idempotent algebra such that $\ob{D}(B)\subset \ob{D}(A)$ is the full subcategory consisting on those $M$ such that $R\iHom_A(D,M)=0$. This defines an uncompleted analytic ring structure on $A^{\triangleright}$. Thus, by  \Cref{TheoCompletedAnRing}, $C\otimes_A B$ is given by the completion of the analytic ring structure on $C^{\triangleright}$ such that a $C^{\triangleright}$-module $M$ is complete if and only if it is $C$ and $B$-complete. This is equivalent to asking $M$ to be $C$-complete and that $R\iHom_{C}(C\otimes_A D, M)=0$. Thus, $C\otimes_A B$ is the open localization complement to the idempotent $C$-algebra $C\otimes_A D$ proving what we wanted. 
\end{proof}

\begin{lemma}\label{LemmaBasicOpenImmersions}
The following maps of analytic rings give rise to open immersions at the level of locales:

\begin{enumerate}

\item $\bb{Z}[T]_{\bb{Z}_{\sol}/} \to \bb{Z}[T]_{\sol}$ with complement idempotent algebra $\bb{Z}((T^{-1}))$.

\item $ \bb{Z}[T]_{\bb{Z}_{\sol}/} \to \bb{Z}[T^{\pm 1}]_{\bb{Z}[T^{-1}]_{\sol}/}$ with complement idempotent algebra $\bb{Z}[[T]]$.

\end{enumerate}
\end{lemma}
\begin{proof}
Part (1) is  \Cref{TheoSolidZT}. Part (2) follows by a similar argument as in \textit{loc. cit.}, we leave the details to the reader. 
\end{proof}

\begin{lemma}\label{LemmaOpenBasicSolid}
Let $R$ be an animated commutative ring and let $f\in R$. Then the map $R_{\sol}\to R[f^{-1}]_{\sol}$ of analytic rings induces an open immersion at the level of locales. 
\end{lemma}
\begin{proof}
By  \Cref{PropositionColimitPreservingSolidRing} we have that 
\[
R[f^{-1}]_{\sol} = R_{\sol}\otimes_{\bb{Z}[T]_{\sol}} \bb{Z}[T^{\pm 1}]_{\sol}
\]
where $T\mapsto f$. Therefore, by  \Cref{LemmaOpenBaseChange}, it suffices to show that the map $\bb{Z}[T]_{\sol}\to \bb{Z}[T^{\pm 1}]_{\sol}$ induces an open immersion of locales, but this follows from  \Cref{LemmaBasicOpenImmersions} (2) by taking the base change of $\Z[T]_{\Z_{\sol}/}\to \Z[T^{\pm 1}]_{Z[T^{-1}]_{\sol}/}$ along  $\bb{Z}[T]_{\bb{Z}_{\sol}/}\to \bb{Z}[T]_{\sol}$.  
\end{proof}

\begin{proposition}\label{PropositionBasicCoversSolid}
Let $R$ be an animated commutative ring and let  $\{U_i=\Spec R[f_i^{-1}]\}_{i=1}^n$ be an open cover of $\Spec R$ by basic affine subspaces. Then the  maps $R_{\sol}\to R[f_i^{-1}]_{\sol}$ induce an open cover in the locale $\Sm(\ob{D}(R_{\sol}))$.
\end{proposition}
\begin{proof}
By  \Cref{LemmaOpenBasicSolid} we know that the pullback maps $\ob{D}(R_{\sol})\to \ob{D}(R[f_i^{-1}]_{\sol})$ are open localizations.  By noetherian approximation we can assume without loss of generality that $R$ is an animated algebra of finite presentation (or for simplicity just a finitely generated algebra over $\bb{Z}$). In such a case, the free solid $R$-module for $S=\varprojlim_n S_n$ profinite is 
\[
R_{\sol}[S]=\varprojlim_n R[S_n]. 
\] 

By construction, the closed complement $Z_i$ of $R_{\sol}\to R[f_i^{-1}]_{\sol}$ is given by the idempotent $R_{\sol}$-algebra
\[
R[[f_i]]:= R_{\sol} \otimes_{\bb{Z}[T]_{\bb{Z}_{\sol}/}} \bb{Z}[[T]]
\]
also given by the derived quotient $R[[X]]/^{\bb{L}}(X-f_i)$. Thus, in order to show the proposition it suffices to prove that $\cap_{i=1}^n Z_i=\emptyset$ in $\Sm(\ob{D}(R_{\sol}))$, namely, that $\bigotimes_{i=1}^n R[[f_i]]=0$. But we have that 
\[
\bigotimes_{i=1}^n R[[f_i]] = R[[X_1,\ldots, X_n]]/^{\bb{L}}(X-f_1,\ldots, X-f_n).
\]
Now let  us write $1=\sum_i a_i f_i$ and consider the map $g:\bb{Z}[T]\to R[[X_1,\ldots, X_n]]$ mapping $T$ to $\sum_i a_i X_i$. Then $g$ extends (uniquely by idempotency) to a map $\bb{Z}[[T]]\to R[[X_1,\ldots, X_n]]$. This shows that the unit $1\in \bigotimes_{i=1}^n R[[f_i]] $ gives rise to a map from $\bb{Z}[[T]] $, but $\bb{Z}[[T]]/(T-1)=0$ proving that $1=0$ in the tensor and so that $\bigotimes_{i=1}^n R[[f_i]] =0$ as wanted. 
\end{proof}

\begin{theorem}\label{theo:ZariskiDescentSolid}
Let $R$ be an animated ring and let $\{U_i=\Spec R_i\}_{i=1}^n$ be an open affine Zariski cover of $\Spec R$. Then the maps $R_{\sol} \to R_{i,\sol}$ induce an open cover of the locale $\Sm(\ob{D}(R_{\sol}))$. 
\end{theorem}
\begin{proof}
First we show that the maps $R_{\sol}\to R_{i,\sol}$ induce open immersion at the level of the locales. For this, let us take open basic affines $\Spec R[f_{i,j}]^{-1}$ covering $U_i$. We have that $R[f_{i,j}^{-1}]=R_i[f^{-1}_{i,j}]$ and by   \Cref{PropositionBasicCoversSolid} we know that $\{ R_{i,\sol}\to R_i[f^{-1}_{i,j}]_{\sol}\}_j$ gives rise an open cover of locales. This implies that $\ob{D}(R_{\sol})\to \ob{D}(R_{i,\sol})$ is an open immersion (being the union of the open immersions $\ob{D}(R_{\sol})\to \ob{D}(R_{i}[f_{i,j}^{-1}]_{\sol})$). Finally the fact that the open immersions $\{\ob{D}(R_{\sol}\to \ob{D}(R_{i,\sol}))\}_i$ cover $\ob{Sm}(\ob{D}(R_{\sol}))$ follows from  \Cref{PropositionBasicCoversSolid} and the fact that we can find a refinement of $\{U_i\}_{i=1}^n$ by basic Zariski open subspaces.  
\end{proof}

\begin{corollary}\label{coroSolidSheavesSchemes}
Let $X$ be a scheme with structural sheaf $\s{O}_X$. Then the functor that maps an open affine Zariski subspace $U\subset X$ to $\ob{D}(\s{O}_X(U)_{\sol})$ is a sheaf. More precisely, let $\ob{D}_{\sol}(X)=\varprojlim_{U\subset X} \ob{D}(\s{O}_X(U)_{\sol})$ be the category of solid quasi-coherent sheaves on $X$, where $U$ runs over the poset of open affine subspaces of $X$. Furthermore,  there is a unique \  surjective morphism of locales
\[
F:\Sm(\ob{D}_{\sol}(X))\to |X|
\]
such that for $U\subset X$ open affinoid  we have $\ob{D}_{\sol}(F^{-1}(U))=\ob{D}(\s{O}_X(U)_{\sol})$. 
\end{corollary}
\begin{proof}
This is a consequence of  \Cref{theo:ZariskiDescentSolid}, the only claim left to verify is the surjectivity of $F$. For this, we can assume that $X=\Spec R$ is affine.  Let $U,U'\subset X$ be two open Zariski subsets such that $F^{-1}(U)=F^{-1}(U')$, we want to show that $U'=U$. By taking intersections with open affines, we can assume without loss of generality that $U'\subset U=X$. Then, to prove the claim it suffices to show that if for $Z\subset X$ a closed subspace the associated idempotent algebra $A_Z$ is zero, then $Z=\emptyset$. Suppose this does not hold, then $A_{Z}=0$ and $Z\neq \emptyset$. There is a closed point $x\in Z$ and a map of analytic rings $R_{\sol}\to \kappa(x)_{\sol}$ that gives rise a commutative map of locales
\[
\begin{tikzcd}
\Sm(\ob{D}(\kappa(x)_{\sol})) \ar[d,"G"] \ar[r,"F_x"] & \mathrm{Spec}( \kappa(x)  )\ar[d] \\ 
\Sm(\ob{D}(R_{\sol})) \ar[r,"F"] & \mathrm{Spec} (R) .
\end{tikzcd}
\]
But then we have that $G^{-1}(F^{-1}(Z))=F_x^{-1}(\Spec (\kappa(x)))$ which shows that $F^{-1}(Z)\neq \emptyset$, a contradiction with the fact that $A_{Z}=0$. 
\end{proof}

\subsection{Discrete Huber pairs}\label{subsec:DiscretePairs}

In  \Cref{ss:SchemesAsstacks} we saw two possible ways to attach categories of quasi-coherent sheaves to a scheme $X$ (later in the notes we will see that this corresponds to realize the scheme $X$ in two different ways as an analytic stack). For the continuation of the theory it will be more convenient to generalize both constructions in the theory of discrete Huber rings and discrete adic spaces. 

\begin{definition}\label{def:DiscreteAdic}
An \textit{animated discrete Huber pair} is a tuple $(R,R^+)$ where $R$ is an  animated discrete ring and  $R^+\subset \pi_0(R)$ is an integrally closed subring. A morphism of discrete Huber pairs $(R,R^+)\to (A,A^+)$ is a map of animated rings $R\to A$ such that $R^+$ is mapped to $A^+$.  We let $\Ring^{\Hub}$ be the $\infty$-category of animated discrete Huber pairs. 
\end{definition}

\begin{lemma}\label{Lem:PropertiesAffDis}
The category $\Ring^{\Hub}$ admits colimits and has compact   generators given by finite colimits of  the discrete pairs $(\bb{Z}[T],\bb{Z})$ and $(\bb{Z}[T], \bb{Z}[T])$.
\end{lemma}
\begin{proof}
Let $\{(R_i,R^+_i)_{i\in I}$ be a diagram in $\Ring^{\Hub}$, then its colimit is given by the Huber pair $(A,A^+)$ where $A=\varinjlim_i R_i$ and $A^+\subset \pi_0(A)$ is the integral closure of the ring generated by the images of $R^+_i\to \pi_0(A)$.   In particular, the functors $(R,R^+)\mapsto R,R^+$ commutes with filtered colimits in $\Ring^{\Hub}$.

 To show that $(\bb{Z}[T], \bb{Z})$ and $(\bb{Z}[T],\bb{Z}[T])$ are compact  generators of  $\Ring^{\Hub}$, it suffices to show that the functor 
\[
F:\Ring^{\Hub}\to \mathrm{Ani}\times \mathrm{Ani}
\]
mapping a ring $(R,R^+)$ to the mapping spaces $\Map_{\Ring^{\Hub}}((\bb{Z}[T], \bb{Z}), (R,R^+))$ and \newline $ \Map_{\Ring^{\Hub}}((\bb{Z}[T], \bb{Z}[T]), (R,R^+))$ is  jointly conservative and commutes with filtered colimits. But we have that 
\[
\Map_{\Ring^{\Hub}}((\bb{Z}[T], \bb{Z}), (R,R^+)) = R
\]
and 
\[
\Map_{\Ring^{\Hub}}((\bb{Z}[T], \bb{Z}[T]), (R,R^+))= R\times_{\pi_0(R)} R^+,
\]
proving the claim. 
\end{proof}

\begin{remark}\label{RemNonProjective}
The Huber pair $(\bb{Z}[T], \bb{Z}[T])$ is \textbf{not} projective. Indeed, consider the quotient $A=\bb{Z}[T]/T^n$. then $(A,A)$ is the  pushout $(\bb{Z}[T], \bb{Z})\otimes_{(\bb{Z}[T], \bb{Z})} (\bb{Z}, \bb{Z})$ mapping $T\mapsto T^n$ and $T\mapsto 0$, writing the pushout as a bar construction one finds that $(R,R^+)\mapsto R^+$ does not commute with geometric realizations (the problem being that colimits on Huber pairs always take the integral closure in the second argument).
\end{remark}

Given a  discrete Huber pairs we can naturally attach a solid analytic ring; it will be convenient to make a more general construction. 

\begin{definition}\label{def:GeneralizedHuberPair}
Let $(A,S)$ be a pair consisting on a solid animated ring and a map of sets $S\to \pi_0(A)(*)$ towards the underlying discrete static ring of $A$. We define the analytic ring $(A,S)_{\sol}$ to be completion of the analytic ring structure on $A$ making a condensed $A$-module complete if and only if it is a solid abelian group and for all $s\in S$ it is $s$-solid as in  \Cref{SolidZT}. 
\end{definition}

\begin{lemma}\label{lem:IntegralClosureSolid}
Let $(R,S)$ be a tuple with $R$ a discrete animated ring and $S\subset \pi_0(R)$ a subset. Let $R^+[S]\subset \pi_0(R)$ be the integral closure of the subalgebra generated by the image of $S$. Then the natural map
\[
(R,S)_{\sol}\to (R,R^+[S])_{\sol}
\]
is an equivalence. 
\end{lemma}
\begin{proof}
 By  \Cref{PropInvarianceHomotopy} we can assume without loss of generality that $R$ is static.   Let $(R,S)$ be a pair with $R$ a discrete ring and $S\subset R$ a set. Let $R^+[S]$ be the integral closure of the subalgebra of $R$ generated by $S$.  We claim that $(R,S)_{\sol}=(R,R^+(S))_{\sol}$. Indeed, let $r\in R^+[S]$, then $r$ is integral over a subalgebra $B$ generated by finitely many elements in $S$. Then, by  \Cref{cor:SolidIntegralInvariance}, we have maps of analytic rings
 \[
 B[r]_{\sol}=B[r]_{B_{\sol}/} \to (R,S)_{\sol},
 \]
 which in turn produces a map of analytic rings 
 \[
 (R,B[r])_{\sol}\to (R,S)_{\sol}.
 \]
 Taking colimits along all $r\in R^+[S]$ we get  maps $(R,S)_{\sol}\to (R,R^+[S])_{\sol}\to (R,S)_{\sol}$ proving that $(R,R^+[S])_{\sol}=(R,S)_{\sol}$ as wanted. 
\end{proof}

\begin{prop}[{\cite[Proposition 3.34]{Andreychev}}]\label{Prop:FUllyFaithDiscreteAnalytic}
The functor $\Ring^{\Hub}\to \AnRing$ mapping $(R,R^+)$ to $(R,R^+)_{\sol}$ is colimit preserving and fully faithful. 
\end{prop}
\begin{proof}

First, note that for any pair $(R,S)$ with $R$ a discrete animated ring, the underlying ring $(R,S)_{\sol}^{\triangleright}=R$ is complete  as any discrete module is solid (resp. $s$-solid for $s\in S$).  The commutativity with colimits follows from the definition of $(R,S)_{\sol}$ by declaring an $R$-module complete if it is solid as abelian group and $s$-solid for all $s\in S$,  the description of colimits of (uncompleted) analytic rings of  \Cref{PropColimitsUncompletedAnalytic}, and \Cref{lem:IntegralClosureSolid}.

 We now prove  fully faithfulness. Since the functor is colimit preserving and $(\bb{Z}[T],\bb{Z})$ and $(\bb{Z}[T], \bb{Z}[T])$ generate $\Ring^{\Hub}$, it suffices to show that  for $(R,R^+)$ a discrete Huber pair the maps
 \begin{equation}\label{eqMapSpace1}
 R\to \Map_{\AnRing}((\bb{Z}[T],\bb{Z})_{\sol} ,(R,R^+)_{\sol})
 \end{equation}
 and 
 \begin{equation}\label{eqMapSpace2}
 R\times_{\pi_0(R)} R^+ \to \Map_{\AnRing}((\bb{Z}[T],\bb{Z}[T])_{\sol} ,(R,R^+)_{\sol})
 \end{equation}
 are equivalences.  By definition, $(\bb{Z}[T],\bb{Z})_{\sol}=\bb{Z}[T]_{\bb{Z}_{\sol}/}$ has the induced solid structure from the integers. Therefore, the mapping space \eqref{eqMapSpace1} is the underlying discrete ring of $(R,R^+)_{\sol}$ which is nothing but $R$ as expected.  For the second claim, note that the map $(\bb{Z}[T],\bb{Z})_{\sol}\to (\bb{Z}[T],\bb{Z}[T])_{\sol}$ is idempotent  and so the mapping space of \eqref{eqMapSpace2} is completely determined by its connected components. But by definition $\pi_0(\Map_{\AnRing}((\bb{Z}[T],\bb{Z}[T])_{\sol},(R,R^+)_{\sol}))$ consists on all the maps $\bb{Z}[T]_{\sol}\to (R,R^+)_{\sol}$.  By  \Cref{PropInvarianceHomotopy} we can assume without loss of generality that $R$ is static, we then have by definition that $(R,R^+)_{\sol}=R_{R^+_{\sol}/}$. Then, for $S=\varprojlim_n S_n$ a light profinite set, we have 
 \[
 (R,R^+)_{\sol}[S]=R\otimes_{R^+} R^+_{\sol}[S].
 \]
 Thus, we can write
 \[
 (R,R^+)_{\sol}[S] =\varinjlim_{\substack{B\subset R^+ \\ M\subset R}} M_{\sol}[S]
 \]
 where $B$ runs over all the finitely generated subalgebras of $R^+$, $M$ runs along all the finite $B$-submodules in $R^+$, and $M_{\sol}[S]=\varinjlim_n M[S_n]$. Now let $T\in R$ be such that we have a map $\bb{Z}[T]_{\sol}\to (R,R^+)_{\sol}$ and let $C\subset R^+$ be the subalgebra generated by the image of $T$. Then $(R,R^+)_{\sol}[S]$ is $C_{\sol}$-complete and we have a map of $C_{\sol}$-modules
 \[
 C_{\sol}[S]\to (R,R^+)_{\sol}[S]
 \]
 Taking $S=\bb{N}\cup \{\infty\}$ and identifying $\bb{Z}_{\sol}[S]\cong \prod_{\bb{N}} \bb{Z}$ we find a map
 \[
 \prod_{\bb{N}} C \to \varinjlim_{\substack{B\subset R^+ \\ M\subset R}} \prod_{\bb{N}} M. 
 \]
In particular, there is a finitely generated ring $B\subset R^+$ and a finite $B$-module $M\subset R$ such that the sequence $(1,T,T^2,\ldots)$ lands in $\prod_{\bb{N}} M$. Thus, all the powers of $T$ are in $M$ which implies that $T$ is integral over $B$ and so that $T\in R^+$. This finishes the proof. 
\end{proof}

\subsection{Discrete Adic spaces}\label{subsec:AdicSpectrum}

Huber's approximation to non-archimedean geometry relies in the formalism of adic spaces \cite{HuberValuations,HuberAdicSpaces,HuberEtaleCohomology}. The traditional approximation to the theory  involves a discussion of topological rings and continuous valuations, then a technical sheafy property appears in order to guarantee descent of the structural sheaf for the analytic topology. In this section we shall follow the approach of \cite[Lexture IX]{ClausenScholzeCondensed2019} and restrict ourselves to the theory of discrete adic spaces where these topological and sheafiness subtleties disappear. 

\subsubsection{Adic spectrum of Huber pairs}

\begin{definition}
Let $(R,R^+)$ be a discrete Huber pair with $R$ a static ring. We let $\Spa(R,R^+)$ be the set of equivalence classes of multiplicative valuations $x: R\to \Gamma$ such that $|f(x)| \leq 1$ for all $f\in R^+$. For a general discrete Huber pair $(R,R^+)$ we define $\Spa(R,R^+):= \Spa(\pi_0(R),R^+)$.  We call $\Spa(R,R^+)$  the \textit{adic spectrum} of $(R,R^+)$.  When $R^+$ is the integral closure of the image of $\bb{Z}$ we simply write $\mathrm{Spv}(R)$; this is the set of equivalent classes of valuations of $R$ and is called the \textit{valuation spectrum} of $R$. We have $\Spa(R,R^+)\subset \mathrm{Spv}(R)$.
\end{definition}

\begin{remark}\label{RemarkEquivalentDescriptionSpa}
An equivalent way to define $\Spa(R,R^+)$ is as follows: it is the set of tuples $(\f{p},V)$ where $\f{p}$ is a prime ideal of $R$ and $V\subset \kappa(\f{p})$ is a valuation ring in the fraction field at $\f{p}$ containing the image of $R^+\to \kappa(\f{p})$. 
\end{remark}

Huber has defined a topology on $\Spa(R,R^+)$ called the \textit{Huber's or analytic topology}. It is defined as follows: let  $(R,R^+)$ be a discrete Huber pair and  let $(f_1,\ldots, f_n,g)$ be a tuple of elements in $R$. Define the rational localization $\Spa(R,R^+)\left(\frac{f_1,\ldots, f_n}{g}\right)$ to be the subset of $\Spa(R,R^+)$ consisting on those equivalence classes of valuations $x:R\to \Gamma$ such that $|g|\neq 0$ and $|f_i(x)|\leq |g(x)|$ for all $i=1,\ldots, n$. The analytic topology of $\Spa(R,R^+)$ is the topology generated by declaring rational subsets to be open subspaces.  We left as an exercise to the reader to prove that for a map $(R,R^+)\to (S,S^+)$ of discrete Huber pairs, the induced map $\Spa(S,S^+)\to \Spa(R,R^+)$ preserves rational localizations (and so it is continuous for Huber's topology), and that intersections of rational localizations are rational localizations.

The following proposition identifies integrally closed subrings of $R$ with suitable subsets of $\mathrm{Spv}(R)$. 

\begin{proposition}[{\cite[Proposition 9.2]{ClausenScholzeCondensed2019}}]\label{PropIntegrallyClosedRings}
Let $R$ be a static discrete ring. There is a bijection between integrally closed subrings $R^+\subset R$ and subsets $U\subset \mathrm{Spv}(R)$ which are intersections of rational localizations $U_{f,1}:=\mathrm{Spv}(R)(\frac{f}{1})$. Explicitly, one has 
\[
R^+ =\{f\in R : \forall x\in U, |f(x)|\leq 1\}
\]
and 
\[
U= \bigcap_{f\in R^+} U_{f,1}=\Spa(R,R^+). 
\]
In particular, 
\[
R^+=\{f\in R: \forall x\in \Spa(R,R^+), |f(x)|\leq 1\}. 
\]
\end{proposition}
\begin{proof}
Let $U=\bigcap_{f\in I} U_{f,1}$ be an intersection as in the lemma, let $R^+\subset R$ be the integral closure in $R$ of the algebra generated by the elements $f\in I$. Then it is easy to see that $U=\Spa(R,R^+)$. Conversely, we want to see that the integrally closed subring $R^+\subset R$ is totally determined by $\Spa(R,R^+)$. For that, let $f\in R$ be such that for all $x\in \Spa(R,R^+)$ one has $|f(x)|\leq 1$. If $f\notin R^+$ then $f$ is not in $R^+[\frac{1}{f}]\subset R[\frac{1}{f}]$ and there is a prime ideal $\f{p}$ of $R^+[\frac{1}{f}]$ that contains $\frac{1}{f}$. Let $\f{q}$ be a minimal prime contained in $\f{p}$. We may then find a valuation ring $V$ with a map $\Spec V\to \Spec R^+[\frac{1}{f}]$ taking the generic point to $\f{q}$ and the special point to $\f{p}$. As the image of $\Spec R[\frac{1}{f}]\to \Spec R^+[\frac{1}{f}]$ contains the minimal prime $\f{q}$ (the map $R^+[\frac{1}{f}]\to R[\frac{1}{f}]$ being injective), we can lift the valuation corresponding to $\Spec V\to \Spec R^+[\frac{1}{f }]$ to $R[\frac{1}{f}]$.  The resulting valuation takes values $\leq 1$ on $R^+$, and value $>1$ on $f$ since $f\in \f{p}$. This gives the contradiction.  
\end{proof}

Next, we want to show that there is a well defined analytic topology on $\Spa(R,R^+)$ and  construct a structural sheaf.

\begin{proposition}\label{PropRatLocDiscrete}
Let $(R,R^+)$ be a discrete Huber ring and let $(f_1,\ldots, f_n,g)$ be a tuple generating the unit ideal of $R$. Let $R(\frac{f_1,\ldots, f_n}{g})^+\subset \pi_0(R[\frac{1}{g}])$ be the integral closure of the subring generated by the image of $R^+$ and $\frac{f_i}{g}$ for all $i=1,\ldots, g$.   Then the natural map 
\[
\Psi:\Spa(R[\frac{1}{g}], R(\frac{f_1,\ldots, f_n}{g})^+) \to \Spa(R,R^+)
\]
induces an homemorphism onto $\Spa(R,R^+)\left(\frac{f_1,\ldots, f_n}{g}\right)$.  Furthermore, if two rational localizations $(R_1,R_1^+)$ and $(R_2,R_2^+)$ are such that $\Spa(R_1,R_1^+)=\Spa(R_2,R_2^+)$ when considered as subspaces in $\Spa(R,R^+)$, one has a unique isomorphism of $R$-algebras $R_1\cong R_2$ such that $R_1^+=R_2^+$.

\end{proposition} 
\begin{proof}[Sketch of the proof]
By the description of the adic spectrum in terms of prime ideals and valuation rings of  \Cref{RemarkEquivalentDescriptionSpa} it is clear that $\Psi$ is injective. To prove that it surjects onto $U= \Spa(R,R^+)\left(\frac{f_1,\ldots, f_n}{g}\right)$, it suffices to see that a point $(\f{p},V)\in \Spa(R,R^+)$ belongs to $U$ if and only if $g$ is non zero  in $\kappa(\f{p})$ and $\frac{f_i}{g}\in V$ for all $i=1,\ldots, n$. But a proofreading of the definitions shows that  these are precisely the elements of $\Spa(R[\frac{1}{g}], R(\frac{f_1,\ldots, f_n}{g})^+)$ after identifying prime ideals of both rings via the inclusion $\Spec R[\frac{1}{g}]\subset \Spec R$, using  that valuation rings are integrally closed to pass to the integral closure.  For proving that it is a homeomorphism for Huber's topology, it suffices to check that a rational localization of $\Spa(R[\frac{1}{g}], R(\frac{f_1,\ldots, f_n}{g})^+)$ is mapped to a rational localization of $\Spa(R,R^+)$, we left this as an exercise to the reader.

Finally, let $(h_1,\ldots, h_k, s)$ be another tuple generating the same rational localization as $(f_1,\ldots, f_n,g)$. First, note that one has an inclusion $F:\Spec R\to \Spec (R,R^+)$ given by mapping a prime ideal $\f{p}$ to the pair $(\f{p},\kappa(\f{p}))$.  It is easy to see that the pre-image along $F$ of a rational localization corresponding to $(f_1,\ldots, f_n,g)$ is $\Spec R[\frac{1}{g}]$.  This shows that  the localizations $R[\frac{1}{s}]$ and $R[\frac{1}{g}]$ induce the same open Zariski subspaces of $\Spec R$ and so that $R[\frac{1}{s}]=R[\frac{1}{g}]$. Then, the equality 
\[
R(\frac{f_1,\ldots, f_n}{g})^+=R(\frac{h_1,\ldots, h_k}{s})^+
\]
follows from  \Cref{PropIntegrallyClosedRings}, proving the last assertion. 
\end{proof}

\begin{definition}
Let $(R,R^+)$ be a discrete Huber pair and denote  $X=\Spa(R,R^+)$. Let $\s{O}_{X}$ and $\s{O}_X^+$ be the presheaf on rational open subsets of $X$ mapping  $\Spa(R',R^{'+})$ to $R'$ and $R^{'+}$ respectively. 
\end{definition}

A final important property of the adic spectrum of discrete Huber rings is that it is immediately sheafy: 

\begin{proposition}\label{propSheafyHuberRings}
Let $(R,R^+)$ be a discrete Huber ring and let $X=\Spa(R,R^+)$ be its adic spectrum. Denote by $X_{\mathrm{rat}}$ the site of rational localizations of $X$. Then the  structural pre-sheaf $\s{O}_X$ on $X_{\mathrm{rat}}$ is a sheaf on   the derived category of discrete $R$-modules. Similarly, $\s{O}_X^+$ is a sheaf on the abelian category of $R^+$-algebras 
\end{proposition}
\begin{proof}
There is a continuous map $F:\Spec R \to \Spa (R,R^+)$ mapping a prime ideal $\f{p}$ to the tuple $(\f{p}, \kappa(\f{p}))$. Then the pre-sheaf $\s{O}_X$ is nothing but the (derived) pushforward of the structural sheaf of $\Spec R$, which shows that is a sheaf. The sheaf property  for $\s{O}_X^+$ follows since for a rational subspace $U\subset X$ one has 
\[
\s{O}_X^+(U)=\{f\in \s{O}_X(U): \forall x\in U, |f(x)|\leq 1\}.
\]
\end{proof}

\subsubsection{Descent for the analytic topology}

In the previous paragraph  we have defined the adic spectrum of a discrete Huber pair, we defined its analytic topology and proved that the structural pre-sheaf is actually a sheaf. Out next goal is to prove that the functor mapping $(R,R^+)\mapsto \ob{D}((R,R^+)_{\sol})$ satisfies descent for the analytic topology (later we will reinterpret this as  $!$-descent on analytic rings). We need a technical lemma.

\begin{lemma}[{\cite[Lemma 2.6]{HuberAdicSpaces}, \cite[Lemma 10.4]{ClausenScholzeCondensed2019}}]\label{LemmaRefinementCover}
Let $(R,R^+)$ be a  discrete Huber pair and $X=\Spa(A,A^+)$. Assume that $U_1,\ldots, U_n\subset X$ are open  rational subsets covering $X$. There there exist $s_1,\ldots, s_N\in A$ generating the unit ideal such that each $X(\frac{s_1,\ldots, s_N}{s_j})$ is contained in some $U_i$; in particular $\{X(\frac{s_1,\ldots, s_N}{s_j})\}_j$ refines $\{U_i\}$.
\end{lemma}

\begin{theorem}\label{theo:DescendAdicSpaces}
Let $(R,R^+)$ be a discrete Huber ring and let $X=\Spa(R,R^+)$ be its adic spectrum. Let $X_{\mathrm{rat}}$ be the site of finite disjoint unions of  rational open subspaces of $X$. Then the functor 
\[
\ob{D}: X_{\mathrm{rat}}^{\op} \to \mathrm{CAlg}({\Pr}^{L,\mathrm{ex}})
\]
on presentably symmetric  monoidal stable $\infty$-categories mapping $U\in X_{\mathrm{rat}}$ to $\ob{D}((\s{O}_X(U), \s{O}_X(U)^+)_{\sol})$ is a sheaf. 
\end{theorem}
\begin{proof}
\textit{Step 1.} Let $\{U_i\}_{i\in I}$ be a finite rational cover of $X$ and let $\{U_J\}_{J\subset I}$ be the poset of finite intersections of the $U_i$. We want to prove that the natural map 
\begin{equation}\label{eqDescentCategories}
\ob{D}((R,R^+)_{\sol})\to \varprojlim_{U_J\subset X} \mathrm{\ob{D}}((\s{O}_X(U_J),\s{O}_X(U_J)^+)_{\sol}) 
\end{equation}
is an equivalence of categories.


 By  \Cref{LemmaRefinementCover} we can find elements $s_1,\ldots, s_N$ generating the unit ideal of $R$ such that the rational cover $\{X(\frac{s_1,\ldots, s_N}{s_j})\}_{j}$ refines $\{U_i\}$. Therefore, it suffices to show descent with respect to the covers of the form $\{U_j=X(\frac{s_1,\ldots, s_N}{s_j})\}_{j}$.

  \textit{Step 2.} Let us first consider the case when  $s_N$ is a unit in $R$. By replacing $s_j$ by $s_{j}/s_N$ we can assume that $s_N=1$. We claim that  the cover $\{U_j=X(\frac{s_1,\ldots, s_{N-1},1}{s_j})\}_{j}$  has a refinement given by composition of   Laurent covers
 \begin{equation}\label{eqLaurentCover}
 X= X(\frac{1}{f})\cup X(\frac{f}{1}).  
 \end{equation}
 Indeed, it is refined by the intersection of the Laurent covers $X=X(\frac{1}{s_j})\cup X(\frac{s_j}{1})$ for the different $j$. Therefore, by an induction argument, it suffices to deal with the case of a single Laurent cover  as in \eqref{eqLaurentCover}. But this cover corresponds to the maps of analytic rings
 \[
 (R,R^+)_{\sol}\to \bb{Z}[T]_{\sol}\otimes_{\bb{Z}[T]}(R,R^+)=(R, R(\frac{f}{1})^+) 
 \]
 and
 \[ (R,R^+)_{\sol}\to (\bb{Z}[T^{\pm 1}], \bb{Z}[T^{-1}])_{\sol}\otimes_{\bb{Z}[T]}(R,R^+)=(R[\frac{1}{f}], R(\frac{1}{f})^+)
 \]
 where $T$ maps to $f$.    \Cref{LemmaBasicOpenImmersions} implies that these maps give rise to open immersions at the level of locales. Therefore, in order to show descent, by base change along $\bb{Z}[T]\to (R,R^+)$ ($T\mapsto f$),  and by  \Cref{theo:openCoversSym}, it suffices to prove that the open localizations $(\bb{Z}[T],\bb{Z})_{\sol}\to \bb{Z}[T]_{\sol}$ and $(\bb{Z}[T],\bb{Z})_{\sol}\to (\bb{Z}[T^{\pm 1}], \bb{Z}[T^{-1}])_{\sol}$ cover the locale $\ob{D}((\bb{Z}[T],\bb{Z})_{\sol})$. This amounts to show that the tensor product of the solid $\bb{Z}[T]$-idempotent algebras $\bb{Z}((T^{-1}))$ and $\bb{Z}[[T]]$ vanish. But we have that 
 \[
\bb{Z}((T^{-1}))\otimes_{\bb{Z}[T]} \bb{Z}[[T]] = \bb{Z}[[X,T]]^{\bb{L}} / (XT-1),
 \]
  but $XT-1$ is a unit in $\bb{Z}[[X,T]]$ proving that the tensor vanishes.

 \textit{Step 3.} We now deal with the general case.  The elements $s_1,\ldots, s_N$  give rise to a Zariski cover $\{\Spec R[\frac{1}{s_j}]\}_j$ of $\Spec R$, we let $V_j=X(\frac{g}{g})$ be the locus where $|g|\neq 0$, and let $\{V_K\}_{K\subset J}$ denote the poset of finite intersections of the $V_j$'s.   Note that $R(\frac{g}{g})^+$ is nothing but the integral closure of $R^+$ in $R[\frac{1}{g}]$. Then, by  \Cref{lem:IntegralClosureSolid} we have that $(R[\frac{1}{g}],R^+)_{\sol}=(R[\frac{1}{g}],R(\frac{g}{g})^+)_{\sol}$.  By Zariski descent  the natural map
 \[
 \ob{D}((R,R^+)_{\sol})= \varprojlim_{V_K \subset X} \ob{D}( (\s{O}_X(V_K), \s{O}_X^+(V_K) )_{\sol})
 \]
is an equivalence. Concretely, the Zariski cover gives rise to a \textbf{closed} cover of the locale of $\ob{D}(R)$, and by  \Cref{theo:openCoversSym} the base change to $\ob{D}((R,R^+)_{\sol})$ gives rise to a closed cover as well.   By Step 2 the rational cover $\{V_K\cap U_j\}_j$ gives rise an open cover of the locale $\Sm(\ob{D}( (\s{O}_X(V_K), \s{O}_X^+(V_K)  )_{\sol}))$ and so we have 
\[
\ob{D}((\s{O}_X(V_K), \s{O}_X^+(V_K)  )_{\sol})=\varprojlim_{U_J\subset X}  \ob{D}( (\s{O}_X(V_K\cap U_J), \s{O}_X^+(V_K\cap U_J)  )_{\sol}),
\]
taking limits with respect to the poset $\{V_K\}$ we get then
\[
\begin{aligned}
\ob{D}((R,R^+)_{\sol}) & = \varprojlim_{V_K\subset X}  \ob{D}((\s{O}_X(V_K), \s{O}_X^+(V_K)  )_{\sol}) \\ 
& = \varprojlim_{V_K\subset X} \varprojlim_{U_J\subset  X} \ob{D}( (\s{O}_X(V_K\cap U_J), \s{O}_X^+(V_K\cap U_J)  )_{\sol}) \\ 
& = \varprojlim_{U_J\subset  X}  \varprojlim_{V_K\subset X}  \ob{D}( (\s{O}_X(V_K\cap U_J), \s{O}_X^+(V_K\cap U_J)  )_{\sol})\\
& =\varprojlim_{U_J\subset  X}  \ob{D}( (\s{O}_X( U_J), \s{O}_X^+( U_J)  )_{\sol})
\end{aligned}
\]
where the first equivalence is Zariski descent for $X$, the second follows from Step 2 applied to the rational subspaces $V_K$, the third is a commutation of  limits, and the last is Zariski descent applied to the rational subspaces $U_J$. This finishes the proof of the theorem. 
\end{proof}

In order to define the correct analogue of  \Cref{Coro:ZariskiDescent,coroSolidSheavesSchemes} for discrete adic spaces we need to modify Huber's topology on $\Spa(R,R^+)$. 

\begin{definition}
Let $(R,R^+)$ be a discrete Huber pair. We let $X=\Spa(R,R^+)^{\mathrm{mod}}$ be the adic spectrum endowed with the coarsest topology making the rational localizations of the form $X(\frac{1}{f})$ and $X(\frac{f}{1})$ open ($f\in R$), and the rational localizations of the form $X(\frac{g}{g})$ closed ($g\in R$).  We call $X=\Spa(R,R^+)^{\mathrm{mod}}$  the \textit{modified adic spectrum}.
\end{definition}

\begin{lemma}\label{LemSpectralModified}
Let $(R,R^+)$ be a discrete Huber pair and let $X=\Spa(R,R^+)^{\ob{mod}}$  be its modified adic spectrum. Then $X$ is a spectral space with qcqs open subspaces given by finite unions of intersections of subspaces of the form $X(\frac{1}{f})$, $X(\frac{f}{1})$ and $X(g=0):=X\backslash X(\frac{g}{g})$. 
\end{lemma}
\begin{proof}
By definition, $\Spa(R,R^+)^{\ob{mod}}$ is the coarsest topology making the two maps 
\[
t\colon \Spa(R,R^+)\to \Spec(R)  \mbox{ and } s\colon \Spa(R,R^+)\to \Spec(R)^{\op}
\]
continuous, where $s$ and $t$ are the maps that send a valuation determined by a pair $x=(V\subset \kappa(\f{p}))$ as in \Cref{RemarkEquivalentDescriptionSpa} to $s(x)=\f{p}$  and $t(x)=\f{m}_V$ (the maximal ideal of $V$). 

 Any subspace given by finite unions of intersections of subspaces as in the lemma is qcqs for the modified topology, namely, it is clopen in the constructible topology of $\Spa(R,R^+)$. Since by definition those objects form a basis for the modified topology and are stable under finite unions, they consist on the qcqs open subspaces of $\Spa(R,R^+)^{\ob{mod}}$. It is easy to check that the space is $T_0$ (or Kolmogorov), this makes $\Spa(R,R^+)^{\ob{mod}}$ a spectral space thanks to Hochster's criteria \cite[Proposition 3.31]{wedhorn2019adicspaces}. 
\end{proof}

\begin{corollary}\label{corodiscreteAdicSpacesAsLocales}
Let $(R,R^+)$ be a discrete Huber pair and let $\Spa(R,R^+)^{\mathrm{mod}}$ be the modified adic spectrum. Then the functor mapping a rational subset $U\subset X$ to $\ob{D}((\s{O}_X(U),\s{O}_X^+(U))_{\sol})$ gives rise to a unique  natural surjective map of locales
\[
F:\Sm(\ob{D}((R,R^+)_{\sol}))\to \Spa(R,R^+)^{\mathrm{mod}}. 
\]
\end{corollary}
\begin{proof}
The existence of the map $F$ follows  formally from  \Cref{theo:DescendAdicSpaces} and \Cref{Coro:ZariskiDescent,coroSolidSheavesSchemes}. 

The  only thing to verify is the surjectivity of $F$.  We want to prove that given $U,U'\subset \Spa(R,R^+)^{\mathrm{mod}}$ two open subspaces for the modified topology, if $F^{-1}(U)=F^{-1}(U')$ then $U=U'$.  For this, let  $x\in \Spa(R,R^+)$  and let $(\kappa(x),\kappa(x)^+)$ be the residue field at $x$. We can assume without loss of generality that $U'\subset U$.  Then by naturality of $F$ and by taking pullbacks along all affinoid points, it suffices to deal with the case where $R=K$ is a field and $K^+\subset K$ an integrally closed subring. In this case, the analytic and the modified topology agree, and the open subspaces of $\Spa (K,K^+)$ form a poset by inclusion. Quasi-compact open subspaces of $\Spa(K,K^+)$ are of the form $\Spa(K,K^{',+})$ for $K^+\subset K^{',+}$ an integrally closed ring.   Thus, by taking intersections of $U$  with open affinoids, to prove surjectivity it suffices to show that for $U\subset \Spa(K,K^+)$ an arbitrary open subspace, if $F^{-1}(U)=F^{-1}(\Spa(K,K^+))$ then $U=\Spa(K,K^+)$. Suppose this does not hold and let $x\in \Spa(K,K^+)$ be the maximal closed point, then $F^{-1}(x)=\emptyset$. Let $A_x$ be the idempotent algebra corresponding to $x$, and let $\{Z_j\}_{j}$ be the poset of closed subspaces with qcqs open complements  with idempotent algebras $A_{Z_j}$. Then we have that 
\[
0=A_x=\varinjlim_{j} A_{Z_j}. 
\]
Since the colimit is filtered, there exists some $j$ such that $1\in A_{Z_j}$ is zero and so $A_{Z_k}=0$ for all $k\geq j$. But this will imply that $F^{-1}(X\backslash Z_j)= F^{-1}(X)$ and $X\backslash Z_j=\Spa(K,K^{',+})$ is affinoid with $K^+$ strictly contained in $K^{',+}$. This contradicts the fully faithfulness of discrete Huber pairs into analytic rings of  \Cref{Prop:FUllyFaithDiscreteAnalytic}, proving what we wanted. 
\end{proof}

\begin{remark}
We have not stated  \Cref{corodiscreteAdicSpacesAsLocales} for general discrete adic spaces since the analytic topology differs from the modified topology, so that discrete adic spaces are not obtained from gluing affinoid spaces along the modified topology. Instead, since the analytic topology contains both closed and open subspaces of the modified topology, discrete adic spaces will be examples of  gluing affinoids via $!$-covers, see \Cref{ss:AnStacksConstruction}.
\end{remark}


\section{Analytic stacks}\label{s:AnStacks}

In  \Cref{SubsecMoreSolidRings} we have constructed categories of classical and solid quasi-coherent sheaves attached to schemes; the solid theory mapped open Zariski subspaces to honest open subspaces in the sense of locales, while the classical theory mapped open Zariski subspaces to closed subspaces. We also saw that discrete adic spaces could be endowed with a theory of solid quasi-coherent sheaves, thought for this last Huber's open subspaces gave rise to locally closed subspaces in the locale. A way to overcome this apparent paradox is to introduce a much larger category of \textit{analytic stacks} where the gluing process can be made for a very large class of ``descendable'' maps\footnote{Descendable algebras in a stable category were introduced by Mathew in \cite{MathewDescent}, the version of descendability appearing in the theory of analytic stacks is a generalization that makes use of the six functor formalisms for analytic rings.}, including open, closed, or locally closed covers.

In order to define categories\footnote{I intentionally say ``categories'' and not ``category'' since in some situations it is better to isolate a convenient subcategory of analytic rings that will suffice to construct the analytic stacks we are interested in.} of analytic stacks, we shall need to recall some facts about six functor formalisms and the category of kernels, see \cite{MannSix,MannSix2},  \cite{SixFunctorsScholze} and \cite{HeyerMannSix}. In fact, the theory of analytic stacks is built up in harmony together with its theory of quasi-coherent sheaves, not being possible to mention one without the other. We then define categories of analytic stacks  and mention some basic facts about their theory of six functor formalisms. We finish this section with a discussion of six functors in algebraic geometry and sheaves on topological spaces via the category of analytic stacks.

\subsection{$6$-functor formalisms}\label{ss:SixFunctors}

The idea of constructing theories of six functor formalisms goes back to Grothendieck where he first constructed such operations for the \'etale cohomology of schemes. Since then, different theories of six functors have appeared in mathematics (such as the theory of $D$-modules, sheaves on topological spaces or motives). However, at least to author's knowledge,  it was not until recently that a formal definition of an abstract six functor formalism was introduced in the work of Mann \cite[Appendix A.5]{MannSix} (building on previous works of Liu-Zheng for \'etale sheaves and Gaitsgory-Rozemblyum for $\ob{IndCoh}$ sheaves). 

The theory of abstract six functor formalisms was then further explored by Scholze  \cite{SixFunctorsScholze}, and its current standard referent is the work of Heyer-Mann \cite{HeyerMannSix}.  In this section we will give a really basic and rough introduction to the theory of abstract six functors. We strongly encourage the interested reader to look at the papers \cite{MannSix,MannSix2,HeyerMannSix,SixFunctorsScholze}
for a complete and better exposition of the subject. Let us higlight that the theory of abstract six functor formalisms is essential in the theory of analytic stacks. 

\subsubsection{Geometric set ups}  The first piece of datum necessary to define a six functor formalism is an ambient category and a marked class of arrows:

\begin{definition}\label{DefGeoSetUp}
A \textit{geometric set up} is a tuple $(\n{C},\n{C}_0)$ consisting on an $\infty$-category and (a non full) subcategory $\n{C}'$ satisfying the following conditions:

\begin{enumerate}

\item $\n{C}_0$ contains all the equivalences of $\n{C}$ (i.e. is a wide subcategory).

\item If $f\colon Y\to X$ is in $\n{C}_0$ and $g\colon X'\to X$ is in $\n{C}$, then the pullback $Y'\to X'$ of $f$ along $g$ exists in $\n{C}$ and it belongs to $\n{C}_0$.

\item The category $\n{C}_0$ admits fiber products and the inclusion $\n{C}_0\subset \n{C}$ preserves them. 
\end{enumerate}

\end{definition}

\begin{remark}\label{RemEquivDefGeoStup}
The datum of a geometric set up is equivalent to the datum of a pair $(\n{C},E)$ where $E$ is a class of arrows in $\n{C}$ satisfying the following conditions (cf. \cite[Lemma 2.1.5]{HeyerMannSix}):
\begin{enumerate}[label = (\roman*)]

\item $E$ contains all isomorphims,

\item $E$ is stable under compositions, and it admits and is stable under pullbacks of edges in $\n{C}$,

\item for every $f\colon Y\to X$ in $E$, the diagonal $f\colon Y\to Y\times_{X} Y$ (which exists by (ii)) is in $E$. 
\end{enumerate}
In this situation, we let $\n{C}_E\subset \n{C}$ be the wide subcategory of $\n{C}$ spanned by the edges in $E$. The pair $(\n{C},\n{C}_E)$ is a geometric set up in the sense of \Cref{DefGeoSetUp}. 
\end{remark}

The idea behind a geometric set up is that the category $\n{C}$ is the underlying category of geometric objects (eg. schemes, adic spaces, complex varieties, analytic stacks, etc.), while $E$ is the class of ``$!$-able maps'', that is, the class of maps that will admit lower $!$-functors in a suitable definition of abstract  $3$-functor formalisms over the geometric set up $(\n{C},E)$. The following lemma shows that the class of arrows in $E$ is right cancellative:

\begin{lemma}
Consider maps $Z\xrightarrow{g} Y \xrightarrow{f} X$ in $\n{C}$. If $f\circ g$ and $f$ are in $E$ then so is $g$. 
\end{lemma}
\begin{proof}
We can write the map $g$ as the composite $Z\xrightarrow{\id_Z\times g} Z\times_X Y \xrightarrow{\pr_2} Y$, so it suffices to show that $\pr_2$ and $\id_Z\times g$ are in $E$. But $\pr_2$ is the base change of $Z\xrightarrow{f\circ g} X$ along $Y\xrightarrow{f} X$ so it is in $E$, and $\id_Z\times g$  is the base change of $Z\times_X Y\xrightarrow{g\times \id_Y} Y\times_X Y$ along the diagonal map $\Delta_f:Y\to Y\times_X Y$ which is also in $E$. 
\end{proof}

\subsubsection{The category of correspondences}

Given a geometric set up $(\n{C},E)$ one constructs an $\infty$-category of correspondences $\Corr(\n{C},E)$
\cite[Definition 3.1]{SixFunctorsScholze}  \cite[Definition 2.2.1]{HeyerMannSix}. Its construction is combinatorially involved and proving that it is an $\infty$-category requires some work \cite[Proposition 2.2.9]{HeyerMannSix}. The intuition behind the category of correspondences is that it is a blueprint for the functors $f^*$ and $g_!$ in a six functor formalism. Informally, the category of correspondences  $\Corr(\n{C},E)$ is the $\infty$-category with objects the same objects as $\n{C}$, a morphism from $X$ to $Y$ is a diagram in $\n{C}$
\[
\begin{tikzcd}
W \ar[d,"f"']\ar[r,"g"] & Y \\
X
\end{tikzcd}
\]
with $g\in E$. A composition of two maps $X \to Y$ and $Y\to Z$ in $\Corr(C,E)$ is the outer correspondence of the diagram with pullback square
\[
\begin{tikzcd}
W\times_Y U  \ar[r] \ar[d] & U  \ar[r] \ar[d] & Z\\ 
W \ar[d] \ar[r] & Y & \\
X & & 
\end{tikzcd}
\]
(note that for the arrow $W\times_Y U\to Z$ to be in $E$ one requires $E$ to be stable under pullbacks and compositions). 

In order to define a six functor formalism from the category of correspondences it is necessary to promote it to a suitable $\infty$-operad denoted by $\Corr(\n{C},E)^{\otimes}$. If $\n{C}$ has finite direct products this $\infty$-operad endows the correspondence category with a symmetric monoidal structure, see \cite[\S 2.3]{HeyerMannSix}. The construction of $\Corr(\mathcal{C},E)^{\otimes}$ is as follows:

 First, given an $\infty$-category $\mathcal{C}$ there is the co-cartesian $\infty$-operad $\mathcal{C}^{\sqcup}$ \cite[Proposition 2.4.3.3]{HigherAlgebra}. If $\mathcal{C}$ has coproducts then $\mathcal{C}^{\sqcup}$ is the symmetric monoidal operad whose objects are pairs $(I,(X_i)_{i\in I})$ with $I$ a finite set and $(X_{i})_{i\in I}$ a collection of elements of $\mathcal{C}$, and morphisms $(I,(X_i)_{i\in I})\to (J,(Y_j)_{j\in J})$ consists on a map of pointed sets
 \[
 f: I\cup \{*\}\to J \cup\{*\}
 \]
 together with maps for all $j\in J$
 \[
 \bigsqcup_{i\in f^{-1}(j)\cap I} X_i\to Y_j. 
 \]
Consider the co-operad $(\mathcal{C}^{\op,\sqcup})^{\op}$. If $\mathcal{C}$ has finite products, one can describe informally this category as follows: its  objects are pairs $(I,(X_{i})_{i\in I})$, and morphisms $(I,(X_{i})_{i\in I})\to (J,(Y_j)_{j\in J})$ are maps of pointed sets
\[
f:J\cup\{*\}\to I\cup\{*\}
\]
and  for all $i\in I$
\[
X_i\to \prod_{j\in f^{-1}(i)\cap J} Y_j.
\]

Now, the category $\Corr(\mathcal{C},E)^{\otimes}$ is itself the correspondence category of $(\mathcal{C}^{\op,\sqcup})^{\op}$ for a suitable choice of arrows $\widetilde{E}$. An arrow $f:(I,(X_i)_{i\in I})\to (J, (Y_j)_{j\in J})$ is in  $\widetilde{E}$ if $f:J\to I$ is an isomorphism and for all $j\in J$ the map $ X_{f(j)}\to Y_j$ lies in $E$. It turns out that $((\mathcal{C}^{\op,\sqcup})^{\op},\widetilde{E})$ forms a geometric set up. One defines the $\infty$-operad
\[
\Corr(\mathcal{C},E)^{\otimes}:= \Corr((\mathcal{C}^{\op,\sqcup})^{\op},\widetilde{E}).
\]

Let's describe more explicitly how $\Corr(\mathcal{C},E)^{\otimes}$ looks like. Objects of $\Corr(\mathcal{C},E)^{\otimes}$ are the same as objects of $(\mathcal{C}^{\op,\sqcup})^{\op}$, i.e.  pairs $(I,(X_i)_{i\in I})$. A morphism $(I,(X_i)_{i\in I})\dashrightarrow (J,(Y_j)_{j\in J})$ in   $\Corr(\mathcal{C},E)^{\otimes}$ is a correspondence in $(\mathcal{C}^{\op,\sqcup})^{\op}$
\[
\begin{tikzcd}
(J, (Z_j)_{j\in J}) \ar[d, "f"] \ar[r, "g" ] & (J, (Y_j)_{j\in J} ) \\ 
(I, (X_{i})_{i\in I}) &
\end{tikzcd},
\]
where $g$ is in $\widetilde{E}$.  More precisely, we have a map of pointed sets $f: I\cup\{*\}\to J \cup \{*\}$, and for all $j\in J$ a correspondence in $\mathcal{C}$
\begin{equation}\label{eqasjdi}
\begin{tikzcd}
Z_j \ar[d, "f"] \ar[r, "g_j"] & Y_j \\ 
 \prod_{i\in I\cap f^{-1}(j)} X_i & 
\end{tikzcd}
\end{equation}
with $g_j\in E$  (cf. \cite[Definition 3.11]{SixFunctorsScholze}).

The category of correspondences $\Corr(\mathcal{C},E)$ admits two fundamental maps:
\begin{equation}\label{eqan9e}
F:\mathcal{C}^{\op}\to \Corr(\mathcal{C},E) 
\end{equation}
and 
\begin{equation}\label{eqndjd}
G:\mathcal{C}_E \to \Corr(\mathcal{C},E).
\end{equation}
Both maps $F$ and $G$ are informally constructed as follows: they are the identity in objects. A map $f^{\op}: X\to Y$ in $\mathcal{C}^{\op}$, corresponding to a map $f:Y\to X$ is sent to the correspondence
\[
\begin{tikzcd}
Y \ar[r, "\id"]  \ar[d,"f"] &  Y\\ 
X &
\end{tikzcd}.
\]
 A map $g:X\to Y$ in $\mathcal{C}_E$ is sent to the correspondence
 \[
 \begin{tikzcd}
 X \ar[d,"\id"] \ar[r,"g"] & Y \\ 
 X &
 \end{tikzcd}.
 \]
 
 The map $F$ can be even upgraded to a map of operads $\mathcal{C}^{\op,\sqcup}\to \Corr(\mathcal{C},E)^{\otimes}$, it is even symmetric monoidal if $\mathcal{C}$ admits finite products. Similarly, the map $G$ can be upgraded to a map of operads $\mathcal{C}^{\times}_E \to \Corr(\mathcal{C},E)^{\otimes}$. 
\subsubsection{Abstract six functor formalisms} 
 
 Let $\Cat_{\infty}^{\times}$ be the $\infty$-category of (small) $\infty$-categories endowed with the catesian symmetric monoidal structure. 
 
 \begin{definition}\label{defSixFunc}
 Let $(\mathcal{C},E)$ be a geometric set up. A \textit{$3$-functor formalism} $\ob{D}$ on $(\mathcal{C},E)$ is a morphism of $\infty$-operads
 \[
 \ob{D}:\Corr(\mathcal{C},E)^{\otimes} \to \Cat_{\infty}^{\times}. 
 \]
 \end{definition}

\begin{remark}\label{RemarkOtherVersions}
Let $(\mathcal{C},E)$ be a geometric set-up. The definition of a $3$-functor formalism is quite flexible in the target. Namely, given $\s{D}$ any symmetric monoidal $2$-category one can define a \textit{$\s{D}$-valued $3$-functor formalism} to be a lax symmetric monoidal functor $\ob{D}\colon \Corr(\mathcal{C},E)\to \s{D}$.  If $\s{C}$ has a final object $*$, then the same argument of  \cite[Lemma 31.9]{HeyerMannSix} shows that $\ob{D}$ factors naturally through a lax symmetric monoidal functor
\[
\ob{D}\colon \Corr(\mathcal{C},E)\to \ob{Mod}_{\ob{D}(*)}(\s{D}).
\]

In case $\s{D}=\Pr^L$ is the $2$-category of presentable categories endowed with Lurie's tensor product, we call a $\Pr^L$-valued $3$-functor formalism a \textit{presentable $3$-functor formalism.}  If $\n{C}$ has a final object $*$,   $\ob{D}$- factors through the category $\Pr^L_{\ob{D}(*)}:=\ob{Mod}_{\ob{D}(*)}(\Pr^L)$ of \textit{$\ob{D}(*)$-linear presentable categories.}
\end{remark}

Let us briefly explain how the functor $\ob{D}$ encodes the important features of a $6$-functor formalism, we send to \cite[Section 3.1]{HeyerMannSix} for a more detailed explanation.

\begin{itemize}
\item[(a)] First, composing $\ob{D}$ with the functor $\mathcal{C}^{\op}\to \Corr(\mathcal{C},E)$ we get a lax symmetric monoidal functor
\[
\mathcal{C}^{\op,\sqcup}\to \Cat_{\infty}^{\times},
\]
this is precisely the same datum as a functor 
\[
\mathcal{C}^{\op}\to \ob{CAlg}(\Cat^{\times}_{\infty})
\]
from $\mathcal{C}^{\op}$ to symmetric monoidal categories (\cite[Theorem 2.4.3.18]{HigherAlgebra}). This could be called the \textit{coefficient functor of $\ob{D}$}. By an abuse of notation we often identify $\ob{D}$ with its restriction to $\n{C}^{\op}$ in case the class of $!$-able arrows and the extenstion of $\ob{D}$ from $\n{C}^{\op}$ to $\Corr(\n{C},E)$ is clear from the context. 

\item[(b)] Any morphism $Z:X\to Y$ in $\Corr(\mathcal{C},E)$ depicted as a diagram
\[
\begin{tikzcd}
Z\ar[r,"g"] \ar[d,"f"'] & Y \\
X 
\end{tikzcd}
\]
can be written as a composite $X\to Z \to Y$, where the first arrow is $X\leftarrow Z = Z$ and the second is $Z=Z\to Y$. Given $f:X\to Y$ a map in $\mathcal{C}$ we define the \textit{pullback functor} $f^*:\ob{D}(Y)\to \ob{D}(X)$ to be $\ob{D}(Y\xleftarrow{f} X=X)$, if $f\in E$ we define $f_!:\ob{D}(X)\to \ob{D}(Y)$ to be $\ob{D}(X=X\to Y)$. Then $\ob{D}$ sends $Z:X\to Y$ to the functor $g_{!}f^*:\ob{D}(X)\to \ob{D}(Y)$.

In this way, the class of arrows $E$ in $\mathcal{C}$ are precisely those arrows that admit $!$-functors, one often call $E$ the class of \textit{$!$-able arrows}. The functoriality of the functors $F$ and $G$ of \eqref{eqan9e} and \eqref{eqndjd} encode the natural compatibility  of the functors $f^*$ and $f_!$ under composition. 

\item[(c)] \textit{Proper base change}. Let $g: Y\to X$ be an arrow in $E$ and let $f:X'\to X$ be a map in $\mathcal{C}$. Consider the cartesian square
\[
\begin{tikzcd}
Y' \ar[r,"{g'}"] \ar[d,"f'"] & X' \ar[d,"f"] \\
Y \ar[r,"g"] & X.
\end{tikzcd}
\]
The proper base change consists on the natural equivalence $f^*g_!\cong g'_{!}f^{'*}$ of functors $\ob{D}(Y)\to \ob{D}(X')$. It arises from the functor $\ob{D}$ as follows: let $h_1=(Y=Y\xrightarrow{g} X)$ and $h_2=(X\xleftarrow{f} X'=X')$ be maps in the correspondence category. Then $h_2\circ h_1= (Y\xleftarrow{f'} Y' \xrightarrow{g'} X')$ and we have 
\[
f^*g_!=\ob{D}(h_2)\circ \ob{D}(h_1) = \ob{D}(h_2\circ h_1) =  g'_{!}f^{'*}.
\]
Note that since this formula only uses the $*$ and $!$-functors,  one  has this identity directly from $\Corr(\mathcal{C},E)$ and its operad promotion is not necessary.

\item[(d)] \textit{Projection formula}. Let $g: X\to Y$  be a map in $E$. The projection formula produces for any $N\in \ob{D}(X)$ and $M\in \ob{D}(Y)$ a natural isomorphism 
\begin{equation}\label{eqnsjs9}
g_!(g^*M \otimes N) \cong M \otimes g_!N.
\end{equation}
It involves all three functors, so it is important to use the operadic version $\Corr(\mathcal{C},E)^{\otimes}$ of the correspondence category.  The projection formula is an equivalence of functors
\begin{equation}\label{eqbshwi}
\ob{D}(Y)\times \ob{D}(X) \to \ob{D}(Y).
\end{equation}
Thus, at the level of the correspondence category we have to find a map 
\[
(\{1,2\}, (Y,X))\to (\{1\}, Y)
\]
that can be written in two different ways yielding the two functors of  \eqref{eqbshwi}.  Consider the commutative diagram in $\Corr(\mathcal{C},E)^{\otimes}$
\[
\begin{tikzcd}
(Y,X) \ar[r] \ar[d] & (Y,Y) \ar[d] \\
X \ar[r] & Y
\end{tikzcd}
\]
where the top horizontal map is $(Y,X)\xrightarrow{(\id_Y,g)} (Y,Y)$, the  lower horizontal map is $g$, the left vertical map is given by $(Y,X) \xleftarrow{(g,\id_X)} X =  X$ and the right vertical map is given by the diagonal $(Y,Y) \xleftarrow{\Delta} Y = Y$. The upper-right composition gives rise the functor 
\[
(M,N)\mapsto (M,g_! N) \mapsto M\otimes g_!N
\]
while the left-lower composition gives rise the functor 
\[
(M,N)\mapsto (g^*M\otimes N) \mapsto g_!(g^*M\otimes N).
\]
As both the upper-right and left-lower compositions give rise to the same functor, we have the desired natural equivalence \eqref{eqnsjs9}.

\end{itemize}

\subsubsection{ $*$ and $!$-descent}

For future reference let us introduce the following notions of descent on a presentable $3$-functor formalism. 

\begin{definition}\label{DefinitionConditionDescent}
Let $(\n{C},E)$ be a geometric set up and $\ob{D}$ a presentable $3$-functor formalism on $(\n{C},E)$.  Suppose that $\n{C}$ has finite products and that the inclusion maps $X\to X\bigsqcup Y$ are $!$-able. 

\begin{enumerate}[label=(\alph*)]

\item A morphism $f:Y\to X$ in $\n{C}$ with \v{C}ech nerve $Y^{\bullet}$ satisfies \textit{$\ob{D}^*$ or $*$-descent} if the natural map  of categories
\[
\ob{D}^*(X)\to \ob{Tot}(\ob{D}^*(Y^{\bullet}))
\]
is an equivalence, where all the transition maps are induced by upper $*$-functors. We say that $f$ satisfies \textit{universal $\ob{D}^*$ or $*$-descent} if the previous holds after pullback along any map $X'\to X$ in $\n{C}$. 

\item A morphism $f:Y\to X$ in $E$ with \v{C}ech nerve $Y^{\bullet}$ satisfies \textit{$\ob{D}^!$ or $!$-descent} if the natural map of categories
\[
\ob{D}^!(X)\to \ob{Tot}(\ob{D}^!(Y^{\bullet}))
\]
is an equivalence, where all the transition maps are induced by upper $!$-functors. We say that $f$ satisfies \textit{universal $\ob{D}^!$ or $!$-descent} if the previous holds after pullback along any  map $X'\to X$ in $\n{C}$.

\end{enumerate}
\end{definition}

\begin{remark}\label{RemarkUniversalShierkDescent}
Universal $!$-descent is a very strong condition. Indeed, by \cite[Theorem 5.12]{SixFunctorsScholze} a map $f\colon Y\to X$ in $E$ satisfies universal $*$ and $!$-descent if and only if $f$ satisfies $!$-descent after pullying back along elements $X'\to X$ in $E$. Furthermore, the morphism $f$ will satisfy universal $*$ and $!$-descent after pullback along an arbitrary presheaf $T\to X$ as consequence of the extension theorem of six functors of \cite[Theorem 5.19]{SixFunctorsScholze}, this was the original definition of \textit{universal $!$-descent} in a previous version of \cite{SixFunctorsScholze}. 

Proving these facts requires the introduction of the presentable categories of kernels that is beyond the scope of these notes. We will content ourselves with an overview of the standard category of kernels in \Cref{ss:KernelsSuavePrim}, and with a special case of \cite[Theorem 5.12]{SixFunctorsScholze} when the six functor formalism is symmetric monoidal, i.e. satisfies categorical K\"unneth, see \Cref{PropDescentShierk}. 
\end{remark}

\subsection{The kernel category,  suave and prim objects}\label{ss:KernelsSuavePrim}

The power of a six functor formalism comes when one looks at its category of kernels. 

\begin{definition}
Let $(\n{C},E)$ be a geometric set up and $\ob{D}$ a $3$-functor formalism on $(\n{C},E)$. Let $S\in \n{C}$ and consider $(\n{C}_{E})_{/S}$ the category of $!$-able maps over $S$. Let $\ob{D}_S$ be the restriction of $\ob{D}$ to $(\n{C}_{E})_{/S}$. The \textit{category of kernels} $\ob{K}_{\ob{D},S}$ over $S$ is the $(\infty,2)$-category obtained by transfer of enrichment along the lax symmetric monoidal functor 
\[
\ob{D}_S: \ob{Corr}((\n{C}_{E})_{/S})^{\otimes}\to \ob{Cat}_{\infty}^{\times}. 
\]
\end{definition}

More informally, the category of kernels over $S$ is the $2$-category described as follows: 

\begin{enumerate}[label = (\roman*)]
\item It has  same objects as $(\n{C}_{E})_{/S}$.

\item  Given two objects $X,Y\in \ob{K}_{\ob{D},S}$ the category of functors $\ob{Fun}_S(X,Y)=\ob{Fun}_{\ob{K}_{\ob{D},S}}(X,Y)$ is given by 
\[
\ob{Fun}_S(X,Y)=\ob{D}(Y\times_S X). 
\]

\item The unit in $\ob{Fun}_S(X,X)=\ob{D}(X\times_S X)$ is given by the object $\Delta_! 1_{X}$ where $\Delta : X\to X\times_S X$ is the diagonal map.

\item The composition law $\ob{Fun}_{S}(Y,Z)\times \ob{Fun}_{S}(X,Y)\to \ob{Fun}_S(X,Z)$ is given by the convolution product
\[
\star: \ob{D}(Z\times_S Y)\times \ob{D}(Y\times_S X)
\]
where for $M\in \ob{D}(Z\times_S Y)$ and $N\in \ob{D}(Y\times_S X)$ one has  
\[
M\star N  = \pi_{13,!}(\pi_{12}^*M\otimes \pi_{23}^* N)
\]
where $\pi_{ij}$ are the projection maps of $Z\times_S Y \times_S X$. 

\end{enumerate}

\begin{remark}\label{RemarkPromotionSymmetricMonoidal2}
The formation of the category of kernels has different functorialities: one arises from changing the  $3$-functor formalism $\ob{D}$, others  srise from pullbacks and lower shierk maps along $Y\to X$ in $\n{C}$ or $\n{C}_E$ respectively. In general, the category of kernels accommodate itself in a $3$-functor formalism of $2$-categories, namely, a lax symmetric monoidal functor 
\[
\ob{K}_{\ob{D},-}:\ob{Corr}(\n{C},E)^{\otimes} \to \ob{Cat}_{(\infty,2)}^{\times},
\]
see \cite[Theorem 4.2.4]{HeyerMannSix}.

In particular, the  category of kernels $\ob{K}_{\ob{D},S}$ is naturally endowed with a symmetric monoidal $2$-categorical structure, and  the $2$-functor 
\[
\Psi: \ob{K}_{\ob{D},S}\to \ob{Cat}_{\infty, /\ob{D}(S)}^{\times}
\]
given by $\ob{Fun}_S(S,-)$ is lax symmetric monoidal. 
\end{remark}

\begin{remark}\label{RemarkPromotionSymmetricMonoidal}
One has a natural  symmetric monoidal functor arising from transfer of enrichement
\[
\Phi: \ob{Corr}((\n{C}_{E})_{/S}) \to \ob{K}_{\ob{D},S}
\]
from the category of correspondences over $S$ (which is nothing but a span category in this case) to the category of kernels. It is the identity on objects and sends a correspondence $X\xleftarrow{f} Z \xrightarrow{g}Y$ to the morphism $(g,f)_! 1_{Z}\in \ob{D}(Y\times_S X)=\ob{Fun}_S(X,Y)$ where $(g,f):Z\to Y\times_S X$ is the graph. In particular, since any object in the correspondence category is self dual, all the objects in $\ob{K}_{\ob{D},S}$ are self dual as well. This produces a natural antiinvolution of $2$-categories $\ob{K}_{\ob{D},S}^{\op}= \ob{K}_{\ob{D},S}$ sending $X\mapsto \underline{\ob{Mor}}_S(X,S)=X$. 

 This technical observation is actually very useful in practice, and it explains why many non-commutative rings  whose categories of modules arise via quasi-coherent sheaves on (analytic) stacks have naturally an antiinvolution up to some twist (eg. distribution algebras of groups, enveloping algebras  or algebras of differential operators). 
\end{remark}

\begin{notation}\label{NotationRealizationFunctors}
Let $(\s{C},E)$ be a geometric set up with fiber products, $\ob{D}$ a $3$-functor formalism on $(\s{C},E)$ and $S\in \s{C}$. We have the following diagram  of symmetric monoidal categories 
\[
\begin{tikzcd}
\s{C}^{\op,\sqcup} \ar[rd] & & \\
& \Corr(\s{C}_{E,/S})  \ar[r,"\Phi"] & \ob{K}_{\ob{D},S} \\ 
\s{C}^{\times} \ar[ur] &  &
\end{tikzcd}
\]
For notational reasons, we shall denote the upper composition by $[-]^*\colon \s{C}^{\op}\to \ob{K}_{\ob{D},S}$ and the lower composition by $[-]_!\colon \s{C}\to \ob{K}_{\ob{D},S}$. By the natural antiinvolution of the correspondence category, the maps $[-]^*$ and $[-]_!$ are duals to each other in $\ob{K}_{\ob{D},S}$. Given $X\in \s{C}_{E,/S}$ we shall write $[X]^*$ and $[X]_!$ for their image in $\ob{K}_{\ob{D},S}$ through $[-]^*$ and $[-]_!$ respectively, in case we want to ignore the difference at the level of objects we simply write $X=[X]^*=[X]_!$. For $f\colon Y\to X$ a map in $\s{C}_{E,/S}$, we shall write $f_!\colon Y\to X$ and $f^*\colon X\to Y$ the image of the map $f$ through $[-]_!$ and $[-]^*$ respectively. 

Finally, if $Y,X\to S$ are $!$-able maps and $M\in \ob{D}(Y\times_S X)=\ob{Fun}_S(X,Y)$ we denote by ${{}_Y M_X}\colon X\to Y$ the associated map in $\ob{K}_{\ob{D},S}$. Under this notation,  $M$ is self dual in the sense that ${{}_Y M_X}\colon X\to Y$ is dual to ${{}_X M_Y}\colon Y\to X$ where now we consider $M\in \ob{D}(X\times_S Y)=\ob{Fun}_S(Y,X)$. 
\end{notation}

\subsubsection{Suave and Prim objects}

One of the main reasons for introducing the $2$-category of kernels is that  the universality in statements such as Poincar\'e and Serre duality can be easily encoded in terms of adjunctions. The main definition around this idea is the following:

\begin{definition}\label{DefPrimSuaveObjects}
Let $(\n{C},E)$ be a geometric set up and $\ob{D}$ a $3$-functor formalism on $(\n{C},E)$. Let $S\in \n{C}$ and let $\ob{K}_{\ob{D},S}$ be the kernel category over $S$. Let $f:X\to S$ be a $!$-able map and $M\in \ob{D}(X)$. 
\begin{enumerate}[label = (\roman*)]
\item The object $M$ is called $f$-\textit{suave} (or just suave if $f$ is clear from the context) if, when seen as a map $M: X\to S$ in $\ob{K}_{\ob{D},S}$, it is a left adjoint.  We let $\ob{SD}_f(M)$ denote its right adjoint and call it the \textit{suave dual of $M$}. 

\item  The object $M$ is called $f$-\textit{prim} (or just prim if $f$ is clear from the context) if, when seen as a map $M:X\to S$ in $\ob{K}_{\ob{D},S}$,  it is a right adjoint. We let $\ob{PD}_f(1)$ denote its left adjoint and call it the \textit{prim dual of $M$} 

\end{enumerate}
\end{definition}

The following proposition explains the relation between suave and prim objects with the six functors (see \cite[Sections 4.4. and 4.5]{HeyerMannSix}): 

\begin{proposition}
Keep the notation of  \Cref{DefPrimSuaveObjects}.  Let $f:X\to S$ be a $!$-able map and $M\in \ob{D}(X)$

\begin{enumerate}

\item  Suppose that $M$ is suave and let $\ob{SD}_f(M)$ be its  right adjoint. Then there are natural equivalences of functors $\ob{D}(S)\to \ob{D}(X)$.
\[
\ob{SD}_f(M) \otimes f^* =  \iHom_X(M, f^!) 
\]
Furthermore, these equivalences hold after arbitrary base change $S'\to S$, the formation of $\ob{SD}_f(M)$ satisfies base change,  and $\ob{SD}_f(M)$ is also suave with suave dual given by $M$. In particular, $\ob{SD}_f(M)=f^! 1_S$.

\item Suppose that $M$ is prim and let $\ob{PD}_{f}(M)$ be its left adjoint. Then there are natural equivalences of functors $\ob{D}(X)\to \ob{D}(S)$
\[
f_!(\ob{PD}_f(M)\otimes -) = f_*(\iHom_X(M, -)) 
\]
Furthermore, these equivalences hold after arbitrary base change $S'\to S$, the formation of $\ob{PD}_f(M)$ satisfies base change,  and $\ob{PD}_f(M)$ is also prim with prim dual given by  $M$. 
\end{enumerate}
\end{proposition}

A special case of suave and prim objects is when $M=1$, it actually deserves its own definition:

\begin{definition}\label{defSuavePrimMap}
Let $\ob{D}$ be a $3$-functor formalism on  a geometric set up $(\n{C},E)$. A $!$-able map $f:X\to S$ is called \textit{suave} if $1$ is suave, we let $\omega_f=\ob{SD}_f(1)$ denote the suave dual of $1$ and call it the \textit{dualizing  sheaf of $f$} . Similarly, $f$ is called \textit{prim} if $1$ is prim, we let $\delta_f=\ob{PD}_f(1)$ be the prim dual of $1$ and call it the \textit{codualizing sheaf of $f$}. 
\end{definition}

Suave and prim objects satisfy many stability properties that can be checked formally from adjunctions (eg. suave objects are stable under prim lower $!$ and suave upper $*$ functors). They also satisfy suitable descent properties: they satisfy universal $*$-descent on targets \cite[Lemma 4.5.6]{HeyerMannSix}, and suitable suave or prim descent on the source \cite[Lemma 4.7.4]{HeyerMannSix}. These general descent properties are very important   since they will allow us to localize the proof of Serre or Poincar\'e duality in different contexts to a much simpler computation, they also help to construct analytic stacks via a very general Grothendieck topology, and to extend the class of $!$-able arrows in $6$-functor formalisms.  

Another very important feature of suave and prim objects is that they are  useful for proving universal $*$ or $!$-descent (cf. \cite[Lemmas 4.7.1 and 4.7.4]{HeyerMannSix}): 

\begin{proposition}
Let $f:Y\to X$ be a $!$-able morphism on a $6$-functor formalism admitting countable limits and colimits.  Let $Y^{\bullet}$ be the \v{C}ech nerve of $f$ and write $f_n: Y^n\to X$.
\begin{enumerate}[label = (\roman*)]

\item Suppose that $f$ is suave, then $f^!$ is conservative if and only if the natural map 
\[
\varinjlim_{[n]\in \Delta^{op}} f_{n,!}f_n^! 1 \to 1
\]
is an isomorphism. If this holds, then so does after any base change $X'\to X$ and $f$ is of universal $*$ and $!$-descent.  We say that $f$ satisfies \textit{suave descent.}

\item  Suppose that $f$ is prim and that the $6$-functor formalism is stable. If $f_* 1_Y\in \ob{D}(X)$ is a descendable algebra then the same holds after any base change $X'\to X$ and $f$ satisfies universal $*$ and $!$-descent.  We say that $f$ satisfies \textit{prim or descendable descent.}

\end{enumerate}
\end{proposition}

Finally, let us state the main theorem that allow us to extend a $6$-functor formalism to stacks. We first need a definition, cf. \cite[Definition 3.4.1]{HeyerMannSix}

\begin{definition}\label{DefSheafy}
Let $(\n{C},E)$ be a geometric set up such that $\n{C}$ is a (essentially small) site. A $6$-functor formalism on $(\n{C},E)$ is said \textit{sheafy} if the map 
\[
\ob{D}^*: \n{C}^{\op}\to \ob{Cat}_{\infty}
\]
satisfies descent. 
\end{definition}

\begin{example}
\begin{enumerate}

\item  A first example is the trivial site $\n{C}$ where covers are generated by isomorphisms. In this case $\ob{Shv}(\n{C})$ is just the category of presheaves on $\n{C}$.

\item A very useful example for the theory of analytic stacks is the $\ob{D}$-topology where covers are generated by canonical covers which are of universal $\ob{D}^*$ and $\ob{D}^!$-descent.  Note that by \Cref{RemarkUniversalShierkDescent} the $\ob{D}$-topology is exactly the same as the \textit{$!$-topology} whose covers are subcanonical  satisfying universal $!$-descent. 

\item The most extreme example of a sheafy $6$-functor formalism is the $\ob{D}^*$-topology where covers are generated by canonical covers which are of universal $\ob{D}^*$-descent. 

\end{enumerate} 
\end{example}

\begin{theorem}[{\cite[Theorem 3.4.11]{HeyerMannSix}}]\label{TheoExtension6}
Let $(\n{C},E)$ be a geometric set up such that $\n{C}$ is a site. Let $\ob{D}$ be a  sheafy $6$-functor formalism on $(\n{C},E)$. Then there is a collection of edges $E'$ in $\n{X}:=\ob{Shv}(\n{C})$ with the following properties:
\begin{enumerate}[label=(\roman*)]
\item The inclusion $\n{C}\hookrightarrow \n{X}$ defines a morphism of geometric set ups $(\n{C},E)\to (\n{X},E')$ and $\ob{D}$ extends uniquely to a $6$-functor formalism on $(\n{X},E')$.

\item $E'$ is $*$-local on the target: Let $f:Y\to X$ be a map in $\n{X}$ whose pullback to every object $X'\to X$ in $\n{C}$ lies in $E'$. Then $f$ lies in $E'$.

\item $E'$ is $!$-local on the target and source: let $f:Y\to X$ be a map in $\n{X}$ that is $!$-locally on the source or target on $E$, then $f$ lies in $E'$.

\item $E'$ is tame: every map $f:Y\to X$ in $E'$ with $X\in \n{C}$ is $!$-locally on the source in $E$. 
\end{enumerate}  

Moreover, there is a minimal choice of $E'$. The same results holds for hypercomplete sheaves in place of sheaves if one assumes that $\n{C}$ is hypersubcanonical\footnote{One can also localize in between sheaves and hypersheaves as long as the representable presheaves  in  $\n{C}$ and the functor $\ob{D}^*$ satisfy the corresponding descent}. 
\end{theorem}

\begin{remark}
By \cite[Lemma 3.4.13]{HeyerMannSix} if the maps  $X\hookrightarrow X\bigsqcup Y$ on $\n{C}$ are $!$-able then the class of $!$-able arrows in  \Cref{TheoExtension6} is also stable under disjoint unions. 
\end{remark}

As stated in  \Cref{TheoExtension6}, the local condition in the target requires to test against \textit{all} the objects in $\n{C}$. Under the $\ob{D}$-topology one only needs to test this against an epimorphism: 

\begin{lemma}
Let $(\n{C},E)$ be a geometric set up and $\ob{D}$ a $6$-functor formalism on $(\n{C},E)$.   Endow $\n{C}$ with the $\ob{D}$-topology and consider the extension of $\ob{D}$ to sheaves $\ob{Shv}_{\ob{D}}(\n{C})$. Let $f:Y\to X$ be a map such that there is an epimorphism $X' \to X$  in $\ob{Shv}_{\ob{D}}(\n{C})$ such that the pullback $f':Y'\to X'$ is $!$-able. Then $f$ is $!$-able. 
\end{lemma}
\begin{proof}
Let $Z\in \n{C}$, we want to show that the pullback $f_Z: Y\times_X Z \to Z$ is $!$-able and apply  \Cref{TheoExtension6} (ii). We know that the map $X'\to X$ is an epimorphism, so there is a $\ob{D}$-cover $Z'\to Z$ in $\n{C}$ and a lift $Z'\to X'$ over $X$. Since $Z'\to Z$ is of universal $!$-descent, by  \Cref{TheoExtension6} (iii) it suffices to show that the pullback $Y\times_{X}Z'\to Z'$ is $!$-able. But this map can also be written as $Y'\times_{X'} Z'\to Z'$ which is $!$-able as the class $E'$ is stable under pullbacks. 
\end{proof}

\subsubsection{\'Etale and proper maps} \label{sss:EtaleProper}

Among \'etale and prim maps there is a suitable class of maps that behave  cohomologically like \'etale or proper maps. We let $(\n{C},E)$ be a geometric set up and let $\ob{D}$ be a six functor formalism on $(\n{C},E)$.

\begin{definition}[{\cite[Definition 4.6.1]{HeyerMannSix}}]\label{DefEtaleProper}
Let $f\colon Y\to X$ be a truncated map in $E$.

\begin{enumerate}

\item We say that $f$ is \textit{\'etale} if $\Delta_f$ is \'etale or an isomorphism and if $f$ is suave.

\item We say that $f$ is \textit{proper} if $\Delta_f$ is proper or an isomorphism and if $f$ is prim. 

\end{enumerate}
\end{definition}

\begin{remark}
Since the map $f$ in \Cref{DefEtaleProper} is truncated, one iterated diagonal will become an isomorphism. Thus, we can say thata truncated map is \'etale (resp. proper) if and only if any iterated diagonal is suave (resp. prim). 

In \cite[Lemma 4.6.3]{HeyerMannSix} it is explain some stability properties of \'etale and proper maps: they are stable under base change, composition and satisfy a suitable right cancellation property (part (i) of \textit{loc. cit.}), they are $\ob{D}^*$-local on the target (part (ii) of \textit{loc. cit.}), and they are $\ob{D}^!$-local on the source by \'etale (resp. proper) maps (part (iii) of \textit{loc. cit.}).  

Finally, \cite[Lemma 4.6.4]{HeyerMannSix} shows that if $f\colon Y\to X$ is an \'etale (resp. proper) map, then there is a natural equivalence of functors $f^!\xrightarrow{\sim} f^*$ (resp. $f_!\xrightarrow{\sim} f_*$) as expected. 
\end{remark}

\subsubsection{Open and closed maps}\label{sss:ClosedOpenMaps}

One can extend the notion of open and closed maps of symmetric monoidal stable categories of \Cref{sss:OpenClosedSymmetricMonoidal} to arbitrary six functor formalisms (see also \cite[Definitionn 3.11]{kesting2025categoricalkunnethformulasanalytic}). Let $(\n{C},E)$ be a geometric set up and $\ob{D}$ a presentable six functor formalism on $(\n{C},E)$. 

\begin{definition}\label{DefOpenClosed}

\begin{enumerate}

\item  A morphism $f\colon Y\to X$ in $E$ is called an \textit{open map} if it is  \'etale and an immersion.

\item A morphism $f\colon Y\to X$ in $E$ is called a \textit{closed map} if it is proper and an immersion. 

\end{enumerate}
\end{definition}

Open and closed maps in a six functor formalism satisfy the expected compatibility properties (see \cite[Proposition 3.12]{kesting2025categoricalkunnethformulasanalytic}).

\begin{proposition}\label{PropStabilityOpenClosedMaps}
The following hold:
\begin{enumerate}

\item Let $Z\xrightarrow{f} Y\xrightarrow{g} X$ be $!$-able maps in $E$ such that $g$ is an open (resp. a closed) map. Then $f$ is an open (resp. a closed) map if and only if $g\circ f$ is an open (resp. a  closed) map. 

\item Let $f\colon Y\to X$ be an open map (resp. a closed map) in $E$, then $f^*\colon \ob{D}(X)\to \ob{D}(Y)$ is an open (resp. closed) map of symmetric monoidal categories.

\item Let $f\colon Y\to X$ be an open (resp. closed) map, let $g\colon X'\to X$ be a map in $\n{C}$ with base change $f'\colon Y'\to X'$. Then $f'$ is an open (resp. closed) map and  the natural map 
\begin{equation}\label{eqKunnethOpenClosedMap}
\ob{D}(X')\otimes_{\ob{D}(X)}\ob{D}(Y)\to \ob{D}(Y')
\end{equation}
is an equivalence in $\cat{Pr}^L$.

\end{enumerate}
\end{proposition}
\begin{proof}

\begin{enumerate}
\item \'Etale and proper maps are stable under base change and satisfy the right cancellation property thanks to \cite[Lemma 4.6.3 (i)]{HeyerMannSix}.  The same holds for immersions.

\item We want to see that if $f$ is \'etale and an immersion then $f^*\colon \ob{D}(X)\to \ob{D}(Y)$ is an open map of symmetric monoidal categories.  Since $f$ is \'etale, $f^*$ has a left adjoint $f_{\natural}$ satisfying base change.  Now, because  $f$ is an immersion, we have that $Y=Y\times_X Y$, and by proper base change we have that  $f^*f_{\natural}=\id$ as endofunctors of  $\ob{D}(Y)$. This shows that $f_{\natural}\colon \ob{D}(Y)\to \ob{D}(X)$ is fully faithful, proving that $f^*$ is an open immersion. A similar argument holds for proper maps, where we use $f_*$ instead of $f_{\natural}$.

\item  Being an immersion and \'etale (resp. proper) is stable under base change. Hence open and closed maps are stable under base change. For the K\"unneth formula of \Cref{eqKunnethOpenClosedMap}, as $f^*$ is an open map of symmetric monoidal stable categories one has that 
\[
\ob{D}(Y)=\ob{coMod}_{f_{\natural}1} \ob{D}(X)
\]
and similarly for $\ob{D}(\widetilde{Y})$. Thus, since $\widetilde{f}_{\natural}1 = g^* f_{\natural} 1$ by base change, one obtains \eqref{eqKunnethOpenClosedMap} by taking Lurie's tensor product of the localization sequence of the coidempotent coalgebra $f_{\natural} 1\in \ob{D}(X)$.

\end{enumerate}

\end{proof}

\subsection{Analytic stacks and quasi-coherent  $6$-functor formalism}\label{ss:AnStacksConstruction}

With the previous preparations on abstract $6$-functor formalisms we can finally introduce the category of analytic stacks.

\subsubsection{Quasi-coherent six functors on affinoids} Let $\ob{AnStk}^{\aff}=\AnRing^{\op}$ be the opposite category of analytic rings that we call \textit{affinoid analytic stacks}. Given $A$ an analytic ring we let  $\AnSpec A\in \ob{AnStk}^{\aff}$ be the object associated to $A$,  called  the \textit{analytic spectrum of $A$}.

To define analytic stacks, we need a suitable Grothendieck topology on analytic rings. First, we need  a $6$-functor formalism for analytic rings. For that, we  apply \cite[Proposition 3.3.3]{HeyerMannSix}, and define classes of arrows $(E,I,P)$ in $\ob{AnStk}^{\aff}$ consisting of $!$-able,  \textit{open imersions}  and \textit{proper maps} of affinoid analytic stacks  respectively.

\begin{remark}\label{RemProperEtaleMaps}
In the generality of a suitable decomposition $(I,P)$ of an abstract geometric set up $(\n{C},E)$, the classes of maps $I$ and $P$ might not be proper and \'etale in the sense of \cite[Definition 4.6.1]{HeyerMannSix} for a given $6$-functor formalism $\ob{D}$.  The only problem being that these maps are not necesarily truncated (only $I\cap P$ is truncated). In practice the class $I$ or $P$ will consist on truncated maps (eg. for analytic rings $I$ will be $(-1)$-truncated, i.e. immersions).

 In the language of Gestalten of Scholze and Stefanich  \cite{ScholzeGestalted}, the arrows in $I$ and $P$ on $\mathcal{C}$ will give rise to $0$-\'etale and $0$-proper morphisms on the Gestalten associated to  the six functor formalism $\ob{D}$ on $(\n{C},E)$. 
\end{remark}

\begin{definition}\label{DefShierkMapsAnring}
Let $f\colon \AnSpec B\to \AnSpec A$ be a map of affinoid analytic stacks. 

\begin{enumerate}

\item  We say that $f$ is in the class $P$ if $B=B_{A/}$ has the induced analytic ring structure.

\item  We say that $f$ is in the class $I$ if the pullback map $f^*\colon \ob{D}(A)\to \ob{D}(B)$ is an open imersion in the smashing spetrum.

\item We say that $f$ is in the class $E$ (or that $f$ is $!$-able), if the map $B_{A/}\to B$ gives rise to an open immersion $\ob{D}(B_{A/})\to \ob{D}(B)$.

\end{enumerate}

\end{definition}

The following lemma proves the required stability properties for the classes of maps $(E,I,P)$:

\begin{lemma}\label{LemStabilityShierkMaps}

Let $\n{C}=\ob{AnStk}^{\aff}$. The following holds: 

\begin{enumerate}

\item The classes of maps $I$, $P$ and $E$ give rise to geometric set ups in $\n{C}$.

\item Any map $f\in E$ is written as the composite $p\circ i$ where $i\in I$ and $p\in P$

\item Every morphism of $I$ is $-1$-truncated, i.e. it is a monomorphism. 

\item Let $j\in I$, then the pullback map $j^*$ has a fully faithful left adjoint $j_{\natural}$ satisfying projection formula and proper base change. 

\item Let $p\in P$, then the right  adjoint $p_*$ satisfies projection formula and proper base change.

\item For every cartesian diagram 
\[
\begin{tikzcd}
U' \ar[r, "j'"] \ar[d,"p'"] & X' \ar[d,"p"] \\ 
U \ar[r, "j"] & X
\end{tikzcd}
\]
in $\n{C}$  with $j\in I$ and $p\in P$, the natural map $j_{\natural} p'_*\xrightarrow{\sim} p_* j'_{\natural}$ is an equivalence. 

\end{enumerate}
\end{lemma}
\begin{proof}
(1). We need to show the stability properties for the arrows $I$, $P$ and $E$ in \Cref{RemEquivDefGeoStup}. By an abuse of notation we will call in the same way the classes of arrows in the  category of analytic rings.

 It is clear that all three classes contain the isomorphisms. The classes $I$ and $P$ are also clearly stable under composition, namely, for $I$ it follows from the fact that the composite of two open immersions of symmetric monoidal categories is an open immersion (consequence of \Cref{theo:openCoversSym}), and for the classes in $P$ it follows from the definition of induced analytic ring structure.  They are also stable under base change of arbitrary arrows in $\n{C}$ thanks to \Cref{PropColimitsAnRingAndMod} (for open immersions) and \Cref{LemmaInducedStructures} (for objects in $P$).

Let us prove stability under composition for the class in $E$. Let $A\to B$ and $B\to C$ be two maps in $E$, then $B_{A/}\to B$ and $C_{/B}\to C$ give rise to open immersions at the level of symmetric monoidal categories. We want to prove that the map $C_{A/}\to C$ gives rise to an open immersion, but it can be written as the composite 
\[
C_{A/}\to C_{/B}  \to C,
\]
 the first map is the base change of $B_{A/}\to  B$ along the map $B_{A/}\to C_{A/}$,  so it is an open immersion. It follows from the previous that $C_{A/}\to C$ is an open immersion, proving what we wanted.

Finally, we are left to show that the classes in $I$, $P$ and $E$ are stable under taking diagonals. We start with $I$. Let $A\to B$ be a map in $I$, we want to show that the multiplication map $m\colon B\otimes_{A} B\to B$ is in $I$. But the map $A\to B$ is idempotent as analytic rings, as $B$ is a uncompleted analytic ring structure on $A^{\triangleright}$. Hence, $m$ is an equivalence and in particular in $I$. Next, for $P$ and $E$, it suffices to see that given $A\to B$ any map of analytic rings, the morphism $B\otimes_A B\to B$ always has the induced structure. This follows from the fact that a $B^{\triangleright}\otimes_A B^{\triangleright}$-module is $B\otimes_A B$-complete if and only if it is $B$-complete when restricted to each term appearing in the tensor, and the composite (for any of the two maps out from $B$) $B\to B\otimes_A B\to B$ is the identity.

(2). This is obvious from the definition of $E$.

(3). This is obvious from the fact that a map $A\to B$ is idempotent as analytic rings. 

(4). This follows from \Cref{PropColimitsAnRingAndMod} and the projection formula and open base change for open immersions in symmetric monoidal categories. 

(5). This follows from \Cref{PropColimitsAnRingAndMod} and the fact that for a symmetric monoidal category $\s{E}$, and a commutative algebra $A\in \ob{CAlg}(\s{E})$, the forgetful map $\ob{Mod}_A(\s{E})\to \s{E}$ is $\s{E}$-linear (i.e. satisfies projection formula), and base change along any $\s{E}$-linear category (\cite[Theorem 4.8.4.6]{HigherAlgebra}). 

(6). Thanks to \Cref{PropColimitsAnRingAndMod}, it suffices to  prove  the following fact: let $\s{E}$ be a stable symmetric monoidal category, $A\in \ob{CAlg}(\s{E})$ a commutative  algebra in $\s{E}$ and $D\in \ob{Idem}(\s{E})$ an idempotent algebra. Let $\s{E}(U)=\s{E}/\ob{Mod}_D(\s{E})$ be its open localization. Consider the following commutative diagram of symmetric monoidal categories
\[
\begin{tikzcd}
\s{E} \ar[r,"j^*"] \ar[d,"f^*"] & \s{E}(U) \ar[d,"f^{\prime,*}"] \\ 
\Mod_A(\s{E}) \ar[r,"j^{\prime,*}"] & \Mod_A(\s{E}(U)).
\end{tikzcd}
\]
Let $f_*$ and $f^{\prime}_*$ be the right adjoints of $f^*$ and $f^{\prime,*}$ respectivelty, similarly let $j_{\natural}$ and $j^{\prime}_{\natural}$  be the left adjoints of $j^*$ and $j^{\prime,*}$ respectively. Then the natural morphism 
\[
j_{\natural} f'_*\xrightarrow{\sim} f_* j'_{\natural}
\] 
is an equivalence. This follows formally from the fact that the morphism $j_{\natural}$ is given by $j_{\natural}j^* M = \ob{fib}(1\to D)\otimes M$ for $M\in \s{E}$, and that the morphism $j^{\prime}_{\natural}$ is given by $j_{\natural}^{\prime} j^{\prime,*} N = (\ob{fib}(A\to A\otimes D))\otimes_A N$ for $N\in \ob{Mod}_A(\s{E})$. Indeed, the associativity of the relative tensor product produces the natural equivalence 
\[
\ob{fib}(1\to D)\otimes_A N\xrightarrow{\sim} \ob{fib}(A\to A\otimes D)\otimes_A N 
\]
for $N\in \ob{Mod}_A(\s{E})$.
\end{proof}

\begin{corollary}\label{CorSixFunctorFormalismAffinoids}
Let $\n{C}=\ob{AnStk}^{\aff}$ be the category of affinoid analytic stacks, and $E$ the class of $!$-able arrows in $\n{C}$ of \Cref{DefShierkMapsAnring}. Then there exists a (unique if done in a $2$-categorical sense) promotion of $\ob{D}\colon \ob{AnStk}^{\aff,\op}\to \Pr^{L}$ to a presentable $6$-functor formalism  
\[
\ob{D}\colon \ob{Corr}(\n{C},E)^{\otimes}\to {\Pr}^{L,\otimes}
\]
where for $j\in I$ one has $j_!=j_{\natural}$ and for $p\in P$ one has $p_!=p_*$. 
\end{corollary}
\begin{proof}
The construction of the $6$-functor formalism follows from \cite[Proposition 3.3.3]{HeyerMannSix} where all the necessary  conditions are verified in \Cref{LemStabilityShierkMaps}, the uniqueness is \cite[Theorem A]{cnossen2025universalityspan2categoriesconstruction} after passing to a $2$-categorical promotion of the category of correspondences. 
\end{proof}

The following lemma identifies the suave maps of affinoid analytic stacks with induced analytic ring structure.

\begin{lemma}\label{LemmaPrimAffine}
Let $f\colon A\to B$ be a morphism of  analytic rings  with the trivial analytic ring structure. Then $f$ is suave if and only if $B$ is dualizable as condensed $A$-module.
\end{lemma}
\begin{proof}
The map $f$ is prim having the induced analytic ring structure. The lemma is a particular case of  \cite[Lemma 4.5.10]{HeyerMannSix}. 
\end{proof}

Thanks to the colimit preserving property of \Cref{PropColimitsAnRingAndMod} we deduce the following fully faithfulness for the category of kernels:

\begin{lemma}\label{LemFullyFaithfulKernelCat}
Let $(\n{C},E)$ be a geometric set up and let $\ob{D}\colon \Corr(\n{C},E)\to \ob{Pr}^L$ be a presentable six functor formalism. Suppose that $\ob{D}\colon \n{C}^{\op}\to \ob{CAlg}(\ob{Pr}^L)$ satisfies categorical K\"unneth formula, i.e, that it commutes with pushouts. Let $S\in \n{C}$, then the natural lax symmetric monoidal functor  
\[
\Psi\colon \ob{K}_{\ob{D},S}\to \ob{Pr}^L_{\ob{D}(S)}
\] 
is symmetric  monoidal and $2$-fully faithful (see  \Cref{RemarkPromotionSymmetricMonoidal2,RemarkPromotionSymmetricMonoidal} for the discussion about symmetric monoidal structures). In particular, this holds for $\n{C}=\AnStk^{\aff}$.
\end{lemma} 
\begin{proof}
Let $X,Y\to S$ be $!$-able objects over $S$. By the categorical K\"unneth formula the functor $\ob{D}\colon \ob{Corr}((\n{C}_{E})_{/S})\to \ob{Pr}^L_{\ob{D}(S)}$ is symmetric monoidal, meaning that $\ob{D}(Y\times_S X) = \ob{D}(Y)\otimes_{\ob{D}(S)} \ob{D}(X)$. On the other hand, the objects of $X\in \ob{Corr}((\n{C}_{E})_{/S})$ are naturally self dual by \cite[Proposition 2.3.9]{HeyerMannSix}, which makes the categories $\ob{D}(Y)$ and $\ob{D}(X)$ naturally self dual in $ \ob{Pr}^L_{\ob{D}(S)}$. This duality excanges the lower $!$ and upper $*$-functor thanks to \cite[Proposition 2.3.9]{HeyerMannSix} and the definition of these functors as those associated to the diagrams $X=X\to S$ and $S \leftarrow X=X$ respectively. Since the functor $\Phi\colon \ob{Corr}((\n{C}_{E})_{/S})\to \ob{K}_{\ob{D},S}$ is symmetric monoidal, essentially surjective and the composite with the lax symmetric monoidal functor $\Psi\colon \ob{K}_{\ob{D},S}\to  \ob{Pr}^L_{\ob{D}(S)}$ is also symmetric monoidal, one deduces that $\Psi$ is symmetric monoidal as well.

On the other hand, by definition, the morphism category from $X$ to $Y$ in $\ob{K}_{\ob{D},S}$ is given by 
\[
\ob{Fun}_{S}(X,Y)= \ob{D}(Y\times_S X).
\]
Then, the morphism of $2$-categories $\Psi$ gives rise to the natural morphism of categories 
\begin{equation}\label{eqFullyFaithKernel}
\ob{D}(Y\times_S X) \to \ob{Fun}_{\ob{D}(S)}(\ob{D}(X), \ob{D}(Y)).
\end{equation}
The equivalence of categories in \eqref{eqFullyFaithKernel} then follows from the natural self duality of $X$ and $\ob{D}(X)$ in $\ob{K}_{\ob{D},S}$ and $\ob{Pr}^L_{\ob{D}(S)}$ respectively, and the categorical K\"unneth formula.
\end{proof}

\subsubsection{$!$-descent}\label{sss:ShierkDescent}

Having constructed six functors for affinoid analytic stacks, we shall study its $\ob{D}$-topology which will be the chosen one for the theory of analytic stacks.  Thanks to \Cref{LemFullyFaithfulKernelCat} one can give a simpler characterization of $\ob{D}$-covers.


Let $(\n{C},E)$ be a  geometric set up with $\n{C}$ a category with finite limits. Let $\ob{D}$ be a presentable  $6$-functor formalism on $(\n{C},E)$. The following proposition is a particular case of \cite[Theorem 5.12]{SixFunctorsScholze} where we assume categorical K\"unneth formula.

\begin{proposition}\label{PropDescentShierk}
Let $(\n{C},E)$ and $\ob{D}$ be as before, suppose that $\ob{D}$ satisfies categorical K\"unneth formula, let $\{Y_i\to X \}_i$ be an arbitrary  family of  $!$-able morphism on $\n{C}$, and let  $Y_{\bullet}\to X$ be the \v{C}ech nerve of $\bigsqcup_i Y_i\to X$ in $\n{P}(\n{C})$.  Then the following are equivalent: 

\begin{enumerate}[label = (\alph*)]

\item The family $\{Y_i \to X\}$ satisfies universal $\ob{D}^*$ and $\ob{D}^!$-descent.

\item The natural map $\ob{D}(X)\to \ob{Tot}(\ob{D}^!(Y_{\bullet}))$ is an equivalence of categories, where the transition maps are given by upper $!$-maps.

\item The natural map $\varinjlim_{[n]\in \Delta^{\op}}\ob{D}_!(Y_{\bullet})\xrightarrow{\sim}  \ob{D}(X)$ is an equivalence in $\ob{Pr}^L_{\ob{D}(X)}$, where the transitions maps are given by lower $!$-maps. 

\item The natural map of $2$-categories $\ob{Pr}^L_{\ob{D}(X)}\to \ob{Tot}(\ob{Pr}^L_{\ob{D}(Y_{\bullet})})$ is an equivalence, here the transition maps are given by base change along the augmented cosimplicial algebra $\ob{D}(X)\to \ob{D}(Y_{\bullet})$ in $\ob{Pr}^L$.

\item The natural map of $2$-categories $\ob{Pr}^L_{\ob{D}(X)}\to \ob{Tot}(\ob{Pr}^L_{\ob{D}(Y_{\bullet})})$ is fully faithful.

\end{enumerate}
In particular, the proposition  holds for $\n{C}=\AnStk^{\aff}$.
\end{proposition}

\begin{remark}\label{RemPowerfullLemmaDescent}
\Cref{PropDescentShierk} is an important instance where the abstract theory of $6$-functor formalism is used in the construction of analytic stacks.  One could think of this lemma as a generalization of the notion of descendable map of algebras of Mathew \cite{MathewDescent}.
\end{remark}

\begin{proof}
For simplicity in the notation we assume that $\{Y_i\to X\}$ consists in a single map $Y\to X$ (the argument in the proof being the same in general). In order to prove \Cref{PropDescentShierk} we will need a key fact about adjunctions of linear categories arising from a $6$-functor formalism satisfying K\"unneth formula. We can assume without loss of generality that $C=(C_{E})_{/S}$. Since the six functor formalism $\ob{D}\colon \Corr(\n{C})\to \ob{Pr}^L_{\ob{D}(S)}$ is symmetric monoidal,  for $f\colon X\to S$ a $!$-able map, the category $\ob{D}(X)\in \Pr^L_{\ob{D}(S)}$ is naturally self dual.  We have a pullback functor $f^{*}\colon \ob{Pr}^L_{\ob{D}(S)} \to \ob{Pr}^L_{\ob{D}(X)}$ which is nothing but the base change along $\ob{D}(S)\to \ob{D}(X)$ with respect to Lurie's tensor product. This functor has a right adjoint $f_*$ which is just the forgetful functor $f_*\colon \ob{Pr}^L_{\ob{D}(X)}=\Mod_{\ob{D}(X)}(\ob{Pr}^L_{\ob{D}(S)})\to \ob{Pr}^L_{\ob{D}(S)}$.  Since $\ob{D}(X)$ is dualizable, $f^*$ also has a left adjoint $f_!$ satisfying 
\[
\begin{aligned}
\Fun_{\ob{D}(S)}(f_!M,N) & =\Fun_{\ob{D}(X)}(M, \ob{D}(X)\otimes_{\ob{D}(S)} N)\\ 
		&  =\Fun_{\ob{D}(X)}(M, \Fun_{\ob{D}(S)}(\ob{D}(X),N)) \\ 
		& = \Fun_{\ob{D}(S)}(\ob{D}(X)\otimes_{\ob{D}(X)} M , N) \\ 
		& = \Fun_{\ob{D}(S)}(M,N),
\end{aligned}
\]
for $M\in \ob{Pr}^L_{\ob{D}(X)}$ and $N\in \ob{Pr}^L_{\ob{D}(S)}$,  where in the second equivalence we used the self dualizability of $\ob{D}(X)$. This shows that $f_!=f_*$, i.e. $f_*$ is both a left and right adjoint of $f^*$.

We now prove the proposition. 

  (a) implies (b) is clear by definition of $!$-descent. For (b) implies (c), note that colimits in $\ob{Pr}^L_{\ob{D}(X)}$ are computed, after passing to right adjoints, as limits in $\ob{Pr}^R_{\ob{D}(X)}$ which are just computed as limits of the uniderlying categories. For (c) implies (d), let $f_{\bullet}\colon Y_{\bullet}\to X$ be the \v{C}ech nerve of $Y\to X$. The pullback functor 
  \[
  \ob{Tot}(f^*_{\bullet}) \colon \ob{Pr}^L_{\ob{D}(X)} \to \ob{Tot}(\ob{Pr}^L_{\ob{D}(Y_{\bullet)}}) 
  \]
has both left and right adjoints as each individual functor $\ob{Pr}^L_{\ob{D}(X)}\to \ob{Pr}^L_{\ob{D}(Y_{\bullet})}$  does. The right adjoint sends a cocartesian section $M_{\bullet}$ to $\ob{Tot}(f_{\bullet,*}M_{\bullet})$, while the left adjoint sends $M_{\bullet}$ to $\varinjlim_{[n]\in \Delta^{\op}}  f_{\bullet, *} M_{\bullet}$. To show that  $\ob{Tot}(f^*_{\bullet})$ is an equivalence it suffices to prove that the unit and co-unit for some of these adjunctions are isomorphisms. 

For fully faithfulness, let $N\in \ob{Pr}^L_{\ob{D}(X)}$, we want to see that the natural map 
\[
N\to \ob{Tot}(f_{\bullet,*} f^{*}_{\bullet} N)
\]
is an equivalence.  The right limit can also be written as 
\[
\begin{aligned}
\ob{Tot}(f_{\bullet,*} f^{*}_{\bullet} N) & = \ob{Tot}(\ob{D}(X_{\bullet})\otimes_{\ob{D}(X)} N) \\
					 & = \ob{Tot}\bigg( \ob{Fun}_{\ob{D}(X)}(\ob{D}(Y_{\bullet}), N) \bigg))\\ 
					 & = \ob{Fun}_{\ob{D}(X)}(\varinjlim_{[n]\in \Delta^{\op}} \ob{D}_!(Y_n), N) \\ 
					 & = \ob{Fun}_{\ob{D}(X)}(\ob{D}(X),N) \\ 
					 & = M
 \end{aligned}
\]
where in the second equivalence we used the self dualizability of $\ob{D}(Y_{\bullet})$ as $\ob{D}(X)$-module,  in the third equivalence we used that the duals of the $*$-pullback maps are precisely the lower $!$ maps, and in the fourth condition (c). This proves that the functor $ \ob{Tot}(f^*_{\bullet})$ is fully faithful. Let us now see that it is an equivalence. For that, let $M_{\bullet}$ be a cocartesian section in  $\ob{Tot}(\ob{Pr}^L_{\ob{D}(Y_{\bullet})}) $, we want to show that the natural map 
\[
\ob{D}(Y_{n})\otimes_{\ob{D}(X)} \ob{Tot}(M_{\bullet})\to M_n
\]
is an equivalence for all $[n]\in \Delta$. To prove this, note that 
\[
\begin{aligned}
\ob{D}(Y_{n})\otimes_{\ob{D}(X)} \ob{Tot}(M_{\bullet}) & =  \ob{Tot}\bigg( \ob{D}(Y_n)\otimes_{\ob{D}(X)} M_{\bullet}  \bigg) \\ 
				&  =   \ob{Tot}( M_{\bullet +n +1}) \\ 
				& = M_n  
\end{aligned}
\]
where in the first equivalence we used that $\ob{D}(Y_n)$ is dualizable over $\ob{D}(X)$, and in the second that the diagram $M_{\bullet}$ is cocartesian and that $\ob{D}(Y_{\bullet})\otimes_{\ob{D}(X)}\ob{D}(Y_n)=\ob{D}(Y_{\bullet+n+1})$ by Kunneth formula, and in the last that $M_{\bullet +n+1}$ is split with limit $M_n$ for all $n\geq 0$. 

The implication $(d)\Rightarrow (e)$ is obvious.  Finally, we prove (e) implies (a). For that, notice that the unit  (arising from seeing $f_*$ as a right adjoint of $f^*$)
\[
\ob{D}(X)\to \ob{Tot}(\ob{D}(X_{\bullet}))
\]
being an equivalence is equivalent to $\ob{D}^*$-descent. Similarly, the co-unit (arising from seeing $f_*$ as left adjoint of $f^*$)
\[
\varinjlim_{[n]\in \Delta}\ob{D}(X_{\bullet}) \to \ob{D}(X)
\]
is equivalent to $\ob{D}^!$-descent (this map arises by taking the dual of $*$-pullbacks maps which are precisely the lower $!$ maps). This condition (d) implies $\ob{D}^*$ and $\ob{D}^!$-descent for $X_{\bullet}\to X$. For universal $*$ and $!$-descent, it suffices to note that if $\ob{D}(X)\to \s{E}$ is a symmetric monoidal map, then $\s{E}\to \ob{D}(X_{\bullet})\otimes_{\ob{D}(X)} \s{E}$ also satisfies descent at the level for linear categories, applying this to $\ob{D}(X)\to \ob{D}(Y)$ for $Y$ any pre-stack on $\n{C}$, and using categorical K\"unneth formula to see that $\ob{D}(X_{\bullet})\otimes_{\ob{D}(X)} \ob{D}(Y)= \ob{D}(X_{\bullet}\times_X Y)$, one deduces the universal $\ob{D}^*$ and $\ob{D}^!$-descent as wanted. 
\end{proof}

\begin{definition}\label{eqCoroDequalShierk}
Let $(\n{C},E)$ be a geometric set up, and $\ob{D}$ a six functor formalism on $(\n{C},E)$ as in \Cref{PropDescentShierk}. By \textit{loc. cit.}  a morphism $f\colon Y\to X$ in $\n{C}$ is a $\ob{D}$-cover if and only if it is a sub-canonical cover and satisfies $!$-descent. The \textit{$!$-topology} on $\n{C}$ is the Grothendieck topology generated by covers  given by \textit{arbitrary} families of $!$-able maps   $\{f_i\colon Y_i\to X\}_{i\in I}$  which are subcanonical covers and satisfy $!$-descent. 
\end{definition}

The following lemma shows that, under some quasi-compact conditions on the categories $\ob{D}(X)$ for $X\in \n{C}$, the $!$-topology is finitary.

\begin{lemma}\label{LemmaEpiInshierkTopology}
Let $(\n{C},E)$  and $\ob{D}$  be as before.  Suppose that for all $X\in \n{C}$ the unit $1_X\in \ob{D}(X)$ is compact.  Then if $\{f_i\colon Y_i\to X\}_I$ is a family of $!$-able morphisms satisfying $!$-descent, there is a finite subfamily $\{f_i\colon Y_i\to X \}_{i\in J}$ (with $I\subset J$ a finite subset) also satisfying $!$-descent. 
\end{lemma}
\begin{proof}
Let $\ob{D}_!\colon \n{C}_{E,/X}\to \Pr^L_{\ob{D}(X)}$ be the functor given by lower $!$-maps, we denote in the same way its extension fo presheaves. Let $\{Y^{\bullet}_{I}\}_{\Delta^{\op}_I}$ be the \v{C}ech nerve of $\{f_i\colon Y_i\to X\}$ (cf. \cite[Definition A.4.5]{HeyerMannSix}), we know that the natural map $f_{I,!}\colon \ob{D}_!(|Y^{\bullet}_I|)\xrightarrow{\sim} \ob{D}(X)$ induced by lower $!$-functors is an equivalence. For a finite subset $J\subset I$ let $Y^{\bullet}_J$ be the \v{C}ech nerve of $\{f_i\colon Y_i\to I\}_{i\in J}$. We want to show that there is some finite subset $J\subset I$ such that the natural map $f_{J,!}\colon \ob{D}_!(|Y_J^{\bullet}|)\to \ob{D}(X)$ is an equivalence. We can write $\ob{D}_I(|Y^{\bullet}_I|)=\varinjlim_{J\subset I} \ob{D}_!(|Y^{\bullet}_J|)$ in $\Pr^L_{\ob{D}(X)}$, thus we have that $1_{X}=\varinjlim_{J} f_{J,!}f_J^! 1_X$ where $f_{J}^!\colon \ob{D}(X)\to \ob{D}_!(|Y^{\bullet}_J|)$ is the right adjoint functor of $f_{J,!}$. Since $1_X$ is compact, there is some finite $J$ such that $1_X$ is a retract of $f_{J,!} f^{!}_{J} 1\in \ob{D}(X)$. Thus, the element $f^!_J 1_X\in \ob{D}_!(|Y^{\bullet}_J|)$ gives rise to a map $h\colon \ob{D}(X)\to \ob{D}_!(|Y^{\bullet}_J|)$ in $\Pr^L_{\ob{D}(X)}$ whose composite with the lower $!$-map $f_{J,!}\circ h$ is a linear endofunctor of $\ob{D}(X)$ that admits the identity as a retract. In other words, $\ob{D}(X)$ is a $2$-retract of the category $\ob{D}_!(|Y^{\bullet}_J|)$.  

Now, the rest of the proof is exactly (and actually a particular case) of that of \cite[Proposition 8.11]{ScholzeGestalted}. By categorical K\"unneth one has that $\ob{D}_!(|Y^{\bullet}_J|)$ is also the \v{C}ech co-nerve of the maps $\{\ob{D}(Y_i)\to \ob{D}(X)\}_{i\in J}$, hence $\ob{D}_!(|Y^{\bullet}_J|)\to \ob{D}(X)$ is a coidempotent coalgebra in $\Pr^L_{\ob{D}(X)}$. Since $\ob{D}(X)$ is a  $2$-retract of $\ob{D}_!(|Y^{\bullet}_J|)$, \cite[Lemma 8.13]{ScholzeGestalted} implies that $\ob{D}(X)$ is in the $\Pr^L$-linear category generated by $\ob{D}(|Y^{\bullet}_J|)$, taking into account the $\ob{D}(X)$-linearity and using that $\ob{D}(|Y^{\bullet}_{J}|)$ is coidempotent, this means that $\ob{D}(X)$ is also a comodule for $\ob{D}(|Y^{\bullet}_{J}|)$, i.e. that $\ob{D}(|Y^{\bullet}_{J}|)\otimes_{\ob{D}(X)} \ob{D}(X)=\ob{D}(X)$, and thus that $\ob{D}(|Y^{\bullet}_{J}|)=\ob{D}(X)$ as wanted. 
\end{proof}

\subsubsection{$!$-equivalences}\label{SectionShierkequivalences}

In \Cref{eqCoroDequalShierk} we saw how the $\ob{D}$ and $!$-topology are the same in a six functor formalism satisfying categorical Kunneth formula.  Thanks to \Cref{TheoExtension6} we can extend $\ob{D}$ to a $6$-functor formalism in the category of sheaves on $\n{C}$ for the $!$-topology.   However, the definition of analytic stacks in \cite{AnalyticStacks} involves an additional localization living in between $!$-sheaves and $!$-hypersheaves. This additional localization is actually subtle, in particular  the definition of \textit{loc. cit.} requires a correction that we will be explain in this paragraph (this correction was explained to the author by Scholze while writing the paper  \textit{de Rham stacks on Fargues-Fontaine curves} \cite{dRFF}).

Let $\n{C}$ be a (small) category with finite limits and terminal object $*$ and let $\n{P}(\n{C})=\ob{PShv}(\n{C},\Ani)$ be the  category of presheaves on anima on $\n{C}$. Let $(\n{C},E)$ be a geometric set up in $\n{C}$ and $\ob{D}\colon \Corr(\n{C},E)\to \Pr^L$ a presentable six functor formalism on $(\n{C},E)$ satisfying categorical K\"unneth formula. Suppose that $\n{C}$ admits finite disjoint unions and that $*\to *\sqcup *$ is $!$-able.  Let $\n{C}_E$ be the non-full wide subcategory of $\n{C}$ spanned by the arrows in $E$.  Consider the functor $\ob{D}^!\colon \n{C}_E^{\op}\to \Pr^R$ and its right Kan extension to $\n{P}(\n{C}_E)$. We also write by $\ob{D}_!\colon \n{P}(\n{C}_E)\to \Pr^L$  for the left Kan extension of the functor on $\n{C}_E$ given by lower $!$-maps (this is nothing but the opposite functor of $\ob{D}^!$ under the equivalence $\Pr^{L,\op}=\Pr^R$).   We make the following assumption on $(\n{C},E)$ that will hold in the case of analytic rings: 

\begin{hypothesis}\label{HypothesisTechnicalShierk}
Let $I$ be a category and let $I^{\triangleright}$ be its right cone with final object $\emptyset$. Let  $X_{\bullet}\colon I^{\triangleright}\to \n{C}$ a diagram in $\n{C}$ such that the natural map 
\[
\ob{D}(X_{\emptyset})\to \varprojlim_{i\in I^{\op}} \ob{D}^*( X_{i})
\]
is an equivalence, where transition maps are given by pullback maps. Then for all $Y\in \n{C}$ the natural map
\[
\Map_{\n{C}}( X_{\emptyset}, Y) \to \varprojlim_{i\in I^{\op}} \Map_{\n{C}}(X_i, Y)
\]
is an equivalence. 
\end{hypothesis}

\begin{remark}
\Cref{HypothesisTechnicalShierk} is helpful for simplifying the condition of a family of maps $\{f_i\colon Y_i\to X\}$ in $\n{C}_E$ to be a $\ob{D}$-cover. Indeed, it is a $\ob{D}$-cover if and only if it satisfies  universal $\ob{D}^*$ and $\ob{D}^!$-descent, namely \textit{loc. cit.} tells us that $\ob{D}^*$-descent already implies that it is a sub-canonical cover. As $\ob{D}$ satisfies categorical Kunneth by construction, \Cref{PropDescentShierk} implies that $\{f_i\}$ is a $\ob{D}$-cover if and only if it satisfies $!$-descent.
\end{remark}

\begin{definition}\label{DefShierkEquivalence}
Suppose that \Cref{HypothesisTechnicalShierk} holds. A map $f\colon Y\to X$ in $\n{P}(\n{C}_E)$ is a \textit{$!$-equivalence} if   any pullback in $\n{P}(\n{C}_E)$ of any iterated diagonal map of $f$ is sent to an equivalence via $\ob{D}^!$. We let $\widetilde{\ob{Shv}}(\n{C}_E)\subset \n{P}(\n{C}_E)$ be the full subcategory of objects $\s{F}$ such that for all $\aleph_1$-compact\footnote{This condition means that the map $f\colon Y\to X$ of presheaves is $\aleph_1$-compact, that is, when pullback along $X'\to X$ with $X'\in \n{C}_E$ the object $Y\times_{X}X'$ is a countable colimit of representable presheaves. The choice of $\aleph_1$ was arbitrary and any uncountable cardinal would do the job, however in practice and in the light set up of condensed mathematics it is a sufficiently large  and very convenient cardinal. } $!$-equivalence $Y\to X$ in $\n{P}(\n{C}_E)$ the map $\s{F}(X)\to \s{F}(Y)$ is an equivalence. 
\end{definition}

\begin{lemma}\label{LemmaStabilityShierkEquivalences}
Let $(\n{C},E)$ be as before. Let $\{f_i\colon Y_i\to X\}$ be a  family of $!$-able maps in $\n{C}$ satisfying $!$-descent and let $X'\subset X$ be the subobject in $\n{P}(\n{C}_E)\subset \n{P}(\n{C})$ defined by the image of the  $f_i$. Then $f\colon X'\to X$ is a $!$-equivalence.  Conversely, any $!$-equivalence in $\n{P}(\n{C}_E)$ is $\infty$-connective for the $!$-topology, in particular  we have  embeddings $\widehat{\ob{Shv}}(\n{C}_E)\subset \widetilde{\ob{Shv}}(\n{C}_E)\subset \ob{Shv}(\n{C}_E)$. 
\end{lemma}
\begin{proof}
For the first claim, note that the iterated diagonals of $f$ are isomorphisms, so it suffices to show that it produces  isomorphisms after any pullback in $\n{P}(\n{C})$ and applying the functor $\ob{D}^!$, this follows from \Cref{PropDescentShierk}. For the second claim, since the diagonal map of a $!$-equivalence is still a $!$-equivalence, to see that it is $\infty$-connective it suffices to show that it is surjective for the $!$-topology (namely, $\infty$-connective maps in a topoi are the same as maps for which all diagonal is an epimorphism). Let $f\colon Y\to X$ be a $!$-equivalence, to prove that it is an epimorphisms we cal localize in $X$, hence we can assume that $X\in \n{C}_E$. Then, we can write $Y=\varinjlim_i Y_i$ as a  colimit of objects in $\n{C}_{E,/X}$. Then, by assumption we have that $\ob{D}(X)=\varprojlim_{i}\ob{D}^!(Y_i)$ along upper $!$-maps, or equivalently, that $\ob{D}(X)=\varinjlim_{i} \ob{D}_!(Y_i)$ in $\cat{Pr}^L$. We claim that the family of maps $\{f_i\colon Y_i\to X\}$ is a $!$-cover, this would imply that $Y\to X$ is an epimorphism in the $!$-topology as wanted. For this last claim, let $\{Y^{\bullet}\}_{\Delta_I^{\op}}$ be the \v{C}ech nerve of $\{f_i\}_i$ (cf. \cite[Definition A.4.5]{HeyerMannSix}). We want to see that the natural map 
\[
\ob{D}_!(|Y^{\bullet}|)=\varinjlim_{\Delta_{I}^{\op}}\ob{D}_!(Y^{\bullet})\to \ob{D}(X)
\]
is an equivalence in $\Pr^L$ (or equivalently, that $\ob{D}(X)=\varprojlim_{\Delta_I} \ob{D}(Y^{\bullet})$). Since $\ob{D}_!(Y)\xrightarrow{\sim}\ob{D}(X)$ and $\ob{D}_!$ satisfies base along maps with target in $\n{C}$ (by K\"unneth), we find that 
\[
\begin{aligned}
\ob{D}_!(|Y^{\bullet}|)=\ob{D}_!(|Y^{\bullet}|)\otimes_{\ob{D}(X)}\ob{D}_!(Y) & =\ob{D}_!(|Y^{\bullet}|\times_X Y) \\ 
				& = \varinjlim_{i} \ob{D}_!(|Y^{\bullet}|\times_X Y_i) \\
				& = \varinjlim_{i} \ob{D}_!(Y_i)\\
				& = \ob{D}(X)
\end{aligned}
\]
where in the third equivalence we have used that $Y^{\bullet}\times_X Y_i$ is split with colimit $Y_i$. This proves the lemma. 
\end{proof}

\begin{corollary}\label{CoroCoTopLocalization}
Let $(\n{C},E)$ be as before. Then $\widetilde{\ob{Shv}}(\n{C}_E)$ is an $\infty$-topoi. Moreorver, $\widetilde{\ob{Shv}}(\n{C}_E)$ is a cotopological localization of $\ob{Shv}(\n{C}_E)$ in the sense of \cite[Definition 6.5.2.17]{HigherTopos}.
\end{corollary}
\begin{proof}
The fact that $\widetilde{\ob{Shv}}(\n{C}_E)$ is an $\infty$-topoi follows from \cite[Theorem 4.3.3]{Anel_2022} and the fact that the class of $!$-equivalences is stable under diagonals and pullbacks by definition. The fact that it is a cotopological localization follows from the fact that it is obtained from sheaves by inverting a class of $\infty$-connective objects. 
\end{proof}

The functor $\n{C}_E\to \n{C}$ sends $!$-covers to $!$-covers, this products a left exact and colimit preserving functor at the level of topoi $\ob{Shv}(\n{C}_E)\to \ob{Shv}(\n{C})$ via left Kan extensions.  

\begin{definition}\label{DefAnStacksFull}
Let $(\n{C},E)$ and $\ob{D}$ be as before. We define the category $\widetilde{\ob{Shv}}(\n{C})\subset \ob{Shv}(\n{C})$ to be the full subcategory of sheaves $\s{F}$ whose restriction to $\n{C}_E$ invert $!$-equivalences. 
\end{definition}

\begin{remark}\label{RemarkHypersheaves}
Let us explain an alternative way to think about $\widetilde{\ob{Shv}}$ in more higher algebraic terms in $\Pr^L$. We have a sequence of geometric maps of $\infty$-topoi
\[
\widehat{\ob{Shv}}(\n{C}_E)\xrightarrow{g_*} \widetilde{\ob{Shv}}(\n{C}_E)\xrightarrow{f_*} \ob{Shv}(\n{C}_E),
\]
with $\widehat{\ob{Shv}}(\n{C}_E)$ the category of hypersheaves on $\n{C}$, whose left adjoints produce idempotent maps of algebras in $\Pr^L$ (with monoidal structure given  by the Cartesian structure)
\[
\ob{Shv}(\n{C}_E)\to \widetilde{\ob{Shv}}(\n{C}_E) \to \widehat{\ob{Shv}}(\n{C}_E). 
\]
We also have a geometric morphism of topoi $h_*\colon \ob{Shv}(\n{C})\to \ob{Shv}(\n{C}_E)$ producing the pullback functor $h^*\colon \ob{Shv}(\n{C}_E)\to \ob{Shv}(\n{C})$. Then, the natural maps of topoi
\[
\ob{Shv}(\n{C})\otimes_{\ob{Shv}(\n{C}_E)}\widetilde{\ob{Shv}}(\n{C}_E)\xrightarrow{\sim} \widetilde{\ob{Shv}}(\n{C}) \mbox{ and } \ob{Shv}(\n{C})\otimes_{\ob{Shv}(\n{C}_E)}\widehat{\ob{Shv}}(\n{C}_E)\xrightarrow{\sim} \widehat{\ob{Shv}}(\n{C})
\]
are isomorphisms. Indeed,   the base change $\ob{Shv}(\n{C})\to \ob{Shv}(\n{C})\otimes_{\ob{Shv}(\n{C}_E)}\widetilde{\ob{Shv}}(\n{C}_E)$ is the localization of $\ob{Shv}(\n{C})$ obtained by inverting the arrows in $\ob{Shv}(\n{C})$ which are inverted in $\widetilde{\ob{Shv}}(\n{C}_E)$, i.e. the localization of $\ob{Shv}(\n{C})$ obtained by inverting the $!$-equivalences in $\ob{Shv}(\n{C}_E)$. This is precisely the category $\widetilde{\ob{Shv}}(\n{C}_E)$ of \Cref{DefAnStacksFull}. In a similar fashion, the tensor $ \ob{Shv}(\n{C})\otimes_{\ob{Shv}(\n{C}_E)}\widehat{\ob{Shv}}(\n{C}_E)$ is the localization of $\ob{Shv}(\n{C})$ by inverting all hypercovers of $\ob{Shv}(\n{C}_E)$, since these hypercovers generate those in $\ob{Shv}(\n{C})$, the tensor product in this case is nothing but $\widehat{\ob{Shv}}(\n{C})$.
\end{remark}


%

The following proposition implies that $\ob{D}^*$ is still a sheaf in $\widetilde{\ob{Shv}}(\n{C})$. 

\begin{proposition}\label{PropDescentShierkEquiv}
Let $Y\to X$ be a $!$-equivalence in $\ob{Shv}(\n{C})$. Then the natural map $\ob{D}^*(X)\to \ob{D}^*(Y)$ is an equivalence.  
\end{proposition}
\begin{proof}
As $\ob{D}^*$ is a sheaf for the $!$-topology, we can work locally on $X$ and assume that $X\in \n{C}$ and that $Y=\varinjlim_i Y_i$ is a colimit of $!$-able maps in $\n{C}_{/X}$. Then $\ob{D}(X)=\ob{D}_!(Y)=\varinjlim_{i} \ob{D}(Y_i)$ in $\Pr^L_{\ob{D}(X)}$. By the fully faithful embedding of \Cref{LemFullyFaithfulKernelCat} one deduces that $\ob{D}_!(Y)$ is a predual of $\ob{D}^*(Y)$, that is, the dual of $\ob{D}_!(Y)$ is the limit $\varprojlim_i \ob{D}(Y_i)^{\vee}$ along the duals of $\ob{D}(Y_i)$, but the categories $\ob{D}(Y_i)$ are self dual in $\Pr^L_{\ob{D}(X)}$ and the dual of lower $!$-functors are upper $*$-functors. The lemma follows.  
\end{proof}

\begin{corollary}\label{CoroSixFunctorAnStk}
Let $(\n{C},E)$ be as before. There is a geometric set up $(\widetilde{\ob{Shv}}(\n{C}),E')$ and a morphism of geometric set ups $(\n{C},E)\to (\widetilde{\ob{Shv}}(\n{C}),E')$ such that the six functor formalism $\ob{D}$ on  $(\n{C},E)$ extends uniquely to $ (\widetilde{\ob{Shv}}(\n{C}),E')$, and such that $E'$ satisfies the properties of \Cref{TheoExtension6}. 
\end{corollary}
\begin{proof}
This follows from \Cref{TheoExtension6} as $\ob{D}^*$  satisfies descent  for the cotopological localization $\ob{Shv}(\n{C})\to \widetilde{\ob{Shv}}(\n{C})$.
\end{proof}

Finally, we discuss how to construct morphisms from topoi towards analytic stacks. 

\begin{lemma}\label{xoajsodasod}
Let $(\n{C},E)$ be as above. Let $\n{D}$ be a site with fiber products, and let  $\ob{Shv}(\n{D})$ be its topos of sheaves. Let $\ob{Shv}(\n{D})\to \n{T}$ be a cotopological localization onto an $\infty$-topoi. Suppose that we are given with a left exact functor $F\colon \n{D}\to \n{P}(\n{C})$ and consider its left Kan extension to a left exact and colimit preserving functor $\n{P}(\n{D})\to \n{P}(\n{C})$. Suppose that any cover $\{Y_i\to X\}$ in $\n{D}$ is mapped to a cover $\{F(Y_i)\to F(X)\}$ in $\n{P}(\n{C})$ for the $!$-topology. Then there is a unique left exact colimit preserving localization $F'\colon \ob{Shv}(\n{D})\to \ob{Shv}(\n{C})$ of $F$. Suppose that in addition all $\infty$-connective object in $\ob{Shv}(\n{D})$ that induces an equivalence in $\n{T}$ is sent to a $!$-equivalence in $\ob{Shv}(\n{C})$ (as $F'$ is colimit preserving and left exact this map is automatically $\infty$-connective). Then there is a unique left exact and colimit preserving localization $\n{T}\to \widetilde{\ob{Shv}}(\n{C})$ of $F'$. 
\end{lemma}
\begin{proof}
This follows from abstract non-sense in topos theory, see \cite[Theorem 4.3.3]{Anel_2022}.
\end{proof}

The previous lemma is helpful to construct natural maps between categories of analytic stacks.


\begin{lemma}\label{LemmaFunctorialityAnStk}
Let $(\n{C},E)$ be a geometric set up with finite limits and $\ob{D}$ a presentable  $6$-functor formalism on $(\n{C},E)$ satisfying categorical K\"unneth formula and \Cref{HypothesisTechnicalShierk}. Let $(\n{C}_1,E_1)$ be another geometric set up with finite limits and consider a left exact morphism of geometric set ups $(\n{C}_1,E_1)\to (\n{C},E)$. Consider $\ob{D}_1$ the restriction of $\ob{D}$ to a six functor formalism on $(\n{C}_1,E_1)$, and suppose that $\ob{D}_1$ also satisfies \Cref{HypothesisTechnicalShierk}. Then, the functor $\n{C}_1\to \n{C}_2$ sends $!$-covers to $!$-covers and it induces a left exact and colimit preserving morphism of sheaves $\ob{Shv}(\n{C}_1)\to \ob{Shv}(\n{C})$ for the $!$-topology. Furthermore, a $!$-equivalence in $\ob{Shv}(\n{C}_1)$ is sent to a $!$-equivalence in $\ob{Shv}(\n{C})$, hence the previous functor localizes to a left exact and colimit preserving functor of topoi $\widetilde{\ob{Shv}}(\n{C}_1)\to \widetilde{\ob{Shv}}(\n{C})$.
 \end{lemma}

 \begin{proof}
The fact that $\n{C}_1\to \n{C}_2$ sends $!$-covers to $!$-covers is clear as the functor preserves pullbacks, and $\ob{D}_1$ is the restriction of $\ob{D}$ to $(\n{C}_1,E_1)$.  The factorization of topoi $F\colon \ob{Shv}(\n{C}_1)\to \ob{Shv}(\n{C}_2)$ is then formal. It is left to show that the previous map sends $!$-equivalences to $!$-equivalences, as it preserves $\infty$-connective objects (being the pullback of a morphism of topoi), this would show that it factors through a left exact and colimit preserving functor $\widetilde{\ob{Shv}}(\n{C}_1)\to \widetilde{\ob{Shv}}(\n{C})$. 


Let $f\colon Y\to X$ be a $!$-equivalence in $\Shv(\n{C}_1)$, we want to show that $F(f)\colon F(Y)\to F(X)$ is a $!$-equivalence in $\ob{Shv}(\n{C})$. Let $Z\in \n{C}$, we want to show that for any iterative diagonal $F(Y)\to  F(Y)^{(n)}$ of $F(f)$, and any map $Z\to F(Y)^{(n)}$, the pullback $F(Y)\times_{F(Y)^{(n)}} Z \to Z$ arises from a $!$-equivalence in $\ob{Shv}(\n{C}_E)$, locally for the $!$-topology on $Z$. Since $F$ is left exact, we have that $F(Y)^{(n)}=F(Y^{(n)})$, and since $!$-equivalences in $\ob{Shv}(\n{C}_1)$ are stable under diagonals, we can assume that $n={0}$, i.e. it suffices to consider a map $Z\to F(X)$ with $Z\in \n{C}$.

 Since $f$ is a $!$-equivalence, we can find an epimorphism $\bigsqcup X_i\to X$ with $X_i\in \n{C}_1$ such that each pullback $Y\times_{X} X_i\to X_i$ arises from a $!$-equivalence in $\ob{Shv}(\n{C}_{1,E_1})$. As $F$ preserves colimits, we have that $\bigsqcup_i F(X_i) \to F(X)$ is an epimorphism, and so there exists a $!$-cover $\widetilde{Z}:=\bigsqcup_{j\in  J} Z_j\to Z$  and a lift $\widetilde{Z}\to \bigsqcup_i F(X_i) $. Hence, we can suppose without loss of generality that $X\in \n{C}_1$, $Y\to X$ arises from a $!$-equivalence in $\ob{Shv}(\n{C}_{1,E_1})$, and $Z\in \n{C}$. In this case, we can write $Y=\varinjlim_{i} Y_i$ as a colimit of $!$-able objects over $X$ such that $\ob{D}(X)=\varinjlim_{i} \ob{D}_!(Y_i)$, and so $F(Y)=\varinjlim_{i} F(Y_i)\to F(X)$.  This implies that $F(Y)\times_{F(X)} Z=\varinjlim_{i} F(Y_i)\times_{F(X)} Z$ and so that $F(Y)\times_{F(X)}Z\to Z$ arises from an object in $\ob{Shv}(\n{C}_E)$, we want to show that this object is actually a $!$-equivalence.   For that, we need to see that $\ob{D}^!$ is identified in any pullback of any iterative diagonal. Since the iterative diagonal arises as pullback of that of $F(f)$, and since   $\ob{D}^!$ satisfies descent for the $!$-topology on $\ob{Shv}(\n{C}_E)$,  by reapeating the previous reduction it suffices to show that $\ob{D}_!(F(Y)\times_{F(X)}Z) =\ob{D}(Z)$, but we have by categorical K\"unneth that 
 \[
 \begin{aligned}
 \ob{D}_!(F(Y)\times_{F(X)}Z) & =\varinjlim_{i} \ob{D}_!(F(Y_i)\times_{F(X)} Z)\\
 &  = \varinjlim_i \ob{D}_!(F(Y_i))\otimes_{\ob{D}(F(X))} \ob{D}(Z) \\ 
 & = \ob{D}(F(X))\otimes_{\ob{D}(F(X))}\ob{D}(Z) \\ 
 & =\ob{D}(Z)
\end{aligned} 
 \]
 proving what we wanted. 
 \end{proof}

\subsubsection{Analytic stacks}

We specialize the constructions of \Cref{SectionShierkequivalences} to analytic stacks.

The following lemma verifies \Cref{HypothesisTechnicalShierk} in the case of analytic rings.

\begin{lemma}\label{LemStarDescentAnStk}
Let $I$ be a category and $I^{\triangleleft}$ its right cone with inital object $-1$. Let $A_{\bullet}\colon I^{\triangleleft}\to \AnRing$ be a diagram. Suppose that the natural map $\ob{D}^*(A_{-1})\to \varprojlim_{i\in I} \ob{D}^*(A_i)$ is an equivalence. Then for all $B\in \AnRing$ the natural map 
\[
\Map_{\AnRing}(B, A_{-1})\xrightarrow{\sim} \varprojlim_{i\in I} \Map_{\AnRing}(B,A_i)
\] 
is an isomorphism. 
\end{lemma}
\begin{proof}
By definition, $\ob{Map}_{\AnRing}(B,A)\subset \ob{Map}_{\CondRing}(B^{\triangleright},A^{\triangleright})$ is the full subanima of maps $f\colon B^{\triangleright}\to A^{\triangleright}$ for which $A$-complete modules are $B$-complete. Suppose that $\ob{D}(A_{-1})=\varprojlim_i (\ob{D}^*(A_{i}))$, then $A^{\triangleright}_{-1}= \varprojlim_i (A_{i}^{\triangleright})$ as condensed animated rings, and we have 
\[
\ob{Map}_{\CondRing}(B^{\triangleright},A^{\triangleright}_{-1}) = \varprojlim_i \bigg(\ob{Map}_{\CondRing}(B^{\triangleright},A^{\triangleright}_{i})\bigg).
\]
Therefore, to finish the proof it suffices to show that if a map $B^{\triangleright}\to A^{\triangleright}_{-1}$ is such that the forgetful sends $A_{i}$-complete modules to $B$-complete modules for all $i\in I$, then it sends $A_{-1}$-complete modules to $B$-complete modules. This is a formal consequence of $\ob{D}(A_{-1})=\varprojlim_i (\ob{D}^*(A_{i}))$.
\end{proof}

With the previous preparations, we give the definition of the category of analytic stacks.

\begin{definition}\label{DefAnStk}
Let $\n{C}=\AnStk^{\aff}$ be the category of affinoid analytic stacks, $E$ the class of $!$-able arrows and $\ob{D}$ the six functor formalism on $(\n{C},E)$ of quasi-coherent sheaves constructed in \Cref{CorSixFunctorFormalismAffinoids}. Let $S\in \n{C}$ and let $\n{C}_0\subset \n{C}_{/S}$ be an essentially small full subcategory stable under finite limits and coproducts. Let $E_0=E\cap \ob{Arr}(\n{C}_0)$ be the class of $!$-able arrows in $\n{C}_0$ and by an abuse of notation denote by  $\ob{D}'$ be the restriction of $\ob{D}$ to $\n{C}_0$. We define the category of \textit{analytic stacks} of $\n{C}_0$, denoted as $\AnStk(\n{C}_0)$, to be $\widetilde{\ob{Shv}}(\n{C}_0)$ as in \Cref{DefAnStacksFull}. By an abuse of notation, we also denote by $\ob{D}$ the extension of the quasi-coherent six functor formalism to analytic stacks $\AnStk(\n{C}_0)$ as in \Cref{CoroSixFunctorAnStk}. 
\end{definition}

\begin{warning}\label{RemarkAnStk}
In \Cref{DefAnStk} we have avoided the case $\n{C}_0=\n{C}$ to construct the category of analytic stacks. In \cite{AnalyticStacks} Clausen and Scholze mentioned that this is actually possible, but it requires to prove some additional technical  properties of the category of analytic rings and the $!$-topology that are unknown to the author (eg. presentability of analytic rings or accesibility of $!$-equivalences).


Nonetheless, in practice one can always choose a large enough essentially small category of analytic rings $\n{C}_0$ that will suffice for applications and perform the construction of \Cref{DefAnStk}. This is the point of view we adopt in these notes for the application to Serre duality in algebraic geometry. With this technical point in mind, and by making an abuse of notation,  we will denote by $\AnStk$ the category of analytic stacks for a large enough  full subcategory $\n{C}_0\subset \AnStk^{\aff}$ (that might depend on the context)  stable under countable limits and finite coproducts (or equivalently a small family of analytic rings stable under countable colimits and finite products).
\end{warning}

\subsection{Algebraic stacks}\label{ssAlgebraicStacks}

In this section we give a realization of a suitable category of algebraic stacks (understood in a broad sense) in analytic stacks. Let $\Ring$ be the category of animated rings and let $\Sch^{\aff}=\Ring^{\op}$ be its opposite category of affine schemes. A standard Grothendieck topology  in $\Ring$ is the flat topology.  Consider the functor of coefficients $R\mapsto \ob{D}(R)=\Mod_R(\ob{D}(\Z))$ given by $R$-modules in derived abelian groups. One can also consider its condensed promotion $\ob{D}^{\cond}(R)= \ob{D}(R)\otimes_{\ob{D}(\Z)} \ob{D}(\Z^{\cond})=\Mod_{R}(\ob{D}(\Z^{\cond}))$.   By taking the class of arrows $E=P=all$ and $I=iso$, \Cref{LemStabilityShierkMaps} constructs a $6$-functor formalism for  both $\ob{D}$ and $\ob{D}^{\cond}$, and we have a morphism of $6$-functor formalisms $\ob{D}\to \ob{D}^{\cond}$ as morphisms of operads $\Corr(\Sch^{\aff})^{\otimes}\to \Cat^{\times}_{\infty}$. 

\begin{remark}\label{RemarkRelationCondensedUsualCategory}
There is a clean relation between the functors $\ob{D}$ and $\ob{D}^{\cond}$. Let $\ob{D}_!$ and $\ob{D}_!^{\cond}$ denote the left Kan extension of $!$-sheaves from $\Sch^{\aff}$ to presheaves. It is clear by construction that $\ob{D}^{\cond}_!=\ob{D}_!\otimes_{\ob{D}(\Z)}\ob{D}(\Z^{\cond})$ as functors $\n{P}(\Sch^{\aff})\to \Pr^L$. Taking $\ob{D}(\Z^{\cond})$-linear duals one gets 
\[
\begin{aligned}
\ob{D}^{\cond} & = \Fun_{\ob{D}(\Z^{\cond})}(\ob{D}_!^{\cond}, \ob{D}(\Z^{\cond})) \\ 
				& = \ob{Fun}_{\ob{D}(\Z^{\cond})}(\ob{D}_!\otimes_{\ob{D}(\Z)} \ob{D}(\Z^{\cond}), \ob{D}(\Z^{\cond})) \\ 
				& = \ob{Fun}_{\ob{D}(\Z)}(\ob{D}_!, \ob{D}(\Z^{\cond})).  
\end{aligned}
\]
More generally, for any analytic ring $\mathcal{A}$ we have a six functor formalism $\ob{D}^{\mathcal{A}}$ on $\Sch^{\aff}_{/\Spec \mathcal{A}(*)}$ constructed in the same way as before,  such that $\ob{D}_!^{\mathcal{A}}=\ob{D}_!\otimes_{\ob{D}(\mathcal{A}(*))} \ob{D}(\mathcal{A})$. If $\ob{D}(\mathcal{A})$ is compactly generated (eg. $\mathcal{A}=\Z_{\sol}$), then it is a dualizable $\ob{D}(\mathcal{A}(*))$-linear category and we have that 
\[
\ob{D}^{\mathcal{A}}= \ob{Fun}_{\ob{D}(\Z)}(\ob{D}_!, \ob{D}(\mathcal{A})) = \ob{D}^*\otimes_{\ob{D}(\mathcal{A}(*))} \ob{D}(\mathcal{A}). 
\]
\end{remark}

\Cref{PropDescentShierk} imply that a map $A\to B$ of animated rings is a $!$-cover for $\ob{D}$ (eq. for $\ob{D}^{\cond}$) if and only if it is a descendable map. Therefore, the $!$-topology on affine schemes with respect to the functor $\ob{D}$  is the same as the descendable topology. Similarly, if $f\colon \s{F}\to \s{G}$ is a pre-stack on $\ob{Sch}^{\ob{aff}}$, the map $f$ is a $!$-equivalence for the functor $\ob{D}$ if and only if it is for the functor $\ob{D}^{\ob{cond}}$ (see \Cref{RemarkShierkTop}).  One formally obtains a left exact colimit preserving functor  of categories 
\begin{equation}\label{eqRealizationAlgInAn}
(-)^{\cond}\colon \widetilde{\ob{Shv}}(\Sch^{\aff})\to \AnStk 
\end{equation}
 such that the composition with the composite $ \widetilde{\ob{Shv}}(\Sch^{\aff})\xrightarrow{\ob{D}}  \AnStk \to \Cat_{\infty}$ being the right  Kan extension of $\ob{D}^{\cond}$ from $\Sch^{\aff}$ to $ \widetilde{\ob{Shv}}(\Sch^{\aff})$. We will denote $\ob{AlgStk}:=  \widetilde{\ob{Shv}}(\Sch^{\aff})$ and call it the category of \textit{algebraic stacks}\footnote{This is a highly non-standard use of the name \textit{algebraic stacks}; classically  one uses the \'etale or flat topologies, and some additional geometric constrains are imposed in the stacks such as being covered by a scheme along a representable \'etale or flat map.  In our framework, this name is convenient to distinguish those analytic stacks that arise from algebraic geometry.} by inverting $!$-equivalences of presheaves as in \Cref{SectionShierkequivalences}.

\begin{remark}\label{RemarkShierkTop}
Thanks to \Cref{PropDescentShierk}, a finite collection of maps $\{A\to B_i\}_{i\in I}$ of animated rings is a $!$-cover if and only if the map $A\to \prod_{i }B_i$ is descendable as in \cite{MathewDescent}.

 Let $\s{F}\in \n{P}(\ob{Sch}^{\ob{aff}})$, by \Cref{RemarkRelationCondensedUsualCategory}  one has that 
\[
\ob{D}_!(\s{F})\otimes_{\ob{D}(\Z)} \ob{D}(\Z^{\ob{cond}}) = \ob{D}_!^{\ob{cond}}(\s{F}),
\]
where $\ob{D}_!$ is the left Kan extension of the functor $\Spec(A)\mapsto \ob{D}_!(A)$ with transition maps given by lower $!$-maps, and $\ob{D}_!^{\ob{cond}}(\s{F})$ is the left Kan extension of $\Spec(A)\mapsto \ob{D}_!(A^{\ob{cond}})$. Thus, a map $f\colon \s{F}\to \s{G}$ of pre-stacks on affine sechemes is a  $!$-equivalence for   $\ob{D}_!$ if and only if it is for $\ob{D}_!^{\ob{cond}}$. Indeed, if   $f$ is a $!$-equivalence for $\ob{D}$ then it is so for $\ob{D}^{\ob{cond}}$ by base change. Conversely, if $f$ is a $!$-equivalence for $\ob{D}^{\ob{cond}}$, one wants to deduce that it is a $!$-equivalence for $\ob{D}$. This follows from the fact that linear categories satisfy descent along $\ob{D}(\Z)\to \ob{D}^{\ob{cond}}(\Z)$ as it has a section as symmetric monoidal categories given by evaluating at the point $\ob{ev}_{*}\colon \ob{D}^{\ob{cond}}(\Z)\to \ob{D}(\Z)$, see \Cref{Lemma:evaluationPointSymmetricMonoidal}.

By \Cref{LemmaFunctorialityAnStk} there is a natural left exact and colimit preserving functor $(-)^{\ob{cond}}\colon \ob{AlgStk}\to \ob{AnStk}$ from algebraic to analytic stacks sending $\Spec A$ to $\AnSpec (A^{\ob{cond}})$.   Examples of $!$-equivalences (or descendable maps), are $\omega_1$-fpqc hypercovers (\Cref{CoroHypersheaves}), or $h$-covers of finite type morphisms of static Noetherian schemes \cite[Theorem 11.26]{BhattScholzeAffineGras}. 
\end{remark}

For the construction \eqref{eqRealizationAlgInAn} to be useful in practice, one needs to produce enough (good) $!$-(hyper)covers of animated rings. In \cite[Proposition 3.31]{MathewDescent} Mathew proves that $\omega_1$-compact faithfully flat maps are descendable:

\begin{lemma}\label{LemmaCountableApproxFlatRing}
Let $f\colon A\to B$ be an $\omega_1$-compact faithfully flat morphism of animated rings. Then $f$ is descendable of index $\leq 2$. 
\end{lemma}
\begin{proof}
Since $A\to B$ is faithfully flat, $Q=\ob{cofib}(A\to B)$ is a flat $A$-module.  Since $Q$ is countably presented, by Lazard's theorem \cite[Theorem 7.2.2.15]{HigherAlgebra}, $Q$ is a countable filtered colimit of finite free $A$-modules. Thus, $Q\otimes_A Q$ is also a countable filtered colimit of  finite free $A$-modules and it has projective dimension $\leq 1$. In particular,  the composition $Q\otimes_A Q[-2]\to Q[-1]\to A$ must vanish. One deduces that $f$ is descendable of index $\leq 2$. 
\end{proof}

We can extend $!$-descent  to faithfully flat hypercovers of countable rings. 

\begin{definition}\label{DefOmega1fpqc}
Let $A$ be an animated ring, an $\omega_1$-fpqc cover of $A$ is a finite collection of flat maps of animated rings $\{A\to A_i\}$ such that each $A_i$ is  an $\omega_1$-compact $A$-algebra,   and $A\to \prod_{i\in I} A_i$ is faithfully flat. 
\end{definition}

\begin{prop}\label{PropFaithFlatdescent}
Let $A\to A_{\bullet}$ be a $\omega_1$-compact faithfully flat hyper-cover of  animated rings. Then the map $|\Spec A_{\bullet}|\to \Spec A$ from the geometric realization is a $!$-equivalence.   
\end{prop}
\begin{proof}
 By \Cref{PropDescentShierk} it suffices to show that the functor $A\to \ob{D}^!(A)$ with pullback maps given by upper $!$-maps is an hypersheaf for the $\omega_1$-compact  flat topology. Since $\omega_1$-compact faithfully flat maps  are descendable by \Cref{LemmaCountableApproxFlatRing}, we know that $A\mapsto \ob{D}^!(A)$ is a sheaf for the  $\omega_1$-compact flat topology. By \cite[Proposition A.3.21]{MannSix} it suffices to show that if $A\to A_{\bullet}$ is an $\omega_1$-compact faithfully flat hypercover in $\Ring_{\omega_1}$ then the natural functor 
\begin{equation}\label{eq013rjo2dwqr}
\ob{D}(A)\to \ob{Tot}(\ob{D}^!(A_{\bullet}))
\end{equation}
is fully faithful. To prove fully faithfulness, we apply the same proof of \cite[Proposition 3.1.26]{MannSix}. The functor \eqref{eq013rjo2dwqr} has a left adjoint sending a cocartesian section $(M_n)_{n\in \Delta}$ to $\varinjlim_{[n]\in \Delta^{\op}} M_n$ where we see all the $M_n$ as $A$-modules by restriction of scalars. Thus, we have to show that for $M\in \ob{D}(A)$ the natural map 
\begin{equation}\label{eq91nkedaped}
\varinjlim_{[n]\in \Delta^{\op}} \Hom_{A}(A_n, M)\to M
\end{equation}
is an equivalence. The ring $A_n$ is a flat $A$-module, by Lazard's theorem \cite[Theorem 7.2.2.15]{HigherAlgebra}, $A_n$ can be written as a filtered colimit of finite free $A$-modules. As $A_n$ is an $\omega_1$-compact $A$-module (being countably generated over $A$ as algebra),  it is a countable filtered  colimits of finite free $A$-modules. In particular, the functor $\Hom_A(A_n, -)$ is left $t$-exact and has cohomological dimension $\leq 1$. 

Writing $M=\varinjlim_{k} \tau_{\geq -k} M$ as a colimit of its canonical truncation, we see that 
\[
\varinjlim_{k} \Hom_{A}(A_n, \tau_{\geq -k} M)\to  \Hom_{A}(A_n,M)
\]
is an isomorphism. Therefore, we can assume without loss of generality that $M$ is $1$-connective. In that case, the simplicial object $(\Hom_A(A_n, M))_{[n]\in \Delta^{\op}}$ is connective. Fix $n\geq 1$,  thanks to \cite[Proposition 1.2.4.5]{HigherAlgebra} for $k<n$ there is an equivalence of homotopy groups
\[
\pi_k(\varinjlim_{ \Delta^{\op}} \Hom_{A}(A_{\bullet}, M) )= \pi_k(\varinjlim_{\Delta^{\op}_{\leq n}} \Hom_A(A_{\bullet},M)). 
\]
On the other hand, let $\Spec(\widetilde{A}_{\bullet})$ be the coskeleton of $(\Spec A_{\bullet})_{\Delta^{\op}_{\leq n}}$. Then, as $\ob{D}^!$ is a flat sheaf, and \cite[Lemma 6.5.3.9]{HigherTopos} we have that 
\[
\ob{D}^!(A)=\ob{Tot}(\ob{D}^!(\widetilde{A}_{\bullet}))
\]
and in particular $M=\varinjlim_{\Delta}\Hom_A(\widetilde{A}, M)$. We deduce that 
\[
\begin{aligned}
\pi_k(\varinjlim_{ \Delta^{\op}} \Hom_{A}(A_{\bullet}, M) ) & =  \pi_k (\varinjlim_{ \Delta^{\op}_{\leq n}} \Hom_{A}(A_{\bullet}, M)) \\ 
& = \pi_k (\varinjlim_{ \Delta^{\op}_{\leq n}} \Hom_{A}(\widetilde{A}_{\bullet}, M)) \\
& = \pi_k (\varinjlim_{ \Delta^{\op}} \Hom_{A}(\widetilde{A}_{\bullet}, M))  \\
& = \pi_k(M). 
\end{aligned}
\]
Since $0\leq k<n$ was arbitrary, we deduce that the natural map 
\eqref{eq91nkedaped} is an equivalence, proving what we wanted. 
\end{proof}

From \Cref{LemmaFunctorialityAnStk} we deduce the following corollary.

\begin{corollary}\label{CoroHypersheaves}
There are  left exact colimit preserving functors
\[
(-)^{\ob{cond}}\colon \widehat{\Shv}_{\omega_1-\ob{fpqc}}(\Sch^{\aff})\to \ob{AlgStk} \to \AnStk
\]
from (accessible) hypersheaves on $\Sch^{\aff}$ in the fppf topology to algebraic stacks (and then to analytic stacks), given by the left Kan extension of the inclusion $\Sch^{\aff}\to \ob{AlgStk}$. 
\end{corollary}

We finish this section proving that qcqs morphisms of schemes are $!$-able and proper in an appropriate sense. 

\begin{prop}\label{PropSixFunctorsOnSchemes}
Let $f\colon Y\to X$ be a qcqs morphism of schemes. Then $f$ is $!$-able for both $\ob{D}$ and $\ob{D}^{\cond}$, it is prim and there is a trivialization of the dualizing sheaf $\delta_f\cong 1$. 
\end{prop}
\begin{proof}
By \cite[Lemmas 4.5.7]{HeyerMannSix} we can prove that $f$ is prim locally in an universal $\ob{D}^*$-cover in the source, so we can assume without lost of generality that $X$ is affine. By \cite[Lemmas 4.5.8 (ii) and 4.7.4]{HeyerMannSix} we can also prove that $f$ is prim after pre-composing with a morphism $Y'\to Y$ of schemes which is prim and descendable. 

Let us first show the proposition if $f$ is affine. If $X$ is affine then so is $Y$, and $f$ is prim by construction of the quasi-coherent six functor formalism. We also have $\delta_f=1$ by construction.  If $f$ is quasi-compact and separated, let   $g\colon \bigsqcup U_i\to Y$  be a finite open affine cover. Then the map $g$ is affine and it is prim by the previous step. By Zariski descent the map $g$ is descendable and therefore $g$ is a prim and descendable cover. Since the composite $\bigsqcup_i U_i\to Y \to X$ is a disjoint union of morphisms of affine schemes, it is $!$-able and prim and then so is the map $f$.  Finally, we deal with the qcqs case. We can find a  finite open affine cover $g\colon \bigsqcup_i U_i\to  Y$ such that $g$ is quasi-compact and separated. By the previous step $g$ is prim, and it is descendable thanks to Zariski descent. As the composite $\bigsqcup_i U_i\to  Y\to X$ is a finite disjoint union of prim maps (being affine again), we deduce that $f$ is prim as wanted. 

Finally, it is left to identify the dualizing sheaf. If $f$ is a morphism of affine schemes, then $\delta_f =1$ by construction of the six functor formalism.  If $f$ is an affine morphism of schemes, we have a trivialization $\delta_* \cong 1$  obtained by glueing  the trivializations of the affine cases  along the target in the Zariski topology. For the quasi-compact and separated case, the diagonal $\Delta f$ is affine so that $\delta_{\Delta f}\cong 1$, and the proof of  \cite[Lemma 4.6.4]{HeyerMannSix} yields an isomorphism $\delta_f \cong 1$ (only depending on the trivialization of the diagonal) as $f$ is prim.  By an inductive argument one proves that $\delta_f \cong 1$ for $f$ a qcqs morphism of schemes.
\end{proof}

In \Cref{SectionSerreDualitySolid} we shall prove Serre duality in this set up,  that is,   if $f\colon Y\to X$ is a proper smooth morphism of schemes, then  it is cohomologically smooth with dualizing sheaf given by  $\det (\Omega^{1}_{Y/X})[\dim(Y/X) ]$.

\subsubsection{Formal schemes} \label{sss:Formal Schemes}

An important class of algebraic stacks that arise in practice are formal schemes.  We begin by discussing descent of the Zariski topology along $!$-covers

\begin{lemma}\label{LemmaDescentSpectrum}
Let $A\to B^{\bullet}$ be a $!$-hypercover of animated rings. Then $|\Spec B^{\bullet}| \to \Spec A$ is an equivalence of topological spaces. 
\end{lemma}
\begin{proof}
It suffices to show that if $A\to B$ is a $!$-cover (i.e. if it is descendable), then the map $f\colon \Spec B\to \Spec A$ is a quotient map of topological spaces. Since base change along $A\to B$ is conservartive, the map $f$ is clearly surjective. To see that $f$ is a quotient map, let $C\subset \Spec A$ be a subspace such that $f^{-1}(C)\subset \Spec B$ is a qcqs open subspace, we want to see that $C$ is itself qcqs open. We will follow the argument of the proof of \cite[Theorem 11.26]{BhattScholzeAffineGras}. Let $A\to B^{\bullet}$ be the \v{C}evh co-nerve of $A\to B$ and let  $B^{\bullet}\to B^{\bullet}_C$ be the idempotent algebra corresponding to the Zariski open subspaces $f^{\bullet,-1}(C)\subset \Spec (B^{\bullet})$ where $f^{\bullet}\colon \Spec(B^{\bullet})\to \Spec(A)$. Since $A\to B$ is descendable, we have an equivalence $\ob{D}(A)= \ob{Tot}(\ob{D}(B^{\bullet}))$ which restricts itself to an equivalence of symmetric monoidal categories $\ob{D}_{+}(A) =\ob{Tot}(\ob{D}_+(B^{\bullet}))$ on eventually connective modules (this follows from the fact that $A\to B$ is already descendable in the stable category $\ob{D}_+(A)$). Hence, $B^{\bullet}\to B^{\bullet}_C$  descends to an idempotent, eventually connective algebra $A\to A_C$ in $\ob{D}_+(A)$. Since $B\to B_C$ is compact as algebras the same holds for $A\to A_C$, and by \cite[Theorem 1.10]{BhattHalpernTannaka} $A_C$ is the idempotent algebra attached to a Zariski open subspace of $\Spec A$ which musth be $C$, proving that it is open as wanted. 
\end{proof}

\begin{proposition}\label{PropFullyFaithfulSchemes}
Let $X$ be a derived scheme considered as a functor $X\colon \ob{Ring}^{\delta}\to \ob{Ani}$ from discrete animated rings to anima. Then $X$ satisfies $!$-descent. Moreover, if $\s{F}\to \s{G}$ is a $!$-equivalence of sheaves for the $!$-topology on $\ob{Sch}^{\ob{aff}}$, then the natural map $X(\s{G})\to X( \s{F})$ is an equivalence. In particular, $X$ is an algebraic stack in the sense of \Cref{RemarkShierkTop}, and we have a fully faithful embedding
\[
\ob{Sch}\hookrightarrow \cat{AlgStk}.
\]
\end{proposition}
\begin{proof}
Let $X$ be a scheme and $A\to B$ be a descendable map with \v{C}ech co-nerve $A\to B^{\bullet}$. We want to show that the natural map $X(A)\to \ob{Tot}(X(B^{\bullet}))$ is an equivalence. By passing to the associated topological spaces, we have a commutative diagram of anima (where the bottom objects are sets)
\[
\begin{tikzcd}
X(A) \ar[r] \ar[d]& \ob{Tot}(X(B^{\bullet})) \ar[d] \\
\ob{Map}_{\cat{Top}}(\Spec A, |X|) \ar[r] & \ob{Tot}( \ob{Map}_{\cat{Top}}(\Spec(B^{\bullet}), |X|)).
\end{tikzcd}
\]
By \Cref{LemmaDescentSpectrum} the lower map is an equivalence, and to show that $X(A)=\ob{Tot}(X(B^{\bullet}))$, it suffices to prove it on the fibers of the vertical arrows.  Hence, we can assume without loss of generality that we have fixed a map $\Spec A\to |X|$ of topological spaces. Hence, by Zariski descent, one reduces to the case when $\Spec A\to |X|$ factors through an open  affine scheme $\Spec C\subset X$ and so to  $X=\Spec C$ itself. The statement then follows from the fact that the $!$-topology on animated rings is subcanonical (eg. as it also satisfies $\ob{D}^*$-descent and by \Cref{LemStarDescentAnStk}). 

Next, we show that $X$ also inverts $!$-equivalences. Let $\s{F}\to \s{G}$ be a $!$-equivalence of $!$-sheaves. We want to see that $X(\s{G})\to X(\s{F})$ is an equivalence of anima. By $!$-descent, we can assume without loss of generality that $\s{F}=\Spec(A)$ is affine itself, and then $\s{G}\to \Spec(A)$  is given  by a colimit $\s{G}=\varinjlim_i \Spec(B_i)$ such that, in particular,  $\ob{D}_!(A)=\varinjlim_i \ob{D}_!(B_i)$, and by taking duals that $\ob{D}^*(A)=\ob{D}^*(B_i)$. Then the statement follows from the same argument as in the affine case:  by \Cref{LemmaDescentSpectrum} the map of topological spaces $|\s{G}|\to \Spec(A)$ is an isomorphism, and we can assume that we have fixed a map $\Spec(A)\to X$. Then, by Zariski descent, we reduce to the case when $X=\Spec(C)$ is affine, in which case it follows from the fact that the $!$-topology is subcanonical and satisfies \Cref{HypothesisTechnicalShierk}. 
\end{proof}

\begin{definition}\label{DefFormalScheme}
Let $X$ be a scheme and $C\subset X$ a locally closed subspace. We define the \textit{formal completion of $X$ at $C$} to be the algebraic substack $X^{\wedge_C}\subset X$ whose values at $\Spec A$ are given by 
\[
X^{\wedge_C}(\Spec A) = \begin{cases}
 X(A) & \mbox{ if } |\Spec A|\to C\subset |X| \\
 \emptyset & \mbox{otherwise}. 
\end{cases}
\]
\end{definition}

\begin{example}\label{ExampleOpen}
Let $X$ be a scheme and $U\subset X$ an open subscheme, then $X^{\wedge_U}=U$.
\end{example}

In the following we will give a more explicit description of the formal completion of a scheme along a Zariski closed subspace, at least under some finite generation condition of the ideal.

\begin{proposition}\label{PropFormalScheme}
Let $A$ be an animated ring with $X=\Spec A$ and let $Z\subset |\Spec A|$ be a  Zariski closed subspace associated to a finitely generated ideal $I\subset \pi_0(A)$, generated by elements $f_1,\ldots, f_n$. Then the natural map 
\begin{equation}\label{eqwpkqpwemfqwef}
\varinjlim_{k} \Spec (A/^{\bb{L}}(f_1^k,\ldots, f_n^k)) \to X^{\wedge_Z} 
\end{equation}
is an equivalence of algebraic stacks.
\end{proposition}
\begin{proof}
 It is clear from the definition that $X^{\wedge_Z}$ is the intersection of the stacks $X^{\wedge_{Z_i}}$ where $Z_i\subset |X_i|$ is the closed subset given by $\{f_i=0\}$ for $i=1,\ldots, n$. Hence, by an inductive argument, we can suppose without loss of generality that $I$ is generated by a single element $f$. In that case, $X^{\wedge_I}$ is the complement of the Zariski localization $\Spec B \to \Spec B[\frac{1}{f}]$. Let $\Spec C\to \Spec B$ be a map of affine schemes that factors through $X^{\wedge_I}$, then $C\otimes_{B} B[\frac{1}{f}]=0$ and there is some $k$ such that $f^k=0$ in $C$. Thus, we find a map  $B/^{\bb{L}} f^k \to C$ proving that \eqref{eqwpkqpwemfqwef} is actually surjective.   Now we prove that it is injective, by base change we can reduce to the universal case 
\[
\varinjlim_{k} \Spec \bb{Z}[X]/X^n \to \widehat{\bb{A}}^1_{\Z}
\]
where $\widehat{\bb{A}}^1_{\Z}$ is the formal completion of the affine line at $0$. It suffices to see that the map  $\varinjlim_{k} \Spec \bb{Z}[X]/X^n\to \Spec (\bb{Z}[X])$ is an immersion of prestacks. By taking diagonals,  one has to show that the  map 
\[
\varinjlim_{n} \Spec \bb{Z}[X]/X^{n} \to \varinjlim_{n} \bb{Z}[X,Y]/^{\bb{L}}(X^n,Y^n,X-Y)
\]
is an equivalence of pre-stacks. This reduces the problem to showing that the map of pro-sequences of rings
\begin{equation}\label{eqpkqpemfqwsf}
(\bb{Z}[X,Y]/^{\bb{L}}(X^n,Y^n,X-Y))_n \to (\bb{Z}[X]/X^n)_n
\end{equation}
is a pro-equivalence. The animated ring $\bb{Z}[X,Y]/^{\bb{L}}(X^n,Y^n,X-Y)$, seen as a $\bb{Z}[X]$-algebra, is a square-zero extension 
\[
\bb{Z}[X]/X^n \oplus     X^n\bb{Z}[X]/X^{2n} [1],
\]
where we identify $X^n\bb{Z}[X]/X^{2n}$ with the conormal sheaf of $\Spec \bb{Z}[X]/X^n\to \Spec \bb{Z}[X]$.  Therefore, the map $\bb{Z}[X,Y]/^{\bb{L}}(X^{2n},Y^{2n},X-Y)\to \bb{Z}[X,Y]/^{\bb{L}}(X^{n},Y^{n},X-Y)$ is given  by the map of square zero extensions 
\begin{equation}
\label{eqapekpasnroqw}
\bb{Z}[X]/X^{2n} \oplus X^{2n} \bb{Z}[X]/X^{4n} [1] \to \bb{Z}[X]/X^{n} \oplus X^{n} \bb{Z}[X]/X^{2n} [1] 
\end{equation}
which sends $ X^{2n} \bb{Z}[X]/X^{4n}$ to zero in $ X^{n} \bb{Z}[X]/X^{2n}$ in the $\pi_1$-groups. Thus, the map \eqref{eqapekpasnroqw} factors through 
\[
\bb{Z}[X]/X^{2n} \oplus X^{2n} \bb{Z}[X]/X^{4n} [1] \to \bb{Z}[X]/X^{2n}\to \bb{Z}[X]/X^n \to \bb{Z}[X]/X^{n} \oplus X^{n} \bb{Z}[X]/X^{2n} [1],
\]
proving that \eqref{eqpkqpemfqwsf} is a pro-equivalence as wanted.  
\end{proof}

The following lemma describes in which situations the map from the closed subscheme to the formal completion is an epimorphism in the $!$-topology.

\begin{lemma}\label{LemmaDescentFormalCompletion}
Let $A\to B$ be a surjective map of animated rings such that $I=\ker(\pi_0(A)\to \pi_0(A))$ is finitely generated by elements $f_1,\ldots, f_n$. Let us denote $X=\Spec(A)$ and $Z=\Spec(B)$.  Consider the composite map $A\to A/^{\bb{L}}(f_1,\ldots, f_n)\to B$. Then $Z\to X^{\wedge_{Z}}$ is an epimorphism in the $!$-topology if and only if $A/^{\bb{L}} (f_1,\ldots, f_n)\to B$ is descendable.
\end{lemma}
\begin{proof}
By \Cref{PropFormalScheme} the map $\Spec(A/^{\bb{L}} (f_1,\ldots, f_n))\to X^{\wedge_Z}$ is an epimorphism for the $!$-topology. Then the map $\Spec(B)\to \Spec(A/^{\bb{L}} (f_1,\ldots, f_n)) \to X^{\wedge_Z}$ is an epimorphism if and only if there is a descendable map $B\to C$ such that $A/^{\bb{L}}(f_1,\ldots, f_n) \to B \to C$ is descendable, and this happens if and only if $A/^{\bb{L}}(f_1,\ldots, f_n)\to B$ is descendable itself. 
\end{proof}

\begin{corollary}\label{CoroDescendable}
Let $A$ be a static  ring and $A\to B$ a quotient map of static rings by a finitely geneated ideal $I$, let  $Z\to X$ be the associated Zariski closed map of affine schemes. Then $Z\to X^{\wedge_Z}$ is an epimorphism of algebraic stacks.  
\end{corollary}
\begin{proof}
Let $f_1,\ldots, f_n\in I$ be a set of generators of $I$, we have a natural map of $A$-algebras $F\colon A/^{\bb{L}}(f_1,\ldots, f_n)\to B$ which is an isomorphism on $\pi_0$. Since $A/^{\bb{L}}(f_1,\ldots, f_n)$ has finitely many homotopy groups, $F$ is descendable by \cite[Proposition 3.34]{MathewDescent}.  By \Cref{LemmaDescentFormalCompletion} the map $Z\to X^{\wedge_I}$ is an epimorphism.
\end{proof}

\begin{lemma}\label{LemmaFormalScheme}
Let $X$ be a scheme and let $Z\subset |X|$ a Zariski closed subspace, locally in the Zariski topology generated by a finitely generated ideal $I$ . Then $j\colon X^{\wedge_Z}\to X$ is $!$-able and it is an open immersion of algebraic (or analytic) stacks (see \Cref{DefOpenClosed}). Furthermore, the fully faithful functor $j_*\colon \ob{D}(X^{\wedge_Z})\to \ob{D}(Z)$ has   essential image those sheaves on $\ob{D}(X)$ which are derived $I$-adically complete. Moreover, if $X=\Spec A$ and $Z=\Spec A/^{\bb{L}}(f_1,\ldots, f_n)$, the map $Z\to X^{\wedge_Z}$ is a suave cover. 
\end{lemma}
\begin{proof}
These statements are local in the Zariski topology of $X$ (by Zariski descent of quasi-coherent sheaves). Thus, we can assume without loss of generality that $X=\Spec A$ and that $I=(f_1,\ldots, f_n)$ with $f_i\in \pi_0(A)$. By taking intersections, we can even assume that $I=(f)$ is generated by a single element and that $Z=\Spec A/^{\bb{L}} (f)$. To show that $X^{\wedge_Z} \to X$ is $!$-able, it suffices to show that $Z\to X^{\wedge_Z}$ is a suave cover (\cite[Lemma 4.7.1]{HeyerMannSix}). For that, as $Z\to X^{\wedge_Z}$ is an epimorphism by \Cref{PropFormalScheme}, it suffices to show that $Z\times_{X^{\wedge_Z}} Z\to Z$ is a  suave cover (\cite[Lemma 4.5.7]{HeyerMannSix}). Since $X^{\wedge_Z}\to X$ is an immersion, one has that 
\[
Z\times_{X^{\wedge_Z}} Z = Z\times_X Z \cong  \Spec  (A/^{\bb{L}} (f)\oplus  (A/^{\bb{L}} (f))[1]). 
\] 
In particular, the projections $\ob{pr}_i\colon  Z \times_{X^{\wedge_Z}} Z\to Z$ are correpresented by a morphism of rings $B\to C$ such that $C$ is a perfect $B$-module. If follows from \Cref{LemmaPrimAffine}   that the projections $\ob{pr}_i$ are suave,  and then so is $Z\to X^{\wedge_Z}$, since this last map is an epimorphism it is in addition a suave cover. By \cite[Lemma 4.5.7]{HeyerMannSix},  the map $Z\to X^{\wedge_Z}$ satisfies universal $\ob{D}^*$ and $\ob{D}^!$-descent, and as the composite $Z\to X$ is $!$-able, the map $j\colon X^{\wedge_Z}\to X$ is also $!$-able and suave. Since $j$ is also an immersion, the diagonal of $j$ is an isomorphism and $j$ is even cohomologically \'etale by \cite[Definition 4.6.1]{HeyerMannSix} and \cite[Lemma 4.6.4]{HeyerMannSix} implies that $j_!$ is the left adjoint of $j^*$. Then, as $j$ is an immersion, proper base change shows that $j_!\colon \ob{D}(X^{\wedge_Z})\to \ob{D}(X)$ is fully faithful, proving that $j$ is an open immersion of algebraic stacks. It is left to see that the  essential image of $j_*$ are the derived $f$-adically complete modules on $A$. This follows easily from the fact that the pullback along $Z\to X^{\wedge_Z}$ is conservative: given $M$ an $A$-module, the previous implies that  $j_*j^* M$ is $f$-adically complete, and that the natural map $M^{\wedge_f}\to j_*j^*M$ is an isomorphism after moding-out by $f$. By derived Nakayama's lemma one has that $M^{\wedge_f}= j_*j^*M$, which yields 
\[
j_*\colon \ob{D}(X^{\wedge_Z})\xrightarrow{\sim} \ob{D}(A)^{\wedge_f} \subset \ob{D}(A)
\]
as wanted. 
\end{proof}

Thanks to \Cref{PropFormalScheme} we can construct six functors for morphisms of formal schemes under suitable finiteness assumptions on the ideals of definition. As a special case that will be useful later we prove the following suaveness for the affine line:

\begin{lemma}\label{LemmaSuaveAffineLine}
Let $f\colon \widehat{\bb{A}}^1_{\Z}\to \Spec \Z$ be the formal completion at zero of the affine line over $\Z$. Then $f$ is suave and there is an equivalence $f^! \Z \cong \s{O}[1]$, where $\s{O}$ is the stuctural sheaf of $\widehat{\bb{A}}^1_{\Z}$. 
\end{lemma}
\begin{proof}
By \Cref{LemmaFormalScheme} the map $\iota\colon \Spec \Z \to \widehat{\bb{A}}^1_{\Z}$ is a suave cover. Since $\Spec \Z$ is suave over itself, it follows that $f$ is suave by \cite[Lemma 4.5.8 (i)]{HeyerMannSix}. To compute the dualizing sheaf, note that 
 \[
 \Z = \iota^! f^! \Z = \iota^! \Z \otimes \iota^* f^! \Z.
 \]
 Hence, $\iota^* f^!\Z = (\iota^! \Z)^{-1}$ and, by derived Nakayama's lemma, it suffices to see that $\iota^! \Z \cong \Z [-1]$. Since $\widehat{\bb{A}}^1_{\Z}\subset \bb{A}^1_{\Z}$ is an open immersion, one has that 
 \[
 \iota^!\Z = \Hom_{\Z[T]}(\Z, \Z[T]) \cong \Z[-1]
 \] 
 proving what we wanted. 
\end{proof}

\subsection{Betti stacks} \label{ssBettiStacks}

We finish this section with the realization of (light) condensed anima in analytic stacks via the so called \textit{Betti stacks}. This realization factors through that of \Cref{ssAlgebraicStacks}, and its $6$-functor formalism  consists on the (sheafification of the)  base change along $\ob{D}(\Z)\to \ob{D}(\Z^{\cond})$ of the $6$-functor formalism of condensed anima of \cite{HeyerMannSix} with $\Z$-coefficients (or equivalently, $\ob{D}^{\cond}=\Fun_{\ob{D}(\Z)}(\ob{D}_!, \ob{D}(\Z^{\cond}))$, see \Cref{RemarkRelationCondensedUsualCategory}).

Let $\Prof^{\light}$ be the site of light profinite sets. Any surjection $S'\to S$ in $\Prof^{\light}$ gives rise to a faithfully flat morphism of rings $C(S,\Z)\to C(S',\Z)$ of locally constant functions (eq. by writing $S'\to S$ as a cofiltered limit of surjections of finite sets). Therefore, the inclusion 
\[
\Prof^{\light}\to \Sch^{\aff}
\]
sends hypercovers to $\omega_1$-fpqc hypercovers, and it extends to a left adjoint in a morphism of topoi
\[
(-)_{\ob{Betti}}\colon \CondAni= \widehat{\Shv}(\Prof^{\light})\to \widehat{\Shv}_{\omega_1-\ob{fpqc}} (\Sch^{\aff})
\]
that we might call the \textit{Betti realization}. By an abuse of notation we denote in the same way the composite  of the Betti realization with the functor $(-)^{\ob{alg}}$ of \Cref{CoroHypersheaves}
\[
(-)_{\ob{Betti}}\colon \CondAni\to \widehat{\Shv}_{\omega_1-\ob{fpqc}} (\Sch^{\aff})\to \AnStk.
\]

Composing the Betti realization with the functor of quasi-coherent sheaves in classical algebraic geometry $\ob{D}\colon \widehat{\Shv}_{\omega_1-\ob{fpqc}} (\Sch^{\aff})\to \Cat_{\infty}$ we recover the $6$-functor formalism of (light) condensed anima of \cite{HeyerMannSix} with $\Z$-coefficients. If instead we compose the Betti realization with the functor $\ob{D}\colon \AnStk\to \Cat_{\infty}$ of quasi-coherent sheaves for analytic stacks, we obtained a $\ob{D}(\Z^{\cond})$-linear variant of the six functors of \cite{HeyerMannSix}. 

\section{Solid geometry and Serre duality}\label{SectionSerreDualitySolid}

In this last section we specialize the theory of analytic stacks to a suitable category of \textit{solid Huber stacks} whose building blocks are generalizations of Huber pairs, the so called \textit{solid Huber rings}. In this framework, we discuss variants of the theory of Huber spaces, a suitable notion of morphisms of finite presentation in the solid realm, \'etale and smooth morphisms in the solid set up, and finally  a general statement of Serre duality that extends that of proper smooth schemes in classical algebraic geometry.

\subsection{Solid Huber rings}\label{ssSolidRing}

In \cite{HuberAdicSpaces}, Huber introduced the theory of adic spaces where the building blocks are the so called \textit{Huber pairs}. In \cite[Proposition 3.34]{Andreychev} Andreychev proved that complete Huber pairs  embed fully faithfully in the category of analytic rings. In  \cite[Definition 2.9.3]{MannSix} Mann introduced the category of (animated) discrete Huber pairs and proved in \cite[Proposition 2.9.6]{MannSix} that it embeds fully faithfully into the category of analytic rings. In this section we generalize these two theories of Huber pairs into a single theory of \textit{solid Huber rings}.

\begin{definition}
\begin{enumerate}

\item A   \textit{solid pre-analytic Huber pair} is a pair $(A^{\triangleright},S)$ consisting of a solid animated ring $A^{\triangleright}\in \Ring_{\bb{Z}_{\sol}}$ and a morphism of anima $S\to A^{\triangleright}(*)$. 

\item  Given $(A^{\triangleright},S)$ a pre-analytic Huber pair, we define  the pre-analytic ring $(A^{\triangleright},S)^{\mathrm{pre}}_{\sol}$ to be the pre-analytic ring structure on $A^{\triangleright}$ making an $A^{\triangleright}$-module $M$ complete if and only if it is $\bb{Z}$-solid, and for all $s\in \pi_0(S)$, the restriction of $M$ as $\bb{Z}[T]$-module (with action of $T$ given by multiplication by the image of $s$ along $S\to A^{\triangleright}(*)$) is $\bb{Z}[T]_{\sol}$-complete. We let $(A^{\triangleright},S)_{\sol}$ denote the completion of $(A^{\triangleright},S)^{\mathrm{pre}}_{\sol}$ to an analytic ring.

\item Let $\mathbb{G}_a=\AnSpec (\bb{Z}[T] )$  be the algebraic affine line with the trivial analytic ring structure. Let $A$ be an analytic ring over $\bb{Z}_{\sol}$, we define the ring of \textit{solid elements of $A$} to be the subring of $A^{\triangleright}(*)= \mathbb{G}_a(A)$ given by 
\[
A^+:=\mathbb{G}_{a,\sol}(A)
\]
where $\mathbb{G}_{a,\sol}=\AnSpec(\bb{Z}[T]_{\sol})$. 

\item We say that a solid analytic ring $A\in \AnRing_{\bb{Z}_{\sol}/}$ is a \textit{solid Huber ring} if the natural map 
\[
(A^{\triangleright}, A^+)_{\sol} \to A
\]
is an isomorphism of analytic rings.   We let $\Ring^{\Hub}_{\bb{Z}_{\sol}}\subset \AnRing_{\bb{Z}_{\sol}/}$ be the full subcategory of solid Huber rings. 

\item Finally, a solid pre-analytic Huber pair $(A^{\triangleright},S)$ (with associated solid ring $A=(A^{\triangleright},S)_{\sol}$) is said a \textit{solid Huber pair} if $A^{\triangleright}$ is $A^+_{\sol}$-complete and the natural map $S\to A^+$ is an isomorphism of anima. 

\end{enumerate}

\end{definition}

\begin{remark}\label{RemRinsStackSolid}
Let $\Ring$ be the category of discrete animated rings and consider the functor $\Ring^{\op} \to \AnStk_{\bb{Z}_{\sol}}$ mapping $R\mapsto \Spec(R_{\sol})$. This functor preserves limits and therefore defines a ring stack $\mathbb{G}_{a,\sol}$ in $\AnStk_{\bb{Z}_{\sol}}$ whose underling stack is $\AnSpec \bb{Z}_{\sol}[T]$. Thus, given $X\in \AnStk_{\bb{Z}_{\sol}}$ a solid analytic stack, the anima of maps $\mathbb{G}_{a,\sol}(X)$ has a natural structure of animated ring. In particular, if $X=\AnSpec A$, $\mathbb{G}_{a,\sol}(A)= A^+$ has a natural structure of animated ring.  
\end{remark}

\begin{remark}\label{RemarkNoDependenceonS}
Let $(A^{\triangleright},S)$ be a pre-analytic Huber pair. The solid analytic ring $(A^{\triangleright}, S)^{\mathrm{pre}}_{\sol}$ only depends on $A^{\triangleright}$ and in the image of $\pi_0(S)\to \pi_0(A^{\triangleright}(*))$.  The reason to use an anima $S$ instead of simply a subset of  $\pi_0(A^{\triangleright}(*))$ is just to allow more flexibility in the notation. 
\end{remark}

\begin{example}\label{ExampleAdicRings}
Let $\Ring^{\Hub}_{\bb{Z}}$ be the category of discrete Huber pairs \cite[Definition 2.9.3]{MannSix},  the fully faithful functor $\Ring^{\Hub}_{\bb{Z}}\to \AnRing$ mapping $(A,A^+)\mapsto (A,A^+)_{\sol}$  factors through $\Ring^{\Hub}_{\bb{Z}_{\sol}}$. Similarly, if $\Ring^{\Hub,\wedge}$ is the category of complete Huber pairs, the natural functor to $\AnRing$ given by $(A,A^+)\mapsto (A,A^+)_{\sol}$ factors through $\Ring^{\Hub}_{\bb{Z}_{\sol}}$ (\cite[Proposition 3.34]{Andreychev}). 
 \end{example}
 
The following proposition proves some basic categorical properties of the category of solid Huber rings. 

\begin{proposition}\label{PropColimits}
The fully faithful functor $\Ring^{\Hub}_{\bb{Z}_{\sol}}\to \AnRing_{\bb{Z}_{\sol}}$ admits a right adjoint that sends a solid analytic ring $A$ to the Huber ring $(A^{\triangleright}, A^+)_{\sol}$. In particular, $\Ring^{\Hub}_{\bb{Z}_{\sol}}$ is stable under all colimits in $\AnRing_{\bb{Z}_{\sol}}$. Furthermore, the following holds: 
\begin{enumerate}

\item Let $S$ be a light profinite set and consider the Huber ring  given by $R=\Sym_{\bb{Z}_{\sol}} \bb{Z}_{\sol}[S]$ with the induced structure from $\bb{Z}_{\sol}$. Then $R$ is compact projective  in $\AnRing_{\bb{Z}_{\sol}}$ and so it is in $\Ring^{\Hub}_{\bb{Z}_{\sol}}$.

\item The rings $\mathbb{Z}[T]_{\sol}$ and $\bb{Z}[[T]]$ are compact (but not projective!) in $\AnRing_{\bb{Z}_{\sol}}$ and  they are in $\Ring^{\Hub}_{\bb{Z}_{\sol}}$.

\item Finite colimits of the rings of (1) and $\bb{Z}[T]_{\sol}$ form a family of compact generators of $\Ring^{\Hub}_{\bb{Z}_{\sol}}$.  In particular, $\Ring^{\Hub}_{\bb{Z}_{\sol}}$ is a presentable category. 

\end{enumerate}
\end{proposition}
\begin{proof}
The functor $\Ring_{\Z}^{\Hub}\to \AnRing_{\Z_{\sol}}$ has a retract given by the functor sending an analytic ring $A$ to the analytic ring $(A^{\triangleright},A^+)_{\sol}$, it is an easy computation that this functor defines a right adjoint for the inclusion. 

The fact that $\Sym_{\Z_{\sol}} \Z_{\sol}[S]$ is compact projective in solid analytic rings follows from the fact that it correpresents the functor $A\mapsto A^{\triangleright}(S)$ (with $S$ light profinite), and this functor commutes with sifted colimits as $\bb{Z}_{\sol}[S]$ are compact projective solid modules, see  \Cref{LemmaSiftedAnalytic}. 

 Next, we show that $\Z_{\sol}[T]$ and $\Z[[T]]$ are compact in animated rings. To see this, by \Cref{PropSolidZTidempotentComplement} a map $(\Z[T],\Z)_{\sol}\to A$ of solid analytic rings factors through $\Z[T]_{\sol}$ (resp. through $\Z[[T]]$) if and only if $\Z((T^{-1}))\otimes_{(\Z[T],\Z)_{\sol}}A =0$ (resp. $(\Z[T^{\pm 1}], \Z[T^{-1}])_{\sol}\otimes_{(\Z[T],\Z)_{\sol}} A =0$).  Indeed, this follows from the fact that $\ob{D}(\Z[T]_{\sol})$ is the open complement of the morphism of idempotent algebras $\Z[T]\to \Z((T^{-1}))$ in $\ob{D}((\Z[T], \Z)_{\sol})$ (resp. $\ob{D}((\Z[T^{\pm 1}], \Z[T^{-1}])_{\sol})$ is the open complement of $\Z[T]\to \Z[[T]]$).  Now, let $A=\varinjlim_i A_i$ be a filtered colimit of analytic rings over $\Z_{\sol}$ and let $\Z[T]\to A$, then there is some $i$ such that there is a lift $\Z[T]\to A_i$. Consider the base change $\Z((T^{-1}))\otimes_{(\Z[T],\Z)_{\sol}} A=\varinjlim_{j\geq i} \Z((T^{-1}))\otimes_{(\Z[T],\Z)_{\sol}} A_j$, then the colimit is zero if and only if there is some $j\geq i$ such that $1\in \Z((T^{-1}))\otimes_{(\Z[T],\Z)_{\sol}} A_j$ is zero, if and only if $\Z((T^{-1}))\otimes_{(\Z[T],\Z)_{\sol}} A_j=0$, if and only if we have an extension $\Z[T]_{\sol}\to A_j$. The same argument holds for $\Z[[T]]$.    Note that neither $\Z_{\sol}[T]$ nor $\Z[[T]]$ are projective by \Cref{RemNonProjective}. 
 
Finally, to see that the rings of (1) and $\Z[T]_{\sol}$ generate $\Ring^{\Hub}_{\Z_{\sol}}$, it suffices to see that the functor sending an Huber ring $A$ to the pair $(A^{\triangleright},A^+)$ is conservative, where $A^{\triangleright}$ is the condensed ring whose $S$-points are correpresnted by $\Sym_{\Z_{\sol}}\Z_{\sol}[S]$, and $A^+$ is the subring of solid elements correpresented by $\Z[T]_{\sol}$. This is clear from the definition of solid Huber rings. 
\end{proof}

\begin{definition}\label{DefTopNilpElements}
Let $A\in \AnRing_{\bb{Z}_{\sol}}$. We define the anima of \textit{topologically nilpotent elements of $A$} to be the full subanima $A^{\circ\circ}(*)\subset A^{\triangleright}(*)$ of elements $a\in A^{\triangleright}$ such that the map $\bb{Z}[T]\to A$ extends to $\bb{Z}[[T]]\to A$. 
\end{definition}

\begin{lemma}\label{LemmaTopNilpIdeal}
Let $A\in \AnRing_{\bb{Z}_{\sol}}$, then   $A^{\circ\circ}(*)\subset A^+$ is an $A^+$-submodule.
\end{lemma}
\begin{proof}
The map $\Z[T]\to \Z[[T]]$ is idempotent, hence $A^{\circ\circ}(*)\subset A^{\triangleright}$ is a full subanima. Thus, we only need to prove the lemma for the connected components of both objects. The inclusion $A^{\circ\circ}(*)\subset A^+$ is clear since we have a factorization of analytic rings $(\Z[T],\Z)_{\sol}\to  \Z[T]_{\sol} \to \Z[[T]]$. The fact that $A^{\circ\circ}(*)$ is an $A^+$-module follows from the following morphisms of rings:
\[
\Z[[X]]\to \Z[[T_1,T_2]], \;\; X\mapsto T_1+T_2 \mbox{ and } \Z[[X]]\to \Z[T][[X]], \;\; X\mapsto TX. 
\]
Indeed, the second shows that $A^{\circ\circ}(*)$ is stable under multiplication by $A^+$, the first shows that it is stable under addition. 
\end{proof}

\begin{remark}\label{RemTopNilpElements}
In \cite[Definition 2.6.1 (2)]{camargo2024analytic} it is constructed a condensed anima $A^{\circ\circ}\subset A^{\triangleright}$ whose values at a point is precisely \Cref{DefTopNilpElements}. Thanks to \cite[Remark 2.6.3]{camargo2024analytic} $A^{\circ\circ}$ is  a solid $A^+$-module. 
\end{remark}

\begin{definition}\label{DefinitionAffineSpace}
Let $A$ be a solid Huber ring. We define the \textit{solid Tate algebra in $n$-variables} to be the solid Huber ring 
\[
A[T_1,\ldots, T_n]_{\sol}:= A\otimes_{\Z_{\sol}} \Z[T_1,\ldots, T_n]_{\sol}.
\]
\end{definition}

\begin{example}\label{ExampleSolidTate}
Let $A$  be a solid Huber ring. The shape of $A[T]_{\sol}$ highly depends on the condensed structure of $A$:

\begin{itemize}

\item Suppose that $A=(R,R^+)_{\sol}$ arises from  a discrete Huber ring, then $A[T]_{\sol}=(R[T], R^+[T])_{\sol}$ is a polynomial algebra over $R$ with solid structure induced from $R^+[T]$.

\item Suppose that $A=(R,\Z)_{\sol}$ has the induced structure and that $R$ is an $I$-adically complete ring for $I\subset \pi_0(R)$ a finitely generated ideal. Then $R[T]_{\sol}=((R[T])^{\wedge_I},\Z[T])_{\sol}$ is the $I$-adic completion of $R[I]$ endowed with the induced structure from $\Z[T]_{\sol}$. Indeed, this follows from the fact that $\Z[T]_{\sol}\otimes_{\Z_{\sol}}-$ commutes with limits and so with $I$-adic completions. 

\item Suppose that $A=(R,R^+)$  is a Tate Huber pair, that is, $R$ is a Banach algebra with pseudo-uniformizer $\pi$. Then $A[T]_{\sol}=A\langle T\rangle_{\sol}=(R\langle T\rangle, R^+\langle T\rangle)_{\sol}$ is the classical Tate algebra over $(R,R^+)$ in one variable. This follows from  \cite[Proposition 3.14]{Andreychev} as $R$ is a filtered colimit of pro-discrete solid abelian groups.

\end{itemize}
\end{example}

As preparation for the proof of Serre duality, let us show the key local computation for the solid Tate algebra. 

\begin{proposition}\label{PropSerreDualityAffine}
The map of solid rings $f\colon \Z_{\sol}\to \Z[T]_{\sol}$ is cohomologically smooth with dualizing sheaf isomorphic to $\Z[T][1]$.
\end{proposition}
\begin{proof}
Let $\bb{G}_{a,\sol}:=\AnSpec \Z[T]_{\sol}$ and let us denote $S=\AnSpec \Z_{\sol}$. By definition, we have to show that the kernel $\Z[T]\in \ob{D}(\Z[T]_{\sol})=\ob{Fun}_{\n{K}_S}(S, \bb{G}_{a,\sol})$ is a right adjoint. By the fully faithfulness of \Cref{LemFullyFaithfulKernelCat}, it suffices to show that the pullback map  $f^*\colon \ob{D}(\Z_{\sol})\to \ob{D}(\Z[T]_{\sol})$ has a $\bb{Z}_{\sol}$-linear left adjoint $f_{\natural}$. We claim that there is a natural isomorphism of functors $\ob{D}(\Z_{\sol})\to \ob{D}(\Z[T]_{\sol})$ 
\[
f^* = \iHom_{\Z_{\sol}}(\Z[T]^{\vee}, - )
\]
where $\Z[T]^{\vee}=\Hom_{\Z}(\Z[T],\Z)$.  Indeed, the dual of the unit $\Z[T]^{\vee}\to \Z$ produces a natural transformation $M\to \iHom_{\Z}(\Z[T]^{\vee}, M)$ for $M\in \ob{D}(\Z_{\sol})$. We claim that $\iHom_{\Z}(\Z[T]^{\vee}, M)$ is $\Z[T]$-solid, and that the natural map 
\[
f^*M\to \iHom_{\Z_{\sol}}(\Z[T]^{\vee}, M)
\]
is an isomorphism. To prove the claim, notice that  $\Z[T]^{\vee}\cong \prod_{\N} \Z$ is a compact projective solid abelian group, so that the functor $\iHom_{\Z_{\sol}}(\Z[T]^{\vee},-)$ commues with limits and colimits. Thus, it suffices to consider the case $M=\Z$, in which case we have that 
\[
\iHom_{\Z_{\sol}}(\Z[T]^{\vee},\Z)= \Z[T]
\]
as wanted.

Now, knowing that $f^*=\iHom_{\Z_{\sol}}(\Z[T]^{\vee},-)$, we see that its has a left adjoint given by $f_{\natural}N= \Z[T]^{\vee}\otimes_{(\Z[T],\Z)_{\sol}} N$  for $N\in \ob{D}(\Z[T]_{\sol})$ proving suaveness of $f$.

Finally, to compute the dualizing sheaf, let us recall the explicit description of $f_!$.  We have a factorization
\[
\bb{G}_{a,\sol}\xrightarrow{j} \bb{G}_{a}\xrightarrow{p} S
\]
where $\bb{G}_a=\AnSpec (\Z[T],\Z)_{\sol}$. Thus, $f_!=p_*j_!$, the functor $p_*$ is a forgetful functor, and since $j$ is an open immersion with complement idempotent algebra $\Z((T^{-1}))$, the functor $j_!$ is given by 
\[
j_!N = \ob{fib}(\Z[T]\to \Z((T^{-1})))\otimes_{(\Z[T],\Z)_{\sol}} N
\]
for $N\in \ob{D}(\Z[T]_{\sol})$. Note that as $\Z[T]$-module, we have an isomorphism $\ob{fib}(\Z[T]\to \Z((T^{-1})))\cong \Z[T]^{\vee}[-1]$. Thus, we have an isomorphism of functors $f_!\cong f_{\natural} [-1]$ from $\ob{D}(\Z[T]_{\sol})\to \ob{D}(\Z_{\sol})$. Passing to right adjoints this implies that $f^* [1]\cong f^!$ proving that $f^!\Z\cong \Z[T][1]$ as wanted. 
\end{proof}

\subsection{The solid  spectrum}\label{ss:SolidSpec}

Attached to any solid Huber ring one can construct a solid  spectrum. Its definition in general is quite inexplicit and only exists thanks to formal reasons. Indeed, the adic spectrum will be defined as the spectral space attached to  the locale  generated by  \textit{rational localizations} of a solid Huber ring via Stone duality \Cref{Theo:StoneDualityLattices}. Later in \Cref{PropComparisonSpectrum} and \Cref{ExamplesSolidHuberSpaces}, we will see how this solid spectrum relates to the topological spectrums appearing in algebraic, formal and adic geometry after considering some concrete examples.

\begin{definition}\label{DefSolidAdicSpec}
A morphism of solid Huber rings $A\to B$ is said a \textit{standard rational localization} if there are elements $f_1,\ldots, f_n,g\in A$ generating the unit ideal, and an isomorphism  of $A$-algebras 
\[
B\cong A\bigg(\frac{f_1,\ldots, f_n}{g} \bigg)_{\sol}:= A[\frac{1}{g}] \otimes_{A[T_1,\ldots, T_n]} A[T_1,\ldots, T_n]_{\sol}
\]
where $T_i\mapsto \frac{f_i}{g}\in A[\frac{1}{g}]$.  A morphism of solid Huber rings $A\to B$ is a \textit{rational localization} if it is a finite composition of iterated standard rational localizations. 
\end{definition}

\begin{remark}\label{RemarkRationalLocalizationsNotStandard}
We highlight that as opposed to the classical theory, the composition of two standard rational localizations of solid Huber rings is not in an obvious way a standard rational localization. For instance, if $A=\big(\Sym_{\Z_{\sol}} \Z_{\sol}[S]\big)[T]$   with $S$ an infinite light profinite set, given any element $f\in \Z[T]_{\sol}[S]\in  A\left(\frac{T}{1}\right)_{\sol}=:A'$ that has not bounded degree on $T$, the composite of rational localizations $A\to A'\to A'\left( \frac{1}{f} \right)_{\sol}$ should not be standard (the heuristic reason behind is that it is impossible to approximate $f$ by elements in $A$ by elements which are more and more topologically nilpotent). Nevertheless, under some reasonable conditions such as $A$ being discrete, a classical Tate Huber ring, or a \textit{Gelfand ring} in the sense of \cite{dRFF}, rational localizations and standard rational localizations agree, see \Cref{PropComparisonSpectrum} and \Cref{RemarkSolidSpectrumBounded}. 
\end{remark}

\begin{lemma}\label{LemmaSolidSpectrumOpen}
Let $f\colon A\to B$ be a  rational localization, then $f$ gives rise to an open map of analytic stacks.  Moreover,  rational localizations are stable under base change and composition. 
\end{lemma}
\begin{proof}
It suffices to show that standard rational localizations give rise to open maps in solid Huber rings. The map $f$ is clearly idempotent being a composition of idempotent maps, namely, the Zariski localization $A\to A[\frac{1}{g}]$ and the solidification of the elements $\frac{f}{g}$. In particular, $f^*\colon \ob{D}(A)\to \ob{D}(B)$ is a localization. Hence, to show that $f$ is an open immersion it suffices to show that $f^*$ admits a $\ob{D}(A)$-linear left adjoint $f_{\natural}$. To prove that, we claim that the natural map 
\begin{equation}\label{eqa9jwwolqweqef}
C:=A\otimes_{\Z_{\sol}} \Z[T_1,\ldots, T_n]_{\sol} /^{\bb{L}}(gT_i-f_i)\to B
\end{equation}
is an isomorphism. Indeed, the ring $C$ is a composition of suave maps over $A$ thanks to \Cref{PropSerreDualityAffine} and \Cref{LemmaPrimAffine}, and the isomorphism $C\cong B$  yields the existence of the linear left adjoint of $f^*$. Next we show the claim. Since $(f_1,\ldots, f_n,g)$ generate the unit ideal, we have a linear combination $\sum_i a_i f_i +b g = 1$. Thus, in the ring $C$, we have the equation  $\sum_i  a_i gT_i +bg=1$ which yields $g(\sum_i a_iT_i+b)=1$. In particular, $g\in C$ is invertible, and therefore $C=C[\frac{1}{g}]=B$. 

The stability of  rational localizations is clear from the definition. 
\end{proof}

\begin{lemma}\label{LemmaCompactOpenLocalization}
Let $A\to B$ be a rational localization of solid Huber rings and let $D\in \ob{D}(A)$ be the idempotent algebra over $A$ complement of $B$. Then $D$ is compact as $\ob{D}(A)$-module.
\end{lemma}
\begin{proof}
This follows from the fact that  the forgetful functor $\ob{D}(B)\to \ob{D}(A)$ commutes with colimits and the unit of $\ob{D}(A)$ is compact. 
\end{proof}

\begin{lemma}\label{LemmaCompact}
Let $A\to B \to C$ be a triangle of analytic rings. The following holds:
\begin{enumerate}

 \item Compact morphisms of analytic rings are stable under finite colimits and  are stable under  base change.

\item   Suppose that $B$ is a compact $A$-algebra. Then $C$ is a compact $A$-algebra if and only if it is a compact $B$-algebra.

\end{enumerate}

\end{lemma}
\begin{proof}
The fact that compact morphisms of analyitc rings are stable under base change and finite colimits is obvious and left to the reader. Suppose that $B$ is compact as $A$-algebra and suppose that $C$ is compact as $B$-algebra. Let $D=\varinjlim_i D_i$ be a filtered colimit of $A$-algebras, we want to see that the natural map of anima
\[
\varinjlim_i \Map_{\AnRing_{A/}}(C, D_i)\to \Map_{\AnRing_{A/}}(C, D)  
\]
is an equivalence. Consider the map of functors $f\colon \Map_{\AnRing_{A/}}(C,-)\to \Map_{\AnRing_{A/}}(B,-)$, the target commute with filtered colimits. Hence, to show that the source commutes with filtered colimits it suffices to show that the fibers do. Hence, we can assume without loss of generality that the $D_i$ are $B$-algebras, in which case the fiber of $f$ at the point $B\to D$ is precisely $\Map_{\AnRing_{B/}}(C,-)$ which commutes with colimits by hypothesis. 

Conversely, suppose that $C$ is compact as $A$-algebra, and let $D=\varinjlim_i D_i$ be a filtered colimit of $B$-algebras. We have a cartesian square
\[
\begin{tikzcd}
\Map_{\AnRing_{B/}}(C, D_i) \ar[r] \ar[d] & \Map_{\AnRing_{A/}}(C,D_i) \ar[d] \\ 
* \ar[r]& \Map_{\AnRing_{A/}}(B,D_i)
\end{tikzcd}
\] 
where the lower horizontal map is the $B$-algebra structure of $D_i$. Since filtered colimits commute with pullbacks and $C$ and $B$ are compact as $A$-algebras, we deduce that $\varinjlim_i \Map_{\AnRing_{B/}}(C, D_i)=\Map_{\AnRing_{B/}}(C, D)$ proving that $C$ is compact as $B$ algebra as wanted. 
\end{proof}

\begin{lemma}\label{LemmaCompactRational}
Let $A\to B$ be a rational localization of solid Huber rings, then $B$ is a compact $A$-algebra.
\end{lemma}
\begin{proof}
By \Cref{LemmaCompact} it suffices to show that a standard rational localization of $A$ is a compact $A$-algebra. This follows from the presentation \eqref{eqa9jwwolqweqef} of a standard rational $A$-algebra,  the compactness of the solid ring $\Z[T]_{\sol}$ of \Cref{PropColimits} (2), and the compactness of the polynomial algebra $\Z[T]$. 
\end{proof}

\begin{definition}\label{DefSolidSpectrum}
Let $A$ be a solid Huber ring and let $\ob{Sm}(\ob{D}(A))$ be its smashing spectrum.  Consider the full subcategory $\s{C}_A\subset\ob{Sm}(\ob{D}(A)) $ generated under finite limits and colimits by rational localizations. Then $\s{C}_A$ is  a distributive lattice as in \Cref{def:DistributiveLattice},  and by Stone duality (\Cref{Theo:StoneDualityLattices}) $\s{C}_A$ is the category of quasi-compact open subspaces  of a  spectral space $\Spec^{\sol}(A)$ (unique up to unique spectral isomorphism) that we call the \textit{solid spectrum of $A$}.  Let $X=\Spec^{\sol} (A)$, given $(f_1,\ldots, f_n,g)$ a sequence in $A$ generating the unit ideal we let $X(\frac{f_1,\ldots, f_n}{g})$ or $\{|f_i|\leq |g|\neq 0\}$ denote the  quasi-compact open subspace of $X$ associated to the standard rational localization $A\to A\big( \frac{f_1,\ldots, f_n}{g} \big)$. 
\end{definition}

\begin{remark}\label{RemBaseChange}
 Given $A\to B$  a morphism of solid rings, there is a natural morphism of smashing spectrum $\ob{Sm}(\ob{D}( B))\to \ob{Sm}(\ob{D}( A))$. Since base change preserves rational localizations, Stone duality produces a natural commutative square of spectral spaces
\[
\begin{tikzcd}
\ob{Sm}(\ob{D}( B))\ar[r] \ar[d] & \Spec^{\sol} B \ar[d] \\
\ob{Sm}(\ob{D}( A))\ar[r] &  \Spec^{\sol} A. 
\end{tikzcd}
\]
Moreover, since the pullback along $F\colon \ob{Sm}(\ob{D}( B))\to \Spec^{\sol} B$ preserves and induces a monomorphism of quasi-compact open subspaces, the map $F$ is a surjection of locales. 
\end{remark}

The formation of the solid spectrum is  compatible with filtered colimits. 

\begin{lemma}\label{LemmaSpectrumColimits}
Let $A=\varinjlim_i A_i$ be a filtered colimit of solid Huber rings. Then the natural map $\Spec^{\sol}(A)\to \varprojlim_i \Spec^{\sol}(A_i)$ is an homeomorphism of spectral spaces.  Furthermore, if $A\to B$ is a rational localization, there is some index $i$ and a rational localization $A_i\to B_i$ such that $B=B_i\otimes_{A_i} A$.
\end{lemma}
\begin{proof}
Let $\ob{Op}^{\ob{rat}}(A_i)$ be the poset of rational open subspaces of $\Spec^{\sol}(A_i)$. Using Stone duality, it suffices to show that we have an equivalence 
\[
\varinjlim_{i} \ob{Op}^{\ob{rat}}(A_i) \xrightarrow{\sim} \ob{Op}^{\ob{rat}}(A)
\]
of posets given by base change. For that, it suffices to show that any rational localization $A\to B$ arises from a rational localization of some $A_i$ (as claimed in the lemma), and that given a ring $A_i$ and two rational localizations $A_i\to B_{i,1}, B_{i,2}$ such that $B_{i,1}\otimes_{A_i} A = B_{i,2}\otimes_{A_i} A$, there is some $j\geq i$ such that $B_{i,1}\otimes_{A_i} A_j=B_{i,2}\otimes_{A_i} A_j$. 

By an inductive argument it suffices to show  the first claim for standard rational localizations, if $A\to A(\frac{f_1,\ldots, f_n}{g})$ is a standard rational localization, there is some index $i$ and elements $(f_1',\ldots, f_n',g')$ in $A_i$ mapping to $(f_1,\ldots, f_n,g)$ and generating the unit ideal, then $A_i\to A_i(\frac{f_1',\ldots, f_n'}{g})$ does the job. 

On the other hand, by \Cref{LemmaCompactRational} rational localizations are compact morphisms. Suppose that  $A_i\to B_{i,1}, B_{i,2}$ are two rational localizations that agree after base change to $A$. Then we have an isomorphism 
\[
\varinjlim_{j\geq i} B_{i,1}\otimes_{A_i} A_j = \varinjlim_{j\geq i} B_{j,2} \otimes_{A_i} A_j. 
\]
By compacity of the $B_{i,k}$, there is some $j\geq i$ and maps of $A_i$-algebras   $B_{i,1}\to B_{i,2}\otimes_{A_{i}} A_{j}$ and $B_{i,2}\to B_{i,1}\otimes_{A_{i}} A_{j}$, since the $B_{i,k}$ are idempotent analytic $A_i$-algebras one has that $ B_{i,1}\otimes_{A_{i}} A_{j}= B_{i,2}\otimes_{A_{i}} A_{j}$ proving what we wanted. 
\end{proof}

The solid spectrum does not depend on the animated structure of the ring and on nilpotent ideals. 

\begin{lemma}\label{LemmaInvarianceNilpotents}
Let $A\to A'$ be a morphism of solid Huber rings such that $\pi_0(A)\to \pi_0(A')$ has nilpotent kernel. Then the natural map $\Spec^{\sol}(A')\to \Spec^{\sol}(A)$ is an homeomorphism identifying (standard) rational subspaces. 
\end{lemma}
\begin{proof}
A rational localization $A\to B$ is induced by an uncompleted analytic ring structure on $A^{\triangleright}$. The lemma follows by the independence of uncompleted analyitc ring structures on higher homotopy groups (\Cref{TheoAnRingStructure}) and nilpotent ideals (\Cref{PropInvarianceNil}).
\end{proof}

Under some adically complete conditions, the solid spectrum only depends on the special fiber.

\begin{lemma}\label{LemmaInvarianceAdicallyComplete}
Let $A$ be a solid analytic ring and let $I\subset \pi_0(A)$ be a finitely generated ideal such that $A$ is $I$-adically complete and either one of the following conditions hold:
\begin{enumerate}

\item[(a)] There is a finite subset $S\subset A^+$ such that $A=(A^{\triangleright}, S)_{\sol}$.

\item[(b)] The quotient $\pi_0(A)/I$ is discrete. 

\end{enumerate}
 Let $(f_1,\ldots, f_n)$ be a set of generators of $I$. Then the natural map of solid spectrum
\[
\Spec^{\sol}(A/^{\bb{L}}(f_1,\ldots, f_n))\xrightarrow{\sim} \Spec^{\sol}(A)
\]
is an equivalence identifying rational localizations.  
\end{lemma}
\begin{proof}
By an induction step we can assume without loss of generality that $I=(f)$ is generated by a single element. Let $B=A/^{\bb{L}}(f)$ and denote $F\colon A\to B$. For a solid Huber ring $C$ let $\ob{Op}^{\ob{rat}}(C)$ be the poset of rational localizations of $\Spec^{\sol}(C)$. We want to show that the natural map 
\[
F^{-1}\colon \ob{Op}^{\ob{\ob{rat}}}(A)\to \ob{Op}^{\ob{rat}}(B)
\]
is an equivalence. We first show that $F^{-1}$ is surjective. By an inductive step, and since rational localizations preserve $I$-adically complete modules being an open immersion, it suffices to show that it preserves standard rational localizations. Let $(f_1,\ldots, f_n,g)$ be elements in $B$  generating the unit ideal, and let $(\widetilde{f}_1,\ldots, \widetilde{f}_n, \widetilde{g})$ be any lift to $A$. Since $A$ is $f$-adically complete, the tuple $(\widetilde{f}_1,\ldots,\widetilde{f}_n, \widetilde{g})$ generates the unit ideal in $A$ and $A\to A\big(\frac{\widetilde{f}_1,\ldots, \widetilde{f}_n}{\widetilde{g}}\big)$ is a rational localization whose base change to $B$ is $B\big(\frac{f_1,\ldots, f_n}{g}\big)$,  proving the surjectivity. 

Now, to prove that $F^{-1}$ is injective, by a devisage to standard rational localizations, it suffices to prove the following facts: 

\begin{itemize}
\item[i.]  An element $g\in A$ is invertible if and only if its image in $B$ is invertible. 

\item[ii.]  An element $g\in A$ is solid if and only if its image in $B$ is solid. 

\end{itemize}

The property (i) is clear since $A$ is $f$-adically complete. We now prove the condition (ii). Let $g\in A$ be such that its image to $B$ is solid. We want to see that $g\in A^+$, or equivalently, that $\Z((T^{-1}))\otimes_{(Z[T], \Z)_{\sol}} A = 0$ where $T\mapsto g$. Note that we have $A\otimes_{(\Z[T],\Z)_{\sol}} \Z((T^{-1}))= (A\otimes_{\Z_{\sol}} \Z[[T]])/^{\bb{L}}(gT-1)$ and that $B\otimes_{(\Z[T],\Z)_{\sol}} \Z((T^{-1}))=(B\otimes_{\Z_{\sol}} \Z[[T]])/^{\bb{L}}(gT-1)=0$ by assumption. 

Suppose that the condition (a) holds, then by \cite[Lemma 2.12.9]{MannSix} the tensor  $A\otimes_{\Z_{\sol}} \Z[[T]]/^{\bb{L}}(gT-1)$ is $I$-adically complete, and by  Nakayama's lemma we get that $A\otimes_{\Z_{\sol}} \Z[[T]]/^{\bb{L}}(gT-1)$ is zero as wanted. If condition (b) holds, then $A\otimes_{\Z_{\sol}} \Z[[T]]/^{\bb{L}}(gT-1)$ is $I$-adically complete by \cite[Proposition 2.12.10 (ii)]{MannSix} and hence zero by Nakayama's lemma. 
\end{proof}

In some special cases the solid spectrum can be understood in terms of more classical spectrum of rings:

\begin{proposition}\label{PropComparisonSpectrum}
The following hold:

\begin{enumerate}

\item Let $(A,A^+)$ be a discrete animated Huber pair. Then we have a natural surjective map of spectral spaces $\Spa(A,A^+)^{\ob{mod}}\to \Spec^{\sol}((A,A^+)_{\sol})$ compatible with the map from the smashing spectrum of \Cref{corodiscreteAdicSpacesAsLocales}.

\item Let $(A,A^+)$ be a sheafy Tate Huber ring. Then there is a natural homeomorphism $\Spec^{\sol}((A,A^+)_{\sol})\xrightarrow{\sim}\Spa(A,A^+)$.

\end{enumerate}
\end{proposition}
\begin{proof}

For part (1), it suffices to notice that rational localizations of $(A,A^+)_{\sol}$ as in \Cref{DefSolidAdicSpec} are particular examples of rational localizations in the usual sense of discrete adic spaces. 

For part (2), since $(A,A^+)$ is Tate, \cite[Proposition 4.3]{Andreychev} says that the family of \textit{nice} rational subspaces of $\Spa(A,A^+)$ (in the sense of \cite[Definition 4.2]{Andreychev}) form a basis for the topology. On the other hand, an iteration of \cite[Proposition 4.4]{Andreychev} implies that if $U\subset \Spa(A,A^+)$ is a nice rational subspace with ring of functions $(A_U,A_U^+)$, then $(A,A^+)_{\sol}\to (A_U,A^+_U)_{\sol}$ is a rational localization as in \Cref{DefSolidAdicSpec}, namely, the proposition shows that  if $f,g$ generate the unit ideal and $U$ is the locus $\{|f|\leq |g|\neq 0\}$ then $(A_U,A_U^+)_{\sol}= (A,A^+)_{\sol}(\frac{f}{g})_{\sol}$ as analytic rings. Thus, if $U\subset \Spa(A,A^+)$ is a general open affinoid with ring of functions $(A_U,A_U^+)$, the map $(A,A^+)_{\sol}\to (A_U,A^+_U)_{\sol}$ is an open immersion of analytic rings as it is so locally along nice rational localizations which are open. Thanks to \cite[Theorem 4.1]{Andreychev}, the functor sending an open affinoid $U\subset \Spa(A,A^+)$ to $\ob{D}((A_U,A_U^+)_{\sol})$ satisfies analytic descent, and has the property that for $U\subset \Spa(A,A^+)$ an arbitrary open the map $\ob{D}(A,A^+)_{\sol}\to \ob{D}(U)$ is an open localization (as it its for the basis of open affinoids). This defines a morphisms of locales $\ob{Sm}(\ob{D}(A,A^+)_{\sol})\to \Spa(A,A^+)$. By construction and the definition of the solid spectrum, this functor factors through a map
\[
\ob{Sm}(\ob{D}(A,A^+)_{\sol})\to \Spec^{\sol}((A,A^+)_{\sol})\to \Spa(A,A^+).
\] 
It is left to see that the map $\Spec^{\sol}((A,A^+)_{\sol})\to \Spa(A,A^+)$ is an homeomorphism identifying rational localizations. But this follows from Stone duality, the fully faithful embedding of complete Huber pairs into analytic rings of \cite[Proposition 3.34]{Andreychev}, and the fact that rational localizations of the solid ring $(A,A^+)_{\sol}$ arise from standard rational localizations of the Tate Huber pair $(A,A^+)$ as these last are stable under composition.
\end{proof}

\begin{remark}\label{RemarkSolidSpectrumBounded}
In \cite{camargo2024analytic} we introduce a suitable subcategory of solid Huber rings over $\Z((\pi))$ called \textit{bounded rings}. A bounded ring $A$ has a $\dagger$-nilradical $\Nil^{\dagger}(A)$ of \textit{elements of norm $0$}, the quotient $A^{\dagger-\red}=A/\Nil^{\dagger}(A)$ is called the $\dagger$-reduction of $A$. Then, the natural map $\Spec^{\sol}(A^{\dagger-\red})\to \Spec^{\sol}(A)$ is an homeomorphism identifying rational localizations.  Bounded rings also admit an uniform completion $A^u$. If the ring $A$ is in addition Gelfand (as introduced in \cite{dRFF}),  $A^u$ is an uniform Tate Huber ring and one has natural homeomorphisms $\Spa(A^u, A^{u,+}) = \Spec^{\sol}(A^u)= \Spec^{\sol}(A)$  identifying rational localizations.  Proving these facts go beyond the scope of these notes, and we leave them for a future update of \cite{camargo2024analytic}. 
\end{remark}

\begin{lemma}\label{LemmaRationalCoversClosedImmersions}
Let $A\to B$ be a morphisms of solid Huber rings with the induced analytic ring structure such that $\pi_0(A)\to \pi_0(B)$ is surjective on underlying rings. Let $B\to B'$ be a rational localization of $B$, then there exists a rational localization $A\to A'$ such that $B'=A'\otimes_A B$. Moreover, if $B\to B'$ is a standard rational localization we can choose $A\to A'$ to be a composite of two standard rational localizations.
\end{lemma}
\begin{proof}
By induction, it suffices to prove the claim for standard rational localization.
 Let $f_1,\ldots, f_n,g\in B$ be elements generating the unit ideal and consider the rational localization $B'=B\left( \frac{f_1,\ldots, f_n}{g} \right)_{\sol}$.
 Consider a linear combination $\sum_i a_i f_i+bg=1$, and take lifts $\widetilde{a_i}$, $\widetilde{f}_i$, $\widetilde{b}$ and $\widetilde{g}$ to $A$. Let $c=\sum_i  \widetilde{a}_i \widetilde{f}_i+ \widetilde{b} \widetilde{g}$.
 We claim that the composite of rational localizations 
\[
A\to A\left(\frac{c^2,1}{c}\right)_{\sol} \to A\left(\frac{c^2,1}{c}\right)_{\sol} \left(\frac{\widetilde{f}_1,\ldots, \widetilde{f}_n}{\widetilde{g}}\right)_{\sol}=:A'
\]
works.  Indeed, the rational localization $A\to A\left(\frac{c^2,1}{c}\right)_{\sol}$ is precisely the base change along $\bb{Z}[T]\to \bb{Z}[T^{\pm 1}]_{\sol}$ of the map $\bb{Z}[T]\to B$ depicted by $c$. Since $c$ is invertible on $B$, we have a natural map surjective on $\pi_0$ with the induced structure
\[
A(\frac{c^2,1}{c})_{\sol}\to B.
\]
It is then clear from this map that $A'\otimes_A B= B'$ as wanted. 
\end{proof}

We finish this section with some lemmas that will come in handy later.

\begin{lemma}\label{LemmaRatLoc}
Let $A$ be a solid Huber ring and $(f_1,\ldots, f_n)\in A$ elements generating the unit ideal. Then the  standard rational localizations
\[
A\to A\left(\frac{f_1,\ldots, f_n}{f_i}\right)_{\sol}
\]
for $i=1,\ldots, n$ form an open cover of $\Spec^{\sol}(A)$.
\end{lemma}
\begin{proof}
Consider a linear combination $\sum_{i=1}^n a_i f_i =1$, and let $R=\bb{Z}[X_1,\ldots X_m,T_1,\ldots, T_n]/(\sum_i X_iT_i-1)$ be endowed with the induced analytic ring structure from $\bb{Z}_{\sol}$. We have a morphism of solid Huber ring 
\[
R\to A
\]
sending $X_i\mapsto a_i$ and $T_i\mapsto f_i$. Then, by base change we can assume without loss of generality that $A=R$. The theorem follows form the fact that the rational localizations of the form $R\left(\frac{f_1,\ldots, f_n}{f_i}\right)$ cover the adic spectrum $\Spa(R,\bb{Z})$ (see \cite[Lemma 10.4]{ClausenScholzeCondensed2019}) and the descent of the category of quasi-coherent sheaves for discrete adic spaces (\cite[Theorem 9.8]{ClausenScholzeCondensed2019}).
\end{proof}

\subsection{Solid Huber stacks}\label{ss:SolidHuberStk}

In this section we introduce the category of solid Huber stacks and solid adic spaces. We also construct some functors from schemes, formal schemes,  discrete adic spaces, and Tate adic spaces to these categories.

\begin{definition}\label{DefSolidStacks}
Let $\ob{Aff}^{\Hub}\subset \ob{AnStk}^{\ob{aff}}_{\Z_{\sol}}$ be the full  subcategory of affinoid analytic stacks over $\Z_{\sol}$ generated by the analytic spectrum of solid Huber rings\footnote{To make the following construction precise we need to take a small subcategory of solid Huber analytic rings, for example, we can fix an uncountable regular  cardinal $\kappa$ and work with $\kappa$-compact solid Huber rings. In the rest of the notes we will ignore these set theoretical issues.}. We let $\ob{SolHubStk}:=\ob{AnStk}(\ob{Aff}^{\Hub})$ be the category of analyitc stacks of \Cref{DefAnStacksFull} generated by $\ob{Aff}^{\Hub}$. We call $\ob{SolHubStk}$ the category of \textit{solid Huber stacks}, given $A\in \Ring^{\Hub}_{\Z_{\sol}}$ a solid Huber ring we let $\ob{AnSpec}(A)\in \ob{SolHubStk}$ be the object corepresented by $A$. 
\end{definition}

Among solid Huber stacks we find the category of solid adic spaces. 

\begin{definition}
Let $\ob{Aff}^{\Hub}$ be endowed with the Grothendieck topology with covers given by rational localizations $\{\AnSpec(B_i)\to \AnSpec(A)\}_i$ such that $\{\Spec^{\sol} B_i\to \Spec^{\sol} A\}$  is jointly surjective, we call this Grothendieck topology the \textit{solid adic topology of $\ob{Aff}^{\Hub}$}. We let $\widetilde{\n{X}}$ be the $\infty$-topos defined by this Grothendieck topology. 

\begin{enumerate}

\item By construction, the functor $\Spec^{\sol}(-)\colon \Ring^{\Hub}_{\Z_{\sol}}\to \ob{Top}$ from solid Huber rings to topological spaces sends rational covers to open covers in $\ob{Top}$.
 Hence, we have a natural left Kan extension $|-|^{\sol} \colon \widetilde{\n{X}}\to \ob{Top}$ that we call the  \textit{underlying solid spectrum} of the solid Huber stack. Concretely, if $X\in \widetilde{\n{X}}$ is written as a colimit $X=\varinjlim_i \AnSpec(A_i)$ with $A_i$ solid Huber rings, then $|X|^{\sol}=\varinjlim_i \Spec^{\sol}(A_i)$, where the last colimit is taking place in topological spaces. 

\item  Let $X\in \widetilde{\n{X}}$ and let $U\subset |X|^{\solid}$ be an open subspace. We define the substack $X_U\subset X$ to be the immersion whose $A$-points for $A$ a solid Huber ring are the maps $\AnSpec(A)\to X$ such that $\Spec^{\sol}(A)\to |X|^{\sol}$ factors through $U$. We say that a morphism $Y\to X$ is an \textit{open immersion} if it is equivalent to a map $X_U\to X$ for $U\subset |X|^{\sol}$ open. 

\item An object $X\in \widetilde{X}$ is called a \textit{solid adic space} if it admits an open cover $U_i\to X$ such that $U_i$ is corepresented by solid Huber rings. We let $\ob{SolAdic}\subset \widetilde{X}$ be the full subcategory of solid adic spaces. 

\end{enumerate} 
\end{definition}

\begin{remark}
Heuristically, we can think of a solid adic space as the datum of a pair $X=(|X|^{\sol}, \s{O}_X)$ where $|X|^{\sol}$ is a topological space, and $\s{O}_X$ is a sheaf valued in solid Huber rings such that there is an open cover $\{U_i\}$ of $X$ such that $U_i = \Spec^{\sol}(\s{O}_X(U_i))$. 
\end{remark}

%
%
%
%

\begin{example}\label{ExamplesSolidHuberSpaces}
In the following we list some examples of solid Huber stacks that occur in practice:

\begin{enumerate}

\item Consider the functor $(-)_{\sol}\colon \Ring^{\delta}\to \Ring^{\Hub}$ from discrete animated rings to solid Huber rings sending $R\mapsto R_{\sol}$. This functor preserves colimits, and sends Zariski open localizations to rational localizations of solid Huber rings.
In particular, for $A$ a discrete ring we have a quotient map of spectral spaces $\Spec^{\sol}(A_{\sol})\to \Spec(A)$. 
 Therefore, it naturally promotes to a functor 
\[
(-)_{\sol}\colon \ob{Sch}\to \ob{SolidAdic}
\]
from schemes to solid adic spaces that we call the \textit{solidification functor}.

\item Let $A\in \ob{Ring}^{\delta}$ be a discrete animated ring which is $I$-adically complete for a finitely generated ideal $I\subset \pi_0(A)$. Let $\ob{Ring}^{\wedge_I}_A$ be the category of $I$-adically complete animated $A$-algebras. Let $\widehat{A}$ be the solid ring obtained as the $I$-adic completion of the discrete condensed ring $A(*)$ We have a natural functor $\widehat{(-)}_{\sol}\colon \ob{Ring}^{\wedge_I}_{A}\to \Ring^{\Hub}_{\widehat{A}_{\sol}}$ sending $B$ to $B(*)^{\wedge_I}_{\sol}$, that is, to the solid ring given by the $I$-adic completion of the discrete  condensed ring $B(*)$ endowed with the solid structure induced from $B(*)_{\sol}$. Then the functor $\widehat{(-)}_{\sol}$ preserves colimits and sends Zariski localizations of the formal spectrum $\Spf(B)=\Spec(B/I)$ to rational localizations of $\widehat{B}_{\sol}$.  By \Cref{LemmaInvarianceAdicallyComplete} (b) we even have a quotient map of spectral spaces
\[
\Spec^{\sol}(\widehat{B}_{\sol})\to \Spf(B).
\]
This functor gives rise to a map 
\[
\widehat{(-)}_{\sol}\colon \ob{FSch}^I_{A}\to \ob{SolAdic}_{\widehat{A}_{\sol}}
\]
from $I$-adic formal schemes over $A$ to solid adic spaces over $\AnSpec(\widehat{A}_{\sol})$.

\item Let $(A,A^+)$ be a Tate sheafy Huber ring. By \Cref{PropComparisonSpectrum} (2) rational localizations of $(A,A^+)$ give rise to rational localizations of the solid Huber ring $(A,A^+)_{\sol}$. By \Cref{PropComparisonSpectrum} we even have an equivalence 
\[
\Spec^{\sol}(A,A^+)\xrightarrow{\sim} \Spa(A,A^+)
\]
identifying rational localizations.  This produces a functor 
\[
\ob{TateAdic}\to \ob{SolAdic}
\]
from the classical category of sheafy Tate adic spaces to the category of solid adic spaces.

\item Finally, let $\ob{Ring}^{\Hub,\delta}\subset \ob{Ring}^{\Hub}$ be the full subcategory of discrete Huber rings, that is, the full subcategory of solid Huber rings $A$ such that $A^{\triangleright}$ is discrete. Thanks to  \cite[Proposition 2.9.6]{MannSix}  $\ob{Ring}^{\Hub,\delta}$ is equivalent to the category of discrete Huber pairs as in \Cref{def:DiscreteAdic}. By \Cref{theo:DescendAdicSpaces}, any classical rational localization $(A,A^+)\to (A_U, A^+_U)$ corresponding to an open subspace $U\subset \Spa(A,A^+)$ in the usual sense of Huber is the composite of a closed localization and an open localization in the sense of analytic rings. Thus, \Cref{theo:DescendAdicSpaces} proves that if $\{U_i\to \Spa(A,A^+)\}$ is a rational cover in Huber's sense, then the maps of solid Huber rings $\{(A,A^+)_{\sol}\to (A_{U_i}, A_{U_i}^+)_{\sol}\}_i$ form a $!$-cover being refined by a composite of closed and open covers. 
 This shows that the natural functor $\ob{Aff}^{\Hub,\delta}\to \ob{SolHubStk}$ from discrete affinoid Huber spaces to solid Huber stacks extends to a functor 
\[
\ob{Adic}^{\delta}\to \ob{SolHubStk}
\]
from the category of discrete adic spaces to the category of solid Huber stacks. 
Note that this functor does not factors through solid adic spaces since not all rational localization in Huber's sense is an open immersion at the level of analytic rings, eg. the localization $(\Z[T],\Z)\to(\Z[T^{\pm 1}], \Z)$ is not.

\end{enumerate}
\end{example}

\subsection{Cotangent complex of analytic rings}

As one might expect, analytic rings also have an analogue theory of cotangent complexes essentially identical to that of animated or connective commutative algebras, eg. as in \cite[Sections 3.2 and 3.3]{LurieDerivedAlgebraic}, \cite[Sections 7.3 and 7.4]{HigherAlgebra} and \cite[Chapters 17 and 25]{LurieSpectralAlg}. In this section we briefly review the theory of the cotangent complex of analytic rings, stating only the main features that will be important for the next sections.

\begin{definition}\label{DefCotangentComplex}
Let $f\colon A\to B$ be a morphism of analytic rings, and let $M\in \ob{D}_{\geq 0}(B)$ be a connective $B$-module. Consider the square-zero extension $B\oplus M$ as an animated $B$-algebra (obtained from animation from the static case) endowed with the induced analytic ring structure from $B$. We define the anima of \textit{$A$-linear derivations of $B$ with values in $M$}, denoted as $\ob{Der}_{A}(B,M)$, as the anima of $A$-algebra morphisms $B\to B\oplus M$ whose composition with the augmentation map $B\oplus M\to B$ is the identity. 

The \textit{cotangent complex} of $f$ is the $B$-module $\bb{L}_{B/A}\in \ob{D}_{\geq 0}(B)$ correpresenting the functor sending $M\in \ob{D}_{\geq 0}(B)$ to $\ob{Der}_A(B,M)$. We let $(\id,d)\colon B\to B\oplus \bb{L}_{B/A}$ denote the universal $A$-linear derivation of $B$ with values in $\bb{L}_{B/A}$. 
\end{definition}

\begin{remark}\label{ExistenceCotangent}
\begin{enumerate}
\item Let $A\to B$ be a morphism of analytic rings. The functor $M\mapsto \ob{Der}_{A}(B,M)$ preserves limits and is accessible, hence by the adjoint functor theorem \cite[Corollary 5.5.2.9]{HigherTopos} it is correpresented by an object in $\ob{D}_{\geq 0}(B)$. This justifies the existence of the cotangent complex in \Cref{DefCotangentComplex}.

\item Let $\Mod_{\geq 0}\to \AnRing$ be the cocartesian fibration associated to the functor sending an analytic ring $A$ to its category $\ob{D}_{\geq 0}(A)$ of connective $A$-modules.   Using the formalism of relative adjunctions of \cite[Section 7.3.2]{HigherAlgebra} one can construct functorially the relative contangent complex as a functor $\ob{Fun}(\Delta^1,\AnRing)\to \Mod$ over $\ob{Fun}(\{1\}, \AnRing)= \AnRing$  that sends a morphism $A\to B$ to the pair $(B, \bb{L}_{B/A})$. 
\end{enumerate}
\end{remark}

The cotangent complex satisfies the expected stability properties:

\begin{lemma}\label{LemmaBasicCotangent}
\begin{enumerate}

\item Let $A\to B$ and $A\to A'$ be morphisms of analytic rings with pushout $B\to B'$. Then the natural map 
\[
B'\otimes_{B} \bb{L}_{B/A}\xrightarrow{\sim} \bb{L}_{B'/A'}
\]
is an equivalence.

\item If $A\to B$ is the colimit of a diagram of morphisms $\{A_i\to B_i\}$ in $\AnRing$, then the natural map 
\[
\varinjlim_i B\otimes_{B_i} \bb{L}_{B_i/A_i}\to \bb{L}_{B/A}
\]
provided by \Cref{ExistenceCotangent} (2) is an equivalence. 

\item Let $A\to B\to C$ be a triangle of analytic rings. Then there is a natural fiber sequence of $C$-modules 
\[
C\otimes_B \bb{L}_{B/A} \to \bb{L}_{C/A}\to \bb{L}_{C/B}
\]

\item Let $A\to A'$ be an idempotent map of analytic rings. Then $\bb{L}_{A'/A}=0$.

\end{enumerate}
\end{lemma}
\begin{proof}
The proofs of (1)-(3) are standard from the definition of the cotangent complex and left to the reader. For (4), notice that $A'=A'\otimes_A A'$ so that 
\[
\bb{L}_{A'/A}=(A'\otimes_{A} A')\otimes_{A'} \bb{L}_{A'/A}= \bb{L}_{(A'\otimes_A A')/A'}=\bb{L}_{A'/A'}=0.
\]
\end{proof}

\begin{example}\label{ExampleCotangentComplex}

In the following we present some examples of cotangent complexes of morphisms of analytic rings. We will see that in many cases we recover what classically is called the \textit{continuous cotangent complex} of a morphism of topological rings. 

\begin{enumerate}

\item Let $A^{\cond}\to B^{\cond}$ be a morphism of discrete animated rings endowed with the trivial analytic ring structure. Then $\bb{L}_{B^{\cond}/A^{\cond}}=\bb{L}_{B/A}$ is also discrete and agrees with the classical cotangent complex of $A\to B$ considered as a discrete condensed $B$-module. Indeed, If $M$ is a conective $B$-module then the anima of condensed $A$-algebra morphisms $B\to B\oplus M$  agrees with the anima of discrete $A$-algebra morphisms $B\to B\oplus M(*)$. 

\item Let $A\to B$ be a morphism of analytic rings, and consider the factorization $A\to B_{A/} \to B$ where $B_{A/}$ has the induced structure from $A$. Since $B_{A/}\to B$ is idempotent, we have that 
\[
\bb{L}_{B/A}= B\otimes_{B_{A/}}\bb{L}_{B_{A/} /A}.
\] 
A similar argument shows that $\bb{L}_{B/A}=\bb{L}_{B/A^{\triangleright,\ob{cond}}}$.

\item One has $\bb{L}_{\Z[T]_{\sol}/\Z_{\sol}}\cong  \bb{Z}[T]$, namely, this follows by base change and (1) and (2) above since  $\bb{L}_{(\bb{Z}[T],\bb{Z})_{\sol}/ \bb{Z}_{\sol}} \cong \bb{Z}[T]$ and $(\Z[T],\Z)_{\sol}\to \Z[T]_{\sol}$ is idempotent. In particular, if $A$ is a solid Huber ring and $A[T]_{\sol}=A\otimes_{\Z_{\sol}}\Z[T]_{\sol}$, one has that $\bb{L}_{A[T]_{\sol}/A}\cong A[T]_{\sol}^{\triangleright}$.

\item Let $R$ be a solid animated ring with induced analytic ring structure from $\Z_{\sol}$. Suppose that $R$ is derived $I$-adically complete for $I\subset \pi_0(R)$ a finitely generated ideal. Let $R\langle T\rangle$ be the derived $I$-completion of $R[T]$. Then, as the solid tensor product preserves $I$-adically complete modules, one has that $R[T]\to R\langle T\rangle$ is an idempotent morphism, and therefore that $\bb{L}_{R\langle T \rangle/R}\cong R\langle T \rangle$. This recovers the \textit{$I$-completed cotangent bundle}.

\end{enumerate}
\end{example}

One has the standard connectivity estimates  of the cotangent complex for analytic rings. 

\begin{lemma}\label{LemmaConnectivitySym}
Let $A$ be an analytic ring and $M$ an $n$-connective $A$-module for $n\geq 1$. If $n\geq 1$ then $\Sym^m_A M$ is $(m+n-1)$-connective. If $n\geq 2$ then $\Sym^m_A M$ is $(n+2m-2)$-connective.   If $A$ is defined over $\Q$ and $n\geq 0$ then $\Sym_A^m M $  is $nm$-connective. 
\end{lemma}
\begin{proof}
By \cite[Proposition 25.2.4.1]{LurieSpectralAlg} we know that the statement holds for discrete animated rings. It is then immediate that it holds for animated rings in the topos of condensed anima. Since $\Sym^m_A$ is the $A$-completion of $\Sym^m_{A^{\triangleright}}$, and $A$-completion preserve connective modules, one deduces the first and second claims of the  lemma. If in addition $A$ is defined over $\Q$, then $\Sym^m_{A}$ is the quotient by the action of $\Sigma_n$ of th $m$-th fold tensor product $M^{\otimes_A m}$ which is $nm$-connective. Since colimits preserve connective objects, we deduce the third claim. 
\end{proof}

\begin{corollary}\label{CoroStimatesWedge}
Let $A$ be an analytic ring and $M$ an $n$-connective $A$-module for $n\geq 0$. Then, for $m\geq 1$, $\bigwedge^m_A M$ is $(n+m-1)$-connective. If $A$ is defined over $\Q$ then it is $nm$-connective.
\end{corollary}
\begin{proof}
By \Cref{LemmaConnectivitySym} we know that $(\bigwedge^m M) [m] = \Sym^m_A (M[1])$ is $(n+1+2m-2)$-connected, thus, $\bigwedge^m M$ is $(n+m-1)$-connected. If $A$ is defined over $\Q$ the estimate follows from the fact that $(\bigwedge^m M)[m]=\Sym^m_A M[1]$ is $(n+1)m$-connected. 
\end{proof}

\begin{proposition}\label{PropConnectivityCotangent}
Let $f\colon A\to B$ be a morphism of analytic rings with the induced analytic ring structure. Let $d\colon B \to  \bb{L}_{B/A}$ be the universal derivation, the composite $A\to B\to \bb{L}_{B/A}$ is zero, so that we have an induced map $\varphi_f\colon B\otimes_A\ob{cofib}(f)\to \bb{L}_{B/A}$.  The map $\varphi_f$ is surjective on $\pi_0$, if $\ob{fib}(f)$ is connective, then $\ob{fib}(\varphi_{f})$ is $2$-connective, if $\ob{fib}(f)$ is $n$-connective for $n>0$, then $\ob{fib}(\varphi_f)$ is $(n+3)$-connective.

 If $A$ is defined over $\Q$ it is even $2n$-connective.
\end{proposition}
\begin{proof}
The statement for classical animated rings is \cite[Proposition 25.3.6.1]{LurieSpectralAlg}, and for $\Q$-algebras follows from \cite[Theorem 7.4.3.1]{HigherAlgebra} (as animated rings and commutative algebras agree). This automatically passes to the  topos of condensed anima, and by \Cref{ExampleCotangentComplex} (2) it extends to analytic rings. 
\end{proof}

\begin{corollary}\label{CoroDetectCotangent}
Let $f\colon A\to B$ be a morphism of analytic rings. Then $f$ is an isomorphism if and only if $\pi_0(A)\to \pi_0(B)$  is an isomorphism (in particular, $f$ has the induced analytic ring structure), and $\bb{L}_{B/A}=0$.
\end{corollary}
\begin{proof}
This is immediate from \cite[25.3.6.6]{LurieSpectralAlg}.
\end{proof}

Square-zero  extensions of analytic rings work in the same way as for classical derived algebraic geometry, see \cite[Section 3.3]{LurieDerivedAlgebraic} and \cite[Section 17.1.3]{LurieSpectralAlg}.  We have the following useful lemma.

\begin{lemma}\label{LemmaApproximationsCotangentComplex}
Let $A\to B$ be a morphism of analytic rings with the induced analytic ring structure. Suppose that $\fib{A^{\triangleright}\to B^{\triangleright}}$ is $k$-connective for $k>0$ and  consider the universal derivation $d\colon B\to  B \oplus \bb{L}_{B/A}$ and let $\widetilde{B}$ be the associated square zero extension obtained as the pullback 
\[
\begin{tikzcd}
\widetilde{B} \ar[d] \ar[r] &  B\ar[d,"d"]\\ 
B \ar[r,"0"] & B\oplus \bb{L}_{B/A}.
\end{tikzcd}
\] Then the fiber of the map $\widetilde{f}\colon A\to \widetilde{B}$ is $(k+1)$-connective. 
\end{lemma}
\begin{proof}
This follows  from \cite[Proposition 3.3.6]{LurieDerivedAlgebraic}.
\end{proof}

Symmetric powers of complete modules remain complete. 

\begin{lemma}\label{LemmaCompletenessSym}
Let $A$ be a solid ring and $I\subset \pi_0(A)$ a finitely generated ideal such that $A$ is derived $I$-adically complete. Let $M\in \ob{D}_{\geq 0}(A)^{\wedge}_I$ be a connective derived $I$-adically complete module. Then for all $k\in \N$, $\Sym^k_{A} M$ is derived $I$-adically complete.  
\end{lemma}
\begin{proof}
As $M$ is connective, it can be written as a geometric realization $M=|P^{\bullet}|$ of direct sums $P^n=\bigoplus_{i\in I_n} A[S_{i,n}]$ of finite projective modules  with $S_{i,n}$ light profinite. Since $M$ is derived $I$-adically complete, and taking $I$-adic completions commutes with geometric realizations of connective objects (thanks to \cite[Proposition 1.2.4.5]{HigherAlgebra} as $I$-adic completions can be tested at the level of cohomology groups), one has that $M=|P^{\bullet, \wedge_I}|$. 
Thus, we get that $\Sym^k_{A} M = | \Sym^{k} P^{\bullet, \wedge_I}|$ and it suffices to show that $\Sym^{k} P^{\bullet, \wedge_I}$  is derived $I$-adically complete. Since derived $I$-adic completions commutes with $\aleph_1$-filtered colimits, we can assume without loss of generality that $P=\bigoplus_{n\in N} A[S_n]$ is a countable direct sum of finite projective modules. By enlarging $S$, eg, to the Cantor set, we can even assume that $P=A[S]^{\oplus_{\N}}$. Hence, we are reduced to the case 
\[
M=A[S]^{\widehat{\oplus}_{\N,I}}.
\]
Consider a map $\Z[[X_1,\ldots, X_n]]\to A$ such that $I$ is generated by the image of the $X_i$, then an $A$-module is $I$-adically complete if and only if its restriction to $\Z[[X_1,\ldots, X_n]]$ is $(X_i)$-adically complete. Since   the solid tensor product preserves complete connective modules (\cite[Proposition 2.12.10]{MannSix}), we have that $M= A\otimes_{\Z[[X_1,\ldots, X_n]]} N$ where $N$ is a completed direct sum of copies of  $\Z[[X_1,\ldots, X_n]][S]$. By functoriality of the symmetric powers, we have that $\Sym^k_A M= A\otimes_{\Z[[X_1,\ldots, X_n]]} \Sym^k_{\Z[[X_1,\ldots, X_n]]} N$.  Thus, using  \cite[Proposition 2.12.10]{MannSix} again, we can assume without loss of generality that $A=\Z[[X_1,\ldots, X_n]]$ and $I$ is generated by $X_1,\ldots, X_n$. Moreover, by an inductive argument we are even reduced to the case of $A=\Z[[X]]$ and $I=(X)$. 

We now follow the standard argument that proves the $I$-adic completeness of the solid tensor product. Let $P=\widehat{\bigoplus}_J \prod_{\N} A$ and $X$-adically completed direct sum of countably many copies of $\prod_{\N} A$. We can assume that $J\cong \N$ as totally ordered set.  We can write 
\[
P=\varinjlim_{f\colon J\to \N} \prod_{J} X^{f(j)} \prod_{N} A
\]
where $f(j)\to \infty$ as $j\to \infty$ and $f\leq g$ if and only if $f(j)\geq g(j)$ for all $j\in J$. Thus, we have that 
\[
\Sym^k_A P = \varinjlim_{f\colon J\to \N} X^{kf(j)} (\Sym^k_A  \left( \prod_{J} \prod_{\N} A) \right), 
\]
and one directly computes 
\[
\Sym^k_A \left( \prod_{K} A \right)= \prod_{K^k/\Sigma_k} A.
\]
Thus, we find that 
\[
\Sym^k_A P = \varinjlim_{f\colon J\to \N} X^{kf(j)} \prod_{(J\times \N)^k /\Sigma_k} A. 
\]
From this presentation one can deduce by a direct computation that $\Sym^k_A P$ is $X$-adically complete. Let us give a more formal argument, from the previous presentation one deduces that $\Sym^k_A P$ is static, and therefore must be equal to the  abelian quotient of $P^{\otimes_A k}$ by the action of the symmetric group $\Sigma_k$. Since $P^{\otimes_A k}$ is derived $X$-adically complete by \cite[Proposition 2.12.10]{MannSix}, its quotient by $\Sigma_k$ is as well (being a cokernel of a map of derived $X$-adically complete modules). 
\end{proof}

\begin{remark}\label{RemCompletenessCotangent}
Let $R$ be a classical $I$-adically complete ring seen as solid ring, with $I\subset \pi_0(R)$ a finitely generated ideal.  From \Cref{ExampleCotangentComplex} (4) one can deduce that the cotangent complex of a  morphism $R\to A$ of classical $I$-complete animated rings of almost finite presentation (that is, such that modulo $I$ the $n$-th truncation is compact as discrete $n$-truncated animated rings for all $n$) is derived $I$-complete. Applying \Cref{LemmaCompletenessSym} one gets  that $\bigwedge^m_A \bb{L}_{A/R}$  are also derived $I$-adically complete for all $m$. We do not know how to prove in general that an $I$-complete solid $R$-algebra $A$ has $I$-adically complete cotangent complex. It was indicated by Akhil Mathew and Peter Scholze that an argument involving the filtration of the Hochschild-Homology of $A\to R$ in terms of the cotangent complex should help (cf. \cite[Theorem 1.2.1]{raksit2020hochschild}). 
\end{remark}

\subsection{Morphisms of solid finite presentation}\label{ss:FinitePresentation}

In classical derived algebraic geometry  there are at least two important finiteness notions for morphisms of rings that appear in practice (see \cite[Section 3.1]{LurieDerivedAlgebraic} and \cite[Section 4]{LurieSpectralAlg}): 

\begin{definition}\label{DefFinitePresentedDiscrete}
Let $f\colon A\to B$ be a morphism of discrete animated rings.

\begin{enumerate}
\item We say that $f$ is of \textit{finite presentation} (called in \cite[Section 4.1.1]{LurieSpectralAlg} of finite generation) if $B$ is compact as animated $A$-algebra.  Equivalently, if $B$ is in the idempotent completion of the full subcategory of $\Ring_{A/}$ generated under finite colimits by  the polynomial algebra $A[T]$.

\item We say that $f$ is of \textit{finite presentation to order $n$} if $\tau_{\leq n} B$ is a compact object in the category of $n$-truncated $A$-algebras. Equivalently (\cite[Proposition 3.1.5]{LurieDerivedAlgebraic}), if there is a morphism of finite presentation $f\colon A\to B'$ and a map of $A$-algebras $B'\to B$ inducing an isomorphism after $n$-truncations. 

\item We say that $f$ is of \textit{almost finite presentation} if it is of finite presentation to order $n$ for all $n\in \N$.

\end{enumerate}

\end{definition}

\begin{remark}\label{RemarkExample}
A morphism $f\colon A\to B$ of discrete animated rings if of finite presentation to order $0$ if and only if $\pi_0(A)\to \pi_0(B)$ is a classical morphism of finite presentation of static rings.  Morphisms of finite presentation to order $0$ are also called \textit{morphisms of finite type}. Note that this differs from the standard convention of morphisms of finite type in classical non-derived algebraic geometry.
\end{remark}

\begin{remark}\label{RemarkCharacterizations}
By \cite[Proposition 4.1.2.1]{LurieSpectralAlg} (or more precisely, its animated analogue \cite[Proposition 3.2.18]{LurieDerivedAlgebraic}) a map of discrete animated rings $f\colon A\to B$ is of (almost) finite presentation if and only if $\pi_0(A)\to \pi_0(B)$ is a classical morphism of finite presentation  of static rings and $\bb{L}_{B/A}$ is a (an almost) perfect module. 
\end{remark}

In the set up of solid Huber rings is not so clear what  the correct  candidate is for a notion  of  finite presentation, eg. there are compact objects such that $\Sym_{\sol} \Z_{\sol}[S]$ for $S$  a light profinite set, which we would like to avoid  as these rings are too big. It turns out that the  differential characterization of morphisms of (almost) finite presentation provides a reasonable generalization.

\begin{definition}
Let $f\colon A\to B$ be a morphism of solid Huber rings. 

\begin{enumerate}

\item We say that $f$ is \textit{of  solid finite   type} if there is a presentation  of static solid Huber rings
\[
\pi_0(B)= \pi_0(A[T_1,\ldots, T_n]_{\sol})/(f_1,\ldots, f_m)
\]
with the induced analytic structure.  Equivalently, if there is a morphism of solid Huber rings $A[T_1,\ldots,T_n]_{\sol}/^{\bb{L}}(f_1,\ldots, f_m)\to B$ which induces an isomorphism on $\pi_0$.

\item We say that $f$ is of \textit{solid finite  presentation} if it is  of solid finite type and $\bb{L}_{B/A}$  arises from is a perfect $B^{\triangleright}(*)$-module. We say that $f$ is of \textit{solid almost finite presentation} if it is of solid finite type and $\bb{L}_{B/A}$ arises from a $B^{\triangleright}(*)$-module of almost finite presentation.

\end{enumerate}

\end{definition}

\begin{example}\label{ExampleMorphismFinitePresentation}

In the following we present different examples of morphisms of solid finite presentation that are related with morphisms topologically of finite type.

\begin{enumerate}

\item  The most basic example is that of schemes. Let $\Ring^{\delta}$ be the category of discrete  animated rings and consider the functor $\Ring^{\delta}\to \Ring^{\Hub}$ sending $R\mapsto R_{\sol}$. Then,  $A\to B$ is a morphism (almost) of finite presentation in $\Ring^{\delta}$ if and only if $A_{\sol}\to B_{\sol}$ is a morphism of  solid (almost)  finite presentation.

\item Another example arises from formal geometry.  Let $A$ be a discrete static ring which is complete for a finitely generated ideal $I\subset A$, eg. take $A=\Z_p$. For simplicity suppose that $A$ is noetherian. A morphism of $I$-complete static rings is said \textit{topologically of finite type} if $B$ is a quotient of the (noetherian!) ring $A\langle T_1,\ldots, T_n \rangle$ given by the $I$-adic completion of the polynomial algebra over $A$. From this, we see that $A\to B$ is topologically of finite type if and only if the morphism of solid rings $A_{\sol}\to B_{\sol}$ endowed with the $I$-adic topology is solid of finite type. As $A$ is noetherian, one can also show that if $A\to B$ is topologically of finite then $A_{\sol}\to B_{\sol}$ is of  solid almost finite presentation, namely, there is a surjective map  $A\langle T_1,\ldots, T_n \rangle\to B$ making $B$  a coherent $A\langle T_1,\ldots, T_n\rangle$-module. This yields that $B_{\sol}= A\langle T_1,\ldots, T_n\rangle_{\sol} \otimes_{A\langle T_1,\ldots, T_n \rangle (*)} B(*)$, so that the solid cotangent complex of $A\langle T_1,\ldots, T_n\rangle \to B$ is the base change of the algebraic cotangent complex which is  almost perfect.

\item Finally, we have examples arising from rigid geometry. Let $K$ be a complete non-archimedean field and  $A$ a Tate algebra over $K$ of finite type. Then $(K,\n{O}_K)_{\sol}\to (A,A^{\circ})_{\sol}$  is a morphism of solid finite type. Since the Tate algebras $K\langle T_1,\ldots, T_n \rangle$ are noetherian, one can prove as in (2) above that  $(K,\n{O}_K)_{\sol}\to (A,A^{\circ})_{\sol}$ is a morphism of  solid almost finite presentation.

\end{enumerate}

\end{example}

\begin{remark}\label{RemarkLocaFinitePresentation}
It is not completely  clear whether the property of being of solid finite presentation for a morphism of solid Huber rings  is local for the  analytic topology of the solid spectrum. For solid schemes, being of solid (almost) finite presentation  is local in the analytic topology  by the standard faithfully flat descent of morphisms  of (almost) finite presentation for discrete rings. We expect that being of (almost) solid finite presentation is local in the analytic adic spectrum of Gelfand Tate algebras, that is, for bounded affinoid rings $A$ such that $A^{\triangleright}/A^{\triangleright,\leq 1}$ is discrete (for the notion of Gelfand ring see \cite{dRFF}). 
\end{remark}

The next lemma proves some basic permanence properties of morphisms of solid finite presentation. 

\begin{lemma}\label{LemmaPermanenceSolidPresentation}
The following hold:

\begin{enumerate}

\item Let $A\to B$ be a morphism of  solid finite type of solid Huber rings and $A\to A'$ a morphism of solid Huber rings with base change $B'=B\otimes_A A'$. Then $A'\to B'$ is of solid finite type.

\item Let $A\xrightarrow{f} B\xrightarrow{g} C$ be morphisms of solid Huber rings. Suppose that $f$ is of solid finite type, then $g$ is of solid finite type if and only if $g\circ f$ is of solid finite type.  In particular, morphisms of solid finite type are stable under composition.

\item Let $A$ be a solid Huber ring and $\Ring^{\sol-\ob{ft}}_{A/}$ be the category of $A$-algebras of solid finite type. Then $\Ring^{\sol\ob{-ft}}_{A/}$ is stable under finite colimits and retracts.

\item Let $f\colon A\to B$ be a rational localization of solid Huber rings, then $f$ is of solid finite presentation.

\end{enumerate}
Moreover, the analogue statements (1)-(3) hold for morphisms of  solid  finite presentation and of solid almost finite presentation. 
\end{lemma}
\begin{proof}
\begin{enumerate}

\item  By definition, $A\to B$ is of finite type if and only if $\pi_0(B)$ is a quotient of a ring $\pi_0(A[T_1,\ldots, T_n]_{\sol})$ by a finitely generated ideal, this condition is obviously stable under base change. The cases of solid finite presentation and  of solid almost finite presentation follows form the fact that the cotangent complex satisfies base change.

\item  Let $f\colon A\to B$ be a morphism of solid finite type. Assume that $g\colon B\to C$ is of solid finite type, we want to show that $g\circ f$ is of solid finite type. Indeed, we have morphisms of solid Huber rings 
\begin{equation}\label{eqFiniteType}
A[T_1,\ldots, T_n]_{\sol}/^{\bb{L}}(a_1,\ldots, a_k)\to B 
\end{equation}
and 
\[
B[X_1,\ldots, X_n]_{\sol}/^{\bb{L}}(b_1,\ldots, b_l)\to C
\] 
which induce isomorphisms on connected components.  Hence, after lifting the elements $b_i$ to elements   $\widetilde{b}_i\in A[T_1,\ldots, T_n]_{\sol}$, we have a map 
\[
A[T_1,\ldots, T_n,X_1,\ldots, X_n]_{\sol}/^{\bb{L}}(a_1,\ldots, a_k,\widetilde{b}_1,\ldots, \widetilde{b}_l)\to C
\]
which is an isomorphism on connected components, proving that $g\circ f$ is of finite type. Conversely, suppose that $g\circ f$ is of finite type, then by (1) above the map $B\to B\otimes_A C$ is of finite type, and by the previous point it suffices to see that the multiplication map $B\otimes_A C\to C$ is of finite type.  

The statement only concerns the connected components of the rings, so we will assume for the rest of this implication that the rings are static, and all the base changes and quotients are as static solid Huber rings.  We have  presentations $B=A[T_1,\ldots, T_n]_{\sol}/I$ and $C=A[X_1,\ldots, X_m]_{\sol}/J$ where $I$ and $J$ are finitely generated ideals.  Lift the images of the $T_i$ in $C$ to elements  $\widetilde{T}_i\in A[X_1,\ldots, X_m]_{\sol}$. By point (4) down bellow, the rational localization $A[X_1,\ldots, X_m]_{\sol}\to A[X_1,\ldots, X_m]_{\sol}\otimes_{A[\widetilde{T}_1,\ldots, \widetilde{T}_n]} A[\widetilde{T}_1,\ldots, \widetilde{T}_n]_{\sol}$ is a morphism of finite presentation. Since the map $A[X_{1},\ldots, X_n]_{\sol}\to C$ factors through the solidification of the elements $\widetilde{T}_i$, and morphisms of finite type are stable under compositions by the first part of this point,  after replacing the given presentation of $C$ we can assume  without loss of generality that the elements $\widetilde{T}_i$ are solid in $A[X_1,\ldots, X_m]_{\sol}$. Hence, we have a commutative diagram 
\[
\begin{tikzcd}
A[T_1,\ldots, T_n]_{\sol} \ar[r] \ar[d] & A[X_1,\ldots, X_m]_{\sol} \ar[d] \\
B \ar[r] & C.
\end{tikzcd}
\]
where the  vertical maps are quotients by finitely generated ideals. We see that (the static) pushout $B\otimes_A C$ is given by $A[T_1,\ldots, T_n,X_1,\ldots, X_m]_{\sol}/(I,J)$, and that the map $B\otimes_A C\to C$ is obtained by taking the quotient of $B\otimes_{A} C$ with the ideal generated by the elements $T_i-\widetilde{T}_i \in A[T_1,\ldots, T_n,X_1,\ldots, X_m]_{\sol}$, proving that it is of finite type.

Finally, the analogue statement for morphisms of (almost) finite presentation follows form the case of finite type and  the fiber sequence of cotangent complexes
\[
B\otimes_A\bb{L}_{B/A}\to \bb{L}_{C/A}\to \bb{L}_{C/B}.
\]

\item We now prove that morphisms of solid finite type over the ring $A$ are stable under finite colimits and under retracts.  Let $B$ be an $A$-algebra of solid finite type and $B'\to B\to B'$ a retract, them  $B'$ is the retract associated to an idempotent morphism $r\colon B\to B$.  Thus, $\pi_0(B')$ is isomorphic to the $\pi_0$ of the pushout  $D$ of the diagram
\[
\begin{tikzcd}
B \otimes_A B  \ar[r, "m"] \ar[d," r\otimes \id"] \ar[d] & B\ar[d] \\ 
B  \ar[r] &  D
\end{tikzcd}
\]
where $m$ is the multiplication map. Thus, it suffices to see that $D$ is of finite type over $A$. The algebras  $B\otimes_A B$ and $B$ are of finite type over $A$ by parts (1) and (2) above, therefore it suffices to show that morphisms of solid finite type  are stable under finite colimits. 

We now prove stability under finite colimits. It suffices to prove that morphisms of solid finite type are stable under pushouts.   Consider a diagram of morphisms of solid finite type over $A$ of the form $C\leftarrow B \rightarrow D$, we want to see that $C\otimes_B D$ is of solid finite type. By (2) above we know that both $B\to C$ and $B\to D$ are morphisms of solid finite type.   By base change of (1), we then know that $C\to C\otimes_B D$ is of solid finite type, by stability under compositions of (2)  we then have that $A\to C\to C\otimes_B D$ is solid of finite type, proving what we wanted.

The statement for morphisms of solid (almost) finite presentation follows from the solid finite type case and  the fact that  (almost) perfect modules are stable under retracts and finite colimits.

\item By  the implication ``$g$ finite presentation implies $g\circ f$ finite presentation'' of part  (2) it suffices to show that a standard rational localization is of finite presentation. Let $f_1,\ldots, f_n,g\in A$ be elements in $A$ generating the unit ideal. By definition, 
\[
A\left(\frac{f_1,\ldots, f_n}{g} \right)_{\sol}= A[T_1,\ldots, T_n]_{\sol}/^{\bb{L}}(gT_i-f_i)
\]
which is clearly of finite presentation.
\end{enumerate}
\end{proof}

\subsection{Solid \'etale and smooth maps}\label{ss:SolidEtaleSmooth}

In this section we define the notions of solid smooth and \'etale maps of morphisms of solid Huber rings. Our definition will be inspired from that in derived algebraic geometry.

\begin{definition}\label{DefinitionEtaleSmooth}
Let $f\colon A\to B$ be a morphism of  solid  finite presentation of solid Huber rings.  We say that $f$ is \textit{solid smooth} if $\bb{L}_{B/A}$ is a finite projective $B$-module. Analogously, we say that $f$ is \textit{solid \'etale} if $\bb{L}_{B/A}=0$. 

More generally, a morphism $f\colon A\to B$ between solid Huber rings is said \textit{locally solid smooth} (resp. \textit{locally solid \'etale}) if, locally in a rational cover of $A$ and $B$, it is of  solid finite presentation and solid smooth (resp. solid \'etale). The \textit{relative dimension of $f$} is the locally constant rank of $\bb{L}_{B/A}$. 
\end{definition}

\begin{remark}\label{RemarkSolidSmoothetale}
In \Cref{DefinitionEtaleSmooth} one is obliged to impose a local condition for solid smooth and \'etale morphisms due to the lack of descent for morphisms of solid almost presentation, see \Cref{RemarkLocaFinitePresentation}. If we are in a context where being of  solid finite presentation is local for the analytic topology, then a morphism is solid smooth (resp. \'etale) if and only if it is locally solid smooth (resp. \'etale). In particular, this holds for discrete Huber rings, Tate algebras of finite type over a non-archimedean field $K$ by classical considerations, and it is expected for Gelfand rings. 
\end{remark}

It is in some occasions convenient to have some precise coordinates when working with smooth and \'etale morphisms of schemes. The following definition is a generalization of the standard charts to the solid case:

\begin{definition}\label{DefStandardEtaleSmooth}
A morphisms $f\colon A\to B$ of solid Huber rings is said \textit{solid standard smooth} if it is of the form 
\[
B=A[T_1,\ldots, T_n]_{\sol}/^{\bb{L}}(a_1,\ldots, a_k)
\]
where $k\leq n$ and $d=\det (\frac{\partial a_i}{\partial T_j })_{1\leq i,j\leq k}$ is invertible on $B$. We say that $f$ is \textit{solid standard \'etale} if $k=n$. 
\end{definition}

\begin{remark}\label{RemarkStandardSolidEtale}
Let $f\colon A\to B$ be a  solid standard smooth map of solid Huber rings. Then $f$ is solid smooth. Indeed, it is by definition of solid finite presentation, and the explicit presentation of $B$ as $A$-algebra shows that $\bb{L}_{B/A}$ is a finite free $B$-module generated by the differentials of the variables $T_j$ for $j>k$. If in addition $f$ is solid standard \'etale then the same presentation shows that $\bb{L}_{B/A}=0$. 
\end{remark}

Our next goal is to characterize locally solid smooth (resp. \'etale) maps of solid Huber rings as those maps which are, locally after rational localizations, solid standard smooth (resp. \'etale).

\begin{proposition}\label{PropStabilityolidEtaleSmooth}
The following holds: 

\begin{enumerate}

\item Let $f\colon A\to B$ be a standard rational localization of solid Huber rings. Then $f$ is solid standard \'etale.

\item Let $f\colon A\to B$ be a solid standard smooth (resp. \'etale) morphism of solid Huber rings. Let $A\to A'$ be a morphism of solid Huber rings, then the base change $A'\to B'=B\otimes_A B'$ is solid standard  smooth (resp. \'etale).   The same holds for (locally) solid smooth  (\'etale) maps.

\item  Let $A\to B$ and $B\to C$ be solid smooth maps of solid Huber rings (resp. locally solid smooth, resp. solid standard smooth), then so is the composite $A\to C$. The same holds for ``\'etale'' instead of ``smooth''.  In particular, rational localizations are solid standard \'etale.

\item Let $A\xrightarrow{f} B\xrightarrow{g} C$ be a morphism of  solid  Huber rings with $f$ and $g\circ f$ (locally) solid \'etale, then $g$ is (locally) solid \'etale.

\item Let $f\colon A\to B$ be a locally solid \'etale morphism  of solid Huber rings that is surjective on $\pi_0$ and has the induced analytic ring structure. Then there is an idempotent element $f\in A$ such that $B=A/(1-f)=A(\frac{1}{f})_{\sol}$. In particular, if $A\to B$ is a locally solid \'etale map, the multiplication map $B\otimes_A B\to B$ is the quotient by an idempotent element in $B\otimes_A B$.

\item Let $A\xrightarrow{f} B \xrightarrow{g} C$ be morphisms of solid Huber rings with $f$ and $g\circ f$  locally solid smooth. Suppose that $g$ has the induced structure and is surjective on $\pi_0$. Then $C$ is, locally in the analytic topology of $B$, a  complete intersection, i.e. a morphism of the form 
\[
C=B/^{\bb{L}}(f_1,\ldots, f_k)
\]
for suitable elements $f_1,\ldots, f_k\in B$.  In particular, $C^{\triangleright}$ is a dualizable $B$-module and $g$ is suave. 
\end{enumerate}
\end{proposition}
\begin{proof}
\begin{enumerate}

\item Let $f\colon A\to B$ be a standard rational localization. By definition, there are elements $(f_1,\ldots, f_n,g)$ in $A$ generating the unit ideal such that $B=A[T_1,\ldots, T_n]_{\sol}/^{\bb{L}}(gT_i-f_i)$. The Jacobian matrix of  the presentation is given by the diagonal matrix with value $g$, since $g$ is invertible in $B$, we deduce that $f$ is standard solid \'etale. 

\item Stability under base change of standard solid  \'etale and smooth maps is obvious from the definition.  The case of (locally) solid smooth  (\'etale) maps is also clear. 

\item Let $A\to B$ and $B\to C$ be standard solid smooth maps. We have presentations 
\[
B=A[T_1,\ldots, T_n]_{\sol}/^{\bb{L}}(f_1,\ldots, f_k)
\]
and 
\[
C=B[X_{1},\ldots, X_m]_{\sol}/^{\bb{L}} (g_1,\ldots, g_l)
\]
such that the Jacobians $b=\det((\frac{\partial f_i}{ \partial T_j})_{1\leq i,j\leq k})$ and $c=\det((\frac{\partial g_i}{ \partial X_j})_{1\leq i,j\leq l})$ are invertible  on $B$ and $C$ respectively. Finding a lift $\widetilde{g}_j$ of $g_j$ from $C$ to $A[T_1,\ldots, T_n,X_1,\ldots, X_m ]_{\sol}$, we see that 
\begin{equation}\label{eqojapemdfawdfqef}
C=A[T_1,\ldots, T_n,X_1,\ldots, X_m]_{\sol}/^{\bb{L}}(f_1,\ldots, f_k,\widetilde{g}_1,\ldots, \widetilde{g}_l).
\end{equation}
Let $\widetilde{c}$ be the Jacobian associated to $(\frac{\partial \widetilde{g}_i}{\partial X_j})_{1\leq i,j\leq l}$, then the Jacobian of the presentation \eqref{eqojapemdfawdfqef} is $a\widetilde{c}$ which becomes invertible on $C$, proving that $A\to C$ is also solid standard smooth. It is clear from the construction that if both $A\to B$ and $B\to C$ are solid standard \'etale then so is $A\to C$. 

The case of solid smooth maps follows from  \Cref{LemmaPermanenceSolidPresentation} (2) and the fiber sequence of cotangent complexes (resp. for solid \'etale). The case of locally solid smooth (resp. locally solid \'etale) follows from the case of solid smooth (resp. solid \'etale) after localizing with rational localizations in both $A$, $B$ and $C$.

\item The statement is local in the analytic topology of $A$, $B$ and $C$. Therefore, we can assume without loss of generality that both $f$ and $g\circ f$ are solid \'etale. By \Cref{LemmaPermanenceSolidPresentation} we know that $g$ is of solid finite type, the fiber sequence of cotangent complexes shows that $\bb{L}_{C/B}=0$ proving that $g$ is solid finite \'etale.

\item Let $f\colon A\to B$ be a locally solid \'etale morphism with the induced structure that surjects on $\pi_0$. By (1)-(3) above, we know that the multiplication map $B\otimes_A B\to B$ is locally solid \'etale, in particular $\bb{L}_{B/B\otimes_A B}=0$. Since the map induces an isomorphism on $\pi_0$, \Cref{CoroDetectCotangent} shows that it is an isomorphism and so $A\to B$ is an idempotent map. We claim that $B$ is dualizable as $A$-module, by \Cref{LemmaIdempotentDualizable} down below this implies that $A\to B$ is also an open immersion and therefore that $A=B\times \ob{fib}(A\to B)$ as algebras. In particular, taking $f=(1,0)$ we get the desired result. To prove that $B$ is dualizable as $A$-module, we can equivalently prove that the morphism $A\to B$ is suave.  Being suave can be proven locally along rational covers of $A$ (and hence of $B$ by \Cref{LemmaRationalCoversClosedImmersions}). In particular, we can assume without loss of generality that $f$ is solid \'etale.  Then, there is a morphism of the form 
\[
A[T_1,\ldots, T_n]_{\sol}/^{\bb{L}}(f_1,\ldots, f_n)\to B
\]
 that is an isomorphism on $\pi_0$. Therefore, the base change
 \begin{equation}\label{eqpakspqmawoqnaeofa}
 A[T_1,\ldots, T_n]_{\sol}/^{\bb{L}}(f_1,\ldots, f_n)\to B[T_1,\ldots, T_n]_{\sol}/^{\bb{L}}(f_1,\ldots, f_n)
 \end{equation}
is solid \'etale and an isomorphism on $\pi_0$ (as $A\to B$ is surjective). Since the relative cotangent complex of $A\to B$ vanishes, the same happens for \eqref{eqpakspqmawoqnaeofa} and it   is an isomorphism by \Cref{CoroDetectCotangent}. On the other hand, the composite map 
\[
A\to A[T_1,\ldots, T_n]_{\sol}/^{\bb{L}}(f_1,\ldots, f_n)=B[T_1,\ldots, T_n]_{\sol}/^{\bb{L}}(f_1,\ldots, f_n)
\] 
is suave being a composition of a  solid polynomial algebra and a derived complete intersection. Let $b_1,\ldots, b_n$ be the images of the $T_i$ in $B$, then we have that 
\[
\bigg(B[T_1,\ldots, T_n]_{\sol}/^{\bb{L}}(f_1,\ldots, f_n) \bigg)/^{\bb{L}}(T_i-b_i)= B/^\bb{L}(0,\ldots, 0)
\] 
is still suave over $A$, where in the right term there are $n$-th occurrences of $0$.  Therefore, as $B$ has the induced structure over $A$, $ B/^\bb{L}(0,\ldots, 0)$  a dualizable as $A$-module, and since $B$ is a retract of $B/^\bb{L}(0,\ldots, 0)$ then so is $B$, finishing the proof of the claim.

The claim about the multiplication map $B\otimes_A B\to B$ for a locally solid \'etale morphism $A\to B$ follows from the fact that $B\otimes_A B\to B$ is locally solid \'etale by (4), it has the induced structure,  and that it surjects on $\pi_0$.

\item By assumption both $f$ and $g\circ f$ are locally solid smooth. Note that since $B\to C$ is surjective on $\pi_0$,  \Cref{LemmaRationalCoversClosedImmersions} shows that any rational localization of $C$ arises from a base change of a rational localization of $B$. Thus, after localizing both on $A$ and $B$, we can assume without loss of generality that both $f$ and $g\circ f$ are of solid finite type. By \Cref{LemmaPermanenceSolidPresentation}   (2) the map $B\to C$ is also of solid finite type. Looking at cotangent complexes, we have a fiber sequence 
\[
C\otimes_B \bb{L}_{B/A}\to \bb{L}_{C/A} \to \bb{L}_{C/B}.
\]

 Since $f$ and $g\circ f$ are solid smooth, $\bb{L}_{B/A}$ and $\bb{L}_{C/A}$ are finite projective, and therefore $\bb{L}_{C/B}$ is a perfect complex proving that $g$ is of solid finite presentation. Furthermore, since $g$ surjects on $\pi_0$, by \Cref{PropConnectivityCotangent} the module $\bb{L}_{B/A}$ is $1$-connective, and so $C\otimes_B \bb{L}_{B/A}\to \bb{L}_{C/A} $ is surjective on $\pi_0$. Since $\bb{L}_{C/A}$ is a vector bundle over $C$,  it is finite projective and we can find a section $s\colon \bb{L}_{C/A}\to C\otimes_B \bb{L}_{B/A}$, and because $C\otimes_B \bb{L}_{B/A}$ is also a $C$-vector bundle we deduce that $\bb{L}_{C/B}[-1]=\ob{fib}(C\otimes_B \bb{L}_{B/A}\to \bb{L}_{C/A})$ is a vector bundle over $C$. 
 
We can find a finite  Zariski cover $\{C\to  C[\frac{1}{g_i}]\}_{i=1}^n$ of $C$,  such that $\bb{L}_{C/B}[-1]$  admits a basis. By \Cref{LemmaRatLoc} the family of rational localizations $\{C\to C\left( \frac{g_1,\ldots, g_n}{g_i} \right)\}_{i=1}^n$ forms a cover of $\Spec^{\sol}(C)$. Thus, lifting the cover of $C$ to rational localizations of $B$ via \Cref{LemmaRationalCoversClosedImmersions}, we can assume without loss of generality that $\bb{L}_{C/B}[-1]$ admits a $C$-basis, say $v_1,\ldots, v_k$.   By taking loops of  \Cref{PropConnectivityCotangent}, we have a surjective map 
\[
\ob{fib}(B\to C)\otimes_B C\to \bb{L}_{C/B}[-1].
\]
If $I=\ker(\pi_0(B)\to \pi_0(C))$, this produces a right exact sequence 
\[
I/I^2 \xrightarrow{d} \pi_0(C\otimes_B \bb{L}_{B/A}) \to \pi_0(\bb{L}_{C/A})\to 0. 
\]
where the first map is just the differential of $B$ restricted to $I$. Let $f_1,\ldots, f_k\in I$ be elements such that $df_i=v_i\in \pi_0(\bb{L}_{C/B}[-1])$. We have an induced map of $A$-algebras
\[
g\colon \widetilde{B}:=B/^{\bb{L}}(f_1,\ldots, f_k)\to C .
\] 
The map $g$ has the induced structure and is surjective on $\pi_0$, it is also of solid finite type (\Cref{LemmaPermanenceSolidPresentation} (2)) and by construction the relative cotangent complex of $g$ vanishes. In particular, $g$ is a solid \'etale morphism, and thanks to part (5) there is an idempotent element $\widetilde{h}\in \widetilde{B}$ such that $C=\widetilde{B}/^{\bb{L}}(1-\widetilde{h})$. Lifting $\widetilde{h}$ to an element $h$ in $B$, we see that 
\[
C=B/^{\bb{L}}(f_1,\ldots, f_k,1-h)
\]
proving what we wanted. 
\end{enumerate}
\end{proof}

The following lemma was used before.

\begin{lemma}\label{LemmaIdempotentDualizable}
Let $\n{C}$ be a  symmetric monoidal stable  category and $1\to A$ an idempotent morphism such that $A$ is dualizable in $\n{C}$. Then $A$ is naturally self dual, and the dual $A\to 1$ of the unit map $1\to A$ endows $A$ with the structure of a co-idempotent coalgebra. In particular $\ob{fib}(1\to A)=\ob{cofib}(A\to 1)$ has the structure of idempotent algebra and $1=A\times \ob{fib}(1\to A)$ as algebras. 
\end{lemma}
\begin{proof}
Given a dualizable object $P\in \n{C}$ we denote by  $P^{\vee}$ its dual. Consider the pullback $\iota^*\colon \n{C}\to \Mod_A(\n{C})$, the map $\iota^*$ has a fully faithful linear right adjoint given by the forgetful functor $\iota_*$. Since $A$ is dualizable, $\iota^*$ also has a fully faithful left adjoint $\iota_!$ given by $A^{\vee}\otimes -=\iHom_{\n{C}}(A,-)$. But then the idempotency of $A$ implies  that 
\[
A^{\vee}\otimes A= \iHom_{\n{C}}(A,A) = \iHom_{A}(A,A) =  A.
\] 
Thus, we see that for $M$ an $A$-module one has  
\[
\iota_! M = A^{\vee}\otimes M = A^{\vee}\otimes A \otimes M = M. 
\]
In other words, there is a natural equivalence $\iota_!=\iota_*$ proving that $A$ is naturally self-dual. In particular, as $A$ is an idempotent algebra, the dual of the unit $1\to A$ gives rise to a trace map $A\to 1$ that endows $A$ is a co-idempotent co-algebra structure. By the previous computation one sees that the categories of modules and comodules of $A$ agree, as full subcategories of $\n{C}$:
\[
\Mod_{A}(\n{C})= \ob{cMod}_{A}(\n{C}).
\]
Passing to the cofiber $B=\ob{cofib}(A\to 1)$, one sees that $B$ is a dualizable idempotent algebra and by the  same argument as before it is self dual,  the dual of the unit map $1\to B$ endows $B$ with a  natural  co-idempotent co-algebra structure so that $\Mod_{B}(\n{C})= \ob{cMod}_{B}(\n{C})$. From the orthogonal decomposition of $\n{C}$ in terms of  $\ob{cMod}_A(\n{C})$ and $\Mod_B(\n{C})$, one deduces a direct product decomposition 
\[
\n{C}=\Mod_{A}(\n{C}) \times  \Mod_B(\n{C})
\]
which in particular shows that $1=A\times B$ as algebras. From this decomposition it is clear that $B=\ob{fib}(1\to A)=\ob{cofib}(A\to 1)$ proving what we wanted. 
\end{proof}

\begin{theorem}\label{TheoremSolidSmoothEtaleMaps}
Let $f\colon A\to B$ be a morphism of solid Huber rings. The following are equivalent: 
\begin{enumerate}

\item  $f$ is locally solid smooth.

\item Locally in the analytic topology of $A$ and $B$, $f$ is solid standard smooth. 

\end{enumerate}
 
The analogue statement holds for ``\'etale'' instead of  ``smooth''. 
\end{theorem}
\begin{proof}
It is obvious form the definitions and \Cref{PropStabilityolidEtaleSmooth} (1) that (2) implies (1) in  the Theorem. Conversely, let $f\colon A\to B$ be a locally solid smooth morphism of solid Huber rings. By taking rational localizations on both $A$ and $B$, we can assume without loss of generality that $f$ is a solid smooth map. In particular, there is a morphism of $A$-algebras 
\[
A[T_1,\ldots, T_n]_{\sol}\to B
\]
with the induced structure that surjects on $\pi_0$. Thanks to \Cref{PropStabilityolidEtaleSmooth} (3) rational localizations of $A[T_1,\ldots, T_n]_{\sol}$ are standard smooth over $A$.  By  \Cref{PropStabilityolidEtaleSmooth} (6), after passing to rational covers  of $A[T_1,\ldots, T_n]$ and $B$, we can find a standard smooth algebra $A\to C$ and a presentation $B=C/^{\bb{L}}(\overline{f}_1,\ldots, \overline{f}_k)$. We will show that we can modify $C$ and the sequence $(\overline{f}_1,\ldots, \overline{f}_k)$ so that $B$ becomes standard \'etale over $A$. For that, write a standard smooth presentation $C=A[T_1,\ldots, T_n]_{\sol}/^{\bb{L}}(g_1,\ldots, g_l)$, and take lifts $f_1,\ldots, f_k$ of  $(\overline{f}_1,\ldots, \overline{f}_k)$ to $A[T_1,\ldots, T_n]_{\sol}$. We have that 
\begin{equation}\label{eqPresentationsolidSMooth}
B=A[T_1,\ldots, T_n]_{\sol}/^\bb{L}(g_1,\ldots, g_l, f_1,\ldots, f_k). 
\end{equation}
Consider the fiber sequence of cotangent complexes
\[
 \bigoplus_{i=1}^l B dg_i \oplus \bigoplus_{j=1}^{k} df_j  \to  \bigoplus_{i=1}^n B dT_i \to \bb{L}_{B/A}.
\]
Since $A\to B$ is solid smooth, $\bb{L}_{B/A}$ is a vector bundle, which forces $l+k\leq n$. We can find a Zariski cover of $B$  of the form $\{B\to B[\frac{1}{h_s}]\}_s$ such that a subset of the $\{dT_i\}_{i=1}^n$ form a basis of $\bb{L}_{B/A}$.  Passing to a rational cover as in \Cref{LemmaRatLoc} with the elements $h_s$,  lifting this to a rational cover of $A[T_1,\ldots, T_n]_{\sol}$, and rearranging the indices of the variables,   we can assume that the elements $dT_{k+l+1},\ldots, dT_{n}$ give rise to a basis of $\bb{L}_{B/A}$. Thus, the composite $\bigoplus_{i=1}^l B dg_i \oplus \bigoplus_{j=1}^{k} df_j  \to \bigoplus_{i=1}^{l+k} B dT_i$ must be an isomorphism, and the Jacobian matrix of the presentation \eqref{eqPresentationsolidSMooth} is invertible on $B$, proving that $B$ is solid standard smooth on $A$ as wanted.  Notice that if $B$ is solid \'etale over $A$, by counting ranks in the contangent complexes one has $k+l=n$, proving that $B$ is solid standard \'etale. 
\end{proof}

\begin{remark}\label{RemarkGlobalizationSolidAdicSpaces}
The discussion of \Cref{ss:FinitePresentation,ss:SolidEtaleSmooth} admits an obvious generalization to solid adic spaces by gluing along rational localizations of solid Huber rings. More precisely, given $f\colon Y\to X$ is a morphisms of solid adic spaces we say that $f$ is \textit{locally of finite type} (resp. \textit{finite presentation} or \textit{almost finite presentation}) if it is locally on affinoids after passing to rational localizations of $Y$ and $X$. The stability properties of these morphisms in \Cref{LemmaPermanenceSolidPresentation} imply that this is a well defined notion stable under base change and composition. Then, we say that $f$ is \textit{locally solid smooth} (resp. \textit{solid \'etale}) if locally after a rational cover in $Y$ and $X$ it is represented by a solid smooth (resp. \'etale) morphism of solid Huber rings.

 If $f$ is a morphism of solid adic spaces, there is a well defined cotangent complex $\bb{L}_f\in \ob{D}(X)$ that locally on affinoids is the cotangent complex of \Cref{DefCotangentComplex}, namely, this follows formally from the fact that rational localizations are formally \'etale by \Cref{PropStabilityolidEtaleSmooth} (1)-(3), and by the base change properties of the cotangent complex.  However, it is a priori not clear that $\bb{L}_f$ satisfies some global deformation theoretic description, one of the problems being that there is not a well defined $t$-structure on quasi-coherent sheaves on analytic stacks. Nevertheless, the map $f$ locally of finite type is locally solid smooth (resp. locally solid \'etale) if and only of $\bb{L}_f$ is, locally in the solid analytic topology, a vector bundle (resp. $0$).
 
 Finally, 
\end{remark}

\subsection{Serre duality}\label{SS:SerreDuality}

In this final section we give a proof of Serre duality for solid smooth morphisms of solid Huber rings following the arguments of Clausen and Scholze \cite{ClausenScholzeCondensed2019,CondensedComplex} (in particular, all the relevant ideas and constructions are due to them).  Of course, the theorem generalizes to locally solid smooth morphisms of solid adic spaces after the obvious modifications as in  \Cref{RemarkGlobalizationSolidAdicSpaces}.

The theorem  splits in two statements: the fact that solid smooth morphisms are cohomological smooth, and the computation of the dualizing sheaf. Since both require different techniques it is natural to split the statement into two different theorems. 

\begin{theorem}\label{TheoLocalSerreDuality}
Let $f\colon A\to B$ be a locally solid smooth morphism of solid Huber rings. Then $f$ is cohomologically smooth and the dualizing sheaf $f^! 1$ is, locally in the analytic topology, a line bundle concentrated in homological degree $\dim(f)$. Moreover, if $f$ is locally solid \'etale then it is cohomologically \'etale.  
\end{theorem}
\begin{proof}
Since rational localizations form an open cover of analytic stacks, the statement is local for the topology of the solid adic spectrum of $A$ and $B$. By \Cref{TheoremSolidSmoothEtaleMaps} we can assume that $B$ is standard solid smooth over $A$:
\[
B=A[T_1,\ldots, T_n]_{\sol}/^{\bb{L}} (f_1,\ldots, f_k). 
\]
Then, by base change the fact that $A\to B$ is cohomologically smooth with invertible sheaf in the desired degree follows from the suaveness of  the map $\bb{Z}_{\sol}\to \bb{Z}[T]_{\sol}$ from \Cref{PropSerreDualityAffine} with dualizing sheaf isomorphic to $\Z[T][1]$, and the suaveness for the  map $\bb{Z}[T]\to \Z[T]/T=\Z$ from \Cref{LemmaPrimAffine}  with dualizing sheaf isomorphic to $\iHom_{\bb{Z}[T]}(\Z,\bb{Z}[T])\cong \bb{Z}[-1]$.  

If in addition $f$ is locally solid \'etale, by \Cref{PropStabilityolidEtaleSmooth} (5) we know that the diagonal of $\AnSpec(B)\to \AnSpec(A)$ is an open immersion, thus cohomologically  \'etale, and therefore $f$ is cohomologically \'etale by \cite[Definition 4.6.1]{HeyerMannSix}.
\end{proof}

In order to compute the dualizing sheaf, we will use a deformation to the normal cone  following Clausen and Scholze \cite{CondensedComplex}. This will reduce the computation of the dualizing sheaf to the case of vector bundles which we handle explicitly in \Cref{PropositionUniversalNormalCone}. 

To warm up, we start with the computation of the dualizing sheaf of vector bundles (or rather of the dualizing sheaf of its zero section). Let $\ob{Vect}_n$ be the algebraic stack classifying vector bundles of rank $n$. Let $\Spec \Z\to \ob{Vect}_n$ be the map associated to the vector bundle $\bb{A}^{n}_{\Z}\to \Z$,  the standard action of $\GL_n$ on $\bb{A}^{n}_{\Z}$ gives rise to an identification $\ob{Vect}_{n}\cong B\GL_n$. 

\begin{prop}\label{PropositionUniversalNormalCone}
Let $V$ be the universal vector bundle over the stack $S:=\ob{Vect}_{n}$ of rank $n$ with zero section $\iota_0\colon S\to V$. Then $\iota_0$ is cohomologically smooth and $\iota_0^{!} 1 [n]$ is a line bundle on $\ob{Vect}_{n}$ giving rise to a map of algebraic stacks 
\[
\delta \colon  \ob{Vect}_n\to \ob{Vect}_1. 
\]
Furthermore, the datum of a commutative triangle 
\begin{equation}\label{eqTrivializationInvertibleSheaf}
\begin{tikzcd}
 & \Spec \Z \ar[rd] \ar[ld] & \\ 
  \ob{Vect}_n \ar[rr,"\delta"] & & \ob{Vect}_1
\end{tikzcd}
\end{equation}
filling the fixed maps $\Spec \Z\to \ob{Vect}_n$ and $\Spec \Z\to \ob{Vect}_{1}$ is  the anima of trivializations of $\widetilde{\iota}_0^![n]$ where $\widetilde{\iota}_0\colon \Spec \Z \to \bb{A}^n_{\Z}$ is the zero section. In particular, the anima of such triangles is a $\Z^{\times}=\{\pm 1\}$-torsor (and therefore a set). Once the triangle \eqref{eqTrivializationInvertibleSheaf} has been fixed, the  map $\delta$  gives rise to a morphism of pointed classifying stacks 
\[
B\GL_n\to B \bb{G}_{m}
\]
which, by delooping,  is equivalent to the morphism of groups given by the determinant $\det\colon \GL_n\to \bb{G}_m$. In particular, if $V$ is the universal vector bundle of $\ob{Vect}_n$, there is a canonical isomorphism 
\[
\iota_0^! 1[n] = \det V= \bigwedge^n V,
\]
depending only on the choice of the triangle \eqref{eqTrivializationInvertibleSheaf} (i.e. up to a sign). 
\end{prop}
\begin{proof}
Suaveness of the map $\iota_0$ is clear: this can be proven after pullback along the map $\Spec \Z \to \ob{Vect}
_n$ where it reduces to cohomological smoothness of the complete intersection $0\in \bb{A}^n_{\Z}$ which then follows form \Cref{LemmaPrimAffine}. The   fact that $\iota_0^! 1[n]$ is a line bundle can be proven after pullback to $\Spec \Z$, and this follows from the fact that $\widetilde{\iota}_{0,*}=\widetilde{\iota}_{0,!}$ being an affine map and the computation  $\iHom_{\Z[T_1,\ldots, T_n]} (\Z, \Z[T_1,\ldots, T_n])\cong \Z[-n]$ via a Koszul complex. 

Now, the object $\iota_0^{!}1[n]$ gives rise to a map of algebraic stacks 
\[
\delta\colon \ob{Vect}_n\to \ob{Vect}_1. 
\]
Since $\ob{Vect}_n\cong B\GL_n$ and $\ob{Vect}_1\cong B\bb{G}_m$ are classifying stacks of groups, the anima of  morphisms $\ob{Vect}_n\to \ob{Vect}_1$ is equivalent to the anima of triangles    \eqref{eqTrivializationInvertibleSheaf}  up to conjugation by $\bb{G}_m(\Z)=\{\pm 1\}$.  On the other hand,  by delooping,   the anima of pointed maps $B\GL_n\to B\bb{G}_m$ is equivalent to  anima of group homomorphisms $\GL_n\to \bb{G}_m$. Since $\bb{G}_m$ is corepresented by a static ring,  the anima of group homomorphisms $\GL_n\to \bb{G}_m$ is a set, and it is given precisely by the group $(\det)^{\Z}$ where $\det\colon \GL_n\to \bb{G}_m$ is the determinant map.  We deduce that the anima of maps $\ob{Vect}_n\to \ob{Vect}_1$ is isomorphic to the groupoid  $(\det)^{\Z}/ \{ \pm 1\}$, in particular, the choice of a triangle \eqref{eqTrivializationInvertibleSheaf}  is a $\{\pm 1\}$-torsor as wanted. 

Let us fix a triangle \eqref{eqTrivializationInvertibleSheaf}, then the map $\delta$ arises from a morphism of groups $(\det)^k \colon \GL_n\to \bb{G}_m$, we are left to show that $k=1$. Let $\bbf{T}\subset \GL_n$ be the diagonal torus, it suffices to show that the composite  $\bbf{T}\to \GL_n \to \bb{G}_m$ is given by the determinant. Furthermore, by  writing $\bbf{T}=\prod_{i=1}^n \bb{G}_m$ as a product of multiplicative groups via the diagonal matrix entries,   over $B\bbf{T}= \prod_{i=1}^n B\bb{G}_m$ we have a decomposition of the universal vector bundle $V= \prod_{i=1}^n L_i$, where $L_i$ is the universal line bundle associated to the $i$-th component. Since the formation of the dualizing sheaf is natural with respect to fiber products (this is just the functoriality of the symmetric monoidal structure after passing to adjoints in the category of kernels), we are reduced to the case $n=1$ and want to show that the map 
\[
\delta\colon B\bb{G}_m\to B\bb{G}_m
\]
is the identity. 

It is now when we do an explicit computation: the universal line bundle over $B\bb{G}_m$ is the quotient $\bb{A}^1_{\Z}/\bb{G}_m\to B\bb{G}_m$ where $\bb{G}_m$ acts by multiplication.  The category $\ob{D}(B\bb{G}_m)$ is naturally equivalent to the category of graded $\Z$-modules, and  $\ob{D}(\bb{A}^1_{\Z}/\bb{G}_m)$ is the category of $\Z[t]$-modules on $\ob{D}(B\bb{G}_m)$  where $t$ is a variable with degree $-1$  (see \Cref{RemarkFileredRingsAsSpec}). Then, using the Koszul resolution 
\[
0\to t\Z[t]\to \Z[t] \to \Z\to 0, 
\]
we have that 
\[
\iota^!_0 1 = \iHom_{\Z[t]}(\Z,\Z[t]) = t^{-1}\Z[-1]
\]
proving that $\iota^!_0 1 [1]$ is given by the line bundle over $B\bb{G}_m$ of degree $1$, this  precisely  corresponds to the identity morphism $B\bb{G}_m\xrightarrow{\id}B\bb{G}_m$, finish the proof of the proposition. 
\end{proof}

\begin{remark}\label{RemarkNaturalityDualizingSheaf}
There are two main contents in  \Cref{PropositionUniversalNormalCone}, the first is that the (shifted) dualizing sheaf of the zero section of a vector bundle is given by the determinant of such, and the second is that this identification is canonical up to a sign (or equivalently, the choice of a triangle \eqref{eqTrivializationInvertibleSheaf}). It is natural to expect that the functor $\delta$ of \Cref{PropositionUniversalNormalCone} can be promoted to a map of commutative monoids $\delta\colon \ob{Vect}\to \ob{Vect}_1$  where $\ob{Vect}$ is the algebraic stack of vector bundles with monoid structure given by direct sum, and that, after fixing a sign, that the map $\delta$ can be naturally identified with the determinant map. 
\end{remark}

We continue with the proof of Serre duality.  For that,  we recall the discussion of pairs and the $I$-adic filtration of \cite[Section 3]{MaoCrystalline}, and adapt it to analytic rings.

\begin{definition}\label{DefinitionPair}
Let $\Ring^{\delta}$ be the category of discrete animated rings.  The category of \textit{pairs} is the full subcategory $\cat{Pair}\subset \ob{Fun}(\Delta^1,\Ring^{\delta})$ consisting on those morphisms $A\to B$ that are surjective on $\pi_0$. Given a pair $f\colon A\to B$, we let $I=\ob{fib}(A\to B)$ be the \textit{ideal of $f$}. 
\end{definition}

We recall the definition of filtered animated rings; for a discussion about filtered and graded categories  we refer to \cite[Section 2.2]{BhattGauges} and the references therein. 

\begin{definition}\label{DefinitionFilteredAnimatedRing}
Let $X$ be an fpqc algebraic stack. Informally, the category of \textit{animated rings over $X$}, denoted by $\Ring_{X}$, is defined as the limit 
\[
\Ring_{X}=\varprojlim_{\Spec R \to X} \Ring^{\delta}_{R/}
\]
of the categories of discrete animated $R$-algebras, where the index diagram  is the category of affine schemes over $X$. Formally, the functor sending $\Spec R\mapsto \Ring^{\delta}_{R/}$ satisfies fqpc descent and $\Ring_{-}$ is the right Kan  extension to an fpqc sheaf on algebraic stacks.  Equivalently, $\Ring_{X}$ is the opposite of the full subcategory of morphisms $Y\to X$ of  fpqc algebraic stacks over $X$ that are represented by affine schemes, i.e., those morphisms $Y\to X$ whose pullback along all map $\Spec R\to X$ from an affine scheme  is an affine scheme. Given $\s{R}\in \Ring_{X}$, we let $\Spec_X(\s{R})\to X$ be its associated algebraic stack over $X$, and call it the \textit{relative spectrum of $X$}.

We define the category of \textit{filtered animated rings} (with decreasing filtration) to be $\Fil(\Ring):=\Ring_{\bb{A}^{1}_{\Z}/\bb{G}_{m}}$, where in $\bb{A}^{1}_{\Z}/\bb{G}_{m}$ the action  of  $\bb{G}_m$ is by multiplication. Similarly, we define the category of \textit{graded animated rings}  to be $\ob{Gr}(\Ring):=\Ring_{B\bb{G}_{m}}$. 
\end{definition}

\begin{remark}\label{RemarkFileredRingsAsSpec}
Let $\Z[t]$ be the polynomial algebra over $B\bb{G}_m$ where $t$ is a  variable  of degree $-1$, then there is a natural equivalence 
\[
\Spec_{B\bb{G}_m}(\Z[t])=\bb{A}^1_{\Z}/\bb{G}_m.
\]
In particular, we have a symmetric monoidal equivalence 
\[
\ob{D}(\bb{A}^1_{\Z}/\bb{G}_m)=\Mod_{\Z[t]}(\ob{D}(B\bb{G}_m)). 
\]
\end{remark}

\begin{remark}\label{RemarkSymmetricAlgebras}
Let $X$ be an fpqc algebraic stack and let $\ob{D}_{\geq 0}(X)$ be the category of connective quasi-coherent sheaves. There is a forgetful functor $G\colon \Ring_{X}\to \ob{D}_{\geq 0}(X)$. By descent from the affine case, the functor $G$ admits a left adjoint that we denote $\Sym_{X}^{\bullet}$; it is nothing but the descent of the symmetric powers from the affine situation. Since $G$ preserves colimits and is conservative (again by descent from the affine case), we  obtain by the monadicity theorem
\[
\Ring_X = \LMod_{\Sym_X^{\bullet}}(\ob{D}_{\geq 0}(X)). 
\]
\end{remark}

\begin{example}\label{ExampleReesConstruction}
Let us explain the Rees construction as in \cite[Proposition 2.2.6]{BhattGauges} in terms of algebraic stacks. Given a filtered $\bb{Z}$-module $\Fil^{\bullet}M\in \ob{D}(\bb{A}^1_{\Z}/\bb{G}_m)$ the \textit{Rees module of $\Fil^{\bullet}M$} is the pullback along $\bb{A}^1_{\Z}\to \bb{A}^1_{\Z}/\bb{G}_m$. Informally, the Rees module is given by 
\[
\ob{Rees}(\Fil^{\bullet} M) = \bigoplus_{n\in \Z} \Fil^n{M} t^{-n}
\]
where $t$ is the coordinate function of $\bb{A}^1_{\Z}$ with degree $-1$. 
\end{example}

The category of  discrete animated rings has compact projective generators given by the full subcategory $\cat{Poly}^{\ob{ft}}\subset \Ring^{\delta}$ of polynomials over $\bb{Z}$ generated by finitely many variables. Colimits on $\cat{Pair}$ are computed pointwise, and the pairs $\bb{Z}[X]\xrightarrow{X\mapsto 0} \bb{Z}$  and $\bb{Z}[X]\xrightarrow{\id} \bb{Z}[X]$ corepresent the functor sending  $A\to B$ to $I$ and $A$ respectively. Thus, $\cat{Pair}$ is a compactly projective generated category, with a family of compact projective generators given by the full subcategory $\cat{Pair}^{\ob{pol}}\subset \cat{Pair}$ consisting of the pairs $\bb{Z}[\underline{X},\underline{Y}]\xrightarrow{\underline{X}\mapsto 0} \bb{Z}[\underline{Y}]$. where $\underline{X}$ and $\underline{Y}$ are finite sets of variables.  

Given a pair $f\colon A\to B$ with ideal $I$ one can construct the \textit{$I$-adic filtration of $f$} as follows.  We say that $f$ is \textit{Koszul regular} if $A$ and $B$ are static, and $I=\ker(A\to B)$ is generated by a finite Koszul regular sequence. All the pairs in $\cat{Pair}^{\ob{pol}}$ are  Koszul regular.  Given a Koszul regular pair $A\to B$ with ideal $I$, one can form the $I$-adic filtration of $A$ given by $(I^n)_{n\in \N}$. The passage to the $I$-adic filtration is functorial on morphisms of Koszul regular pairs, and gives rise to a functor 
\[
\Fil_{\ob{ad}}^{\bullet}\colon \cat{Pair}^{\ob{pol}}\to \Fil(\Ring^{\delta})
\]
on filtered animated rings. An easy computation shows that $\Fil^{\bullet}_{\ob{ad}}$ preserves pushouts.  Extending by sifted colimits, i.e. passing to the animation, we have constructed the \textit{$I$-adic filtration functor}
\[
\Fil_{\ob{ad}}^{\bullet} \colon \cat{Pair} \to  \Fil(\Ring^{\delta}). 
\]
By construction, the functor $\Fil_{\ob{ad}}^{\bullet}$ commutes with sifted colimits, and thanks to \cite[Proposition 5.5.8.15 (3)]{HigherTopos} it also commutes with all colimits.

\begin{example}\label{ExampleIadicFiltration}
Ler $R$ be a discrete ring. The datum of a map $f\colon \Spec R\to \bb{A}^1_{\Z}/\bb{G}_{m}$  is the same as a \textit{virtual Cartier divisor}, that is a morphism  $I\to R$ where $I$ is line bundle over $R$. The map $f$ is represented in affine schemes and therefore  gives rise to a filtered animated ring which is nothing but $\bigoplus_{n\in \Z} I^nt^{-n}$ (with $t$ the coordinate of $\bb{A}^1_{\Z}$ with degree $-1$). 
\end{example}

\begin{lemma}\label{LemmaPositivelyFiltered}
Let $\Fil^{\geq 0}(\Ring^{\delta})\subset \Fil(\Ring^{\delta})$ the full subcategory of \textit{positively filtered animated rings}, that is, those filtered animated rings $\Fil^{\bullet} A$ such that for $n\geq  0$ the map $\Fil^{0}A\to \Fil^{-n}A$ is an isomorphism.  Then the inclusion $\Fil^{\geq 0}(\Ring^{\delta})\subset \Fil(\Ring^{\delta})$ admits a right adjoint sending a filtered animated ring $\Fil^{\bullet} A$ to its positive truncation $\Fil^{\geq 0} A$.  
\end{lemma}
\begin{proof}
Let $\ob{D}_{\geq 0}(\bb{A}^1_{\Z}/\bb{G}_{m})$ be the category of connective filtered $\Z$-modules, considered as graded $\Z[t]$-modules where $t$ has degree $-1$. The category  $\ob{D}_{\geq 0}(\bb{A}^1_{\Z}/\bb{G}_{m})$ is compactly projective generated by the objects $\bigoplus_{i=1}^s \Z[t] t^{k_i}$ where $k_i\in \Z$. Hence, $\Fil(\Ring^{\delta})$ is also compactly projectively generated by the graded $\Z[t]$-algebras $\Z[t][X_1,\ldots,X_s]$ where $X_i$ has degree  $k_i\in \Z$. Similarly, the full subcategory  $\Fil^{\geq 0}(\ob{D}_{\geq 0}(\bb{A}^1_{\Z}/\bb{G}_{m}))$ of modules $M$ such that $\Fil^0 M=\Fil^{-n} M$ for $n\geq 0$ is compactly projective generated with generators $\bigoplus_{i=1}^n \Z[t] t^{k_i}$ with $k_i \leq 0$. Thus, $\Fil^{\geq 0}(\Ring^{\delta})$ is also compactly generated with generators given by the graded $\Z[t]$-algebras $\Z[t][X_1,\ldots,X_s]$ as above with $|X_i|=k_i\leq 0$. The description of the right adjoint of the inclusion $\Fil^{\geq 0}(\Ring^{\delta})\subset \Fil(\Ring^{\delta})$ follows from animation, we leave the details to the reader. 
\end{proof}

\begin{example}\label{ExampleAdicFiltration}
Let us give a more geometric interpretation of the $I$-adic filtration in the case of a line bundle. Let $f\colon \Spec R\to \bb{A}^1_{\Z}/\bb{G}_{m}$ be a map of algebraic stacks associated to the Cartier divisor $I\to R$. Let $I^{\bullet}$ be the  filtered ring  associated to the map $f$ as in \Cref{ExampleIadicFiltration}. Then, the $I$-adic filtration of $I\to R$ is the positive truncation 
\[
\Fil_{\ob{ad}}^{\bullet}(R\to R/^{\bb{L}}I) = \Fil^{\geq 0}(I^{\bullet}).
\]
\end{example}

Given a Koszul regular pair $f\colon A\to B$ with  adic filtration $\Fil^{\bullet}_{\ob{ad}}(A\to B)= (I^{n})_{n\in \N}$, the  \textit{associated graded ring}  is given by $\gr^{\bullet}_{\ob{ad}}(A\to B)= \bigoplus_{n\in \N} I^n/I^{n+1}$. The Koszul regularity of $f$ yields that 
\[
\gr^{\bullet}_{\ob{ad}}(A\to B)= \Sym_{B} (I/I^2).
\]
Furthermore, as $I$ is generated by a Koszul regular sequence, a standard computation of cotangent complexes implies that the natural map $I\to \bb{L}_{B/A}[-1]$ induces an equivalence of $B$-modules  $I/I^2\xrightarrow{\sim} \bb{L}_{B/A}[-1]$.  Animating this construction we obtain the following result.

\begin{proposition}[{\cite[Corollary 3.61]{MaoCrystalline}}]\label{PropGraded}
The composite functor $\cat{Pair}\to \Fil(\Ring^{\delta})\xrightarrow{\gr} \ob{Gr}(\Ring^{\delta})$ sends a pair $A\to B$  to the graded animated ring $\Sym_{B}^{\bullet} (\bb{L}_{B/A}[-1])$. 
\end{proposition}

\begin{remark}\label{RemarkMapFilteredRings}
Let $A\to B$ be a pair of discrete animated rings. There is a natural morphism of filtered animated rings $A\to \Fil_{\ob{ad}}(A\to B)$, where $A$ is endowed with the trivial filtration given by $\Fil^{n}(A)=A$ for $n\leq 0$ and $0$ for $n>0$. Indeed, this can be proven by animation from the construction of the adic filtration of $\cat{Pair}^{\ob{pol}}$. Similarly, one has a map of filtered $A$-algebras  $\Fil_{\ob{ad}}(A\to B) \to B$.
\end{remark}

In the construction of the deformation to the normal cone we shall need to compute the right adjoint of the adic filtration functor.

\begin{lemma}\label{LemmaAdjointAdicFiltration}
Consider the functor $\Fil_{\ob{ad}}^{\bullet}\colon \cat{Pair}\to \Fil(\Ring^{\delta})$ from pairs to filtered animated rings given by the  adic filtration. Then $\Fil_{\ob{ad}}^{\bullet}$ factors through $\Fil^{\geq 0}(\Ring^{\delta})$,  its restriction to a functor 
\[
\Psi\colon \cat{Pair}\to \Fil^{\geq 0}(\Ring^{\delta})
\]
is fully faithful, and it  admits a right adjoint given by sending a positively filtered animated ring $\Fil^{\geq 0} A$ to the pair $\Fil^0 A \to \gr^0 (A)$.  
\end{lemma}
\begin{proof}
In this proof we will see $\Fil(\Ring^{\delta})$ as graded $\Z[t]$-algebras with $|t|=-1$.  Recall that $\cat{Pair}$ is compactly projectively generated by the pairs $\Z[\underline{x}, \underline{y}] \to \Z[\underline{y}]$ where $\underline{x}$ and $\underline{y}$ are finite sets of variables. The functor $\Psi$  sends $\Z[\underline{x}, \underline{y}] \to \Z[\underline{y}]$ to the graded $\Z[t]$-algebra $\Z[t][\underline{x}, \frac{y}{t}]$, i.e. the free animated $\Z[t]$-algebra generated by elements $\underline{Y}=\frac{y}{t}$ of degree $-1$ and $\underline{X}=\underline{x}$ of degree $0$. Then, from the description of $\Fil^{\geq 0}(\Ring^{\delta})$ as the animation of the graded $\Z[t]$-algebras $\Z[t][X_1,\ldots, X_n]$ with $|X_i|\leq 0$, it follows that the functor $\Psi$ has by right adjoint the functor sending a positively filtered animated ring $\Fil^{\geq 0} A$ to the pair $\Fil^{0} A\to \ob{gr}^0(A)$. The functor $\Psi$ is then fully faithful as for a pair $A\to B$, one has that $\Fil^0_{\ob{ad}}(A\to B)=A$ and $\gr^0_{\ob{ad}}(A\to B)=B$ (with the obvious map between them being the original $A\to B$). 
\end{proof}

\begin{construction}\label{DefinitionNormalCone}
Let $f\colon A\to B$ be a pair of discrete animated rings and $\Fil_{\ob{ad}}(A\to B)$ the filtered animated ring associated to $f$ seen as an animated ring over the stack $\bb{A}^1_{\bb{Z}}/\bb{G}_m$. The \textit{deformation to the normal cone of $f$}, denoted by $\widetilde{\n{N}}_f$, is the relative spectrum $\Spec_{\bb{A}^1_{\bb{Z}}/\bb{G}_m}(\Fil_{\ob{ad}}(A\to B))\to \bb{A}^1_{\Z}/\bb{G}_{m}$. By construction, the fiber of $\widetilde{\n{N}}_f$ at $\Spec \bb{Z}=\bb{G}_{m}/\bb{G}_{m}\to \bb{A}^1_{\Z}/\bb{G}_{m}$ is given by $\Spec(A)$ while the fiber at $B\bb{G}_{m}$ is given by the \textit{normal cone}, i.e.  the spectrum of the graded ring $\n{N}_f=\Spec_{B\bb{G}_{m}} (\Sym^{\bullet}_{B} (\bb{L}_{B/A}[-1]))\to B\bb{G}_{m}$. 

There are  natural maps $A\to \Fil_{\ob{ad}}(A\to B)\to B$ of filtered rings giving rise to maps of stacks over $\bb{A}^1_{\Z}/\bb{G}_{m}$ 
\[
\Spec(B)\times \bb{A}^1_{\Z}/\bb{G}_{m} \to \widetilde{\n{N}}_f \to \Spec(A)\times \bb{A}^1_{\Z}/\bb{G}_{m}
\]
functorial on the pair $f$. 
\end{construction}

\begin{remark}\label{DescentSchemes}
Let us explain a more geometric construction of the deformation to the normal cone. Let $X$ be an algebraic stack. The \textit{cone of $f$} is the stack 
\[
X^{\ob{cone}}\to \bb{A}^1_{\Z}/\bb{G}_{m}
\]
whose values at a ring $R$ with a Cartier divisor $I\to R$ is the anima 
\[
X^{\ob{cone}}(R,I)=X(R/I).
\]
The natural map $R\to R/I$ gives rise to a natural transformation  of stacks over  $\bb{A}^1_{\Z}/\bb{G}_{m}$
\[
X\times \bb{A}^1_{\Z}/\bb{G}_{m}\to X^{\ob{cone}}.
\]
Given $f\colon Y\to X$  a morphism of algebraic stacks, the \textit{cone of $f$} is the pullback 
\[
Y^{\ob{cone}/X}:= Y^{\ob{cone}}\times_{X^{\ob{cone}}} (X\times \bb{A}^1_{\Z}/\bb{G}_{m}).
\]
The pullback of $Y^{\ob{cone}/X}$ along $\bb{G}_{m}/\bb{G}_{m}$ sends a map $\Spec(R)\to X$ to the anima $Y(R/0)=Y(0)=*$. Hence, 
\[
Y^{\ob{cone}/X}\times_{\bb{A}^1_{\Z}/\bb{G}_{m}}  \bb{G}_{m}/\bb{G}_{m}=X.
\]
On the other hand, the pullback along $B\bb{G}_{m}\to \bb{A}^1_{\Z}/\bb{G}_{m}$ is the stack sending a map $\Spec(R)\to X \times B\bb{G}_{m}$ with line bundle $\n{L}$ over $R$ to the anima of maps $Y^{\ob{cone}/X}(R)=Y(R\oplus \n{L}[1])$. In particular, if $f$ admits a cotangent complex, one has that 
\[
Y(R\oplus \n{L}[1]) = \tau_{\geq 0}(\Hom_{Y}(\bb{L}_f, \n{L}[1]). 
\]
In particular, if $f\colon Y\to X$ is a Zariski closed immersion of schemes, one has that
\[
Y^{\ob{cone}/X}\times_{\bb{A}^1_{\Z}/\bb{G}_{m}} B\bb{G}_m = \Spec_Y(\Sym_{Y} (\n{L}^{-1}\otimes \bb{L}_{Y/X}[-1]))
\]
is the normal cone of $Y$ over $X$. 
\end{remark}

\begin{remark}\label{RemarkComparisonNormalConeConstructions}
Let us compare the construction of \Cref{DescentSchemes} with that of \Cref{DefinitionNormalCone}. Let $Y\to X$ be the Zariski closed immersion associated to a pair $A\to B$. Given $\Spec(R)\to X\times \bb{A}^1_{\Z}/\bb{G}_{m}$, the space $Y(R/I)$ is the anima of commutative squares of animated rings 
\[
\begin{tikzcd}
R \ar[r]  \ar[d]&  A \ar[d] \\ 
R/I \ar[r] & B.
\end{tikzcd}
\]
By  the fully faithfulness of \Cref{LemmaAdjointAdicFiltration}, this space  is the same as the anima  of maps $\Fil_I^{\bullet} R\to \Fil_{\ob{ad}}(A\to B)$ of filtered animated rings, where $\Fil_I^{\bullet} R$ is the $I$-adic filtration of $R$. This proves that 
\[
Y^{\ob{cone}/X}=\Spec_{\bb{A}^1/\bb{G}_m}(\Fil_{\ob{ad}}(A\to B))
\]
as wanted. 
\end{remark}

\begin{remark}\label{RemarkNormalConeCondensed}
Finally, let us explain how   \Cref{DefinitionNormalCone} naturally promotes to analytic rings.  Let $\n{X}$ be an $\infty$-topos, and $\mathcal{C}$ a presentable $\infty$-category. By \cite[Proposition 4.8.1.17]{HigherAlgebra} Lurie's tensor product $\n{X}\otimes \n{C}$ is naturally equivalent to the category of sheaves $\ob{Shv}(\n{X}, \n{C})$ of $\n{X}$ valued in $\n{C}$. Applying this to $\n{C}$ as the categories of discrete animated rings, pairs, filtered rings, etc,  one obtains the categories of animated rings on $\n{X}$, pairs on $\n{X}$, filtered animated rings on $\n{X}$, etc. Furthermore, the functor $\Fil_{\ob{ad}}\colon \cat{Pair}\to \Fil(\Ring^{\delta})$ preserves colimits, and thus it defines a colimit preserving functor 
\[
\Fil_{\ob{ad}}\colon  \ob{Shv}(\n{X}, \cat{Pair})\to \ob{Shv}(\n{X}, \Fil(\Ring^{\delta}))
\]
which is the \textit{adic filtration} from the category of  pairs over to filtered animated rings over $\n{X}$. Since passing to the graded ring also commutes with colimits, we also have an induced functor 
\[
\ob{gr}\colon \ob{Shv}(\n{X}, \Fil(\Ring^{\delta}))\to \ob{Shv}(\n{X}, \ob{Gr}(\Ring^{\delta}))
\]
whose composite with $\Fil_{\ob{ad}}$ sends a pair $\s{A}\to \s{B}$ of animated rings over $\n{X}$ to the symmetric algebra $\Sym_{\s{B}}(\bb{L}_{\s{B}/\s{A}}[-1])$, where  $\bb{L}_{\s{B}/\s{A}}[-1])$ is  defined as in the animated variant of the cotangent complex of  \cite[Section 17.1]{LurieSpectralAlg}.

 Applying this construction to the topos of condensed anima $\CondAni$, we obtain a condensed variant of the adic filtration.  Now, a \textit{pair of analytic rings}  is a morphism of analytic rings  $A\to B$  with the induced analytic ring structure and surjective on $\pi_0$. The associated  morphism of condensed rings $A^{\triangleright}\to B^{\triangleright}$ is a pair of animated condensed rings, and so it has an adic filtration $\Fil_{\ob{ad}}(A^{\triangleright}\to B^{\triangleright})$ together with a morphism of filtered condensed rings $A^{\triangleright}\to \Fil_{\ob{ad}}(A^{\triangleright}\to B^{\triangleright})$ where $A^{\triangleright}$ has the trivial filtration as in \Cref{DefinitionNormalCone}.  Taking $A$-completions, we obtain the \textit{adic completion of analytic rings} as the morphism of filtered animated $A$-algebras (with induced structure from $A$)  $A\to \Fil_{\ob{ad}}(A\to B)$. The graduation of  $\Fil_{\ob{ad}}(A\to B)$ is then isomorphic to the animated graded $B$-algebra $\Sym^{\bullet}_{B} (\bb{L}_{B/A}[-1])$, where $\bb{L}_{B/A}$ is the cotangent complex of animated rings. Given a pair $f\colon A\to B$ of analytic rings, the adic filtration gives rise to  a morphism of analytic stacks 
 \[
 \widetilde{\n{N}}_{f}\colon= \AnSpec(\Fil_{\ob{ad}}(A\to B)) \to  \AnSpec(A)\times \bb{A}^{1,\ob{cond}}_{\Z}/\bb{G}^{\ob{cond}}_{m}
 \]
  where the left hand side is a relative analytic spectrum, and  $\bb{A}^{1,\ob{cond}}_{\Z}/\bb{G}^{\ob{cond}}_{m}$ is the  analytic stack obtained from the fpqc algebraic stack as in  \Cref{ssAlgebraicStacks}.  We call $\widetilde{\n{N}}_{f}$ the \textit{deformation to the normal cone of $f$}.  In total, we have a diagram of analytic stacks over $\bb{A}^{1,\ob{cond}}_{\Z}/\bb{G}^{\ob{cond}}_{m}$
 \[
 \AnSpec(B)\times \bb{A}^{1,\ob{cond}}_{\Z}/\bb{G}^{\ob{cond}}_{m}\to \widetilde{\n{N}}_{f} \to \AnSpec(A)\times \bb{A}^{1,\ob{cond}}_{\Z}/\bb{G}^{\ob{cond}}_{m}. 
 \]
One can also define a general cone construction for analytic stacks as in \Cref{DescentSchemes} that agrees with the previous definition in the affine case. 
\end{remark}

In the next lemmas we study the some basic geometric properties of the deformation to the normal cone. To light notation, given an algebraic stack $S$, we will simply write $S$ for its realization $S^{\ob{cond}}$ as an analytic stack. 

\begin{lemma}\label{LemmaNormalConeEtale}
 The following hold:

\begin{enumerate}

\item The functor sending an analytic stack $X$ to the ($!$-sheafification of the) cone $X^{\ob{cone}}$ commutes with pullbacks.

\item  Consider a pushout diagram of analytic rings 
\[
\begin{tikzcd}
A \ar[r] \ar[d] & B \ar[d] \\ 
A' \ar[r] & B'
\end{tikzcd}
\]
such that $A\to B$ is a pair (and so it is $A'\to B'$). Let  us denote $X=\AnSpec A$, $Y=\AnSpec B$ (and similar for $A'$ and $B'$). Then the following diagram has cartesian squares
\[
\begin{tikzcd}
Y' \times \bb{A}^1_{\Z}/\bb{G}_m \ar[r] \ar[d] & Y^{\prime,\ob{cone}/X'} \ar[r] \ar[d] & X'\times \bb{A}^1_{\Z}/\bb{G}_m \ar[d] \\ 
Y \times \bb{A}^1_{\Z}/\bb{G}_m \ar[r] &  Y^{\ob{cone}/X} \ar[r] & X \times \bb{A}^1_{\Z}/\bb{G}_m
\end{tikzcd}
\]

\item Consider a diagram of analytic rings $A\to B\to  C$ with corresponding diagram of analytic stacks $Z\to Y \to X$. Suppose $A\to B$ and $B\to C$ are pairs and that   $\bb{L}_{C/B}=0$.  Then the map $B\to C$  is idempotent, the following  is a pullback diagram of analytic stacks 
\[
\begin{tikzcd}
Z^{\ob{cone}/Y} \ar[r] \ar[d] & Z^{\ob{cone}/X} \ar[d] \\
Y\times \bb{A}^1_{\Z}/\bb{G}_m  \ar[r] & Y^{\ob{cone}/X},
\end{tikzcd}
\]
and the map $ Z\times \bb{A}^1_{\Z}/\bb{G}_m \xrightarrow{\sim} Z^{\ob{cone}/Y}\times_{Y\times \bb{A}^1_{\Z}/\bb{G}_m} (Z\times \bb{A}^1_{\Z}/\bb{G}_m)$ is an equivalence. In particular, $ Z\times \bb{A}^1_{\Z}/\bb{G}_m \to Z^{\ob{cone}/Y}$ is a closed immersion   of analytic stacks \Cref{DefOpenClosed}, and if  $C$ is a dualizable $B$-module it is a clopen immersion.

\item Let $\bb{G}^n_{a}=\AnSpec \Z[T_1,\ldots, T_n]$ seen as a vector bundle and consider the zero section $\iota_0\colon \AnSpec \Z \to \bb{G}_{a}^n$. Let $t$ be the   coordinate of $S:=\bb{A}^1_{\Z}/\bb{G}_{m}$ with degree $-1$. Then,  we have the following presentations:
\begin{itemize}
\item[(a)] $\bb{G}^{\ob{cone}}_{a}= \bb{G}_{a,S}/ t \bb{G}_{a,S}$ where $t\in \bb{G}_{a,S}(S)$ is seen as a Cartier divisor of the ring stack $\bb{G}_{a,S}$. 
\item[(b)] $\widetilde{\n{N}}_{\iota_0}= t\bb{G}^n_{a,S}$ and the map $\widetilde{\iota}_0\colon S\to \widetilde{\n{N}}_{\iota_0}= t\bb{G}^n_{a,S}$  is the zero section of $t\bb{G}^n_{a,S}$. Moreover, we have a natural identification (up to an universal sign as in \Cref{PropositionUniversalNormalCone}) 
\[
\widetilde{\iota}_0^! 1 [n]= \det (t\bb{G}^n_{a,S}) = t^{n} \bb{G}_a=\s{O}(-n). 
\]
is the line bundle of degree $-n$ (equivalently, the $n$-th power of the universal Cartier divisor of $S$). 
\end{itemize}

\end{enumerate}
\end{lemma}
\begin{proof}
\begin{enumerate}

\item  This follows from the functor of points description of $X\mapsto X^{\ob{cone}}$. 

\item This follows from \Cref{RemarkComparisonNormalConeConstructions}, and the fact that $\ob{Fil}_{\ob{ad}}(A\to B)\otimes_{A} A' = \ob{Fil}_{\ob{ad}}(A'\to B')$ by the colimit preserving property of $\ob{Fil}_{\ob{ad}}$. 

\item The fact that $B\to C$ is idempotent follows from the fact that $C\otimes_B C\to C$ is a formally \'etale morphism of connective rings with which induces an isomorphism on connected components (see \Cref{CoroDetectCotangent}). To prove the second claim,  similarly as in (2), by translating in terms of filtered animated rings we have to show that the natural map 
\[
\ob{Fil}_{\ob{ad}}(A\to C)\otimes_{\ob{Fil}_{\ob{ad}}(A\to B)} B\to C
\]
is an isomorphism. By translating the problem to pairs via the fully faithful embedding of \Cref{LemmaAdjointAdicFiltration}, we have to show that the natural map  of pairs
\[
(A\to C)\otimes_{(A\to B)} (B\xrightarrow{\id} B) \to (B\to C)
\]
is an isomorphism, this follows from the fact that pushouts in $\cat{Pair}$ are computed pointwise.

 To see that $Z\times \bb{A}^1_{\Z}/\bb{G}_m \xrightarrow{\sim} Z^{\ob{cone}/Y}\times_{Y\times \bb{A}^1_{\Z}/\bb{G}_m} (Z\times \bb{A}^1_{m}/\bb{G}_m)$ is  an equivalence, it suffices to verify the analogue statement at the level of pairs. This translates to the fact that $(B\to C)\otimes_{(B\to B)} (C\to C) = (C\to C)$ as $B\to C$ is idempotent. If in addition $C$ is a dualizable $B$-module, then $Z\to Y$ is a clopen immersion by \Cref{LemmaIdempotentDualizable}, and the same holds for $Z\times \bb{A}^1_{\Z}/\bb{G}_m\to Z^{\ob{cone}/Y}$ by base change.

\item The description of $\bb{G}_a^{\ob{cone}}$ follows by definition: it is the $!$-sheafification of the stack sending a ring $R$ with Cartier divisor $I\to R$ to the quotient $R/I$. It is also clear that the map $S\to \bb{G}_a^{\ob{cone}}$ arising from the map of pairs $(\Z[t]\to \Z) \to (\Z\to \Z)$ gives rise to the zero section of the ring stack $\bb{G}_a^{\ob{cone}}$. It follows that $\widetilde{\n{N}}_{\iota_0}=\ob{fib}( \bb{G}_{a,S}^n\to \bb{G}_a^{\ob{cone},n})= t\bb{G}_{a,S}^n$ as wanted. The computation of $\widetilde{\iota}_0^! 1$ follows from \Cref{PropositionUniversalNormalCone}.  

\end{enumerate}
\end{proof}

\begin{lemma}\label{LemmaNormalConeAffineLine}
Let $A\to B$  be a pair of analytic rings  such that $B$ is of the form $B=A/^{\bb{L}} (a_1,\ldots, a_n)$ for certain elements $a_i\in \pi_0(A)$. Let $f\colon Y\to X$ be the associated map of analytic stacks.  The following hold: 
\begin{enumerate}

\item The map $\widetilde{f}\colon Y \times \bb{A}^1_{\bb{Z}}/ \bb{G}_{m}\to \widetilde{\n{N}}_{f}$ is a  complete intersection of codimension $n$. In particular, it is cohomologically smooth. 

\item Consider the map $g\colon \bb{A}^{1}_{\Z}/\bb{G}_{m}\to \AnSpec \Z$. Then $g$ is prim and the natural map $\Z\to g_*1$ is an equivalence. In particular, $g^*$ is fully faithful and the same holds after base change along any morphism of analytic stacks $S\to \AnSpec \Z$.

\item For $k\in \Z$ let  $\s{O}(k)$ be the line bundle over $B\bb{G}_{m}$ of degree $k$. Then $\widetilde{f}^! 1 \otimes \s{O}(n)$ belongs to the essential image of the pullback along $g\colon Y \times \bb{A}^1_{\bb{Z}}/ \bb{G}_{m} \to Y$. In particular, if $\iota_0$ and $\iota_1$ denote the sections $\AnSpec \Z \to \bb{A}^1_{\Z}/\bb{G}_{m}$ at $0$ and $1$ respectively, there is a natural equivalence of sheaves over $Y$ 
\[
\iota^*_0(\widetilde{f}^! 1 \otimes \s{O}(n)) = \iota^*_1(\widetilde{f}^! 1 \otimes \s{O}(n)). 
\]

\item In the situation of (3), we have a natural equivalence (up to an universal sign)
\[
\iota^*_0(\widetilde{f}^{!}1\otimes \s{O}(n))= \det(N_f)[-n]
\]
where $\n{N}_f=(\bb{L}_{Y/X}[-1])^{\vee}$ is the normal cone of $f$. 
\end{enumerate}
\end{lemma}
\begin{proof}
\begin{enumerate}

\item We can find a morphism of analytic rings $\Z[T_1,\ldots, T_n]\to A$ sending $T_i\mapsto a_i$ such that $B=A\otimes_{\Z[T_1,\ldots, T_n]} \Z$. This yields a cartesian square of analytic stacks 
\[
\begin{tikzcd}
Y \ar[r] \ar[d] & X \ar[d] \\ 
\Spec \Z \ar[r] & \bb{A}^{n}_{\Z}.
\end{tikzcd}
\]
By \Cref{LemmaNormalConeEtale} (2) we have a Cartesian square in the corresponding deformations to the normal cones, and by base change it suffices to consider the case of the zero section $\iota\colon \Spec \Z\to \bb{A}^n_{\Z}$. This case follows from the explicit computation of \Cref{LemmaNormalConeEtale} (4.b).

\item Let us see that $g$ is $!$-able and prim. We rite $g$ as the composite $\bb{A}^1_{\Z}/\bb{G}_m\to B \bb{G}_m \to \Spec \Z$, the first map is prim being represented by affine schemes, it is left to show that $B\bb{G}_m\to \Spec \Z$ is prim. For that,  by \cite[Lemma 4.7.4]{HeyerMannSix} it suffices to show that $h\colon \Spec\Z \to B \bb{G}_m$ is prim and descendable; it is prim being represented in affine schemes, it is descendable since $h_* 1 =\bigoplus_{n\in \Z} \mathcal{O}(n)$ and the unit of $B\bb{G}_m$ is a retract of it.  Finally, the fact that $\Z\to g_*1$ is an equivalence follows from the fact that $g_*$ is the functor sending a $\Z$-graded module to its degree $0$-part. The fully faithfulness of $g^*$ (even after base change) follows from linearity over the base and the previous computation.

\item  By \Cref{LemmaNormalConeEtale} (2) and the discussion of (1) above, it suffices to prove the statement for the case of the zero section $\Spec \Z\to \bb{A}^n_{\Z}$. The statement then follows from \Cref{LemmaNormalConeEtale} (4.b). The equivalence  along the pullbacks of $\iota_0$ and $\iota_1$ follow from the fact that the composite with $Y\times \bb{A}^1_{\Z}/\bb{G}_m\to Y$ agree.

\item It is left to compute $\iota^*_0(\widetilde{f}^! 1\otimes \s{O}(n))$. For that, consider the pullback along the zero section $\Spec \Z \to \bb{A}^1_{\Z}/\bb{G}_m $. Then $\widetilde{\n{N}}_f\times_{ \bb{A}^1_{\Z}/\bb{G}_m} \AnSpec \Z= \n{N}_f = \AnSpec (\Sym_B(\bb{L}_{B/A}[-1]))$ is the normal cone of $f$. In particular, $\n{N}_f\to \AnSpec B$ is a vector bundle and the pullback  $h\colon \AnSpec B\to \n{N}_f$ of $\widetilde{f}$ is the zero section. By suave base change we have that 
\[
\iota^*_0(\widetilde{f}^! 1 \otimes \s{O}(n))= h^! 1
\]
where we have used a fixed trivialization of $\s{O}(n)$ (only depending on $n$) after pullback along $\AnSpec \Z\to \bb{A}^1_{\Z}/\bb{G}_m$. The identification 
\[
h^! 1 = \det(\n{N}_f)[-n]
\]
follows from \Cref{PropositionUniversalNormalCone}.
\end{enumerate}
\end{proof}

\begin{theorem}\label{TheoSerreDualityDualizing}
Let $f\colon A\to B$ be a locally solid smooth morphism of solid Huber rings of dimension $n$.  Then the deformation to the normal cone of $f$ give rise an isomorphism $f^! 1 = \det( \bb{L}_{B/A}) [\dim f]$ (natural up to the choice of a sign only depending on the dimension of $f$). 
\end{theorem}
\begin{proof}
Let $X=\AnSpec B$, $S=\AnSpec A$ and $f\colon X\to S$. Consider the diagonal map $X\xrightarrow{\Delta_S} X\times_S X$, and let $\ob{pr}_i\colon X\times_S X\to X$ be the $i$-th projection map. By suave base change, we have that
\[
f^! 1 = \Delta^*(\ob{pr}_1^* f^! 1 ) = \Delta^* \ob{pr}_2^! 1. 
\]
Hence, by the diagonal trick, we can assume without loss of generality that the map $f\colon X\to S$ has a section $s\colon S\to X$, and we want to produce a natural isomorphism 
\[
s^*f^! 1 = s^*(\det \bb{L}_{X/S})[n].  
\]  
Since $s\colon S\to X$ is a Zariski closed immersion of locally solid smooth affinoid stacks over $S$, by \Cref{PropStabilityolidEtaleSmooth} (6) the map $s$ is, locally in the solid analytic topology of $X$, a  complete intersection of codimension $n$. In particular, $s$ is suave and we have that 
\[
1= s^! f^! 1 = s^!1 \otimes s^* f^! 1. 
\]
Therefore, it suffices to produce a natural equivalence 
\[
s^! 1 = s^*(\det \bb{L}_{X/S}^{\vee})[-n]
\]
where $ \bb{L}_{X/S}^{\vee}=N_s=\bb{L}_{S/X}[-1]$ is the normal cone of $s$. To do that, consider the deformation to the normal cone fitting in the diagram 
\[
\begin{tikzcd}
 \widetilde{\n{N}}_f \ar[r,"\widetilde{f}"] & X\times\bb{A}^1_{\Z}/\bb{G}_m \ar[d] \\
& S\times \bb{A}^1_{\Z}/\bb{G}_m \ar[lu,"\widetilde{s}"].
\end{tikzcd}
\]
Since $s\colon S\to X$ is locally in the solid analytic topology of $X$ a  complete intersection of codimension $n$, and the formation of the deformation to the normal cone is stable under pullbacks as in \Cref{LemmaNormalConeEtale} (2) and (3), one sees by \Cref{LemmaNormalConeAffineLine} that the map $\widetilde{s}$ is, locally in the analytic topology of $X$, a  complete intersection of codimension $n$,  and hence that it is suave.  On the other hand, let $\iota_0,\iota_1\colon S\to S\times \bb{A}^1_{\Z}/\bb{G}_m$  denote the section at $0$ and $1$ respectively, by  \Cref{LemmaNormalConeAffineLine} (3) the sheaf $\widetilde{s}^!1\otimes  \s{O}(n)$ lies in the essential image of the fully faithful functor $\ob{D}(S)\to \ob{D}(S\times \bb{A}^1_{\Z}/\bb{G}_m)$, and hence we have natural equivalences
\[
\iota^*_0(\widetilde{s}^! 1 \otimes \s{O}(n))=\iota^*_1(\widetilde{s}^! 1\otimes \s{O}(n)).
\]
The pullback of $\widetilde{f}\colon \widetilde{N}_f\to X\times \bb{A}_{\Z}/\bb{G}_m$ along the section $1$ is an isomorphism, this produces a natural equivalence 
\[
\iota^*_1(\widetilde{s}^! 1\otimes \s{O}(n)) = s^! 1.
\]
On the other hand,  \Cref{LemmaNormalConeAffineLine} (4) produces a natural equivalence $\iota^*_0(\widetilde{s}^! 1 \otimes \s{O}(n))= \det(N_s)[-n]$, which yields a natural equivalence 
\[
s^!1=\det(N_s)[-n]
\]
as wanted. 
\end{proof}

\begin{remark}\label{RemarkGlobalization}
Both \Cref{TheoLocalSerreDuality,TheoSerreDualityDualizing} have natural extensions to solid adic spaces by descent along the analytic topology thanks to the functoriality of the  cone construction as in \Cref{DescentSchemes} and \Cref{LemmaNormalConeEtale}. Moreover, since \Cref{LemmaNormalConeEtale,LemmaNormalConeAffineLine} hold for general maps of analytic rings, in order to prove Serre duality in other contexts (eg. for the gaseous theory  in complex analytic geometry), one only needs to have an analogue of  \Cref{PropStabilityolidEtaleSmooth} and  \Cref{TheoLocalSerreDuality}. 
\end{remark}

For future reference we make explicit Serre duality in the case of schemes and rigid varieties.

\begin{corollary}\label{CorollarySerreDuality}
\begin{enumerate}
\item Let $f\colon Y\to X$ be a smooth morphism of schemes adic spaces of relative dimension $n$, and let $Y_{\sol}\to X_{\sol}$ be the associated morphism of solid stacks. Then, $f$ is cohomologically smooth and there is a natural identification $\omega_f= \det \Omega^1_{Y/X}[n]$.

\item Let $f\colon Y\to X$ be a smooth morphism of rigid varieties over a complete non-archimedean field $K$. Then $f$ is cohomologically smooth and there is a natural identification $\omega_d=\det \Omega^1_{Y/X}[n]$.

\end{enumerate}

\end{corollary}
\begin{proof}
To prove cohomological smoothness, we can localize in the adic topology on both $Y$ and $X$ and assume without loss of generality that $Y$ and $X$ are affinoid, and that the map $Y\to X$ is  smooth. By further localizing we can assume that $f$ factors through an \'etale map $Y\to \bb{A}^n_X\to X$ where in the case of (b) we  consider an analytic affine line. In the case of (a), the affine line is cohomologically smooth by \Cref{PropSerreDualityAffine}, in the case of (b) it is smooth since it is a union of open affinoid discs and $(K,K^+)_{\sol}\otimes_{\Z_{\sol}}\Z[T]_{\sol}= (K\langle T\rangle, K^+\langle T\rangle)_{\sol}$ is smooth over $(K,K^+)_{\sol}$.  Then, as $Y$ is smooth over $X$, the map $Y\to \bb{A}^n_X$ is a local complete intersection, and so it is cohomologically smooth by \Cref{LemmaPrimAffine}.

 The computation of the dualizing sheaf follows from a deformation to the normal cone argument as in \Cref{TheoSerreDualityDualizing}.
\end{proof}

\begin{remark}\label{RemarkClassicalSerreDuality}
Let $X$ be a proper smooth scheme over a ring $A$. \Cref{CorollarySerreDuality} implies that the map of algebraic stacks $f\colon X\to \Spec A$ is cohomological smooth for the classical six functor formalism of quasi-coherent sheaves. Indeed, we want to show that $f_*$ is naturally a left adjoint to $\omega_f\otimes f^*$ with $\omega_f=\det \Omega_{X/A} [\dim f]$, compatibly with pullback along $\Spec A'\to \Spec A$. By \Cref{CorollarySerreDuality} we know that this holds for the solid six functor formalism. Denote $f_{\sol}\colon X_{\sol}\to \AnSpec A_{\sol}$, since $f$ is proper, the functor $f_{\sol,*}=f_{\sol,!}$ preserves discrete objects. Hence, the adjunction $f_{\sol,*}\rightleftharpoons \omega_f\otimes f_{\sol}^*$ restricts to an adjunction at the level of classical quasi-coherent sheaves. Since pullbacks of solid sheaves along maps of schemes preserve discrete objects, the adjunction $f_{*}\rightleftharpoons \omega_f\otimes f^*$ is preserved under base change of classical quasi-coherent sheaves along $\Spec A'\to \Spec A$, proving that $f$ is cohomologically smooth as wanted. 
\end{remark}


\bibliographystyle{alpha}
\bibliography{BiblioSolidGeometryI}
\end{document}